\documentclass[A4,reqno]{amsart}
\usepackage[dvips]{graphicx}
\usepackage{euscript}
\usepackage{amscd}
\usepackage{amssymb,latexsym}
\usepackage{mathrsfs}
\usepackage[all]{xy}
\usepackage{amssymb}
\usepackage[latin1]{inputenc}
\xyoption{all}

\newtheorem{thm}{Theorem}[section]

\newtheorem{prop}[thm]{Proposition}
\newtheorem{cor}[thm]{Corollary}
\newtheorem{lma}[thm]{Lemma}
\newtheorem{rem}[thm]{Remark}
\newenvironment{rmk}{\begin{rem}\rm}{\end{rem}}

\numberwithin{equation}{section}

\newcommand{\R}{{\mathbb{R}}}
\newcommand{\Z}{{\mathbb{Z}}}
\newcommand{\C}{{\mathbb{C}}}
\newcommand{\Q}{{\mathbb{Q}}}

\begin{document}
\title[Symplectic topology of representation varieties, I]{Symplectic topology of SU(2)-representation varieties and link homology, I:\\ Symplectic braid action \\ and the first Chern class}

\author[Magnus Jacobsson and  Ryszard Rubinsztein ]{Magnus Jacobsson and  Ryszard L. Rubinsztein}
\address{Dept. of Mathematics, Uppsala University, Box 480,
Se-751 06 Uppsala, Sweden;}
\email{Jacobsson: Magnus.Jacobsson{\@@}math.uu.se}
\email{Rubinsztein: ryszard{\@@}math.uu.se}

\begin{abstract}
There are some similarities between cohomology of $SU(2)$-represen- tation varieties of the fundamental group of some link complements and the Khovanov homology of the links. We start here a program to explain a possible source of these similarities. We introduce a symplectic manifold ${\mathscr M}$ with an action of the braid group $B_{2n}$ preserving the symplectic structure. The action allows to associate a Lagrangian submanifold of  ${\mathscr M}$ to every braid. The representation variety of a link can then be described as the intersection of such Lagrangian submanifolds, given a braid presentation of the link.  We expect this to go some way in explaining the similarities mentioned above.  
\end{abstract}



\maketitle

\section{Introduction}\label{intro}

\bigskip

 Let $C$ be the conjugacy class of traceless matrices in the Lie group $SU(2)$.
The origin of this paper is the following observation. 

Let $L$ be a link in $\R ^3$. Consider 
the space $\widetilde{J}_C(L)$ of representations of the fundamental group of the link complement
in $SU(2)$ satisfying the condition that they send meridians to the conjugacy class $C$. Typically $\widetilde{J}_C(L)$ has several connected components 
$J_C^r(L)$. Let $$H^*(\widetilde{J}_C(L);\mathbb{Z}) = \bigoplus_r H^*(J_C^r(L);\mathbb{Z})$$ be the singular 
cohomology of this space.

Let, on the other hand,  $Kh^{i,j}(L)$ denote the integral Khovanov homology of the link $L$ 
and let $Kh^k(L)$ denote the singly graded homology theory obtained by collapsing the grading along $k = i - j$.

\subsection{Observation} For every prime knot $L$ with seven crossings or less, 
and for every $(2,n)$-torus link, there are integers $N_r = N_r(L)$ such that 
\begin{equation}\label{intro1}
Kh^*(L) = \bigoplus_r H^*(J_C^r(L);\mathbb{Z})\{N_r\}.
\end{equation}

That is, the Khovanov homology consists of pieces, which are isomorphic to the 
cohomology of the components of the representation variety $\widetilde{J}_C(L)$, but with 
each such component shifted by some integer.

To prove this observation let us recall that the conjugacy class $C$ is an example of a topological quandle. In  \cite{R1}, given a topological quandle $Q$ and an oriented link $L$ in $\R ^3$, the second author constructed an invariant space $J_Q(L)$ of the link. When the quandle $Q$ is a conjugacy class $K$ in a topological group $G$ (or a union of conjugacy classes), the invariant space $J_K(L)$ can be identified with the space of homomorphisms from the fundamental group of the link complement $\pi _1(\R ^3-L)$ to the group $G$ mapping the positively oriented meridians into $K$ (see \cite{R1}, Lemma 4.6). In particular, if $K$ is the conjugacy class $C$ of traceless matrices in $SU(2),$ the representation variety $\widetilde{J}_C(L)$ can be identified with the invariant space  $J_C(L)$ of \cite{R1}. From now on we shall use the symbol  $J_C(L)$ for both of them. 
 The invariant spaces   $J_C(L)$ of all prime knots with up to seven crossings and of all $(2,n)$-torus links can be easily determined by the techniques used in Section 5 of \cite{R1}. The Khovanov homology of such knots, on the other hand, can be computed by using the Mathematica package JavaKh \cite{KT}. The isomorphism (\ref{intro1}) then follows.

We leave details to the reader but we illustrate with 

\subsection{Example} The singly graded Khovanov homology $Kh^*(T_{2,n})$ for 
$(2,n)$-torus knots can be written as

$$Kh^*(T_{2,n}) = H^*(S^2)\{n-2\} \oplus \bigoplus_{i=1}^{\frac{n-1}{2}} 
H^*(\mathbb{R}P^3)\{2i - 2 + n\},$$ see \cite{K}.  And, as we compute in Section \ref{ex},  
$$J_C(T_{2,n}) = S^2 \cup \bigcup_{j=1}^{\frac{n-1}{2}} \mathbb{R}P^3.$$ 

For more examples, see Section \ref{ex}. 

\medskip

One of the future aims of our work is to enhance the link invariant $J_C(L)$ by introducing a 
grading on its connected components. The grading should explain the shifts in our examples. 

Note that  the group $SU(2)$ acts on  $J_C(L)$ by conjugation. If the space $J_C(L)$ 
is factored out by this action it becomes the  moduli space of flat $SU(2)$-connections 
on the link complement with holonomies along meridians belonging to the conjugacy class $C$. 

In \cite{GHJW}, Guruprasad, Huebschmann, Jeffreys and Weinstein study 
the moduli space $\mathcal{M}$ of flat $SU(2)$-connections on a surface with $n$ marked 
points. They express $\mathcal{M}$ by symplectic reduction of a certain finite-dimensional symplectic 
manifold $\mathscr{M}$ called an ``extended moduli space''. The space $\mathcal{M}$ plays a crucial role in the quantized 
Chern-Simons gauge theory, which is the 3-dimensional origin of the Jones polynomial, 
which in turn is the Euler characteristic of Khovanov homology.

The main idea of our paper is to express the representation variety $J_C(L)$ of the link $L$ as an intersection of two Lagrangian submanifolds of a symplectic manifold. We achieve this in the following way. The link $L$ can be represented as a closure of a braid $\sigma$ on $n$ strands. The braid group ${\mathcal B}_n$ on $n$ strands acts on the product $P_n$ of $n$ copies of $C$. The invariant space $J_C(L)$ is defined as the subspace of fix-points of $\sigma$ acting on $P_n$  \cite{R1}. We can also describe $J_C(L)$ as an intersection of a ``twisted'' diagonal of $P_n$ in $\, P_{2n}=P_n\times P_n$ with a ``twisted'' graph of $\sigma$. We shall use this description of $J_C(L)$. Thus we consider the product $P_{2n}$ of $2n$ copies of the conjugacy class $C$. Following Guruprasad,  Huebschmann, Jeffreys and Weinstein \cite{GHJW} we interpret  $P_{2n}$ as a space naturally associated to the moduli space of flat $SU(2)$-connections on the $2$-sphere $S^2$ with $2n$ punctures and with monodromies around  the punctures belonging to $\, C$. As shown in  \cite{GHJW},  $P_{2n}$ contains the extended moduli space  ${\mathscr M}$ mentioned above as an open submanifold with a symplectic structure. We show that the action of the braid group ${\mathcal B}_{2n}$ on  $P_{2n}$ maps the submanifold   ${\mathscr M}$ to itself and preserves its symplectic form. (It should, perhaps, be stressed that the symplectic structure on  ${\mathscr M}$ is {\it not} the restriction of the {\it product} symplectic structure on $P_{2n}=\prod C$. That one would not have been preserved by the action of ${\mathcal B}_{2n}$.) Furthermore, to every braid $\sigma$ on $n$ strands we associate a Lagrangian submanifold $\Gamma _{\sigma}$ in  ${\mathscr M}$. In particular, the identity braid defines a Lagrangian submanifold $\Lambda$. Finally,  $\Gamma _{\sigma}$ and $\Lambda$ intersect in a space homeomorphic to $J_C(L)$.

This picture fits the general scheme of a classical field theory. One may think 
of $L$ as composed of two braids; one inside a 2-sphere, piercing it in $2n$ points, 
and one outside (usually, but not necessarily, the identity braid). The theory 
associates a symplectic manifold to the 2-sphere with marked points and Lagrangian 
submanifolds to the 3-discs with braids which are bounded by it. 

In good cases (when the intersection is clean) we hope to be able to use the Lagrangians to 
define an index associated to the components of $J_C(L)$. This will be the 
topic of our next paper. In later work we plan to extend this to a Floer-type Lagrangian
homology theory based on the invariant spaces  $J_C(L)$. 

(The idea of describing $J_C(L)$ as the intersection of  $\Lambda$ and   $\Gamma _{\sigma}$ goes back to X-S.Lin, \cite{L2}. However, no symplectic topology was was involved in his considerations.)

It should be emphasized that there are examples where the cohomology of $J_C(L)$ does not give back 
Khovanov homology groups of $L$. The first such that we know of is the knot 
$9_{42}$ in Rolfsen's table. In section $10$ we determine the representation variety $\, J_C(L)\,$ for $\, L= 9_{42}\,$ and point out the difference between its cohomology and the collapsed Khovanov homology of the knot.  
 
\bigskip

As already indicated, this paper is the first one in a series. It is devoted to the construction of the symplectic structure on the manifold $\, {\mathscr M}\,$, to the  proof of its invariance under the braid group action, to the  definition of the Lagrangian manifolds $\, \Gamma _{\sigma}\,$ and to the study of the first Chern class of  $\, {\mathscr M}\,$. 

The contents of the present paper are as follows.

In Section 2 we study the topology of certain spaces $\, K_{2n}\,$. These spaces are closely related to the moduli spaces of flat $\, SU(2)$-connections on the $2$-dimensional sphere $\, S^2\,$ with $2n$ holes with holonomy around the holes belonging to the conjugacy class $\, C\,$ of traceless matrices in  $\, SU(2)\,$. To be exact,  $\, K_{2n}\,$ is the space $\, Hom (\pi , \,\, SU(2))_{\, C}\,$  of representations  of the fundamental group $\, \pi \,$ of the punctured sphere with values in  $\, SU(2)\,$ such that the standard generators of  $\, \pi \,$ are mapped to  $\, C\,$. The space  $\, K_{2n}\,$ is a manifold with singularities. Observe that we do not remove singularities. The space  $\, K_{2n}\,$ is simply connected  and we determine its second homotopy group $\, \pi _2(K_{2n})\,$ (Proposition 2.9). In order to do that we study the structure of the singularities of  $\, K_{2n}\,$. These are completely described in Appendix A.

The space  $\, K_{2n}\,$ is a subspace of the product $\, P_{2n}\,$ of 2n copies of $\, S^2\,$, $\, P_{2n} = \prod\limits_{i=1}^{2n} S^2\,$.  (Observe that the conjugacy class $\, C\,$ can be identified with $\, S^2\,$.) In Section $3$ we recall  the construction of Guruprasad, Huebschmann, Jeffrey and Weinstein, \cite{GHJW}, of a symplectic structure on an open neighbourhood $\, {\mathscr M}\,$ of   $\, K_{2n}\,$ in $\, P_{2n}\,$.

The space  $\, P_{2n}\,$ admits an action of the braid group $\, {\mathcal B}_{2n}\,$ on $2n$ strands. The manifold  $\, {\mathscr M}\,$ can be chosen to be invariant under the action of  $\, {\mathcal B}_{2n}\,$. In Section 4 we show that the braid group $\, {\mathcal B}_{2n}\,$ acts  on $\, {\mathscr M}\,$ by symplectomorphisms  (see Theorem 4.1).

In Section 5 we introduce Lagrangian submanifolds $\, \sigma (\Lambda ) , \,\, \sigma \in  {\mathcal B}_{2n}\,$, of the symplectic manifold $\, {\mathscr M}\,$. Some of these are then used to express the invariant spaces $\, J_C(L)\,$ of links as intersections of two Lagrangian submanifolds in  $\, {\mathscr M}\,$. 

In Section \ref{almCh} we study the first Chern class $\, c_1( {\mathscr M})\,$  of   $\, {\mathscr M}\,$. We consider first the evaluation of  $\, c_1( {\mathscr M})\,$ on the homology classes of maps $\, \gamma _{k, \epsilon} : S^2\rightarrow K_{2n}\,$ given by $\, \gamma _{k, \epsilon}(A)= (J, ... ,J, A, \epsilon A, J, ..., (-1)^n\epsilon J)\,$, where $\, A\in S^2\,$ and $\, J\,$ is a fixed base point of $\, S^2\,$ (here $\, A\,$ appears on the right hand side in the $\, k$-th factor). We formulate Theorem \ref{almChTh1} which says that the evaluation $\, \langle \,\, c_1( {\mathscr M})  \, \vert \,  [  \gamma _{k, \epsilon} ]      \,\, \rangle \,$ is equal to $\, 0\,$.  Theorem \ref{almChTh1} is proven in Appendix B.

We construct then some special mapping $\, f_n :S^2 \rightarrow K_{2n}  \,\, (n\ge 2)\,$ and formulate Theorem \ref{almChTh2} which says that the evaluation  $\, \langle \,\,  c_1( {\mathscr M})  \, \vert \,     [ f_n ]   \,\, \rangle\,$ is equal to $\, -2  \,$ for all $\, n\ge 2 \,$. We prove  Theorem \ref{almChTh2} in Appendix C.

Finally, we show that the subgroup $\, {\mathcal G}\,$ of $\, \pi _2(K_{2n})\,$ generated by the homotopy classes $\, [f_n]\,$ and $\, [  \gamma _{k, \epsilon} ], \,\, \epsilon = \pm 1,\,  k=1, ... , 2n-1,\,$ is of index at most 2.

In Section \ref{mono} we apply the results of Section \ref{almCh} to show that the symplectic manifold ${\mathscr M}$ is monotone (Theorem \ref{monoTh4}). It is perhaps worth mentioning that our proof of this fact is carried out completely in a finite-dimensional context.

In Section \ref{quandle} we point to the connection between the Lagrangian submanifolds constructed in Section 5 and the invariant spaces (representation varieties) $\, J_C(L)\,$ of links and indicate the subject of study which will be described in the second part of this paper.

In Section \ref{ex} we discuss several examples of links $L$ whose singly graded Khovanovs homology $Kh^*( L ) $ is isomorphic to the cohomology $H^*( J_C(L))$ of its invariant space. In particular, we show this to be true for all two-bridge knots, at least over $\, \Q\,$, and point to the arguments which make us strongly believe this to be true even over the integers $\Z\,$.

In Section \ref{counter} we determine the invariant space $\, J_C(L)\,$ for the knot $\, L=9_{42}\,$. We show that it is a disjoint union of a $2$-sphere $S^2$ and of seven copies of $\, \R P^3\,$ and point out the difference between its cohomology  and the collapsed Khovanov homology $\, Kh^*(9_{42})\,$  of the knot.

There are three appendices.
In Appendix A we analyse the singularities of the space $\, K_{2n}\,$. There are several components of the singular locus, all isomorphic to each other. Every component has a neighbourhood which is identified with a fibre bundle over the $2$-dimensional sphere $\, S^2\,$ with the fiber being the cone $\, C(S^{2n-3} \times S^{2n-3} )\,$ over a product of two $\, (2n-3)$-dimensional spheres (see Theorem A.7).

In Appendix B we prove Theorem \ref{almChTh1} and in Appendix C  Theorem \ref{almChTh2}.

\bigskip

The research described in this paper was initiated by the observation of similarities between the cohomology of the $SU(2)$-representation varieties $\, J_C(L)\,$ and the Khovanov homology $\, Kh^*(L)\,$  of certain links $\,L,\,$ the first instances of which we made in September 2006. 

 The recently published paper of  P.B.Kronheimer and T.S.Mrowka, \cite{KM}, has also the same observation as the starting point.  The approach and the setting of our paper is however different as we work with Lagrangian intersections while the authors of \cite{KM} construct a link homology based on the instanton Floer homology.  

\bigskip  

\bigskip

{\it Acknowledgments}: The authors are grateful to Oleg Viro and to Tobias Ekholm for many illuminating discussions on the subject of the paper and for constant encouragement in the course of this work. They also want to thank Charles Frohman for directing them to the paper \cite{L2} by X-S.Lin and the paper \cite{Kl1} by Klassen, and for some very fruitful discussions on representation varieties. The first author was partially supported by the Center for the Topology 
and Quantization of Moduli Spaces (CTQM) at the Århus University and by a grant from the Carlsberg Foundation. The second author gratefully acknowledges 
hospitality of CTQM where a part of the work on the paper has been conducted.

\bigskip

\section{Topology of the space of representations}\label{repr}

\bigskip

\qquad  2.1.\quad {\bf The group action and the differential of $\, p\,$.}

\bigskip

\bigskip

Let $\, G=SU(2)\,$ be the Lie group of $2\times2$ unitary matrices with determinant $1$. Elements of  $\, SU(2)\,$ are matrices $\, \left( \begin{array}{rr} a&b\\-\overline{b}&\overline{a} 
\end{array}\right), \,\, a, b \in \C , \,\, |a|^2+|b|^2=1\,$. The center of  $\, SU(2)\,$ is $\, Z=\{I, \, -I\}\,$, where $\, I\,$ is the identity matrix.

Let us identify the 2-dimensional  sphere $\, S^2\,$ with the subset $\, {\mathcal S}\, $ of  $\, SU(2)\,$ consisting  of matrices of trace $0$\, i.e. with the set  $\, {\mathcal S}\,$ of matrices of the form  $\, \left( \begin{array}{rr} it&z\\-\overline{z}&-it 
\end{array}\right), \,\, t\in\R, \, \, z \in \C , \,\, t^2+|z|^2=1\,$.The group  $\, SU(2)\,$ acts transitively on the set  $\, {\mathcal S}\,$ by conjugation and  $\, {\mathcal S}\,$ is a conjugacy class of  $\, SU(2)\,$.
Observe that for any matrix $\, A\in  {\mathcal S}\,$ one has $\, A^2=-I\,$.

For $\, A\in {\mathcal S}\,$ denote by $\, G_A\,$ the isotropy subgroup of $\, A\,$ under that action of $\, G\,$. The matrix  $\, J=  \left( 
\begin{array}{rr} i&0\\ 0&-i 
\end{array}\right)\,$ belongs to  $\, {\mathcal S}\,$. The isotropy subgroup $\, G_J\,$ of $\, J\,$ consists of matrices  $\, \left( 
\begin{array}{ll} e^{it}&0\\ 0&e^{-it} 
\end{array}\right), \,\, t\in\R, \,$ and is the unique maximal torus of $\, G\,$ containing $\, J\,$. The isotropy subgroups $\, G_A\,$ of other matrices $\, A\in {\mathcal S}\,$ are conjugates of  $\, G_J\,$ and are the unique maximal tori of $\, G\,$ containing $\, A\,$. Observe that for $\, A, B \in {\mathcal S}\,$ with $\, A\ne \pm B\,$ one has $\, G_A=G_{-A}\,$ and $\, G_A \cap G_B = \{I, \, -I\}\,$.

Let $\, P_{2n}\,$ be the product of $\, 2n\,$ copies of $\, S^2, \,\, P_{2n}=S^2\times ...\times S^2\,$.  Denote by $\, p_{2n}: P_{2n}\longrightarrow SU(2)\,$ the map defined by $\, p_{2n}(A_1,...,A_{2n})=A_1\cdot A_2\cdot ...\cdot A_{2n}\,$ for $\, A_j\in SU(2),\,\, \,  tr(A_j)=0, \,\, j=1, ...,2n\,$. We shall denote the mapping $\, p_{2n}\,$ by  $\, p\,$ whenever it will be clear from the context which $\, n\,$ is meant.

Let $\, K_{2n} = p^{-1}_{2n}(I), \,\,  \widetilde {K}_{2n} = p^{-1}_{2n}(-I)\,$.

The group $\, SU(2)\,$ acts diagonally (by conjugation) on $\, P_{2n}\,$ and on itself (again by conjugation) and the map $\, p_{2n}\,$ is equivariant. (Both actions are from the left.)
 Since $\, I,-I \in SU(2)\,$ are fixed points of the action on $\, SU(2)\,$, the spaces $\, K_{2n} = p^{-1}_{2n}(I)\,$ and $\,  \widetilde {K}_{2n} = p^{-1}_{2n}(-I)\,$ are equivariant subsets of $\, P_{2n}\,$ and inherit an action of $\, SU(2)\,$.

Let $\, \epsilon _1, ... ,  \epsilon _{2n} \in \{ 1, -1\}\,$. Denote by $\, \Delta _{ \epsilon _1, ... ,  \epsilon _{2n}}\,$ the closed equivariant submanifold of $\, P_{2n}\,$ defined by 
\begin{displaymath}
 \Delta _{ \epsilon _1, ... ,  \epsilon _{2n}} = \{ (\epsilon _1 A, ... , \epsilon _{2n} A ) \in P_{2n} \,\, | \,\, A\in   {\mathcal S} \, \, \} .
\end{displaymath}
Each submanifold  $\, \Delta _{ \epsilon _1, ... ,  \epsilon _{2n}}\,$ is diffeomorphic to $\, S^2\,$.
The isotropy subgroup of a point $\, (\epsilon _1 A, ... , \epsilon _{2n} A )\,$ in $\, \Delta _{ \epsilon _1, ... ,  \epsilon _{2n}} \,$ is equal to the isotropy subgroup $\, G_A\,$ of the matrix $\, A\,$ in  $\, SU(2)\,$. Two subsets  $\, \Delta _{ \epsilon _1, ... ,  \epsilon _{2n}}\,$ and  $\, \Delta _{ \delta _1, ... ,  \delta _{2n}}\,$ are equal if $\, \delta _i = \epsilon \epsilon _i\,$ for some $\, \epsilon = \pm 1\,$ and all $\, i\,$, and disjoint otherwise.

Every subset  $\, \Delta _{ \epsilon _1, ... ,  \epsilon _{2n}} \,$ is contained either in  $\, K_{2n}\,$ or in $\, \widetilde {K}_{2n}\,$ depending on whether $\,(-1)^n  (\prod \epsilon _i)\,$ is equal to $\, 1\,$ or to $\, -1\,$. Let us denote by $\, \Sigma\,$ the union of these subsets  $\, \Delta _{ \epsilon _1, ... ,  \epsilon _{2n}} \,$ which are contained in $\, K_{2n}\,$ and by $\, \widetilde{\Sigma}\,$ the union of those  $\, \Delta _{ \epsilon _1, ... ,  \epsilon _{2n}} \,$ which are contained in  $\, \widetilde {K}_{2n}\,$.

Observe that the multiplication by $\, -1\,$ on, say, the first coordinate of $\, P_{2n}\,$ gives a diffeomorphism of $\, P_{2n}\,$ onto itself which maps  $\, K_{2n}\,$ onto  $\, \widetilde {K}_{2n}\,$ and vice-versa, and maps  $\, \Sigma\,$ onto  $\, \widetilde{\Sigma}\,$.

\smallskip

Let $\, W_{2n}= 
P_{2n}- (\Sigma \cup \widetilde{\Sigma})
\,$. It is an open equivariant subset of $\, P_{2n}\,$.

\bigskip

\begin{lma}
For every point $\, (A_1, ... , A_{2n}\,) \in W_{2n},\, \, n\ge 1\,$, the isotropy subgroup of the action of $\, SU(2)\,$ is equal to $\, Z=\{ I, -I\, \}\,$. 
\end{lma}

\begin{proof}
The center $\, Z\,$ is contained in isotropy subgroups of all points of $\, P_{2n}\,$. Since $\,  (A_1, ... , A_{2n}\,) \in W_{2n}\,$ there is a pair of indices $\, 1\le i, j \le 2n\,$ such that $\, A_j\ne \pm A_i\,$. Let $\, G_{A_i}\,$ and  $\, G_{A_j}\,$ be the isotropy subgroups of the matrices $\, A_i\,$ and $\, A_j\,$ in $\, SU(2)\,$. The isotropy subgroup of  $\,  (A_1, ... , A_{2n}\,)\,$ is contained in   $\, G_{A_i}\cap G_{A_j}= \{ I, \, -I\,\}\,$ and thus equal to it.
\end{proof}

Thus  the group $\, SO(3)=SU(2)/\{ \pm I\}\,$ acts freely on $\, W_{2n}\,$.

\bigskip

\begin{lma}
The mapping $\, p_{2n}: P_{2n}\longrightarrow SU(2)\,$ is a submersion at every point of $\, W_{2n}\,$.
\end{lma}

\begin{proof}
It is enough to prove the lemma in the case when $\, n=1\,$. Indeed, if $\, (A_1, ... , A_{2n})\in W_{2n}\,$ then there is an index $\, i\,$ such that $\, A_i\ne \pm A_{i+1}\,$. Consider the mapping $\, j:W_2 \longrightarrow W_{2n}\,$ given by $\, j((B_1,B_2))= (A_1, ... , A_{i-1}, B_1, B_2, A_{i+2}, ... ,A_{2n})\,$ and the diffeomorphism $\, f: SU(2) \longrightarrow SU(2)\,$ given by $\, f(X)=A_1\cdot ... \cdot A_{i-1}\cdot X \cdot A_{i+2}\cdot ... \cdot A_{2n}\,$. The diagram 
\begin{displaymath}
\begin{CD}
W_2  @>j>>  W_{2n}\\
@VV{p_2}V   @VV{p_{2n}}V\\
SU(2)  @>f>> SU(2)
\end{CD}
\end{displaymath}
commutes. Hence, if $\, p_2\,$ is a submersion at the point $\, (A_i, \, A_{i+1} )\,$ in $\, W_2\,$ then $\, p_{2n}\,$ is a submersion at the point $\,j((A_i, \,A_{i+1}))= (A_1, ... , A_{i-1}, A_i, A_{i+1}, A_{i+2}, ... ,A_{2n})\,$ in $\, W_{2n}\,$.

We assume now that $\, n=1\,$. Let $\, (A_1,\,  A_2 )\,$ be a point in $\, W_2\,$.  We have $\, A_1\ne \pm A_2\,$. Conjugating by en element of $\, SU(2)\,$ if necessary, we can assume that $\, A_1=J= \left( 
\begin{array}{rr} i&0\\ 0&-i 
\end{array}\right) \,$. Then $\, A_2=  \left( \begin{array}{rr} is&v\\-\overline{v}&-is 
\end{array}\right)\,$ with $\, \, s\in\R, \, \, v \in \C , \,\, s^2+|v|^2=1\,$ and $\, v\ne 0\,$.

For $\, t\in\R , \,\, -1<t<1, \,$ define matrices $\, A_1(t)=  \left( 
\begin{array}{cc} i\sqrt{1-t^2}&-t\\ t&-i\sqrt{1-t^2} 
\end{array}\right)\,$,
\,\,$ \, \widetilde{A}_1(t)= \left( 
\begin{array}{cc} i\sqrt{1-t^2}&-it\\ -it&-i\sqrt{1-t^2} 
\end{array}\right) \,$\,\,\, and \,\,\, $\, A_2(t)= \left( \begin{array}{cc} is&ve^{-it}\\-\overline{v}e^{it} &-is 
\end{array}\right)\,$. 

Observe that $\, A_1(t),\, \widetilde{A}_1(t), \, A_2(t)\in {\mathcal S}\,$ for $\, -1<t<1\,$ and $\, A_1(0)=\widetilde{A}_1(0) = A_1\,$ while $ \,  A_2(0)=A_2\,$. Consider the smooth curves $\, h_1(t)= (A_1(t),A_2),\, \, h_2(t) = (\widetilde{A}_1(t),A_2)\,$ and $\, h_3(t) = (A_1, A_2(t))\,$ in $\, P_2\,$. We have $\, h_i(0)=(A_1, A_2)\,$ for $\, 1\le i \le 3\,$ and 

\begin{equation*}
\begin{split}
\frac d{dt} (p_2(h_1(t)))| _{t=0}& =  \frac d{dt} (A_1(t)\cdot A_2 )|_{t=0}=\\ 
&=\left( 
\begin{array}{cc} -i\frac t{\sqrt{1-t^2}}&-1\\ 1& i \frac t{\sqrt{1-t^2}}
\end{array}\right)_{t=0}\cdot  \left( \begin{array}{rr} is&v\\-\overline{v}&-is  \end{array}\right) =\\
&= \left( \begin{array}{rr} 0&-1\\ 1&0  \end{array}\right) \cdot  \left( \begin{array}{rr} is&v\\-\overline{v}&-is 
\end{array}\right) =\left( \begin{array}{rr}\overline{v}& is\\ is&v 
\end{array}\right),\\
\frac d{dt} (p_2(h_2(t)))| _{t=0}& =  \frac d{dt} (\widetilde{A}_1(t)\cdot A_2 )|_{t=0}=\\ 
&=\left( 
\begin{array}{cc} -i\frac t{\sqrt{1-t^2}}&-i\\ -i& i \frac t{\sqrt{1-t^2}}
\end{array}\right)_{t=0}\cdot  \left( \begin{array}{rr} is&v\\-\overline{v}&-is  \end{array}\right) =\\
&= \left( \begin{array}{rr} 0&-i\\ -i&0  \end{array}\right) \cdot  \left( \begin{array}{rr} is&v\\-\overline{v}&-is 
\end{array}\right) =\left( \begin{array}{rr}i\overline{v}& -s\\ s&-iv 
\end{array}\right)
\end{split}
\end{equation*}
and
\begin{equation*}
\begin{split}
\frac d{dt} (p_2(h_3(t)))| _{t=0}& =  \frac d{dt} (A_1\cdot A_2(t) )|_{t=0}=\\
&= \left( \begin{array}{rr} i&0\\ 0&-i  \end{array}\right) \cdot  \left( \begin{array}{rr} 0&-iv\\-i \, \overline{v}&0 
\end{array}\right) =\left( \begin{array}{rr} 0&v  \\ -\overline{v}& 0 
\end{array}\right). 
\end{split}
\end{equation*}
\smallskip

\noindent
The three matrices $\, \left( \begin{array}{rr}\overline{v}& is\\ is&v 
\end{array}\right),\,\, \left( \begin{array}{rr}i \, \overline{v}& -s\\ s&-iv 
\end{array}\right), \,\, \left( \begin{array}{rr} 0&v  \\ -\overline{v}& 0 
\end{array}\right)\,$ with $\, s^2+|v|^2=1,\,\,  v\ne 0\,$ are linearly independent over real numbers. Indeed, if
\begin{displaymath}
\alpha\cdot \left( \begin{array}{rr}\overline{v}& is\\ is&v 
\end{array}\right) + \beta \cdot\left( \begin{array}{rr}i\overline{v}& -s\\ s&-iv 
\end{array}\right) + \gamma \cdot \left( \begin{array}{rr} 0&v  \\ -\overline{v}& 0 
\end{array}\right) =0 \,\, , \quad  \alpha, \beta, \gamma \in \R ,
\end{displaymath}
then $\, ( \alpha + i\beta )\overline{v} = 0\,$ and $\, (i\alpha - \beta )s + \gamma v =0\,$. Since $\, v\ne 0\,$ we get $\, \alpha =\beta =\gamma =0\,$. 

\smallskip

Thus the rank of the differential $\, d(p_2)_{(A_1,A_2)}\,$ is equal to 3 and the mapping $\, p_2\,$ is a submersion at the point $\, (A_1,A_2)\in W_2\,$.
\end{proof}

\medskip

Let $\, \widetilde{B}=SU(2)-\{\pm  I \}\,$ and $\, \widetilde{E}=p^{-1}_{2n}(B) = P_{2n}-(K_{2n}\cup \widetilde{K}_{2n})\,$. Denote  by $\,\widetilde{p}\,$  restriction of the mapping $\, p_{2n}\,$ to $\,\widetilde{E}\,$. Then $\, \widetilde{p}: \widetilde{E} \longrightarrow \widetilde{B}\,$ is a proper map. Since 
$\, \widetilde{E} \subset W_{2n}\,$, Lemma 2.2 implies that $\,\widetilde{p} : \widetilde{E} \longrightarrow \widetilde{B}\,$ is a submersion. Hence, we get

\smallskip 
\begin{cor}
The mapping $\,\widetilde{p} : \widetilde{E} \longrightarrow \widetilde{B}\,$ is a locally trivial fibration. 
\end{cor} 

\smallskip

The subsets $\, K_{2n}\,$ and $\,  \widetilde{K}_{2n}\,$ are real algebraic subvarieties of $\, P_{2n}\,$. Let $\, E = P_{2n}-\widetilde{K}_{2n}\,$. We have $\, K_{2n} \subset E\,$. 
\smallskip

\begin{cor}
 The subspace $\, K_{2n}\,$ is a strong deformation retract of 
$\, E\,$.
\end{cor} 

\begin{proof}
Let $\, B = SU(2)-\{-I\}\,$. We identify $\, B\,$ with a 3-dimensional disc with center at $\, I\,$. Since  $\, K_{2n}\,$ is a  real algebraic subvariety of $\, P_{2n}\,$ there exists an open neighbourhood $\, U\,$ of  $\, K_{2n}\,$ in  $\, P_{2n}\,$ such that  $\, K_{2n}\,$ is  a strong deformation retract of $\, U\,$, (see \cite{L1}). There exists a closed concentric subdisc $\, D\,$ of $\, B\,$ with center at $\, I\,$ such that $\, p_{2n}^{-1} (D) \subset U\,$. Since $\, \widetilde{p} : \widetilde{E} \longrightarrow \widetilde{B}\,$ is a locally trivial fibration (Corollary 2.3) the radial strong deformation retraction of  $\, B\,$ onto $\, D\,$ can be lifted to $\, E\,$. That deformation followed by the deformation retraction of $\, U\,$ onto  $\, K_{2n}\,$ gives a strong deformation retraction of  $\,  E\,$ onto  $\, K_{2n}\,$. 
\end{proof}

\medskip

\begin{rmk}
It follows from Lemma 2.2 that the singularities of the real algebraic varieties  $\, K_{2n}\,$ and $\,  \widetilde{K}_{2n}\,$        (if any) are contained in the subsets   $\,\Sigma \,$ respectively $\,\widetilde{\Sigma} \,$  .
\end{rmk}

\smallskip

\begin{cor}
The spaces $\,  K_{2n}\,$  and $\,  \widetilde{K}_{2n}\,$  are path-connected.
\end{cor}

\begin{proof}
Let $\, x\,$ and $\,  y\,$ be two points in  $\,  K_{2n}\,$. Choose a smooth path $\, \sigma \,$  in $\, P_{2n}\,$ from   $\, x\,$ to $\,  y\,$. We can choose $\, \sigma\,$ to be transversal to, hence disjoint from, the 2-dimensional submanifold $\,\widetilde{\Sigma}\,$. Once that has been done, we can farther deform $\, \sigma\,$ to be transversal to, hence disjoint from, the submanifold  $\,  \widetilde{K}_{2n} -\widetilde{\Sigma} \,$ which is of codimension 3. As a result the path  $\, \sigma\,$ lies in $\, E\,$. Applying the strong deformation retraction of $\, E\,$ onto  $\,  K_{2n}\,$ to  $\, \sigma\,$ we get a path  $\, \sigma '\,$ in  $\,  K_{2n}\,$ from  $\, x\,$ to $\,  y\,$. Hence  $\,  K_{2n}\,$ is path-connected. Since $\, \widetilde{K}_{2n}\,$ is homeomorphic to  $\,  K_{2n}\,$, it is also path-connected.
\end{proof}

\bigskip

\qquad  2.2.\quad {\bf Singularities of  $\, K_{2n}\,$ and  $\,  \widetilde{K}_{2n}\,$  .}

\bigskip

The main aim of this subsection is to prove that the space  $\,  \widetilde{K}_{2n} - \widetilde{\Sigma}\,$ is path-connected. This claim will be a consequence of a much stromger result which we formulate as Proposition 2.7 and which in turn will follow from a detailed analysis of the singularities of  $\, K_{2n}\,$ and  $\,  \widetilde{K}_{2n}\,$.  
\bigskip

It follows from Lemma 2.2 that the singularities of the real algebraic varieties  $\, K_{2n}\,$ and $\,  \widetilde{K}_{2n}\,$        (if any) are contained in the subsets   $\,\Sigma \,$ respectively $\,\widetilde{\Sigma} \,$  . We shall now study the nature of these singularities.

It is enough to study the structure of  $\, K_{2n}\,$ in a neighbourhood of one subset  $\, \Delta _{ \epsilon _1, ... ,  \epsilon _{2n}}\,$ with $\,(-1)^n  (\prod \epsilon _i) = 1\,$.
Indeed, let $\,f_j:P_{2n}\longrightarrow P_{2n}\,$ be the diffeomorphism $\, f_j(A_1, ... ,A_j, ... ,A_{2n})=(A_1, ... ,-A_j, ... ,A_{2n})\,$ of multiplication of the $\, j$-th coordinate of $\, P_{2n}\,$ by $\, -I\,$. Then $\, f_j\,$ exchanges  $\, K_{2n}\,$ with $\,  \widetilde{K}_{2n}\,$ and exchanges different subsets  $\, \Delta _{ \epsilon _1, ... ,  \epsilon _{2n}}\,$. Actually, the mappings $\, f_j\,$ define a smooth action of the group $\, (\Z /2\Z)^{2n}\,$ on $\, P_{2n}\,$ which acts transitively on the family of all subsets  $\, \Delta _{ \epsilon _1, ... ,  \epsilon _{2n}}\,$. It follows that the structure of  $\, K_{2n}\,$ in a neighbourhood of every subset  $\, \Delta _{ \epsilon _1, ... ,  \epsilon _{2n}}\,$ contained in $\, K_{2n}\,$ is the same.

Consider a sequence $\,  (\epsilon _1, ... ,  \epsilon _{2n})\in \{ \pm 1\}^{2n}\,$ with $\, (-1)\prod\limits_{j=1}^{2n}\epsilon _j = 1\,$. Then $\, \Delta =  \Delta _{ \epsilon _1, ... ,  \epsilon _{2n}} \subset  K_{2n}\,$ and $\, \Delta\,$ is diffeomorphic to $\, {\mathcal S} = S^2\,$. 

\begin{prop}
There is an open neighbourhood $\, U\,$ of $\, \Delta \,$ in  $\, K_{2n}\,$ and a continuous mapping $\, \xi :\overline{U} \rightarrow \Delta\,$  whose restriction to $\, \Delta \,$ is the identity and which is a locally trivial fibration with a fibre homeomorphic to the cone $\, C(S^{2n-3}\times S^{2n-3})\,$ over the product of two $\, (2n-3)$-dimensional spheres.
\end{prop}

\smallskip

Proof of Proposition 2.7 is deferred to the Appendix A. (See Theorem A.7.)

\smallskip

By the remarks above Proposition 2.7 is also valid for $\,\Delta =  \Delta _{ \epsilon _1, ... ,  \epsilon _{2n}}\,$ with $\,(-1)^n  (\prod \epsilon _i) = -1\,$ and with  $\, K_{2n}\,$ replaced by  $\,  \widetilde{K}_{2n}\,$.

\smallskip

\begin{cor} If $\, n\ge 2\,$ then
the spaces $\, K_{2n}-\Sigma \,$ and  $\, \widetilde{K}_{2n}-\widetilde{\Sigma} \,$ are path-connected.
\end{cor}

\begin{proof}
It is enough to prove the statement for  $\, K_{2n}-\Sigma \,$. Let $\, x\,$ and $\, y\,$ be points in $\, K_{2n}-\Sigma \,$. By Corollary 2.6 there is a path $\, \sigma (t), \, 0\le t \le 1, \,$ in  $\, K_{2n}\,$ from  $\, x\,$ to $\, y\,$. Let $\, k\,$ be the number of components $\,  \Delta _{ \epsilon _1, ... ,  \epsilon _{2n}}\,$ which the path  $\, \sigma \,$ meets. The proof is by induction on $\, k\,$. Suppose $\, \Delta \,$ is a component of $\, \Sigma\,$ met by  $\, \sigma \,$. Choose a neighbourhood $\ U\,$ of $\, \Delta \,$ in $\, K_{2n}\,$  as in Proposition 2.7. We can assume that $\,  \overline{U}\,$ does not meet other components of $\, \Sigma \,$.  Let $\, t_1\,$ be the smallest value of $\, t\,$ such that $\, \sigma (t)\in \overline{U}\,$ and let $\, t_2\,$ be the largest one. Since $\, n\ge 2\,$, the base $\,S^{2n-3}\times S^{2n-3}\,$ of the cone $\, C(S^{2n-3}\times S^{2n-3})\,$ is path-connected and, therefore, there exists a path $\, \tau \,$ in $\, \overline{U} - \Delta \,$ from $\, \sigma (t_1)\,$ to  $\, \sigma (t_2)\,$. Then the path $\,\sigma (t)\,$ from $\,\sigma (0)\,$ to   $\,\sigma (t_1)\,$ followed by the path  $\, \tau \,$ and then followed by path  $\,\sigma (t)\,$ from $\,\sigma (t_2)\,$ to   $\,\sigma (1)\,$ gives a path in  $\, K_{2n}\,$ from  $\, x\,$ to $\, y\,$ with a smaller number $\, k\,$. That proves Corollary.
\end{proof}

\bigskip

\qquad  2.2.\quad {\bf The first and the second homotopy groups of $\, K_{2n}\,$.} 

\bigskip

The first two homotopy groups of the space $\,P_{2n}\,$ are $\, \pi _1( P_{2n} ) = 0\,$ and  $\, \pi _2( P_{2n} ) \cong  \bigoplus\limits_{i=1}^{2n} \Z _i\,$, where all $\,\Z _i \cong \Z\,$. Let $\, j:  K_{2n}\longrightarrow P_{2n}\,$ be the inclusion.  The first two homotopy groups of $\, K_{2n}\,$ are given by 

\bigskip

\begin{prop}
If $\, n\ge 2\,$, then

\smallskip

\qquad \, (i) \,\, $\, \pi _1( K_{2n} ) = 0\,$.

\smallskip

\qquad (ii) \, \,  There is a short split exact sequence 
\begin{displaymath}
\begin{CD}
0 @>>>\Z @>>>  \pi _2( K_{2n} ) @>{j_*}>>  \pi _2( P_{2n} ) @>>> 0 .
\end{CD}
\end{displaymath}
\end{prop}

\medskip

\begin{proof}
Recall that $\, K_{2n}\,$ is a strong deformation retract of $\, E= P_{2n}- \widetilde{K}_{2n}\,$ (see Corollary 2.4). Therefore it is enough to prove the Proposition for the space $\, E\,$ instead of $\, K_{2n}\,$. Let $\, \widetilde{j} : E \longrightarrow P_{2n}\,$ be the inclusion.

(i) \,\, Let $\, \gamma :S^1 \longrightarrow E\,$ be an arbitrary loop in $\, E\,$. Then $\,  \widetilde{j} \circ \gamma \,$ is a loop in $\,  P_{2n}\,$. Since the space $\, P_{2n}\,$ is 1-connected there is an extension of  $\,  \widetilde{j} \circ \gamma \,$ to a map $\, \Gamma :D^2 \longrightarrow  P_{2n}\,$. The submanifold  $\, \widetilde{\Sigma}\,$   is of codimension $\, 4n-2\ge 6\,$ in  $\, P_{2n}\,$. Therefore the mapping $\, \Gamma \,$ can be chosen to be disjoint from $\,  \widetilde{\Sigma}\,$. One gets $\, \Gamma :D^2 \longrightarrow  P_{2n}-  \widetilde{\Sigma} \,$.  It follows from Lemma 2.2 that $\,  \widetilde{K}_{2n} - \widetilde{\Sigma}\,$ is a submanifold of $\, P_{2n} -  \widetilde{\Sigma}\,$ of codimension 3. Therefore  $\, \Gamma \,$ can deformed rel $\partial D^2\,$ even further and chosen to be disjoint from  $\,  \widetilde{K}_{2n} - \widetilde{\Sigma}\,$ and hence from  $\,  \widetilde{K}_{2n}\,$. In that way we get a map  $\, \Gamma :D^2 \longrightarrow  P_{2n}-  \widetilde{K}_{2n} = E \,$ extending the loop  $\, \gamma\,$. Hence the homotopy class of  $\, \gamma\,$ in $\, \pi _1(E )\,$ is trivial, showing that $\, \pi _1(E ) = 0\,$.

\smallskip

(ii) Let $\, \zeta : S^2 \longrightarrow P_{2n}\,$ be arbitrary. The same dimensional argument as in the proof of part (i) shows that $\, \zeta \,$ can be deformed to a map $\, \zeta ': S^2 \longrightarrow E\,$. This proves that $\, \widetilde{j}_*:\pi _2( E )\longrightarrow  \pi _2( P_{2n} )\,$ is surjective.

Let $\, H \subset \pi _2( E )\,$ be the kernel\, of $\, \widetilde{j}_*:\pi _2( E )\longrightarrow  \pi _2( P_{2n} )\,$. We shall now construct an isomorphism $\, \chi: H \longrightarrow \Z\,$.

We choose an orientation of $\, P_{2n}\,$ and of $\, SU(2)\,$. That gives an orientation of $\, P_{2n} -  \widetilde{\Sigma}\,$. By Lemma 2.2 the subset  $\,  \widetilde{K}_{2n} - \widetilde{\Sigma} = p_{2n}^{-1}(-I) - \widetilde{\Sigma}\,$ is a submanifold of $\, P_{2n} -  \widetilde{\Sigma}\,$ and we give it an orientation such that locally at its points the orientation of $\, P_{2n}\,$ is equal to that of  $\,  \widetilde{K}_{2n} - \widetilde{\Sigma}\,$ followed by that of $\, SU(2)\,$. 

 Let $\, \varphi :S^2 \longrightarrow E\,$ be a smooth map representing an element of $\, H\,$. Since $\, \widetilde{j}_* ([\varphi])=0\,$, there exists a smooth extension $\, \Phi : D^3 \longrightarrow P_{2n}\,$ of $\, \varphi\,$. Again, since the codimension of the submanifold $\, \widetilde{\Sigma}\,$ of $\, P_{2n}\,$ is $\, 4n-2\ge 6\,$ we can choose  $\, \Phi\,$ whose image is disjoint from  $\, \widetilde{\Sigma}\,$. Moreover, we can choose  $\, \Phi\,$ to be transversal to 
the oriented submanifold   $\,  \widetilde{K}_{2n} - \widetilde{\Sigma}\,$. We define $\, \chi (\varphi ) \in \Z\,$ to be the oriented intersection number of  $\,  \widetilde{K}_{2n} - \widetilde{\Sigma}\,$ with  $\, \Phi : D^3 \longrightarrow P_{2n} -\widetilde{\Sigma} \,$ in $\, P_{2n} -\widetilde{\Sigma} \,$. This number is independent of the choice of $\, \varphi\,$ in its homotopy class and of the choice of its extension $\, \Phi: D^3 \longrightarrow P_{2n} -\widetilde{\Sigma}  \,$. 
Indeed, let  $\, \varphi ':S^2 \longrightarrow E\,$ be homotopic to $\, \varphi \,$ with $\, F:S^2\times I\longrightarrow E\,$  being a homotopy between them. Choose an extension $\, \Phi ' :  D^3 \longrightarrow P_{2n}\,$ of $\, \varphi '\,$ with properties as those of  $\, \Phi \,$. Identify the oriented 3-dimensional sphere $\, S^3\,$ with the union  $\, D^3\,  \cup \, S^2\times I \,\cup\, (-D^3)\,$, where $\, -D^3\,$ denotes the standard 3-dimensional disc with the reversed orientation and where one component of the oriented boundary of $\, S^2 \times I\,$ has been identified with the boundary of $\, D^3\,$ while the other component with the boundary of $\, -D^3\,$. The mappings $\, \Phi\,$ on $\, D^3\,$,$\, F\,$ on $\, S^2\times I\,$ and $\, \Phi '\,$ on $\, -D^3\,$  fit together and define a mapping $\, G:S^3\longrightarrow P_{2n}\,$ which is transversal to $\,\widetilde{K}_{2n}\,$ and does not meet $\,  \widetilde{\Sigma}\,$.\, (Observe that the restriction of $\, G\,$ to $\, S^2\times I\,$ is given by $\, F\,$ whose image lies in $\, E\,$ and hence does not intersect  $\,\widetilde{K}_{2n}\,$).
 The oriented intersection number $\,m\,$ of $\,\widetilde{K}_{2n} -  \widetilde{\Sigma}\,$ with $\, G\,$ is equal to $\, \chi (\varphi ) -  \chi (\varphi ')\,$. Here $\, \chi (\varphi )\,$ is defined by the choice of the extension $\, \Phi\,$ and $\, \chi (\varphi ' )\,$  by the choice of the extension $\, \Phi '\,$. 

On the other hand, the oriented real subvariety  $\,\widetilde{K}_{2n}\,$ of $\, P_{2n}\,$ with singularities of codimension at least 6 defines a homology class in $\, H_{2n-3}(P_{2n} ; \Z )\,$ and the mapping $\, G\,$ defines a homology class in $\, H_3(P_{2n} ; \Z )\,$. Both these classes are equal to $\, 0\,$ since $\, P_{2n}\,$ has only trivial homology in odd degrees. At the same time, the intersection number $\, m\,$ is equal to the homological intersection number of those two homology classes and hence equal to $\, 0\,$.  Therefore  $\, \chi (\varphi ) -  \chi (\varphi ') = m = 0\,$ and, finally,  $\, \chi (\varphi ) =  \chi (\varphi ')\,$.

Thus, we get a well defined integer $\, \chi ([ \varphi ] ) \in \Z\,$ depending only on the homotopy class of $\,\varphi:  S^2\longrightarrow E\,$ for $\,[ \varphi ] \in H = ker (\widetilde{j}_*)\,$. It is obvoius from the construction that in that way we get a group homomorphism $\, \chi :H \longrightarrow \Z\,$. 

Let us choose a point $\, x\,$ in   $\,  \widetilde{K}_{2n} - \widetilde{\Sigma}\, $\ and a normal space $\, N_x\,$ to  $\,\widetilde{K}_{2n}\,$ in $\, P_{2n}\,$ at $\, x\,$ (a fiber of a regular neighbourhood) . We have $\, \dim N_x = 3\,$.  A sphere $\, S_x\,$ in $\, N_x\,$ with centrum at  $\, x\,$ and a small radius gives us a mapping $\,\varphi:  S^2\longrightarrow E\,$ with   $\, \chi ([ \varphi ] ) = \pm 1\,$. Thus  $\, \chi \,$ is an epimorphism.

To show that  $\, \chi :H \longrightarrow \Z\,$ is a monomorphism let us
 assume that  $\,\varphi:  S^2\longrightarrow E\,$ is a map  such that 
$\, [\varphi ]\in H\,$ and $\, \chi ([ \varphi ] ) = 0\,$. It follows that  
$\, \varphi \,$ extends to a map   $\,\Phi:  D^3\longrightarrow P_{2n}\,$ whose 
image does not meet $\, \widetilde{\Sigma}\,$ and intersects 
$\, \widetilde{K}_{2n}\,$ transversally in an even number of points, half of 
them with index $\, +1\,$ and half with index $\, -1\,$. Accoring to Corollary 2.8,  the manifold 
$\, \widetilde{K}_{2n} - \widetilde{\Sigma}  \,$ is path-connected. Thus,
 applying the same 
arguments on cancellation of intersection points in pairs as in the proof of 
the theorem of Whitney (see \cite{M1}, Theorem 6.6.) one proves that the map  
$\,\Phi \,$ can be deformed rel. $\partial D^3\,$ to a map whose image does not 
intersect   $\, \widetilde{K}_{2n}\,$  at all. This new deformed map is of the 
form  $\,\Psi:  D^3\longrightarrow E\,$ with 
$\,\Psi\,  | \, \partial D^3 = \varphi \,$. We conclude that the homotopy 
class of 
$\, [\varphi ]\in \pi _2 (E)\,$ is trivial. Hence 
$\, \chi:H \longrightarrow \Z \,$ is a monomorphism and therefore 
an isomorphism.

That proves the second part of Proposition 2.9 since $\, H\,$ was the kernel of $\, \widetilde{j}_*:\pi _2( E )\longrightarrow  \pi _2( P_{2n} )\,$. The short exact sequence is split because $\, \pi _2( P_{2n} )\,$ is a free $\, \Z$-module.
\end{proof}

\medskip

Let $\, h_n : P_{2n} \rightarrow  P_{2n+2}\,$ be the map defined by 
\begin{equation}\label{2.1}
h_n(A_1, ... , A_{2n}) =(A_1, ... , A_{2n}, -J,\, J\, ) \, .
\end{equation}
The diagram 
\begin{equation}\label{2.2}
\xymatrix{P_{2n}\ar[rd]_{p_{2n}}\ar[rr]^{h_n}&& P_{2n+2}\ar[ld]^{p_{2n+2}}\\ 
& SU(2) }
\end{equation}
commutes.
Hence $\, h_n\,$ maps $\, K_{2n}=p_{2n}^{-1}(I)\,$ into $\, K_{2n+2}=p_{2n+2}^{-1}(I)\,$. We denote the restriction of $\, h_n\,$ to $\,  K_{2n}\,$ by the same symbol, $\, h_n:  K_{2n}\rightarrow  K_{2n+2}\,$.

Consider the homomorphism induced by $\, h_n\,$  on the second homotopy groups 
\begin{displaymath}
 h_{n *}: \pi _2( K_{2n})\rightarrow \pi _2 (K_{2n+2})\, .
\end{displaymath}
(We do not specify the base points as it is inessential. Whenever necessary we can take $\, x_0=(-J, \, J, ... , -J, \, J)\in K_{2n}\,$ as a base point.)

\begin{prop} If $\,n\ge 2\,$ then the homomorphism $\, h_{n *}: \pi _2( K_{2n})\rightarrow \pi _2 (K_{2n+2})\,$ is a monomorphism mapping the kernel of $\, j_*:\pi _2(K_{2n})\rightarrow \pi _2(P_{2n})\,$ isomorphically onto the kernel of  $\, j_*:\pi _2(K_{2n+2})\rightarrow \pi _2(P_{2n+2})\,$.
\end{prop}

\begin{proof} The diagram
\begin{equation}\label{2.3}
\xymatrix{
0\ar[r]& \Z \ar[r]& \pi _2(K_{2n})\ar[r]^{j_*}\ar[d]_{ h_{n *}}&  \pi _2(P_{2n})\ar[r]\ar[d]_{h_{n *}} &0\\
0\ar[r]& \Z \ar[r]& \pi _2(K_{2n+2})\ar[r]^{j_*}&  \pi _2(P_{2n+2})\ar[r]&0 }
\end{equation}
commutes. Hence the homomorphism  $\, h_{n *}: \pi _2( K_{2n})\rightarrow \pi _2 (K_{2n+2})\,$ maps the kernel $\,\Z\,$ of $\, j_*:\pi _2(K_{2n})\rightarrow \pi _2(P_{2n})\,$ into the kernel $\, \Z\,$ of   $\, j_*:\pi _2(K_{2n+2})\rightarrow \pi _2(P_{2n+2})\,$.

Let $\, \gamma :D^3\rightarrow P_{2n}\,$ be an embedding of a small $3$-dimensional disc which intersects $\, \widetilde{K}_{2n}=p_{2n}^{-1}(-I)\,$ transversally in one non-singular point with the intersection number $1$. Then the homotopy class of the restriction of $\, \gamma\,$ to the boundary $\, \partial D^3\,$ of the disc, $\,\gamma \,\vert _{\partial D^3}:\partial D^3 \rightarrow (P_{2n}-  \widetilde{K}_{2n}) \simeq K_{2n}\,$, gives a generator $\,\xi\,$ of the kernel of  $\, j_*:\pi _2(K_{2n})\rightarrow \pi _2(P_{2n})\,$. 

Since the diagram (2.2) commutes, it follows that the composition $\, h_n\circ \gamma\,$ is an embedding of $\, D^3\,$ into $\, P_{2n+2}\,$ which again intersects $\,  \widetilde{K}_{2n+2}=p_{2n+2}^{-1}(-I)\,$ transversally in one non-singular point with the intersection number $1$. The homotopy class of its restriction $\, (h_n\circ \gamma )\, \vert _{\partial D^3}\,$  to the boundary of the disc gives a generator $\, \xi '\,$ of the kernel of $\, j_*:\pi _2(K_{2n+2})= \pi _2(P_{2n+2}-  \widetilde{K}_{2n+2}) \rightarrow \pi _2(P_{2n+2})\,$.

Thus $\, h_{n *}: \pi _2( K_{2n})\rightarrow \pi _2 (K_{2n+2})\,$ maps the generator $\, \xi\,$ of the kernel of $\, j_*:\pi _2(K_{2n})\rightarrow \pi _2(P_{2n})\,$ to the generator $\, \xi '\,$ of the kernel of  $\, j_*:\pi _2(K_{2n+2})\rightarrow \pi _2(P_{2n+2})\,$ and, hence, maps the first of the kernels isomorphically onto the second one. (Actually, the restriction of  $\, h_{n *}\,$ to the kernels is the identity homomorphism, when both are identified with $\, \Z\,$ in the way indicated in the proof of Lemma 2.9). As the right-hand vertical 
arrow in the diagram (2.3) is obviously a monomorphism, so is the homomorphism    $\, h_{n *}: \pi _2( K_{2n})\rightarrow \pi _2 (K_{2n+2})\,$.
\end{proof}

\bigskip

\section{The symplectic structure }\label{sympl}

\bigskip

Let $\, G\,$ be a compact Lie group and $\, \mathfrak{g}\,$ its Lie algebra.  We choose an invariant $\, \R$-valued positive definite inner product $\, \bullet \,$ on  $\, \mathfrak{g}\,$. 

Let $\,\mathcal C = ( C_1, ... , C_m) \,$ be a sequence of $\, m\,$ conjugacy classes in $\, G\,$ (not necessarily distinct).

Denote by $\, F\,$ a free group on $\, m\,$ generators $\, z_1, ... z_m\,$. 

Let $\,\text{Hom} ( F, \, G )_{\mathcal C}\,$ be the set of all group homomorphisms $\,\varphi :F \rightarrow G\,$ such that $\, \varphi (z_j)\in C_j, \,\, 1\le j \le m\,$. We identify $\,\text{Hom} ( F, \, G )_{\mathcal C}\,$ with the smooth manifold $\, C_1\times ... \times C_m\,$ by identifying $\, \varphi \,$ with $\, (g_1, ... , g_m) = ( \varphi (z_1), ... , \varphi (z_m) )\,$. The group $\, G\,$ acts on  $\,\text{Hom} ( F, \, G )_{\mathcal C}\,$, on the conjugacy classes $\, C_j\,$ and on itself by conjugation.

Denote by $\, f_j:\text{Hom} ( F, \, G )_{\mathcal C}\rightarrow C_j, \,\, j=1, ... ,m\,$, the projection on the $\, j$th factor and by $\, r: \text{Hom} ( F, \, G )_{\mathcal C} \rightarrow G\,$ the mapping $\, r(\varphi )=  \varphi(z_1)\cdot ... \cdot \varphi (z_m)\,$. All these are $\, G$-equivariant maps.

 Let $\, K=r^{-1}(e)\subset  \text{Hom} ( F, \, G )_{\mathcal C} \,$, where $\, e\in G\,$ is the unit element.

According to \cite{GHJW} there exists a neighbourhood $\, \mathscr{M}\,$ of $\, K\,$ in 
$\,\text{Hom} ( F, \, G )_{\mathcal C}\,$ and a closed 2-form $\, \omega _{\mathcal C}\,$ on  $\, \mathscr{M}\,$ which is non-degenerate and hence gives a symplectic structure on  $\, \mathscr{M}\,$. The form $\, \omega _{\mathcal C}\,$ is denoted by $\, \omega _{c, \mathscr P,\mathcal C}\,$ in \cite{GHJW}. We shall now recall its definition.
 
The form  $\, \omega _{\mathcal C}\,$ is a sum of three parts.

\noindent

{\it The first part:} \,\, 
Let $\, C_{*}(F)\,$ be the chain complex of the nonhomogeneous reduced normalised bar resolution of $\, F$, see \cite{ML1}. (We use the notation of that reference.)  

The form $\, \omega _{\mathcal C}\,$ is not uniquely defined.  Its definition depends on a choice of a 2-chain $\, c\,$ in $\, C_2(F)\,$ satisfying 
\begin{displaymath}
\partial c = [z_1\cdot z_2\cdot ... \cdot z_m ]-[z_1]- ... -[z_m] \,\, ,
\end{displaymath}  
where, for any $\, x\in F\,$, the symbol $\, [x]\,$ denotes the corresponding generator of $\, C_1(F)\,$. As pointed out in \cite{GHJW}, page 391, one such possible choice is 
\begin{equation}
c = -[z_1\cdot ... \cdot z_{m-1} | z_m]-[z_1\cdot ... \cdot z_{m-2} | z_{m-1}] - ... - [z_1 | z_2] \,\, . 
\end{equation} 
We shall stay with this choice of $\, c\,$ for the rest of the paper.

Let us first recall some definitions from \cite{W1} and \cite{H1}. 
(We follow the convention used in \cite{H1} that if $\, \alpha\,$ is a $p$-form and  $\, \beta\,$ is a $q$-form then 
\begin{displaymath}
(\alpha \wedge \beta) \, (v_1, ... , v_{p+q}) = \sum\limits _{p,q\, \text{shuffles}} \, (\text{sgn}\,  \pi ) \alpha (v_{\pi (1)}, ... , v_{\pi (p)})\, \beta (v_{\pi (p+1)}, ... , v_{\pi (p+q)}), 
\end{displaymath}
for tangent vectors $\, v_1, ... , v_{p+q}\,$.)

We start by defining a $3$-form $\, \lambda\,$ on $\, G\,$ and a $2$-form $\, \Omega\,$ on $\, G\times G\,$.

Denote by $\, \omega\,$ the $\, {\mathfrak g}$-valued, left-invariant 1-form on $\, G\,$ which maps each tangent vector to the left-invariant vector field having that value. The corresponding right-invariant form will be denoted by $\, \bar{\omega}\,$. 

The tripple product 
\begin{displaymath}
\tau (x,y,z)= \frac 12 \, [x, y]\bullet z ,  \qquad  x, y, z \in {\mathfrak g}
\end{displaymath}
yields an alternating trilinear form on  $\, {\mathfrak g}$. Denote by $\, \lambda \,$ the left translate of $\, \tau\,$. It is a closed invariant $3$-form on $G$ and it satisfies
\begin{displaymath}
\lambda = \frac 1{12} \, [\, \omega, \omega ] \bullet \omega \,\, .
\end{displaymath}
(Note that by the convention used $\, \, [\, \omega, \omega ](X, Y ) = 2 [ X, Y ]\,$ for arbitrary left-invariant vector fields $\, X\,$ and $\, Y\,$ on $\, G\,$.)

For any differential form $\, \alpha \,$ on $\, G\,$ denote by $\, \alpha _j\,$ the pulback of  $\, \alpha \,$ to $\, G\times G\,$ by the projection $\, p_j\,$ to the $\, j$th factor. Let 
\begin{displaymath}
\Omega = \frac 12 \,\, \omega _1 \bullet  \bar{\omega} _2 \,\,\, .
\end{displaymath}
This is a real-valued 2 form on  $\, G\times G\,$. According to the convention 
\begin{equation}
( \omega _1 \bullet  \bar{\omega} _2 )(U, V) =  \omega _1 (U) \bullet  \bar{\omega} _2 (V) -  \bar{\omega} _2 (U) \bullet  \omega _1 (V) \,\,\, .
\end{equation}

Let us now consider the evaluation map 
\begin{displaymath}
E : F^2 \times \text{Hom} (F, \, G)_{\mathcal C} \rightarrow G^2 \,\, , 
\end{displaymath}
and, for every $\, x, y \in F\,$, let us denote by $\, E_{[x | y]}:\text{Hom} (F, \, G)_{\mathcal C}\rightarrow G^2\,$ the map given by 
$\, E_{[x | y]} (\varphi ) = ( \varphi (x), \, \varphi (y) )\,$. 
Denote by $\, \omega_{[x | y]}\,$ the $2$-form on $\,  \text{Hom} (F, \, G)_{\mathcal C}\,$ which is the pullback of $\, \Omega\,$ by the map $\, E_{[x | y]}\,$
\begin{displaymath}
 \omega_{[x | y]}= E_{[x | y]}^* (\Omega ) \,\, .
\end{displaymath}  
Finally, for $\, c= - \sum\limits _{j=1}^{m-1} \, [z_1\cdot ... \cdot z_j | z_{j+1} ] \in C_2(F)\,$, let 
\begin{equation}
\omega _c = - \sum\limits _{j=1}^{m-1} \, \omega _{ [z_1\cdot ... \cdot z_j | z_{j+1} ]} \,\, .
\end{equation}

\medskip

{\it The second part:} \,\, Let $\, \widetilde{\mathscr{O}}\,$ be an open 
$\, G$-equivariant subset of the Lie algebra  $\, {\mathfrak g}$ containing $\, 0\,$ and such that (i) the exponential mapping is a diffeomorphism of $\, \widetilde{\mathscr{O}}\,$ onto $\, \text{exp} ( \widetilde{\mathscr{O}})\,$, and (ii) $\, \widetilde{\mathscr{O}}\,$ is star-shaped w.r.t. $\, 0\,$.

Let $\, \mathscr{H} = r^{-1} ( \text{exp} (\widetilde{\mathscr{O}} ) )\,$. It is an open subset of  $\,\text{Hom} ( F, \, G )_{\mathcal C}\,$ containing $\, K\,$. Denote by $\,\eta \,$ the embedding of $\,\mathscr{H}\,$ into $\, \text{Hom} ( F, \, G )_{\mathcal C}  \,$ and by 
$\,
\tilde{r}: \mathscr{H}\rightarrow \widetilde{\mathscr{O}}
\,$
the composition of the map $\, r:\text{Hom} ( F, \, G )_{\mathcal C}\rightarrow G\,$ restricted to  $\,\mathscr{H}\,$ with the map  $\, \text{exp} ^{-1}: \text{exp} ( \widetilde{\mathscr{O}}) \rightarrow \widetilde{\mathscr{O}}\,$. Also denote   by
$\, 
\tilde{f} _j: \mathscr{H}\rightarrow C_j \, , \, j=1, ... ,m,
\,$
the restrictions of $\, f_j:\text{Hom} ( F, \, G )_{\mathcal C}\rightarrow C_j\,$ to $\,\mathscr{H}\,$.

Let $\, h:\Omega ^*( {\mathfrak g}  )\rightarrow \Omega ^{*-1}( {\mathfrak g}  ) \,$ be  the standard homotopy operator given by integration of forms along linear paths in  $\, {\mathfrak g}$,  see \cite{N1}, Lemma 2.13.1. Then $\, \beta = h( \text{exp}^*(\lambda ))\,$ is a $2$-form on $\,{\mathfrak g}\,$ and we define a $2$-form  $\, \omega _{c, \mathscr P}\,$ on
$\,\mathscr{H}   \,$ by
\begin{equation}
 \omega _{c, \mathscr P} = \eta ^* \omega _c - \tilde{r}^* \beta \,\, .
\end{equation}

\medskip

{\it The third part:} \,\, Let $\, C\,$ be a conjugacy class in $\, G\,$. For a point $\, p\,$ of  $\, C\,$, an arbitrary tangent vector is of the form 
\begin{equation}
Xp -pX = (X-Ad(p)X)p \in T_pC,
\end{equation} 
where $\, X\,$ is an element of the Lie algebra $\, \mathfrak{g}\,$ identified with the tangent space $\, T_eG\,$ of $\, G\,$ at $\, e\,$ and where $\, Xp\,$ and $\, pX\,$ are the right and the left translation of $\, X\,$ by $\, p\,$ respectively. The formula 
\begin{equation}
\tau (\,  Xp -pX, \,  Yp -pY\, ) =  \frac 12 (\, X \bullet Ad(p)\, Y - Y \bullet Ad(p)\, X \,), \qquad p\in C, 
\end{equation}
yields an equivariant 2-form $\, \tau\,$ on $\, C\,$, \cite{GHJW}, Sec.6.

For the conjugacy classes $\, C_1, ... , C_m\,$ in $\, G\,$ denote by $\, \tau _1, ... , \tau _m\,$ the corresponding 2-forms on them.
 
 Finally define a $2$-form   $\, \omega _{\mathcal C}\,$ on
$\,\mathscr{H}   \,$ by 
\begin{equation}
\omega _{\mathcal C} = \omega _{c, \mathscr{P}} + f_1^*\tau _1 + ... + f_m^*\tau _m \,\, .
\end{equation}
(The form $\, \omega _{\mathcal C}\,$ is denoted by  
$\, \omega _{c, \mathscr{P}, \mathcal{C}}\,$ in \cite{GHJW}.)

\medskip

One of the main results of \cite{GHJW}, in the special case considered in this Section (this is the case of a surface of genus $\, g=0\,$), is

\begin{thm}
[\cite{GHJW}; Thm.8.12] There is a $\, G$-invariant neighbourhood $\, \mathscr{M}\,$ of $\, K\,$ in $\, \mathscr{H}\,$ where the $2$-form  $\, \omega _{\mathcal C}\,$ is symplectic.
\end{thm} 

For the proof of the Theorem we refer the reader to \cite{GHJW}.

\bigskip

We conclude this Section with a straightforward observation.

For some special conjugacy classes $\, C_1, ... , C_m\,$ the $2$-form $\,  \omega _{\mathcal C}\,$ simplifies somewhat. That holds in particular for the conjugacy classes of interest in this paper. 

Let us suppose that $\, C\,$ is a conjugacy class in $\, G\,$ such that for every $\, p\in C\,$ one has
\begin{equation}
p^2 \in Z(G)\, ,
\end{equation} 
where $\, Z(G)\,$ is the center of $\, G\,$.

We have

\begin{lma}
If $\, C\,$ is a conjugacy class in $\, G\,$ satisfying the condition (3.8) then the $2$-form $\, \tau \,$ vanishes on $\, C\,$.
\end{lma}

\begin{proof}
Since the inner product $\, \bullet \,$ on  $\, \mathfrak{g}\,$ is invariant with respect to the adjoint action of $\, G\,$ and since, as a consequence of the condition (3.8), $\, Ad(p^2)= Id\,$ for every $\, p\in C\,$, it follows that 
\begin{equation*}
\begin{split}
 X \bullet Ad(p)\, Y &=  Ad(p)\, X \bullet  Ad(p)^2 \, Y =  Ad(p)\, X \bullet  Ad(p^2) \, Y =\\ & =  Ad(p)\, X \bullet \, Y = \\ & =  Y \,  \bullet \, Ad(p)\, X
\end{split}
\end{equation*}
for all $\, X, Y \in   \mathfrak{g}\,$.

Therefore the expression on the right hand side of (3.6) vanishes and so does the $2$-form $\, \tau\,$ on $\, C\,$.
\end{proof}

\begin{cor}
If the conjugacy classes in $\,\mathcal{C} = ( C_1, ... , C_m )\,$ all satisfy the condition (3.8) then the $2$-form $\, \omega _{\mathcal C}\,$ on $\, \mathscr{H}\,$ is given by 
\begin{displaymath}
 \omega _{\mathcal C} = \omega _{c, \mathscr P} = \eta ^* \omega _c - \tilde{r}^* \beta \,\, .
\end{displaymath}
\end{cor}

\medskip

Observe that the conjugacy class $\, \mathcal{S}\,$ in $\, SU(2)\,$ considered in the previous Sections does satisfy the condition (3.8) as 
$\, A^2= -I\,$ for every $\, A\in \mathcal{S}\,$.

\medskip

We shall abuse the notation and write $\,  \omega _c\,$ for the form $\, \eta ^* \omega _c\,$ on $\,  \mathscr{H}\,$.

\begin{cor}
If the conjugacy classes in $\,\mathcal{C} = ( C_1, ... , C_m )\,$ all satisfy the condition (3.8) then      at all points of $\, K=r^{-1}(e)=\tilde{r}^{-1}(0)\,$ the 2-form $\, \omega _{\mathcal C}   \,$ is given by 
\begin{displaymath}
 \omega _{\mathcal C} = \omega _c\quad  .
\end{displaymath}
\end{cor}

\begin{proof}
The mapping $\, \tilde{r}: \mathscr{H}\rightarrow \widetilde{\mathscr{O}}\,$ maps $\, K\,$ to $\, 0\in \mathscr{O}\subset  {\mathfrak g}\,$. Since the 2-form $\, \beta\,$ on $\,  {\mathfrak g}\,$ vanishes at $\, 0$, the 2-form $\,  \tilde{r}^*\beta\,$ vanishes over  $\, K\,$. Hence $\,  \omega _{\mathcal C} = \eta ^* \omega _c   = \omega _c\,$.
\end{proof}

\bigskip

\section{The action of the braid group}\label{action}

\bigskip

We continue with the notation of Section 3.

Let $\, {\mathcal B}_m\,$ be the braid group on $m$ strands. It is generated by $m-1$ elementary braids $\, \sigma _k, \,\, k= 1, ... ,m-1\,$ subject to relations 
\begin{equation*}
\begin{split}
&\sigma _k \, \sigma _j =\sigma _j \, \sigma _k \qquad \text{if}\,\,\,  |k-j|>1,\\
&\sigma _k \, \sigma _{k+1} \, \sigma _k = \sigma _{k+1} \, \sigma _k \, \sigma _{k+1} \qquad \text{for} \, \,\, k=1, ... , m-2.
\end{split} 
\end{equation*}

Let $\, F_m\,$ be the free group on $m$ generators $\, z_1, ... , z_m\,$. The braid group  $\, {\mathcal B}_m\,$ acts from the right on $\, F_m\,$. The action of the $k$-th elementary braid $\, \sigma _k\,$ is given by the automorphism $\, \sigma _k:F_m\rightarrow F_m\,$ defined by 
\begin{equation}
\sigma _k (z_j) = \left\{\begin{array}{lll}
z_j & \qquad & j\ne k, k+1,\\
z_k z_{k+1} z_k^{-1} & & j=k,\\
z_k & & j=k+1,
\end{array}   \right.
\end{equation}
see \cite{B1}, Cor. 1.8.3. 

The right action of  $\, {\mathcal B}_m\,$ on  $\, F_m\,$ induces a left action of  $\, {\mathcal B}_m\,$ on the space of all group homomorphisms $\, \text{Hom} (F_m, \, G )\,$ given by $\, \sigma (\varphi )(f)=\varphi (\sigma (f))\,$ for any $\, \sigma \in  {\mathcal B}_m, \,\, \varphi \in  \text{Hom} (F_m, \, G )\,$ and $\, f\in F_m\,$.

Let $\, C\,$ be a conjugacy class in $\, G\,$ and let $\, {\mathcal C}= (C, ... , C )\,$ be the sequaence of $m$ copies of $\, C\,$. Recall that  $\, \text{Hom} (F_m, \, G )_{\mathcal C}\,$ is the subspace of $\, \text{Hom} (F_m, \, G )\,$ consisting of those homomorphisms $\, \varphi :F_m\rightarrow G\,$ which map all generators $\, z_j\,$ of $\, F_m\,$ to $\, C\,$. As a consequence of the definition (4.1) we get that for any braid $\, \sigma \in {\mathcal B}_m\,$ and any generator $\ z_j\in F_m\,$ the element $\, \sigma (z_j)\,$ is conjugate to some generator $\, z_i\,$ of $\, F_m\,$. It follows  that if $\, \varphi \in  \text{Hom} (F_m, \, G )_{\mathcal C}\,$ then, for any braid $\, \sigma\,$, the homomorphism  $\,  \sigma (\varphi )\,$ also belongs to  $\, \text{Hom} (F_m, \, G )_{\mathcal C}\,$ i.e. that  $\, \text{Hom} (F_m, \, G )_{\mathcal C}\,$ is a  $\, {\mathcal B}_m$-invariant  subspace of  $\, \text{Hom} (F_m, \, G )\,$. Therefore the braid  group  $\, {\mathcal B}_m\,$ acts also on   $\, \text{Hom} (F_m, \, G )_{\mathcal C}\,$. (This action is a special case of actions of  $\, {\mathcal B}_m\,$ studied in \cite{R1}.)

If we identify  $\, \text{Hom} (F_m, \, G )_{\mathcal C}\,$ with the product $\,\prod\limits _{j=1}^m C\,$ of $m$ copies of $\, C\,$ by identifying $\, \varphi \in \text{Hom} (F_m, \, G )_{\mathcal C}\,$ with $\, (\varphi (z_1), ... ,\varphi (z_m)) \in \prod\limits _{j=1}^m C\,$ then the action of the $k$-th elementary braid $\, \sigma _k\,$  on $\, g=(g_1, ... , g_m)\in \prod C\,$ is given by 
\begin{equation}
(\sigma _k(g))_j = \left\{\begin{array}{lll}
g_j& \qquad & j\ne k, k+1,\\
g_k\, g_{k+1}\, g_k^{-1}&& j=k,\\
g_k && j=k+1,
\end{array}  \right.
\end{equation}
where $\, (.)_j\,$  stands for the projection on the $j$-th factor of $\,  \prod C\,$.

Recall the mapping $\, r:\text{Hom} (F_m, \, G )_{\mathcal C}\rightarrow G\,$ 
defined by $\, r(g_1, ... , g_m)=g_1\cdot ...\cdot  g_m\,$.  It follows from (4.2) that for every elementary braid $\, \sigma _k\,$  and $\, g=(g_1, ... , g_m)\in  \prod C\,$
\begin{displaymath}
r(\sigma _k(g)) =g_1\cdot ...\cdot  g_m = r(g) \, . 
\end{displaymath} 
Therefore, for every braid $\, \sigma \in {\mathcal B}_m\,$, we have 
\begin{equation}
r \circ \sigma = r  \,.
\end{equation}  

It follows from (4.3) that the subset $\,{\mathscr H} = r^{-1}( \text{exp} ( \widetilde{\mathscr{O}}))\,$ of   $\, \text{Hom} (F_m, \, G )_{\mathcal C}\,$ is invariant with respect to the action of $\, {\mathcal B}_m\,$. Moreover, since $\, \widetilde{r}:{\mathscr H}\rightarrow \widetilde{{\mathscr{O}}}\,$ is the composition of $\, r\,$ and $\,  \text{exp} ^{-1}\,$, we have 
\begin{equation}
\widetilde{r} \circ \sigma = \widetilde{r}  \,
\end{equation} 
and, consequently, 
\begin{equation}
\sigma ^* (\widetilde{r}^* \beta ) = \widetilde{r}^* \beta
\end{equation} 
for every braid  $\, \sigma \in {\mathcal B}_m\,$. 

Finally, consider the open subset $\, {\mathscr M} \,$ of $\,  {\mathscr H} \,$ the existence of which is ascertained in Theorem 3.1. Since the space $\, \prod C\,$ is compact and the group $\, G \,$ is Hausdorff, we can assume (possibly choosing a smaller $\, {\mathscr M} \,$) that $\, {\mathscr M} \,$  is of the form  $\, {\mathscr M}= r^{-1}(U) \,$ for some open $\, Ad$-invariant neighbourhood U of  the identity  element $\, e\,$ of the group $\ G\,$. It follows then from (4.3) that such a  subset  $\, {\mathscr M} \,$ of $\,  {\mathscr H} \,$ is not only $\, G$-invariant but also $\, {\mathcal B}_m$-invariant and inherits the action of $\, {\mathcal B}_m\,$. It follows also from (4.3) that $\, K=r^{-1}(e)\,$ is a  $\, {\mathcal B}_m$-invariant subset of $\, \prod C\,$ and of  $\, {\mathscr M} \,$.

The main aim of this section is to prove

\begin{thm}
If the conjugacy class $\, C\,$ satisfies the condition (3.8) then the symplectic $2$-form $\, \omega _{\mathcal C}\,$ on  $\, {\mathscr M} \,$ satisfies 
\begin{displaymath}
\sigma ^* (  \omega _{\mathcal C}  ) =  \omega _{\mathcal C}
\end{displaymath}
for every braid  $\, \sigma \in {\mathcal B}_m\,$ i.e. the braid group $\, {\mathcal B}_m\,$ acts on  $\, {\mathscr M} \,$ through symplectomorphisms.
\end{thm}  

\begin{proof}
According to Corollary 3.3 the $2$-form $\, \omega _{\mathcal C}\,$ on  $\, {\mathscr H} \,$ is given by
\begin{displaymath}
 \omega _{\mathcal C} = \eta ^* \omega _c - \widetilde{r}^* \beta \, .
\end{displaymath} 
By (4.5) we know already that $\, \sigma ^* (\widetilde{r}^* \beta ) = \widetilde{r}^* \beta \,$ for any braid $\, \sigma \in  {\mathcal B}_m\,$. Hence it is enough to show that $\, \sigma ^*( \eta ^* \omega _c) = \eta ^* \omega _c\,$.

The embedding  $\,\eta \,$ of $\, {\mathscr H}\,$ into $\, \text{Hom} (F_m, \, G )_{\mathcal C} = \prod C \,$ commutes, by the definition, with the action of 
$\, {\mathcal B}_m\,$. Therefore, it will be enough to prove that
\begin{equation}
 \sigma ^*( \omega _c) =  \omega _c
\end{equation}
for any braid  $\, \sigma \in  {\mathcal B}_m\,$.

According to (3.3) we have
\begin{equation}
\omega _c = - \sum\limits _{j=1}^{m-1} \, \omega _{ [z_1\cdot ... \cdot z_j | z_{j+1} ]} \,\, .
\end{equation}
Recall that $\, \omega _{ [z_1\cdot ... \cdot z_j | z_{j+1} ]} = E_j^*(\Omega )\,$, where $\, E_j: \prod\limits _{j=1}^m C \rightarrow G\times G\,$ is given by $\, E_j(g_1, ... , g_m)= (g_1\cdot ... \cdot g_j, \, g_{j+1})\,$ and where $\, \Omega \,$ is the $2$-form on $\, G\times G\,$ defined by (3.2).

Denote by $\, F_k: \prod C \rightarrow \prod C\,$ the diffeomorphism of $\, \prod C\,$ given by the action of the $k$-th elementary braid $\, \sigma _k\,$ and described in (4.2). 

Thus $\, \sigma _k^*(\omega _{ [z_1\cdot ... \cdot z_j | z_{j+1} ]}\, ) =
 F_k^*(\omega _{ [z_1\cdot ... \cdot z_j | z_{j+1} ]}\, ) = F_k^* E_j^* (\Omega ) = (E_j \circ F_k)^*(\Omega )\,$. 

Observe that if $\, j+1 \ne k, k+1\,$ then  the mapping $\, E_j \circ F_k : \prod C \rightarrow G\times G\,$ satisfies $\, 
E_j\circ F_k = E_j\,$.
Hence 
\begin{equation}
\sigma _k^*(\omega _{ [z_1\cdot ... \cdot z_j | z_{j+1} ]}\, ) = (E_j \circ F_k)^*(\Omega ) = E_j^* ( \Omega ) = \omega _{ [z_1\cdot ... \cdot z_j | z_{j+1} ]}
\end{equation}
if $\, j+1 \ne k, k+1\,$.

We shall now show that 
\begin{equation}
\sigma _k^*(\omega _{ [z_1\cdot ... \cdot z_{k-1} | z_k ]}+ \omega _{ [z_1\cdot ... \cdot z_k | z_{k+1} ]}      \, ) = \omega _{ [z_1\cdot ... \cdot z_{k-1} | z_k ]}+ \omega _{ [z_1\cdot ... \cdot z_k | z_{k+1} ]}  
\end{equation} 
on $\, \prod C =\text{Hom} (F_m, \, G )_{\mathcal C}\,$. 

Before we start let us observe that the $2$-forms $\, \omega _{[x|y]} = E_{[x|y]}^*(\Omega )\,$ are defined not only on $\, \prod C =\text{Hom} (F_m, \, G )_{\mathcal C}\,$ but also on the whole manifold $\, \prod\limits _{i=1}^m G =\text{Hom} (F_m, \, G )\,$. 

Let $\, g=(g_1, ... , g_m)\in \prod\limits _{i=1}^m C\,$. We represent tangent vectors to the conjugacy class  $\, C\,$ at the point $\, g_i\in C\,$ in the form $\, X_i\, g_i\,$ with $\, X_i\in {\mathfrak g}\,$. Let $\, X = ( X_1\, g_1, ... ,  X_m\, g_m )\,$ and  $\, Y = ( Y_1\, g_1, ... ,  Y_m\, g_m )\,$ be two arbitrary tangent vectors to $\ \prod C\,$ at the point $\, g\,$. Here $\, X_1, ... , X_m,  Y_1, ..., Y_m \in {\mathfrak g}\,$.

As $\, E_{k-1}: \prod C \rightarrow G\times G ,\,\,  E_{k-1}(g_1, ... , g_m)= (g_1\cdot ... \cdot g_{k-1}, \, g_k )\,$ and  $\, E_k: \prod C \rightarrow G\times G ,\,\,  E_k(g_1, ... , g_m)= (g_1\cdot ... \cdot g_k, \, g_{k+1} )\,$, we have 
\begin{equation*}
\begin{split}
dE_{k-1}(X)& =dE_{k-1}(  X_1\, g_1, ... ,  X_m\, g_m ) = \\ & =\bigg( \sum\limits _{j=1}^{k-1} \, g_1...g_{j-1} X_j g_j...g_{k-1},\, \, X_k g_k \bigg) = \\
& = \bigg( g_1... g_{k-1} \, \sum\limits _{j=1}^{k-1} \, Ad(g_j ...g_{k-1} )^{-1}(X_j), \,\, X_k g_k \bigg)
\end{split}
\end{equation*}
and  
\begin{equation*}
dE_k (X) = \bigg( g_1... g_k \, \sum\limits _{j=1}^k \, Ad(g_j ...g_k )^{-1}(X_j), \,\, X_{k+1}\, g_{k+1} \bigg)
\end{equation*}
and similarily for the tangent vector $\, Y\,$.

Hence, according to (3.2), 
\begin{equation}
\begin{split}
\omega _{ [z_1\cdot ... \cdot z_{k-1} | z_k ]} (X, \, Y )& = E_{k-1}^* (\Omega ) ( X, \, Y ) = \Omega (dE_{k-1} (X), \, dE_{k-1} (Y)) =\\
& = \frac 12 \big[ \omega _1(dE_{k-1} (X))\bullet \bar{\omega }_2 (dE_{k-1} (Y)) - \\ 
& \qquad \qquad \qquad -\omega _1(dE_{k-1} (Y))\bullet \bar{\omega }_2 (dE_{k-1} (X)) \big]=\\
& = \frac 12 \bigg[ \bigg( \sum\limits _{j=1}^{k-1} \, Ad(g_j ...g_{k-1} )^{-1}(X_j)\bigg) \bullet Y_k - \,\, \text{Symm}\,\,  \bigg]
\end{split}
\end{equation}
and
\begin{equation}
\begin{split}
\omega _{ [z_1\cdot ... \cdot z_k | z_{k+1} ]} (X, \, Y )& = \frac 12 \bigg[ \bigg( \sum\limits _{j=1}^k \, Ad(g_j ...g_k )^{-1}(X_j)\bigg) \bullet Y_{k+1} - \,\, \text{Symm}\,\,  \bigg] \, ,
\end{split}
\end{equation}
where the term ``Symm''  means ``the same terms with the r$\hat o$les of $\, X\,$ and $\, Y\,$ exchanged''.

On the other hand 
\begin{equation*}
\begin{split}
dF_k (X) &= dF_k(X_1 g_1, ... , X_m g_m) = \\
& = ((X_1 g_1, ... , X_{k-1} g_{k-1},\, X_k g_k g_{k+1} g_k^{-1} + g_k X_{k+1} g_{k+1} g_k^{-1} -\\ 
&\qquad \qquad -  g_k g _{k+1} g_k^{-1} X_k, \, X_k g_k,\,  X_{k+2} g_{k+2}, ... , X_m g_m ) = \\
& = ((X_1 g_1, ... , X_{k-1} g_{k-1},\, ( X_k + Ad(g_k) (X_{k+1}) -\\ 
&\qquad \qquad - Ad( g_k g _{k+1} g_k^{-1}) (X_k) ) g_k g_{k+1} g_k^{-1} , \, X_k g_k,\,  X_{k+2} g_{k+2}, ... , X_m g_m )
\end{split}
\end{equation*}
and similarily for $\, Y\,$.

Therefore, by the same argument  as in the derivation of (4.10) and (4.11), we get
\begin{equation}
\begin{split}
(\sigma _k^* \, \omega _{ [z_1\cdot ... \cdot z_{k-1} | z_k ]}\, ) (X, \, Y )& = 
 \omega _{ [z_1\cdot ... \cdot z_{k-1} | z_k ]} (dF_k (X), \, dF_k(Y) )
=\\
& = \frac 12 \bigg[ \bigg( \sum\limits _{j=1}^{k-1} \, Ad(g_j ...g_{k-1} )^{-1}(X_j)\bigg) \bullet \bigg( Y_k +\\ 
& \qquad \qquad + Ad(g_k)(Y_{k+1}) - Ad( g_k g _{k+1} g_k^{-1}) (Y_k)\bigg)- \\
& \qquad \qquad \qquad - \text{Symm} \,\,\, \bigg]\, .
\end{split}
\end{equation}

Similarily,
\begin{equation}
\begin{split}
(&\sigma _k^* \, \omega _{ [z_1\cdot  ... \cdot z_k | z_{k+1} ]} \, ) (X, \, Y ) = 
 \omega _{ [z_1\cdot ... \cdot z_k | z_{k+1} ]} (dF_k (X), \, dF_k(Y) )
=\\
& = \frac 12 \bigg[ \bigg( \sum\limits _{j=1}^{k-1} \, Ad(g_j ...g_{k-1} g_k g_{k+1} g_k^{-1}  )^{-1}(X_j) + \\ 
&\qquad +  Ad( g_k g_{k+1} g_k^{-1}  )^{-1} \big( X_k + Ad(g_k)(X_{k+1}) - Ad( g_k g_{k+1} g_k^{-1}  )(X_k)\big)     \bigg) \bullet Y_k -\\
& \qquad \qquad \qquad  - \text{Symm} \,\,\, \bigg] = \\
& = \frac 12 \bigg[ \bigg( \sum\limits _{j=1}^{k-1} \, Ad(g_j ...  g_k )^{-1}(X_j) +  \\
& \qquad +  Ad( g_k  )^{-1} \big( X_k + Ad(g_k)(X_{k+1}) - Ad( g_k g_{k+1} g_k^{-1}  )(X_k)\big)     \bigg) \bullet Ad( g_{k+1} g_k^{-1}) (Y_k) -\\
& \qquad \qquad\qquad  - \text{Symm} \,\,\, \bigg] =\\
& = \frac 12 \bigg[ \bigg( \sum\limits _{j=1}^k \, Ad(g_j ...  g_k )^{-1}(X_j) + X_{k+1}  - Ad( g_{k+1} g_k^{-1}  )(X_k) \bigg) \bullet Ad( g_{k+1} g_k^{-1}) (Y_k) -\\
& \qquad \qquad\qquad  - \text{Symm} \,\,\, \bigg]\, .
\end{split}
\end{equation}
The next to the last equality follows from the invariance of the inner product $\bullet \,$ with respect to $\, Ad( g_{k+1} g_k^{-1})\,$ .

It follows from (4.10) and (4.12) that 
\begin{equation}
\begin{split}
A:& = (\sigma _k^* \, \omega _{ [z_1\cdot  ... \cdot z_{k-1} | z_k ]} \, ) (X, \, Y ) - \omega _{ [z_1\cdot  ... \cdot z_{k-1} | z_k ]} \,  (X, \, Y ) =\\
&= \frac 12 \bigg[ \bigg( \sum\limits _{j=1}^{k-1} \, Ad(g_j ...g_{k-1} )^{-1}(X_j)\bigg)\bullet  \bigg( Ad(g_k)(Y_{k+1}) - Ad( g_k g _{k+1} g_k^{-1}) (Y_k)\bigg)- \\
& \qquad \qquad \qquad - \text{Symm} \,\,\, \bigg] =\\
&= \frac 12 \bigg[ \bigg( \sum\limits _{j=1}^{k-1} \, Ad(g_j ...g_{k-1} g_k )^{-1}(X_j)\bigg)\bullet  \bigg( Y_{k+1} - Ad( g _{k+1} g_k^{-1}) (Y_k)\bigg) - \text{Symm} \,\,\, \bigg]\,.
\end{split}
\end{equation}
The last equality holds by the invariance of $\, \bullet \,$ with respect to $\, Ad(g_k)^{-1}\,$.

Similarily, it follows from (4.11) and (4.13) that 
\begin{equation}
\begin{split}
B:&=  (\sigma _k^* \, \omega _{ [z_1\cdot  ... \cdot z_k | z_{k+1} ]} \, ) (X, \, Y ) - \omega _{ [z_1\cdot  ... \cdot z_k | z_{k+1} ]} \,  (X, \, Y ) =\\
& = \frac 12 \bigg[ \bigg( \sum\limits _{j=1}^k \, Ad(g_j ...  g_k )^{-1}(X_j)\bigg) \bullet \bigg( Ad( g_{k+1} g_k^{-1}) (Y_k) - Y_{k+1} \bigg) + \\
& \qquad \qquad  + \bigg( X_{k+1}  - Ad( g_{k+1} g_k^{-1}  )(X_k) \bigg) \bullet \bigg( Ad( g_{k+1} g_k^{-1}) (Y_k)\bigg) \,\,\,  - \,\, \, \text{Symm} \,\,\, \bigg]\, .
\end{split}
\end{equation}

Consequently,
\begin{equation*}
\begin{split}
A+B&= \frac 12 \bigg[ \bigg(Ad(g_k)^{-1} (X_k)\bigg) \bullet   \bigg( Ad( g_{k+1} g_k^{-1}) (Y_k) - Y_{k+1} \bigg) + \\
& \qquad \qquad  + \bigg( X_{k+1}  - Ad( g_{k+1} g_k^{-1}  )(X_k) \bigg) \bullet \bigg( Ad( g_{k+1} g_k^{-1}) (Y_k)\bigg) \,\,\,  - \,\, \, \text{Symm} \,\,\, \bigg] = \\
& = \frac 12 \bigg[ \bigg((Ad(g_k)^{-1} (X_k)) \bullet ( Ad( g_{k+1} g_k^{-1}) (Y_k)) - \\ 
&\qquad \qquad \qquad  - (Ad(g_k)^{-1} (Y_k)) \bullet ( Ad( g_{k+1} g_k^{-1}) (X_k))\bigg) - \\
&\quad \quad  -  \bigg( (Ad(g_k)^{-1} (X_k)) \bullet  Y_{k+1} -(Ad(g_k)^{-1} (Y_k)) \bullet  X_{k+1}\bigg) +\\
&\quad \quad  + \bigg( X_{k+1} \bullet ( Ad( g_{k+1} g_k^{-1}) (Y_k))-  Y_{k+1} \bullet ( Ad( g_{k+1} g_k^{-1}) (X_k)) \bigg) -\\
&\quad \quad  - X_k \bullet Y_k +  Y_k \bullet X_k \,\, \bigg] \, .
\end{split}
\end{equation*}
Rearranging the terms and using the invariance of the inner product one gets 
\begin{equation}
\begin{split}
A+B& = \frac 12 \bigg[ \bigg((Ad(g_{k+1}^{-1}g_k^{-1}) (X_k)) \bullet ( Ad( g_k^{-1}) (Y_k)) - \\ 
&\qquad \qquad \qquad  - (Ad(g_k)^{-1} (Y_k)) \bullet ( Ad( g_{k+1} g_k^{-1}) (X_k))\bigg) - \\
& \quad  \quad -  \bigg( (Ad(g_k)^{-1} (X_k)) \bullet  Y_{k+1} + (Ad(g_{k+1}^{-1}) (Y_{k+1})) \bullet ( Ad( g_k^{-1}) (X_k)) \bigg) +\\
& \quad \quad  + \bigg((Ad(g_{k+1}^{-1}) (X_{k+1}) \bullet ( Ad(g_k^{-1}) (Y_k))+ (Ad(g_k)^{-1} (Y_k)) \bullet  X_{k+1}  \bigg) \,\, \bigg] =\\
& =\frac 12 \bigg[ \bigg((Ad(g_{k+1}^{-1}g_k^{-1}) (X_k))-  ( Ad( g_{k+1} g_k^{-1}) (X_k))\bigg)     \bullet ( Ad( g_k^{-1}) (Y_k)) - \\ 
& \quad \quad  -  \bigg(  Y_{k+1} + Ad(g_{k+1}^{-1}) (Y_{k+1})\bigg) \bullet ( Ad( g_k^{-1}) (X_k))  +\\
& \quad \quad  + \bigg( X_{k+1} + Ad(g_{k+1}^{-1}) (X_{k+1})\bigg) \bullet ( Ad(g_k^{-1}) (Y_k)) \,\, \bigg] \, .
\end{split}
\end{equation}  
Since the conjugacy class $\, C\,$ satisfies the condition (3.8) we have $\, g_{k+1}^2 = d \,$ with $\, d\,$ belonging to the center $\, Z(G)\,$ of $\, G\,$. Hence $\, g_{k+1}^{-1} = d^{-1} g_{k+1}\,$ and $\, Ad( g_{k+1}^{-1}) = Ad( d^{-1} g_{k+1}) = Ad( g_{k+1})\,$. Therefore 
\begin{equation}
Ad(g_{k+1}^{-1}g_k^{-1}) (X_k))=  Ad( g_{k+1} g_k^{-1}) (X_k)
\end{equation} 
and the first term in the last expression of (4.16) vanishes.

Moreover, if $\, X_j  g_j\,$ is a tangent vector to the conjugacy class $\, C\,$ at the point $\, g_j\,$ then there exists $\, \widetilde{X}_j\in {\mathfrak g}\,$ such that $\,  X_j  g_j =  \widetilde{X}_j g_j - g_j  \widetilde{X}_j = ( \widetilde{X}_j - Ad(g_j)(\widetilde{X}_j) ) g_j\,$.  Therefore $\, X_j =  \widetilde{X}_j - Ad(g_j)(\widetilde{X}_j) \,$. Applying this observation with $\, j=k+1\,$ we get 
\begin{equation}
\begin{split}
X_{k+1} + Ad( g_{k+1}^{-1})(X_{k+1})& =  \widetilde{X}_{k+1} - Ad(g_{k+1})(\widetilde{X}_{k+1} ) + \\
& \qquad   + Ad( g_{k+1}^{-1}) \big(  \widetilde{X}_{k+1}  - Ad(g_{k+1} )(\widetilde{X}_{k+1} )\big) = \\
& =  - Ad(g_{k+1})(\widetilde{X}_{k+1} ) + Ad( g_{k+1}^{-1}) (  \widetilde{X}_{k+1})=\\
& =- Ad(g_{k+1}^{-1})(\widetilde{X}_{k+1} ) + Ad( g_{k+1}^{-1}) (  \widetilde{X}_{k+1}) =\\
& = \,\,\, 0
\end{split}
\end{equation} 
and, similarily, 
\begin{equation}
Y_{k+1} + Ad( g_{k+1}^{-1})(Y_{k+1}) = 0 \, .
\end{equation}
Hence the second and the third term in the last expression of (4.16) vanish.

\medskip

The equality (4.16) together with (4.17), (4.18) and (4.19) gives us now
\begin{equation*}
\begin{split}
\big(\sigma _k^* \,(& \omega _{ [z_1\cdot  ... \cdot z_{k-1} | z_k ]} +  \omega _{ [z_1\cdot  ... \cdot z_k | z_{k+1} ]} \, )\big) (X, \, Y ) - \\
& \qquad \qquad - ( \omega _{ [z_1\cdot  ... \cdot z_{k-1} | z_k ]} +  \omega _{ [z_1\cdot  ... \cdot z_k | z_{k+1} ]} \, ) (X, \, Y ) = \\
& =  A+B =\\
& = 0
\end{split}
\end{equation*}
for all vectors $\, X\,$ and $\, Y\,$ tangent to $\, \prod C\,$ at a point $\, g = ( g_1, ... , g_m )\,$. Therefore 
\begin{equation*}
\begin{split}
\sigma _k^* \,& \omega _{ [z_1\cdot  ... \cdot z_{k-1} | z_k ]} + \sigma _k^* \, \omega _{ [z_1\cdot  ... \cdot z_k | z_{k+1} ]}  = 
  \omega _{ [z_1\cdot  ... \cdot z_{k-1} | z_k ]} +  \omega _{ [z_1\cdot  ... \cdot z_k | z_{k+1} ]} \, .
\end{split}
\end{equation*}
Together with the identities (4.8) that proves 
\begin{equation*}
\begin{split}
\sigma _k^* \,\omega _c & =\sigma _k^* \bigg(- \sum\limits _{j=1}^{m-1}  \omega _{ [z_1\cdot  ... \cdot z_j | z_{j+1} ]} \bigg) = - \sum\limits _{j=1}^{m-1}  \omega _{ [z_1\cdot  ... \cdot z_j | z_{j+1} ]} =\\
& = \omega _c
\end{split}
\end{equation*}
for all elementary braids  $\, \sigma _k , \,\, k=1, ... , m-1\,$. Therefore 
\begin{displaymath}
\sigma ^* \omega _c =\omega _c 
\end{displaymath}
for all braids $\, \sigma \in {\mathcal B}_m\,$.  According to (4.6) that concludes the proof of Theorem 4.1.
\end{proof}

\bigskip

\section{The Lagrangian submanifolds}\label{lagr}

\bigskip

We continue with the notation of Section 3.

We shall now consider the case when $\, m=2n\,$.

Let $\, C_1, ... , C_n\,$ be conjugacy classes in the compact Lie group $\, G\,$. For $\, i=1, ... , n\,$ define conjugacy classes $\, C_{2n-i+1}\,$ 
in $\, G\,$ by
\begin{displaymath}
 C_{2n-i+1} = C_i^{-1} = \{g\in G \,\vert \, g^{-1}\in C_i\,\}
\end{displaymath}
and let $\,{\mathcal C}=( C_1, ... , C_n, C_{n+1},  ... , C_{2n})\,$. Consider a submanifold $\, \Lambda \,$ of $\, \text{Hom}(F_{2n}, \, G)_{\mathcal C} = C_1 \times ...\times C_{2n}\,$,
\begin{displaymath}
\Lambda = \{ (p_1, ... p_{2n}) \in  C_1 \times ...\times C_{2n} \, \vert \, p_{2n-i+1} = p_i^{-1}, \,\,\, i=1, ... , n \,\}.
\end{displaymath}
Observe that the projection $\, C_1 \times ...\times C_{2n}\rightarrow  C_1 \times ...\times C_n\,$  onto the first $\, n\,$ factors gives a diffeomorphism between  $\, \Lambda \,$ and the product 
$\,  C_1 \times ...\times C_n\,$ with the inverse mapping given by  $\,(p_1, ... , p_n)\mapsto (p_1, ... , p_n, p_n^{-1}, ... , p_1^{-1})\,$.

Recall the mapping $\, r: C_1 \times ...\times C_{2n} \rightarrow G, \,\, r(g_1, ... , g_{2n})=g_1\cdot ... \cdot g_{2n}\,$. Note that $\, r\,$ maps all points of  $\, \Lambda \,$ onto the identity $\, I\in G\,$. Hence  $\, \Lambda  \subset K = r^{-1}(I) \subset {\mathscr M}\,$. According to Theorem 3.1,  $\,  {\mathscr M}\,$ equipped with the $2$-form $\,\omega _{\mathcal C}\,$ is a symplectic manifold.

We assume for the rest of this section that all the conjugacy classes $\, C_1, ... , C_n\, $ satisfy the condition (3.8). 

\begin{prop}\label{P5.1} $\, \Lambda \,$ is a Lagrangian submanifold of $\,  {\mathscr M}\,$.
\end{prop} 

\begin{proof}
Let us consider a point $\,x = (p_1,... , p_n, p_n^{-1}, ... , p_1^{-1}) \in \Lambda \,$. Let $\, i\,$ and $\, j\,$ be a pair of integers such that $\, 1\le i,j\le n\,$ and let $\, X_i, \, Y_j\,$ be elements of the Lie algebra $\, {\mathfrak g}\,$ such that $\, X_i\cdot p_i\in T_{p_i}G\,$ and  $\, Y_j\cdot p_j\in T_{p_j}G  \,$  are vectors tangent to $\, C_i\,$ at $\, p_i\,$  and  to $\, C_j\,$ at $\, p_j\,$ respectively. Then  $\, U=\oplus\, U_k, \,\,\,  V=\oplus\, V_k \in T_x(  C_1 \times ...\times C_{2n})\,$ given by 
\begin{equation}\label{5.1}
U_k =\begin{cases} X_i \cdot p_i & \qquad \text{if} \quad k=i,\\
-p_i^{-1} \cdot X_i = (-Ad(p_i^{-1})(X_i))\cdot p_i^{-1} & \qquad \text{if} \quad k=2n-i+1,\\ 
0 & \qquad \text{otherwise}
\end{cases} 
\end{equation}
and
\begin{equation}\label{5.2}
V_k =\begin{cases} Y_j \cdot p_j & \qquad \text{if} \quad k=j,\\
-p_j^{-1} \cdot Y_j = (-Ad(p_j^{-1})(Y_j))\cdot p_j^{-1} & \qquad \text{if} \quad k=2n-j+1,\\ 
0 & \qquad \text{otherwise}
\end{cases} 
\end{equation}
are tangent to the submanifold  $\, \Lambda \,$ of $\, C_1 \times ...\times C_{2n}\,$ at the point $\, x\,$ (and the whole tangent space $\, T_x(\Lambda )\,$ is spanned by such vectors).

We shall now show that  $\,\omega _{\mathcal C} (U, \, V )= 0\,$.
Since the conjugacy classes $\, C_1, ... ,C_n\,$ (and, hence, their inverses  $\, C_{n+1}, ... ,C_{2n}\,$) satisfy the condition (3.8) and since $\, \Lambda \subset K\,$, we have  $\,\omega _{\mathcal C} (U, \, V )=\omega _c (U, \, V ) \,$ (see Corollary 3.4).

Let us recall that $\, \omega _c =  - \sum\limits _{k=1}^{2n-1} \, \omega _{ [z_1\cdot ... \cdot z_k | z_{k+1} ]}\,$.

We shall consider two cases.

(i) Let us suppose that $\, i=j\,$. Then, by the definition of $\, \omega _{ [z_1\cdot ... \cdot z_k | z_{k+1} ]}\,$,  we have $\, \omega _{ [z_1\cdot ... \cdot z_k | z_{k+1} ]}(U, \, V)=0\,$  for all $\, k\ne 2n-i\,$. For $\, k=2n-i\,$ we get 
\begin{equation*}
\begin{split}
\omega _{ [z_1\cdot ... \cdot z_{2n-i} | z_{2n-i+1} ]}&(U, \, V)= \tfrac 12 \, \bigl\{ 
\omega(df_{z_1...z_{2n-i}}(U))\bullet \overline{\omega}(df_{z_{2n-i+1}}(V)) - \\
&\qquad \qquad \qquad\qquad   
-\omega(df_{z_1...z_{2n-i}}(V))\bullet \overline{\omega}(df_{z_{2n-i+1}}(U))
\bigr\}=\\
&=\tfrac 12 \, \bigl\{ \omega(p_1...p_{i-1}X_ip_i...p_{2n-i})\bullet \overline{\omega}(-Ad(p_i^{-1})(Y_i)\cdot p_i^{-1}  ) - \\&\qquad \qquad \quad  
 -\omega( p_1...p_{i-1}Y_ip_i...p_{2n-i} )\bullet \overline{\omega}(-Ad(p_i^{-1})(X_i)\cdot p_i^{-1}  )
\bigr\}=\\
&=\tfrac 12 \, \bigl\{( Ad((p_i...p_{2n-i})^{-1})(X_i))\bullet (-Ad(p_i^{-1})(Y_i))) - \\& \qquad \qquad \quad 
 -( Ad((p_i...p_{2n-i})^{-1})(Y_i))\bullet (-Ad(p_i^{-1})(X_i))
\bigr\}=\\
&=\tfrac 12 \, \bigl\{( Ad((p_i...p_{2n-i}p_i^{-1})^{-1})(X_i))\bullet (-Y_i)) - \\& \qquad \qquad \quad 
 -( Ad((p_i...p_{2n-i}p_i^{-1} )^{-1})(Y_i))\bullet (-X_i)
\bigr\}
\end{split}
\end{equation*}
Since for every point $\, x=(p_1, ... , p_n, p_n^{-1}, ... ,p_1^{-1})\in \Lambda\,$ one has $\, p_i...p_{2n-i}p_i^{-1}=I\,$, we get 
\begin{equation}\label{5.3}
\omega _{ [z_1\cdot ... \cdot z_{2n-i} | z_{2n-i+1} ]}(U, \, V)     = \tfrac 12 \, \bigl\{
X_i\bullet (-Y_i)
 -Y_i\bullet (-X_i)
\bigr\} = 0\, .
\end{equation}
Thus  $\,\omega _{\mathcal C} (U, \, V )=  \omega _c (U, \, V )= 0\,$.
\medskip

(ii) Let us now suppose that $\, i<j\,$. Observe first of all that $\, \omega _{ [z_1\cdot ... \cdot z_k | z_{k+1} ]}(U, \, V)=0\,$  for all $\, k\ne j-1,\,\,  2n-j\,$ and $\, 2n-i\,$. That follows from  the definition of $\, \omega _{ [z_1\cdot ... \cdot z_k | z_{k+1} ]}\,$. Then we have  
\begin{equation}\label{5.4}
\begin{split}
\omega _{ [z_1\cdot ... \cdot z_{j-1} | z_j ]}(U, \, V)&= \tfrac 12 \, \bigl\{
\omega (df_{z_1...z_{j-1}}(U))\bullet \overline{\omega}(df_{z_j}(V) -\\
&\qquad\qquad \qquad 
-\omega (df_{z_1...z_{j-1}}(V))\bullet \overline{\omega}(df_{z_j}(U))
\bigr\}=\\
&=\tfrac 12 \, \bigl\{
\omega ( p_1 ... p_{i-1}X_ip_i...p_{j-1})\bullet \overline{\omega}(Y_j\cdot p_j)
\bigr\}=\\
&=\tfrac 12 \, \bigl\{( 
Ad((p_i...p_{j-1})^{-1})(X_i))\bullet Y_j
\bigr\}\,,
\end{split}
\end{equation}
where the second equality follows from $\, df_{z_j}(U)=0\,$. Similarily
\begin{equation}\label{5.5}
\begin{split}
\omega _{ [z_1\cdot ... \cdot z_{2n-j} | z_{2n-j+1} ]}&(U, \, V)= \tfrac 12 \, \bigl\{
\omega (df_{z_1...z_{2n-j}}(U))\bullet \overline{\omega}(df_{z_{2n-j+1}}(V)) -\\
&\qquad\qquad \qquad 
-\omega (df_{z_1...z_{2n-j}}(V))\bullet \overline{\omega}(df_{z_{2n-j+1}}(U))
\bigr\}=\\
&=\tfrac 12 \, \bigl\{
\omega ( p_1 ... p_{i-1}X_ip_i...p_{2n-j})\bullet \overline{\omega}(-p_j^{-1}\cdot Y_j)
\bigr\}=\\
&=\tfrac 12 \, \bigl\{
\omega ( p_1 ... p_{i-1}X_ip_i...p_{2n-j})\bullet \overline{\omega}(-Ad(p_j^{-1})(Y_j)\cdot p_j^{-1}  )
\bigr\}=\\
&=\tfrac 12 \, \bigl\{
Ad((p_i...p_{2n-j})^{-1})(X_i))\bullet (-Ad(p_j^{-1})(Y_j))
\bigr\}=\\
&=- \tfrac 12 \, \bigl\{
( Ad((p_i...p_{2n-j}p_j^{-1} )^{-1})(X_i))\bullet Y_j
\bigr\}  \,.
\end{split}
\end{equation}
Now, for every point $\, x=(p_1, ... , p_n, p_n^{-1}, ... ,p_1^{-1})\in \Lambda\,$, one has $\, p_i...p_{2n-j}p_j^{-1}= p_i...p_{j-1}  \,$, provided $\, i<j\,$. Hence
\begin{equation}\label{5.6}
\omega _{ [z_1\cdot ... \cdot z_{2n-j} | z_{2n-j+1} ]}(U, \, V)=- \tfrac 12 \, \bigl\{
( Ad((p_i...p_{j-1} )^{-1})(X_i))\bullet Y_j
\bigr\}  \,.
\end{equation} 

\medskip

Finally, observe that 
\begin{equation*}
\begin{split}
 df_{z_1...z_{2n-i}}(V)&=p_1...p_{j-1}Y_jp_j...p_{2n-i}
 - p_1...p_{2n-j}p_j^{-1}Y_jp_{2n-j+2}...p_{2n-i}=\\
&= p_1...p_{j-1}\bigl( Y_jp_j...p_{2n-j+1} - p_j...p_{2n-j+1}Y_j \bigr) p_{2n-j+2}...p_{2n-i}=\\
&= 0\, ,
\end{split}
\end{equation*}
because for all points $\,  x= (p_1, ... ,p_n, p_n^{-1}, ... , p_1^{-1})\in \Lambda\,$ one has $\, p_j\cdot ... \cdot p_{2n-j+1}=I\in G\,$.  Since, according to  (\ref{5.2}), we also have $\, df_{z_{2n-i+1}}(V) = 0\,$, it follows that  
\begin{equation}\label{5.7}
\begin{split}
\omega _{ [z_1\cdot ... \cdot z_{2n-i} | z_{2n-i+1} ]}(U, \, V)&= \tfrac 12 \, \bigl\{
\omega (df_{z_1...z_{2n-i}}(U))\bullet \overline{\omega}(df_{z_{2n-i+1}}(V)) -\\
&\qquad\qquad \qquad 
-\omega (df_{z_1...z_{2n-i}}(V))\bullet \overline{\omega}(df_{z_{2n-i+1}}(U))
\bigr\}=\\
& = 0.
\end{split}
\end{equation}

Therefore, 
\begin{equation*}
\begin{split}
\omega _{\mathcal C}(U, \, V ) &= -\sum\limits _{k=1}^{2n-1} \omega _{ [z_1\cdot ... \cdot z_k | z_{k+1} ]}(U, \, V) =\\
&= -\omega _{ [z_1\cdot ... \cdot z_{j-1} | z_j ]}(U, \, V) - \omega _{ [z_1\cdot ... \cdot z_{2n-j} | z_{2n-j+1} ]}(U, \, V)- \\
&\qquad \qquad \qquad\qquad \qquad \qquad \qquad  - \omega _{ [z_1\cdot ... \cdot z_{2n-i} | z_{2n-i+1} ]}(U, \, V)=\\
& = 0 \,.
\end{split}
\end{equation*}
That proves Proposition 5.1.
\end{proof}

Let now $\, C\,$ be a conjugacy class in $\, G\,$ satisfying the condition (3.8) and such that 
\begin{equation}\label{5.8}
 C=C^{-1}= \{g\in G \,\vert \, g^{-1}\in C\,\}  \, .
\end{equation}
(An example of such a class is the conjugacy class $\, {\mathcal S}\,$ in $\, G= SU(2)\,$ consisting of trace-free matrices considered in Section 2.)

Let $\, C_1= ... =C_n=C\,$ and let $\,{\mathcal C}= (C, ...,C)\,$ ($2n$ copies of $\, C$). Under these assumptions the braid group $\,{\mathcal B}_{2n}\,$ on $2n$ strands acts on $\, {\mathscr M}\,$ and on $\, K=r^{-1}(I)\,$ (see Section 4). According to Theorem 4.1, the action of  $\,{\mathcal B}_{2n}\,$ on $\, {\mathscr M}\,$ preserves the symplectic structure. The following corollary is then an immediate consequence of Proposition 5.1.

\begin{cor}
For every $\, \tau \in {\mathcal B}_{2n}\,$, $\, \tau (\Lambda )\,$ is a Lagrangian submanifold of $\, {\mathscr M}\,$.
\end{cor}  

Observe that both $\, \Lambda \,$ and  $\, \tau (\Lambda )\,$ are contained in $\, K=r^{-1}(I)\,$.

\medskip


Let us now consider the embedding of the braid group $\,{\mathcal B}_n\,$  
on $n$ strands into $\,  {\mathcal B}_{2n}\,$ extending a braid $\, \sigma \in {\mathcal B}_n\,$ to 
a braid $\, \widetilde{\sigma} \in  {\mathcal B}_{2n}\,$ which is trivial on {\it the  first $n$  strands}. Let us denote by $\, \Gamma _{\sigma}\,$ the submanifold  $\,  \widetilde{\sigma} (\Lambda )\,$.

\begin{cor}
For every $\, \sigma \in {\mathcal B}_n\,$, 
 $\, \Gamma _{\sigma} =  \widetilde{\sigma} (\Lambda ) \,$
 is a Lagrangian submanifold of $\, {\mathscr M}\,$.
\end{cor} 

Again, for every $\, \sigma \in {\mathcal B}_n\,$, the manifold  $\, \Gamma _{\sigma}\,$ is contained in the subspace $\, K=r^{-1}(I)\,$ of $\, {\mathscr M} \,$. 

The manifolds  $\, \Gamma _{\sigma}\,$ will play an important r{\^ o}le in later parts of the paper.
\bigskip

\section{The almost complex structure and the first Chern class}\label{almCh}

\bigskip

We shall now apply the constructions of Section 3 to the situation considered in Section 2. 

We assume that the Lie group $\, G\,$ is equal to $\, SU(2)\,$ and that the conjugacy classes $\, C_j\,$ are all equal to $\, {\mathcal S} = \{ A\in SU(2) \, | \, \text{trace} (A) = 0 \,\}$. 
We consider the constructions of Section 3 in that special case.

Let  $\, F\,$ be the free group on $\, 2n\,$ generators $\, z_1, ... , z_{2n}\,$, $\, n\ge 1\,$, and let $\, {\mathcal C}=(C_1, ... ,C_{2n}) = ( {\mathcal S}, ... , {\mathcal S} )\,$($2n$ elements). 

Consider the manifold  $\,\text{Hom} ( F, \, G )_{\mathcal C}\,$ of all group homomorphisms $\,\varphi :F \rightarrow G\,$ such that $\, \varphi (z_j)\in  {\mathcal S}, \,\, 1\le j \le m\,$, and the map $\, r:\text{Hom} ( F, \, G )_{\mathcal C}\rightarrow G , \,\,  r(\varphi )=  \varphi(z_1)\cdot ... \cdot \varphi (z_m)\,$.  In Section 2 the manifold   $\,\text{Hom} ( F, \, G )_{\mathcal C}\,$ was denoted by $\, P_{2n}\,$ and the map $\, r\,$ by $\, p_{2n}\,$.

Choose the open subset  $\, \widetilde{\mathscr{O}}\,$ in  $\, {\mathfrak g}= su(2)\,$  to be the connected component of $\,  \text{exp}^{-1}( SU(2)-\{-I\})\,$ containing $\, 0\,$ and let  $\, \mathscr{H} = r^{-1} ( \text{exp} (\widetilde{\mathscr{O}} ) ) = r^{-1} (  SU(2)-\{-I\})\,$.

The subspace  $\, \mathscr{H}\,$ is an open neighbourhood of  $\, K = r^{-1}(I)\,$ in $\, \text{Hom} ( F, \, G )_{\mathcal C}\,$.
Again, recall that  in  Section 2 the subset  $\, \mathscr{H}\,$ was denoted by $\, E\,$ and   the subset  $\, K\,$ by  $\, K_{2n}\,$.

According to Theorem 3.1 there exists an open neighbourhood  $\, \mathscr{M}\,$ of $\, K = K_{2n}\,$ in $\, \mathscr{H}=E \,$ such that the $2$-form $\, \omega _{\mathcal C}\,$ is symplectic on  $\, \mathscr{M}\,$. We can assume that  $\, {\mathscr M}\,$ has been choosen so that it is homotopy equivalent to $\, K_{2n}\,$.  From now on we  shall consider  $\, \mathscr{M}\,$ as a symplectic manifold equipped with the form   $\, \omega _{\mathcal C}\,$. 

Let $\, T {\mathscr M}\,$ be the tangent bundle of $\, \mathscr{M}\,$. We choose a complex structure on  $\, T {\mathscr M} \,$ compatible with  $\, \omega _{\mathcal C}\,$. Let $\, c_1( {\mathscr M}) =c_1( T {\mathscr M} )\,$ be the first Chern class of $\, T {\mathscr M} \,$,  $\, c_1( {\mathscr M})  \in \text{H}^2( \mathscr{M} ; \Z )\,$. 

\bigskip

The main aim of this Section is to  determine how $\, c_1 ( {\mathscr M})  \,$ evaluates on some elements of $\, \pi _2( K )\,$. As a consequence, we shall prove in the next Section that the symplectic manifold  $\, \mathscr{M}\,$ is {\it monotone}.

\bigskip
For every integer $\, k\,$ such that $ \, 1\le k\le 2n-1\,$ and for $\, \epsilon =\pm 1\,$ let us 
consider  elements of $\, \pi _2( K ) \,$ represented by the maps 
$\, \gamma _{k, \epsilon}: S^2 \rightarrow K \,$ given by 
\begin{equation}\label{almCh1}
\gamma _{k, \epsilon }(A)= (J, J, ... ,J, A, \epsilon A, J, ...,(-1)^n \epsilon J ), \qquad \text{for} \,\,\,  A\in S^2,
\end{equation}
where, as in Section 2, $\, S^2\,$ has been identified with  
$\,{\mathcal S}\,$,  $\,  J= \left(\begin{smallmatrix}i&0\\0&-i\end{smallmatrix}\right) \,$  and the first factor of  $\, A\,$ on the RHS is in the 
$k$-th place. In the case when $\, k=2n-1\,$ the sign $\, (-1)^n \epsilon \,$ is placed at the first factor of $\, J\,$. The maps $\,\gamma _{k, \epsilon } \,$ are embeddings of $\, S^2\,$  into $\, K\,$ and, hence, \, into $\,  \mathscr{M}\,$ and $\, P_{2n}\,$. 

Let $\ [ \gamma _{k, \epsilon } ]  \in \text{H}_2( \mathscr{M} ; \Z )\,$ be the homology classes represented by  the corresponding mappings.

\medskip

\begin{thm}\label{almChTh1}
For all integers $\, k\,$ such that $ \, 1\le k \le 2n-1\,$ and $\, \epsilon =\pm 1\,$, the evaluation of the first Chern class $\, c_1(\mathscr{M})\,$ on the homology classes
 $\ [ \gamma _{k, \epsilon } ]\,$ is equal to $\, 0\,$,
\begin{displaymath}
\langle\, c_1(\mathscr{M})  \, | \,  [ \gamma _{k, \epsilon } ]\, \rangle = 0. 
\end{displaymath}
\end{thm}

\medskip

Proof of Theorem \ref{almChTh1} is deferred to the Appendix B. (See Theorem B.1.)

\bigskip

We shall now,  for every integer $n\ge 2$, define a mapping $f_n:S^2 \rightarrow K_{2n}$.

We identify the $2$-dimensional sphere $S^2$ with a manifold  $\, \widetilde{{\mathcal S}}\,$ which is a union of two hemispheres $U_1$ and $U_2$ and of the cylinder $U_3$ joining the boundaries of the hemispheres as in Figure 6.1.

\vskip1truecm
$$\includegraphics[width=7cm, height=2cm]{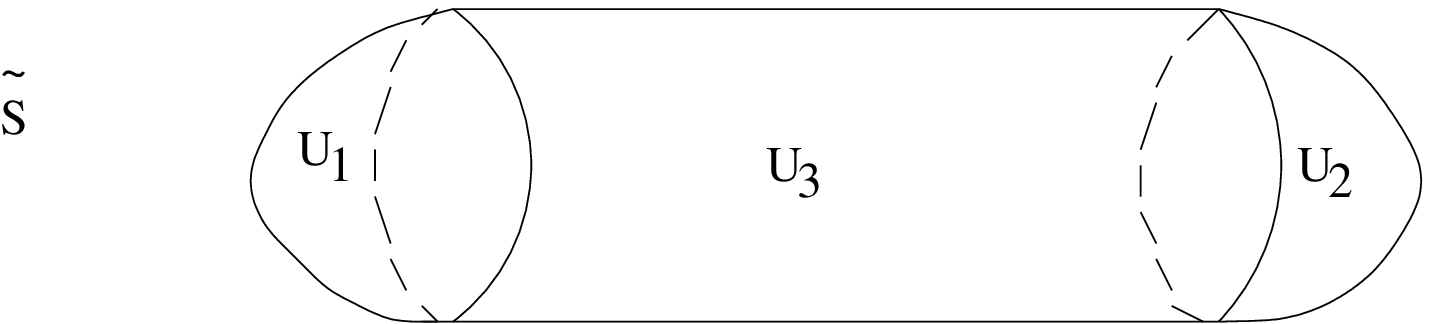}$$

\bigskip
\begin{center}
{Fig. \ref{almCh}.1}
\end{center}

\medskip

\noindent
We look upon the cylinder $\, U_3\,$ as the product $\, [0,\, \pi ] \times S^1\,$ of the interval $\, [0,\, \pi ]\,$ and of the circle $\, S^1=\R  / 2\pi \Z\,$. The identification of $\,S^2\,$ with $\,  \widetilde{\mathcal S}\,$ is such that the orientation of $\,S^2\,$ corresponds to the product orientation on the cylinder $\, U_3=[0, \, \pi ]\times S^1= [0, \, \pi ]\times \R/2\pi\Z\,$.

Moreover, we identify the hemispheres $\, U_1\,$ and $\, U_2\,$, each one separately, with subspaces of the conjugacy class $\, {\mathcal S}\,$. The hemisphere $\, U_1\,$ is identified with the subset of  $\,{\mathcal S} \,$ containing matrices $\, A=\left(\begin{smallmatrix}it&z\\-\overline{z}&-it\end{smallmatrix}\right)\,$ with $\, \text{Im}(z)\ge 0\,$ and the hemisphere $\, U_2\,$ with the subset containing matrices $\,A\,$ with   $\, \text{Im}(z)\le 0.\,$  

The points $\, (0,\theta _2)\,$ in $\, U_3\,$ are identified with matrices 
$\, \left( \begin{array}{rr}
i \cos \theta _2 & \sin \theta _2 \\
- \sin \theta _2 & -i \cos \theta _2
\end{array}\right)\,$ in $\, U_1\,$. The points  $\, (\pi ,\theta _2)\,$ in $\, U_3\,$ are identified with matrices 
$\, \left( \begin{array}{rr}
i \cos \theta _2 & \sin \theta _2 \\
- \sin \theta _2 & -i \cos \theta _2
\end{array}\right)\,$ in $\, U_2\,$.

We define the mapping $\, f_n: \widetilde{\mathcal S}\rightarrow K_{2n}\,$  by 
\begin{equation}\label{almCh2} 
f_n(A)=\begin{cases}(J, J, A, A, -J, J, ..., -J, J)    \qquad \qquad& \text{if}\qquad  A\in U_1,\\
(-J, J, A, -A, -J, J, ..., -J, J)  \qquad \qquad& \text{if}\qquad  A\in U_2.
\end{cases} 
\end{equation}
To define the restriction of $\, f_n\,$ to the cylinder $\, U_3=[0, \, \pi ]\times S^1= [0, \, \pi ]\times \R/2\pi\Z   \,$, let us denote by $\, A_{\theta}\,$ the matrix $\, \left( \begin{array}{rr}
i \cos \theta & \sin \theta \\
- \sin \theta & -i \cos \theta
\end{array}\right)\in  {\mathcal S}  \,$. Then, for $\, {\bf x}= [\theta_1, \, \theta_2 ] \in U_3= [0, \, \pi ]\times \R/2\pi\Z\,$, we set
\begin{equation}\label{almCh3} 
f_n({\bf x})=f_n([ \theta_1, \, \theta_2 ] ) =(A_{\theta_1}, J, A_{\theta_2}, A_{\theta_1+\theta_2}, -J, J, ..., -J, J)  \, .
\end{equation} 

As defined above, the mapping $\, f_n\,$ apriori takes its values in $\, P_{2n}\,$.
Let $\, p_j:P_{2n}\rightarrow {\mathcal S}\,$ be the projection on 
the $j$-th factor, $\, j=1, ... , 2n,$ and let 
$\, f_{n,j}=p_j\circ f_n:\widetilde{\mathcal S}\rightarrow {\mathcal S}\,$ be 
the coordinate maps of $\, f_n\,$. By the definition 
\begin{displaymath}
f_{n,j}(x) = (-1)^j J \qquad \text{for all} \,\, x\in  
\widetilde{{\mathcal S}} \,\, \text{and} \,\,\, 5\le j\le 2n. 
\end{displaymath}
Let $\, r:P_{2n} \rightarrow SU(2)\,$ be the mapping $\, r(A_1, ... ,  A_{2n}) = A_1\cdot ... \cdot A_{2n}\,$. One can check directly  
\begin{equation*}
f_{n,1}(x) \cdot f_{n,2}(x) \cdot f_{n,3}(x) \cdot f_{n,4}(x) = I \qquad 
\text{for all}\,\, x\in  \widetilde{{\mathcal S}}.
\end{equation*}
\noindent
It follows that $\, r(f(x)) = f_{n,1}(x) \cdot f_{n,2}(x) \cdot ... 
\cdot f_{n,2n}(x) = 
f_{n,1}(x) \cdot f_{n,2}(x) \cdot f_{n,3}(x) \cdot f_{n,4}(x)\cdot (-JJ)^{n-2} = 
I\cdot I^{n-2} = I\,$ for all $\,  x\in  \widetilde{{\mathcal S}}\,$. 
 Hence the image of $\, f_n\,$ lies in $\, K_{2n}= r^{-1}(I) \subset P_{2n}\,$ and we get the mapping $\, f_n:  \widetilde{{\mathcal S}}\rightarrow K_{2n}\,$.

The mapping $\, f_n\,$ is smooth away from the boundaries of the subsets $\, U_i, \,\, i=1,2,3\,$. One can see directly that $\, f_n\,$ can be smoothed in small collars of $\, \partial U_i\,$ without changing the image $\, f_n( \widetilde{{\mathcal S}})\,$. Whenever necessary we shall assume, without further mentioning it, that such a smoothing has been done.

We choose an orientation of  $\, \widetilde{{\mathcal S}}\,$ so that it coincides  on $\, U_3=[0, \pi ]\times \R /2\pi \Z\,$ with the product of the standard orientations of $\, [0, \, \pi ]\,$ and of $\, \R /2\pi \Z\,$. 
Let us then choose an orientation of $\, {\mathcal S}\,$ so that under 
the identification of $\, U_1\,$ and $\, U_2\,$ with subsets of  
$\, {\mathcal S}\,$ the orientations of  $\, {\mathcal S}\,$ and of  
$\, \widetilde{{\mathcal S}}\,$ coincide.

The degrees of the coordinate maps $\, f_{n,j}=p_j \circ f_n : 
\widetilde{{\mathcal S}}\rightarrow {\mathcal S}\,$ are 
\begin{equation}\label{almCh4}
\deg( f_{n,j}) = \deg (p_j \circ f_n ) = \begin{cases}
0& \qquad j\ne 3, \,\, 1\le j \le 2n,\\
1& \qquad j=3.
\end{cases}
\end{equation}
Indeed, the mappings $\, f_{n,j}: \widetilde{{\mathcal S}}\rightarrow 
{\mathcal S}\,$ are not surjective for $\, j\ne 3\,$, while 
$\, \deg( f_{n,3})=1\,$ follows immediately from the construction.

Let $\, f_{n\, *}: H_2( \widetilde{{\mathcal S}}, \, \Z )\rightarrow  H_2( K_{2n}, \, \Z )\,$ be the homomorphism induced by $\, f_n\,$ on the second homology groups and let  $\, [ \widetilde{{\mathcal S}}]\in H_2( \widetilde{{\mathcal S}}, \, \Z )\,$ be given by the orientation of $\,  \widetilde{{\mathcal S}}\,$.  Consider the element $\, f_{n\, *}  [ \widetilde{{\mathcal S}}] \in  H_2( K_{2n}, \, \Z )\,$.

\medskip

\begin{thm}\label{almChTh2}
For every $\, n\ge 2\,$
\begin{displaymath}
\langle\, c_1(\mathscr{M})  \, | \,\, f_{n\, *}  [ \widetilde{{\mathcal S}}] \, \, \rangle = -2. 
\end{displaymath}

\end{thm}

\medskip

\noindent
Proof of Theorem \ref{almChTh2} is deferred to  the Appendix C. (See Theorem C.1.)

\bigskip

We shall now consider the subgroup of $\, \pi _2(K_{2n})\,$ generated by the elements $\, [f_n]\,$ and  $\,[\gamma _{k,\epsilon}]\,$ with $\, 1\le k \le 2n-1, \,\, \epsilon = \pm 1.\,$  We shall show that this subgroup is of index $2$ in  $\, \pi _2(K_{2n})\,$.

 Let $\, j:K_{2n}\rightarrow P_{2n}\,$ be the inclusion.
The product structure of $\, P_{2n}=\prod\limits _{i=1}^{2n} {\mathcal S}\,$ gives an identification $\, \pi _2(P_{2n})=\bigoplus\limits _{i=1}^{2n} \Z\,$.  It follows directly from the definition of the mappings  $\, \gamma _{k,\epsilon}\,$ that under  the induced homomorphism  $\, j_*: \pi _2(K_{2n})\rightarrow \pi _2 (P_{2n})\,$ the homotopy class $\, [\gamma _{k,\epsilon}]\in  \pi _2(K_{2n})\,$ is mapped to $\, j_*([\gamma _{k,\epsilon}])= \oplus _{i=1}^{2n}  (j_*([\gamma _{k,\epsilon}]))_i\in  \pi _2(P_{2n})\,$ with
\begin{equation}\label{almCh5}
 (j_*([\gamma _{k,\epsilon}]))_i=\begin{cases}
1& \qquad \qquad \text{if} \quad i=k,\\
\epsilon &\qquad \qquad \text{if} \quad i=k+1,\\
0 &\qquad \qquad \text{otherwise}.
\end{cases}
\end{equation}

It follows from (\ref{almCh5}) that the homotopy class $\, [f_n]\in \pi _2 (K_{2n})\,$ is mapped by $\, j_*\,$ to  $\, j_*([f_n])= \oplus _{i=1}^{2n}  (j_*([f_n]))_i\in  \pi _2(P_{2n})\,$ with 
\begin{equation}\label{almCh6}
 (j_*([f_n]))_i=\begin{cases}
1& \qquad \qquad \text{if} \quad i=3,\\
0 &\qquad \qquad \text{otherwise}.
\end{cases}
\end{equation}

 Observe that if $\, h_n:K_{2n}\rightarrow K_{2n+2}\,$ is the mapping $\, h_n(A_1,..., A_{2n}) = (A_1,...,A_{2n},$ $ -J, J)\,$ as in (2.1) then 
\begin{equation}\label{almCh7}
h_n\circ f_n = f_{n+1}
\end{equation}
for $\, n\ge 2\,$.

Let us consider the element 
\begin{equation}\label{almCh8}
\alpha _n= 2[f_n]-[\gamma _{3,1}] -[\gamma _{3,-1}] \in \pi _2 (K_{2n})\, .
\end{equation}
Because of (\ref{almCh7}) and of the definition of $\, \gamma _{k,\epsilon }\,$ one has 
\begin{equation}\label{almCh9}
{h_n}_*(\alpha _n ) = \alpha _{n+1}\, .
\end{equation}
It follows from (\ref{almCh5}) and (\ref{almCh6}) that 
\begin{equation}\label{almCh10}
 j_*(\alpha _n ) = 0  \, .
\end{equation}
Let $\, \xi _n \,$ be a generator of $\, \text{Ker}\bigl(\, j_*: \pi _2(K_{2n})\rightarrow \pi _2 (P_{2n})\bigr)\cong \Z\,$. 

\begin{lma}\label{almChL3}
For $\, n\ge 2\,$,
\begin{displaymath}
\alpha _n = k \,\xi _n \quad with \quad k\in \{\pm 1, \pm 2\}
\end{displaymath}
in $\, \pi _2(K_{2n})\,$.
\end{lma}

\begin{proof}
According to Lemma 2.10, the homomorphism $\, {h_n}_*:\pi _2 (K_{2n}) \rightarrow \pi _2(K_{2n+2})\,$ maps  $\, \text{Ker}\bigl(\, j_*: \pi _2(K_{2n})\rightarrow \pi _2 (P_{2n})\bigr)\,$ isomorphically onto   $\, \text{Ker}\bigl(\, j_*: \pi _2(K_{2n+2})\rightarrow \pi _2 (P_{2n+2})\bigr)\,$. Thus $\,  {h_n}_*(\xi _n ) = \pm \xi _{n+1}\,$. Since, by (\ref{almCh9}), $\, {h_n}_*(\alpha _n ) = \alpha _{n+1}\,$, in order to prove the Lemma, it is enough to show that 
\begin{equation}\label{almCh11}
 \alpha _2 = k \, \xi _2\quad \text{with}\quad  k\in \{\pm 1, \pm 2\}  \, 
\end{equation} 
in $\, \pi _2 (K_4)\,$.

Let us now consider the quotient projection $\, q:K_4 \rightarrow K_4/G\,$.  
The space $\,  K_4/G\,$ is known as the ``pillow-case'', Figure \ref{almCh}.2, see \cite{L2}. 

It can be described as the space $\, [0, \, \pi ] \times S^1/\sim\,$,\,\,\, where $\, S^1=\R / 2\pi \Z\,$, 

\vskip1truecm

$$\includegraphics[width=7cm, height=4cm]{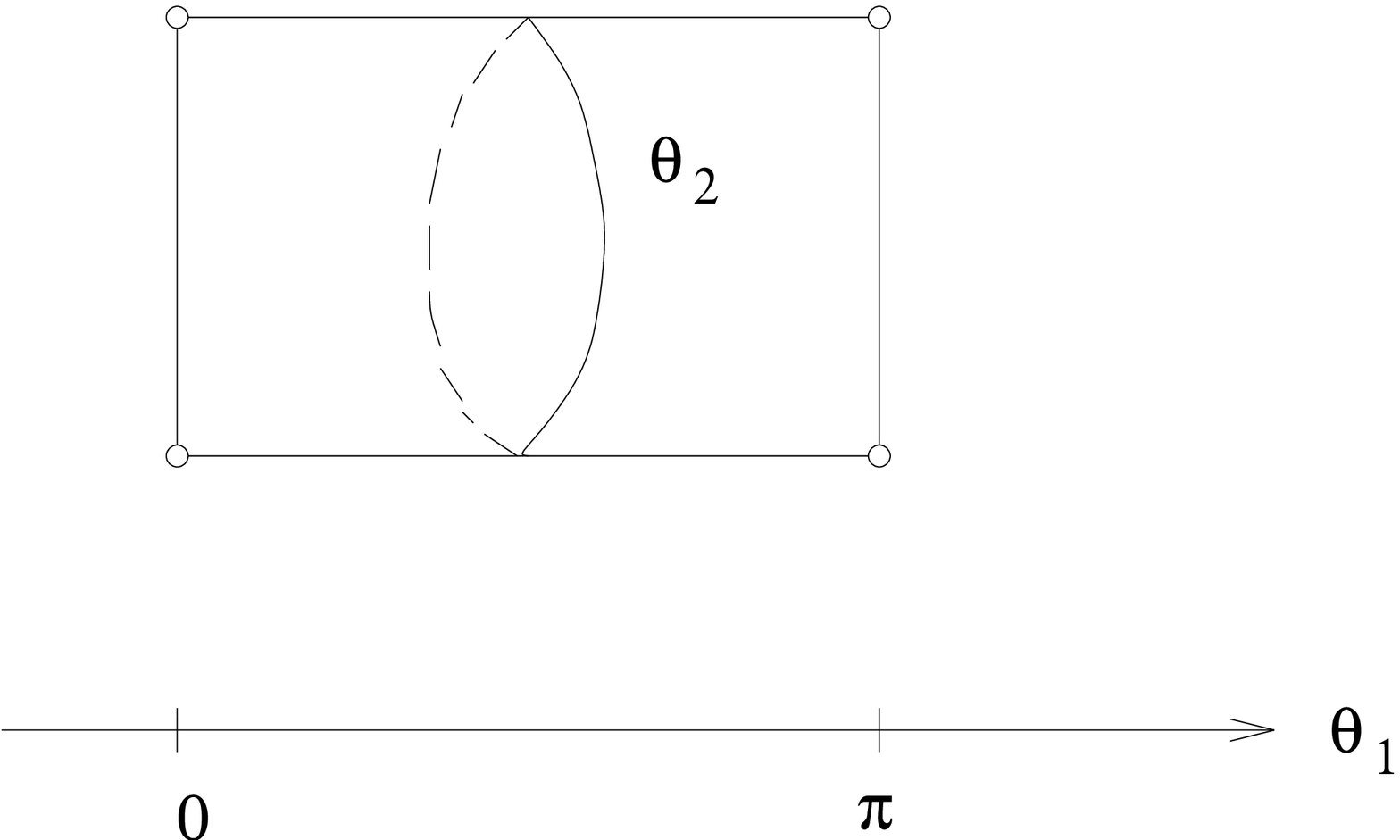}$$

\bigskip

\begin{center}
Fig. \ref{almCh}.2
\end{center}
with the identifications $\, (0, \theta _2) \sim (0, -\theta _2)\,$ and 
 $\, (\pi , \theta _2) \sim (\pi , -\theta _2)\,$ for $\, \theta _2\in \R / 2\pi \Z\,$. A point $\, [\, \theta _1, \, \theta _2 ]\,$ in  $\, K_4/G = [0, \, \pi ] \times S^1/\sim  \,$ is represented  by the point
\begin{multline*}
\left( \left(\begin{array}{rr}
i\cos \theta _1& \sin \theta _1\\
-\sin \theta _1&-i\cos \theta _1
\end{array}
\right), \,\, J, \left(\begin{array}{rr}
i\cos \theta _2& \sin \theta _2\\
-\sin \theta _2&-i\cos \theta _2
\end{array}
\right), \,\, \right.\\
\left.\left(\begin{array}{rr}
i\cos (\theta _1 + \theta _2)& \sin (\theta _1+\theta _2)\\
-\sin (\theta _1+\theta _2)&-i\cos (\theta _1 + \theta _2)
\end{array}
\right)\right)
\end{multline*} 
in $\, K_4\,$, (see \cite{L2}, p.347).

Therefore $\, K_4/G\,$ is homeomorphic to a $2$-dimensional sphere and 
we orient it by the product of the natural orientations of $\, [0, \, \pi ]\,$ and of $\, S^1=\R /2\pi \Z\,$. The choice of orientation gives us an isomorphism $\, \pi _2(K_4/G ) =\Z\,$.

Let us consider the induced homomorphism 
\begin{displaymath}
q_*:\pi _2(K_4)\rightarrow \pi _2( K_4/G)=\Z \,.
\end{displaymath} 
According to Lemma 2.10, $\, q_*([\,\gamma _{3,1}]) = q_*([\,\gamma _{3,-1}]) = 0\,$. Thus $\, q_*(\alpha _2 ) = 2\, q_*([f_2]) = 2[\, q\circ f_2]\,$.  

Now, the mapping $\, f_2:S^2 = \widetilde{\mathcal S}\rightarrow K_4\,$ was devised in such a way that the composition $\, q\circ f_2 : S^2\rightarrow K_4/G\cong S^2\,$ has degree $1$. Indeed, $\,\widetilde{\mathcal S}=U_1\cup U_2\cup U_3\,$ and $\, q\circ f_2\,$ maps the interior of $\, U_3\,$ homeomorphically and preserving the orientations onto the open subset of the ``pillow-case'' corresponding to $\, 0< \theta _1<\pi\,$, while it maps  $\, U_1\,$ into the closed subset  given by $\, \theta _1=0\,$  and maps $\, U_2\,$ into the closed subset  given by $\, \theta _1=\pi\,$. 

Therefore $\, [\, q\circ f_2]=1\,$  in $\, \pi _2( K_4/G) = \Z\,$ and 
\begin{equation}\label{almCh12}
q_*(\alpha _2) = 2 \, .
\end{equation}
As $\, \alpha _2\,$ belongs to the kernel of $\,  j_*: \pi _2(K_4)\rightarrow \pi _2 (P_4)\,$ one has 
\begin{displaymath}
\alpha _2 = k \, \xi _2
\end{displaymath}
for some $\, k\in \Z\,$. Applying $\,\, q_*\,$ we get
\begin{displaymath}
2 = q_*(\alpha _2) = k \cdot q_*(\xi _2)\, 
\end{displaymath} 
in $\, \pi _2(K_4/G) = \Z\,$.
Thus $\, k\,$ divides $\, 2\,$. That proves Lemma \ref{almChL3}.
\end{proof}

Let $\, {\mathcal G}\,$ be the subgroup of $\, \pi _2(K_{2n})\,$ generated by the homotopy classes $\, [f_n]\,$ and $\, [\gamma _{k,\epsilon }],\,\, k=1, ... , 2n-1, \, \, \epsilon = \pm 1\,$.

\begin{lma}\label{almChL4}
Let $\, n\ge 2\,$.  $\, {\mathcal G}\,$ is a subgroup of  index at most $2$ in  $\, \pi _2(K_{2n})\,$.
\end{lma} 

\begin{proof}
It follows directly from (\ref{almCh5}) and (\ref{almCh6}) that the homomorphism $\, j_*:\pi _2(K_{2n})\rightarrow \pi _2(P_{2n}) \,$ maps the subgroup  $\, {\mathcal G}\,$ surjectively onto $\, \pi _2(P_{2n}) \,$. Moreover, since $\, \alpha _n= 2[f_n]-[\gamma _{3,1}] -[\gamma _{3,-1}] \in  {\mathcal G}\,$, it follows from Lemma \ref{almChL3} that $\,  {\mathcal G}\,$ intersects the kernel of $\, j_*\,$ in a subgroup of index $1$ or $2$.  Hence  $\, {\mathcal G}\,$ is a subgroup of index $1$ or $2$  in  $\, \pi _2(K_{2n})\,$.
\end{proof}

\medskip

\section{The manifold $\, {\mathscr M}\,$ is monotone}\label{mono}

\bigskip

The aim of this Section is to prove that the symplectic manifold  $\, {\mathscr M}\subset P_{2n} \,$ is monotone.
Denote $\, K=K_{2n}\subset{\mathscr M}\subset P_{2n} \,$. 

Let $\, \omega_C\,$ be the symplectic $2$-form on $\, {\mathscr M}\,$.  According to Corollary \ref{sympl}.4, the restriction of the $2$-form  $\, \omega _C\,$ to the bundle $\, T {\mathscr M} |_{K_{2n}}\,$ is equal to the $2$-form $\, \omega _c\,$, the definition of which we shall now recall.

The $2$-form  $\, \omega _c\,$ is defined on the whole manifold  $\, P_{2n}\,$. Let $\, {\mathfrak g}\,$ be the Lie algebra of the Lie group $\, G=SU(2)\,$ equipped with an invariant positive definite inner product $\, \bullet\,$. 
Denote by $\, \omega\,$ the $\, {\mathfrak g}$-valued, left-invariant 1-form on $\, G\,$ which maps each tangent vector to the left-invariant vector field having that value. The corresponding right-invariant form will be denoted by $\, \bar{\omega}\,$. 

For any differential form $\, \alpha \,$ on $\, G\,$ denote by $\, \alpha _j\,$ the pulback of  $\, \alpha \,$ to $\, G\times G\,$ by the projection to the $\, j$th factor. Let 

\begin{displaymath}
\Omega = \tfrac 12 \,\, \omega _1 \bullet  \bar{\omega} _2 \,\,\, .
\end{displaymath}
This is a real-valued $2$-form on  $\, G\times G\,$. 
According to the convention 
\begin{equation}\label{mono5}
( \omega _1 \bullet  \bar{\omega} _2 )(U, V) =  \omega _1 (U) \bullet  \bar{\omega} _2 (V) -  \bar{\omega} _2 (U) \bullet  \omega _1 (V) \,\,\, .
\end{equation}

For any integer $\, j\,$ such that $\, 1\le j\le 2n-1,\,$ let $\, E_j:P_{2n}\rightarrow G\times G \,$ be the mapping $\, E_j(A_1, ... ,A_{2n})= (A_1 \cdot ... \cdot A_j, \,\,  A_{j+1})\,$ and let 
\begin{displaymath}
\omega _{[z_1\cdot ... \cdot z_j | z_{j+1}]} = E^*_j(\Omega) \,\, .
\end{displaymath} 
The $2$-form $\, \omega _c\,$ on $\, P_{2n}\,$ is defined by 
\begin{equation*}
\omega _c= -\sum\limits^{2n-1}_{j=1} \omega _{[z_1\cdot ... \cdot z_j |  z_{j+1}]} \,\, .
\end{equation*}

\medskip

Let us again consider the  mappings $\, \gamma_{k, \epsilon}:S^2\rightarrow K, \,\, 1\le k\le 2n-1,\,\, \epsilon =\pm 1, \,$ defined in (\ref{almCh1}).

\begin{lma}\label{monoL1}
For all $\, 1\le k\le 2n-1,\,\, \epsilon =\pm 1,   \,$ the pull-back of the form  $\, \omega_C\,$ by the mapping $\,\gamma_{k, \epsilon}\,$ vanishes,
\begin{displaymath}
 \gamma_{k, \epsilon}^*(\omega_C)=0\,.
\end{displaymath} 
\end{lma}

\begin{proof} The restriction of the symplectic $2$-form $\, \omega_C\,$ to the bundle $\, T{\mathscr M}\,\big|_K\,$ is equal to $\, \omega_c=-\sum\limits_{j=1}^{2n-1}\omega_{[z_1...z_j\,|\, z_{j+1}]}\,$ (Corollary \ref{sympl}.4). Since the image of  $\,\gamma_{k, \epsilon}\,$ lies in $\, K\,$, we have to show that $\,  \gamma_{k, \epsilon}^*(\omega_c)=0\,$.

Let us denote $\, p=A\in {\mathcal S}\,$ and $\, (p_1, ... , p_{2n})=(J, J, ... ,J, A, \epsilon A, J, ... , (-1)^n\epsilon J), \break \,\, p_j\in {\mathcal S} , \,\, p=p_k=\epsilon p_{k+1}.\,$ 
Every tangent vector $\,v\in  T_p{\mathcal S} \,$ is of the form $\, v= Xp-pX =\widetilde{X} p\,$  for some  $\, X\in {\mathfrak su}(2)\,$ and  $\, \widetilde{X}=X-Ad(p)(X)\in {\mathfrak su}(2).\,$

Let $\, X,Y\in  {\mathfrak su}(2).\,$ We have $\,\widetilde{X} p, \,\widetilde{Y}p\in T_p{\mathcal S} \,$. We shall now calculate the value $\, \omega_c( d\gamma_{k, \epsilon}(\widetilde{X} p), \,  d\gamma_{k, \epsilon}(\widetilde{Y} p))\,$ of the $2$-form $\, \gamma_{k, \epsilon}^*(\omega_c)\,$ on the pair $\, (\widetilde{X} p, \,\widetilde{Y}p)\,$.  

Let $\, \pi_1,\pi_2:G\times G\rightarrow G\,$ be the projections onto the first and the second factor.   
Since for $\,j\ne k\,$ at least one of the compositions $\,\pi_1\circ E_j\circ  \gamma_{k, \epsilon}\,$ or $\,\pi_2\circ E_j\circ  \gamma_{k, \epsilon}\,$ is a constant mapping into a point, it follows from the definition (\ref{mono5}) of $\, \Omega\,$ that 
\begin{displaymath}
 \gamma_{k, \epsilon}^*(\omega_{[z_1...z_j\,|\, z_{j+1}]})= \gamma_{k, \epsilon}^*(E_j^*(\Omega))=0    \qquad \text{for} \quad j\ne k.
\end{displaymath}
For $\, j=k\,$ we have $\, d(E_k\circ\gamma_{k, \epsilon})(\widetilde{X} p)=(p_1\cdot ... \cdot p_{k-1}\cdot \widetilde{X} p_k, \, \widetilde{X} p_{k+1})\in T_{E_k(\gamma_{k, \epsilon}(p))}(G\times G).\,$  It follows that 
\begin{equation}\label{mono1}
\begin{split}
 \omega_c( d&\gamma_{k, \epsilon}(\widetilde{X} p), \,  d\gamma_{k, \epsilon}(\widetilde{Y} p))=-\omega_{[z_1...z_k\,|\, z_{k+1}]}( d\gamma_{k, \epsilon}(\widetilde{X} p), \,  d\gamma_{k, \epsilon}(\widetilde{Y} p))=    \\
&=-(E_k^*(\Omega))( d\gamma_{k, \epsilon}(\widetilde{X} p), \,  d\gamma_{k, \epsilon}(\widetilde{Y} p))=    \\
&=-\Omega( d(E_k\circ\gamma_{k, \epsilon})(\widetilde{X} p), \,  d(E_k\circ\gamma_{k, \epsilon})(\widetilde{Y} p))=    \\
&=-\Omega\left((p_1\cdot ... \cdot p_{k-1}\cdot \widetilde{X} p_k, \, \widetilde{X} p_{k+1}),\,\, (p_1\cdot ... \cdot p_{k-1}\cdot \widetilde{Y} p_k, \, \widetilde{Y} p_{k+1})\right)=\\
&=-\Omega\left((p_1\cdot ... \cdot p_{k-1}\cdot p_k \cdot Ad(p_k^{-1})(\widetilde{X} ), \, \widetilde{X} p_{k+1}),\right.\\
& \left. \qquad \qquad  \qquad \qquad \qquad \,\, (p_1\cdot ... \cdot p_{k-1}\cdot p_k\cdot Ad(p_k^{-1})(\widetilde{Y}), \, \widetilde{Y} p_{k+1})\right)=\\
&=-\tfrac 12 \left(\omega\left( p_1\cdot ... \cdot p_{k-1}\cdot p_k \cdot Ad(p_k^{-1})(\widetilde{X} )\right)\bullet \overline{\omega}\left( \widetilde{Y} p_{k+1}\right) - \right. \\
&\left.  \qquad \qquad  \qquad   \,\,- \omega\left( p_1\cdot ... \cdot p_{k-1}\cdot p_k \cdot Ad(p_k^{-1})(\widetilde{Y}) \right)\bullet \overline{\omega}\left( \widetilde{X} p_{k+1}\right) \right)=\\
&=-\tfrac 12 \left(( Ad(p_k^{-1})(\widetilde{X} ))\bullet \widetilde{Y} - ( Ad(p_k^{-1})(\widetilde{Y}) )\bullet  \widetilde{X}  \right)=\\
&=-\tfrac 12 \left(\widetilde{X}\bullet (Ad(p_k)(\widetilde{Y}))-  ( Ad(p_k)(\widetilde{Y}) )\bullet  \widetilde{X}    \right)=\\
&=\,\,\, 0\, .
\end{split}
\end{equation} 
The next to the last equality follows from the invariance of the inner product $\bullet$ and from the fact that $\, Ad(p_k^{-1})=Ad(-p_k)=Ad(p_k)\,$.

The equality (\ref{mono1}) means that 
$\, \gamma_{k, \epsilon}^*(\omega_c)=0, \,$
as claimed.
\end{proof}

Let us denote by $\, \nu _k :P_{2n}\rightarrow P_{2n}\,$ the multiplication by $\, -I\,$ in the $k$-th factor, $\,\nu _k (A_1, ..., A_{k-1}, A_k, A_{k+1},..., A_{2n}) =(A_1, ..., A_{k-1},-A_k, A_{k+1},..., A_{2n}), \,\, 1\le k \le 2n\,$.   Let  $\, \nu _{k l} = \nu _k \circ  \nu _l : P_{2n}\rightarrow P_{2n}\,$, $\, 1\le k, l\le 2n\,$.   

\begin{lma}\label{monoL2}
For any pair of integers $\, k, l\,$ such that $\, 1\le k, l\le 2n\,$,
\begin{displaymath}
 \nu _{k l}^*\, (\omega _c) = \omega _c \,\, .
\end{displaymath}
\end{lma}

\begin{proof}
Let $\, \eta :G\rightarrow G\,$ be the multiplication by $\, -I\,$, $\, \eta (A)=-A\,$. Since $\, -I\,$ belongs to the center of $\, G,\,$ we have $\, \eta ^* (\omega)= \omega\,$ and $\,\eta ^* ({\bar\omega})= {\bar\omega}\,$. Therefore both maps $\, \eta \times id :G\times G \rightarrow G\times G\,$ and  
 $\, id \times \eta :G\times G \rightarrow G\times G\,$ satisfy
\begin{equation}\label{mono2}
 (\eta \times id)^*(\Omega )=\Omega \qquad \text{and}\qquad (  id \times \eta   )^*(\Omega )=\Omega \,\, .
\end{equation} 
Now, $\, E_j \circ  \nu _{k l} =(\psi _1 \times \psi _2)\circ E_j\,$, where $\, \psi _1,\psi _2 : G\rightarrow G\,$ are equal to $\, \eta\,$ or to $\, id\,$. It follows from (\ref{mono2}) that $\,(\psi _1 \times \psi _2)^*(\Omega)=\Omega\,$. Therefore
\begin{displaymath}
\begin{split}
 \nu _{k l}^*\, (\omega _{[z_1\cdot ... \cdot z_j |  z_{j+1}]})&= \nu _{k l}^*\, E^*_j(\Omega) = ( E_j \circ  \nu _{k l})^*(\Omega) =  ((\psi _1 \times \psi _2)\circ E_j)^*\,(\Omega)=  \\
&= E_J^*\,(\psi _1 \times \psi _2)^*(\Omega)=  E_J^*\,(\Omega)=\\
&= \omega _{[z_1\cdot ... \cdot z_j |  z_{j+1}]} 
\end{split}
\end{displaymath} 
for $\, 1\le j\le 2n-1.\,$ Consequently
\begin{displaymath}
 \nu _{k l}^*\, (\omega _c) = \nu _{k l}^*\, (  -\sum\limits^{2n-1}_{j=1} \omega _{[z_1\cdot ... \cdot z_j |  z_{j+1}]}  ) =  -\sum\limits^{2n-1}_{j=1} \omega _{[z_1\cdot ... \cdot z_j |  z_{j+1}]} = \omega _c \,\, . 
\end{displaymath}
 
\end{proof}

Observe that $\, \nu _{k l}(K_{2n})=K_{2n}\,$  for all pairs $\, (k,l)\,$. We denote by $\, {\mathcal V}\,$  the finite group of diffeomorphisms of $\, P_{2n}\,$ generated by $\,  \nu _{k l},\,\,  1\le k,l\le 2n\,$. $\, {\mathcal V}\,$  is isomorphic to $\,\Z _2^{2n-1}\,$.

Consider again the mapping $\, f_n:S^2\rightarrow K_{2n}\,$ defined in (\ref{almCh2}) - (\ref{almCh3}). Recall that for the purpose of this definition the $2$-sphere $\, S^2\,$ has been identified  with the manifold $\, \widetilde{\mathcal S} = U_1 \cup U_2 \cup U_3\,$ (see Figure \ref{mono}.1),

\medskip
 
$$\includegraphics[width=7cm, height=2cm]{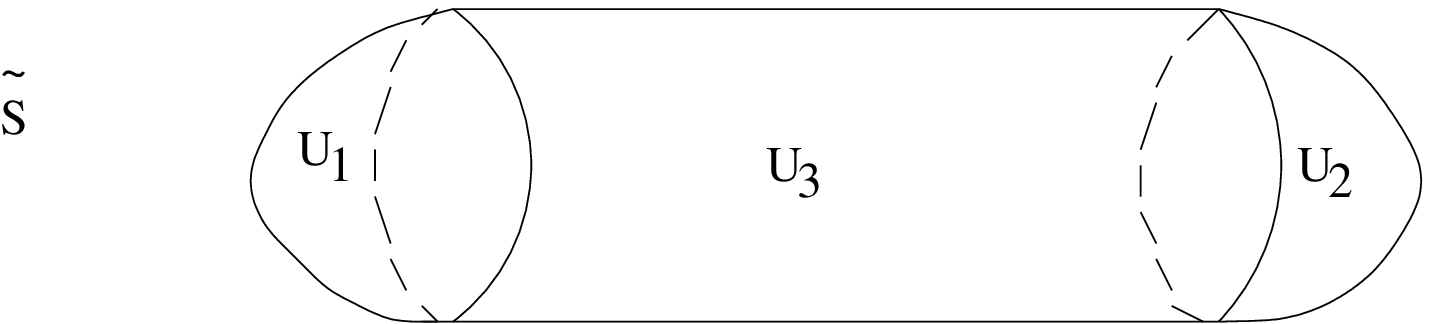}$$

\bigskip
\begin{center}
{Fig. \ref{mono}.1}
\end{center}
where $\, U_1, U_2, U_3\,$ are submanifolds (with boundary)  of $\, \widetilde{\mathcal S}\,$, intersecting only along boundaries. The identification of $\, S^2\,$ with  $\, \widetilde{\mathcal S}\,$ maps the standard orientation of  $\, S^2\,$ to the product orientation of $\, U_3= [\, 0, \, \pi ] \times ( \R /2\pi \Z )\,$.

Let $\, [S^2]\in H_2( S^2, \, \Z)\,$ be the fundamental class given by the standard orientation of $\, S^2\,$ and let $\, [\,\omega_C]\in H^2({\mathscr M},\, \R )\,$ be the deRham cohomology class of the closed form $\,\omega_C\,$. Let $\, f_{n*}:H_2(  S^2, \, \R)\rightarrow H_2({\mathscr M},\, \R )\,$ be the homomorphism induced by $\, f_n:S^2\rightarrow K_{2n}\subset {\mathscr M}\,$. 

The symplectic $2$-form $\, \omega_C\,$ depends on a choice of the invariant symmetric inner product $\, \bullet\,$ on the Lie algebra $\, {\mathfrak su}(2)\,$, which we have so far not specified in this paper. From now on we choose $\, \bullet \,$ to be given by $\, X\bullet Y=-\frac 12 \, \text{Trace}(XY)\,$ for $\, X,Y\in {\mathfrak su}(2)\,$.  

\begin{lma}\label{monoL3}
The evaluation $\, \langle \, [\, \omega_C] \, \, \big|\, \, f_{n*}( [S^2])   \, \rangle = -\pi ^2\,$.
\end{lma}

\begin{proof} It is clear from (\ref{almCh2}) and (\ref{almCh3}) that the mapping   $\, f_n:\widetilde{\mathcal S} \rightarrow {\mathscr M}\,$ is piece-wise smooth, being smooth on each of the pieces $\, U_1, U_2\,$ and $\, U_3\,$ separately. We have then 
\begin{equation}\label{mono3}
 \langle \, [\omega_C] \, \, \big|\, \, f_{n*}( [S^2])   \, \rangle  = \sum\limits_{j=1}^3 \,\, \int_{U_j} (f_n\,\big|_{U_j})^*( \omega_C) \, .
\end{equation}
As the image of $\, f_n\,$ lies in $\, K_{2n},\,$ we have  $\, (f_n\,\big|_{U_j})^*( \omega_C)=  (f_n\,\big|_{U_j})^*( \omega_c)\,$ for $\, j=1,2,3.\,$

Now, we have $\, (f_n\,\big|_{U_1})^*( \omega_C)= (f_n\,\big|_{U_2})^*( \omega_C)=0\,$. Indeed, this can be proven in exactly the same way as Lemma \ref{monoL1} or can be shown as follows.  We see from  (\ref{almCh2}) that $\, f_n\,\big|_{U_1}=\nu'\circ\gamma_{3,\, 1}\,\big|_{U_1}\,$ and  $\, f_n\,\big|_{U_2}=\nu''\circ\gamma_{3,\, -1}\,\big|_{U_2}\,$ for some diffeomorphisms $\, \nu', \nu''\in {\mathcal V}\,$, where $\, {\mathcal V}\,$ is the group of diffeomorphisms of $\, P_{2n}\,$ introduced above. According to Lemma \ref{monoL2} we have $\, \nu'(\omega_c)=  \nu''(\omega_c) = \omega_c\,$.  Thus $\, (f_n\,\big|_{U_1})^*( \omega_C)=  (f_n\,\big|_{U_1})^*( \omega_c)= (\gamma_{3,\, 1}\,\big|_{U_1})^*(\nu'^*(\omega_c))= (\gamma_{3,\, 1}\,\big|_{U_1})^*(\omega_c)= (\gamma_{3,\, 1}\,\big|_{U_1})^*(\omega_C)= 0,\,$ the last equality following from  Lemma \ref{monoL1}. The proof that $\, (f_n\,\big|_{U_2})^*( \omega_C)=0\,$ is similar.

We shall now calculate $\,(f_n\,\big|_{U_3})^*( \omega_C)\,$. Let $\, \partial_{\theta_1}= \frac {\partial}{\partial \theta_1}\,$ and  $\, \partial_{\theta_2}= \frac {\partial}{\partial \theta_2}\,$ be the coordinate vector fields given by the coordinates $\, (\theta_1, \, \theta_2)\,$ on $\,U_3= [0, \, \pi ]\times (\R/2\pi\Z)\,$. Let  $\, B=\left( \begin{array}{rr}
0&i\\
i&0
\end{array}\right)  \,$. 
Recall that $\, A_{\theta}=
\left( \begin{array}{rr}
i \cos \theta & \sin \theta \\
- \sin \theta & -i \cos \theta
\end{array}\right)  \,$ 
and observe that $\, \dfrac {d A_{\theta}}{d\theta}= 
B\cdot A_{\theta}\,$. 
\smallskip

Let us denote 
\begin{displaymath}
f_n([\theta_1, \, \theta_2 ])=(p_1, ... ,p_{2n})
\end{displaymath}
with $\, p_j= p_j(\theta_1, \, \theta_2 )\in {\mathcal S}, \,\, j=1, ..., 2n.\,$ Thus $\, p_1=A_{\theta_1}, \,\,\,  p_2=J, \,\,\, p_3=A_{\theta_2},\,\, \, p_4=A_{\theta_1+\theta_2}\,$ and $\, p_j=(-1)^jJ\,$ for $\, j\ge 5\,$. We have 
\begin{equation*}
\begin{split}
df_n(\partial_{\theta_1})&=(  B\cdot A_{\theta_1},\, 0\cdot J, \, 0\cdot  A_{\theta_2}, \,  B\cdot A_{\theta_1 + \theta_2}, \, 0\cdot p_5, ... , 0\cdot p_{2n})=\\
&=(  B\cdot p_1,\, 0\cdot p_2, \, 0\cdot p_3, \,  B\cdot p_4, \, 0\cdot p_5, ... , 0\cdot p_{2n})
\end{split}  
\end{equation*}
and
\begin{equation*}
\begin{split}
df_n(\partial_{\theta_2})&=(  0\cdot A_{\theta_1},\, 0\cdot J, \, B\cdot  A_{\theta_2}, \,  B\cdot A_{\theta_1 + \theta_2}, \, 0\cdot p_5, ... , 0\cdot p_{2n})=\\
&=(  0\cdot p_1,\, 0\cdot p_2, \, B\cdot p_3, \,  B\cdot p_4, \, 0\cdot p_5, ... , 0\cdot p_{2n}) .  
\end{split}
\end{equation*}
It follows that $\, \omega_{[z_1...z_j\,|\, z_{j+1}]}( df_n(\partial_{\theta_1}), \,\,  df_n(\partial_{\theta_2}))=0\,$ for $\, j=1\,$ and $\,j\ge 4\,$. Furthermore, 
\begin{equation}\label{mono4}
\begin{split}
\omega_C(&df_n(\partial_{\theta_1}), \, df_n(\partial_{\theta_2}))= \omega_c(df_n(\partial_{\theta_1}), \, df_n(\partial_{\theta_2}))=\\
&=-\sum\limits_{j=1}^{2n-1}\omega_{[z_1...z_j\,|\, z_{j+1}]}( df_n(\partial_{\theta_1}), \,\,  df_n(\partial_{\theta_2}))=\\
&=-\sum\limits_{j=2}^3\omega_{[z_1...z_j\,|\, z_{j+1}]}( df_n(\partial_{\theta_1}), \,\,  df_n(\partial_{\theta_2}))=\\
&=- \big(\Omega ((B\cdot p_1p_2, \, 0\cdot p_3), \,\, ( 0\cdot p_1p_2, \, B\cdot p_3))+ \\ 
& \qquad \qquad + \Omega ((B\cdot p_1p_2p_3, \, B\cdot p_4), \,\, ( p_1p_2\cdot B \cdot p_3, \, B\cdot p_4))   \big)=\\
&=-\tfrac 12 \big(\omega(B\cdot p_1p_2)\bullet \overline{\omega}(B\cdot p_3)\,\, +\\
&\qquad \quad + \omega(B\cdot p_1p_2p_3)\bullet \overline{\omega}(B\cdot p_4)-\omega(p_1p_2\cdot B\cdot p_3)\bullet \overline{\omega}(B\cdot p_4)\big)=\\
&=-\tfrac 12 \big( \big(Ad((p_1p_2)^{-1})(B)\big)\bullet B\,+ \big(Ad((p_1p_2p_3)^{-1})(B)-Ad(p_3^{-1})(B)\big)\bullet B \big).
\end{split}
\end{equation}
Since $\, A_{\theta} B A_{\theta}^{-1}=-B\,$ for all $\, \theta \in \R,\,$ we have $\, Ad(p_j)(B)=-B\,$ for $\, j=1,2,3\,$. From (\ref{mono4}) we get then that
\begin{equation*}
\begin{split}
\omega_C(&df_n(\partial_{\theta_1}), \, df_n(\partial_{\theta_2}))= 
-\tfrac 12 \big( B\bullet B\big)=\tfrac 14\, \text{Trace}(BB)= 
-\tfrac 12\,.
\end{split}
\end{equation*}
Therefore $\, (f_n\,\big|_{U_3})^*(\omega_C)=-\frac 12\,  d\theta_1\wedge d\theta_2\,$ on $\, U_3\,$ and, consequently,
\begin{displaymath}
\int_{U_3}(f_n\,\big|_{U_3})^*(\omega_C)= -\pi ^2\, .
\end{displaymath}

Lemma \ref{monoL3} follows now from (\ref{mono3}) due to the fact that  $\, (f_n\,\big|_{U_1})^*( \omega_C)= \break =(f_n\,\big|_{U_2})^*( \omega_C)=0\,$.
\end{proof}

\medskip

Let us recall that a symplectic manifold $\, (M, \omega)\,$ is called {\it monotone}\, if the cohomology class $\,[\omega]\in H^2(M, \, \R)\,$ is a positive multiple of the first Chern class $\, c_1(M)\,$.

\begin{thm}\label{monoTh4}
The symplectic manifold $\, {\mathscr M}\subset P_{2n}\,$ is monotone,
\begin{displaymath}
[\omega_C]= \frac {\pi ^2}2\,\, c_1({\mathscr M})\,.
\end{displaymath}
\end{thm} 

\begin{proof} Recall that $\, K_{2n}\,$ is a strong deformation retract of $\,{\mathscr M}\,$ and that it is $1$-connected (Proposition 2.9).  
According to Lemma \ref{almChL4}, the subgroup $\, {\mathcal G} \,$ of $\, \pi_2({\mathscr M})=H_2({\mathscr M}, \,\Z)\,$ generated by the homotopy classes $\, [\gamma_{j, \, \epsilon}], \,\, j=1,...,2n-1, \, \epsilon =\pm 1,\,$ and by $\, [f_n]\,$ is of index at most $2$ in $\, \pi_2({\mathscr M}).\,$  According to  Theorems \ref{almChTh1} and \ref{almChTh2},  one has 
\begin{equation*}
 \langle \, c_1({\mathscr M})\, \,\big| \,\, [f_n]\,\rangle = -2 \qquad \text{and} \qquad
\langle \, c_1({\mathscr M})\, \,\big|\,\, [\gamma_{j, \, \epsilon}]\,\rangle =0\end{equation*}  
for $\, j=1,...,2n-1, \, \epsilon =\pm 1.\,$ 

Now, Lemma \ref{monoL1} shows that 
$\, \langle \, [\omega_C]\, \,\big|\,\, [\gamma_{j, \, \epsilon}]\,\rangle =0 \,$ for $\, j=1,...,2n-1, \, \epsilon =\pm 1.\,$  Lemma  \ref{monoL3} shows that 
$\, \langle \, [\omega_C]\, \,\big|\,\, [f_n]\,\rangle = -\pi^2\,$. Thus the identity 
\begin{equation*}
[\omega_C]= \frac {\pi ^2}2\,\, c_1({\mathscr M})
\end{equation*}
holds when evaluated on the elements of the subgroup  $\, {\mathcal G} \,$ of index at most $2$ in $\,H_2({\mathscr M}, \,\Z)\,$. Therefore it holds when evaluated on all elements of  $\,H_2({\mathscr M}, \,\Z)\,$ and thus holds in  $\,H^2({\mathscr M}, \,\R)\,$.  
\end{proof}

\bigskip

\section{ Topological quandles and link invariants}\label{quandle}

\bigskip

Every conjugacy class in a topological group is a topological quandle.  In  \cite{R1} the second author introduced invariants of oriented links in $\, \R ^3\,$ derived from topological quandles. If $\, L\,$ is an oriented link in  $\, \R ^3\,$  and $\, Q\,$ is a topological quandle one defines an invariant $\, J_Q(L)\,$ which is a topological space well-defined up to a homeomorphism. If $\, Q\,$ is a conjugacy class $\, C\,$ in a Lie group G and  $\, L\,$ is the closure of a braid $\, \sigma\,$ on $\,n\,$ strands then the invariant space $\, J_Q(L)\,$ 
is, by definition, equal to the subspace of fixed points $\, Fix(\sigma )\,$ of the action of the braid $\,\sigma\,$ on the product $\, \prod\limits_{i=1}^n C\,$ of $n$ copies of $\, C\,$.

However, it is immediate that the space  $\, Fix(\sigma )\,$ and, hence, the invariant space  $\, J_Q(L)\,$ is also homeomorphic  
to the intersection $\, \Lambda \cap \Gamma _{\sigma}\,$ in $\, \prod\limits _{i=1}^{n} C \times  \prod\limits _{i=1}^{n} C^{-1}\,$ (see Section \ref{lagr}) of the ``twisted diagonal'' $\, \Lambda\,$ with the ``twisted graph'' $\, \Gamma _{\sigma}\,$ of $\, \sigma\,$. 

In the present paper we concentrate on the case when $\, Q=C\,$ is the conjugacy class $\, {\mathcal S}\,$ in the group $\, G=SU(2)\,$. Note that the quandle structure on  $\, {\mathcal S}\,$ given by the conjugation in  $\, SU(2)\,$ is the same as the quandle structure obtained from looking at  $\, {\mathcal S}\,$ as  the symmetric Riemannian  manifold isometric to the standard  2-dimensional sphere $\, S^2\,$.

\bigskip

\section{Examples}\label{ex}

\bigskip
Let $G$ be a Lie group. We shall now recall the  action of the braid group on $n$ strands ${\mathcal B}_n$ on 
$G^{\times n}$ used in \cite{R1}.

\bigskip

$$\includegraphics[width=7cm, height=3cm]{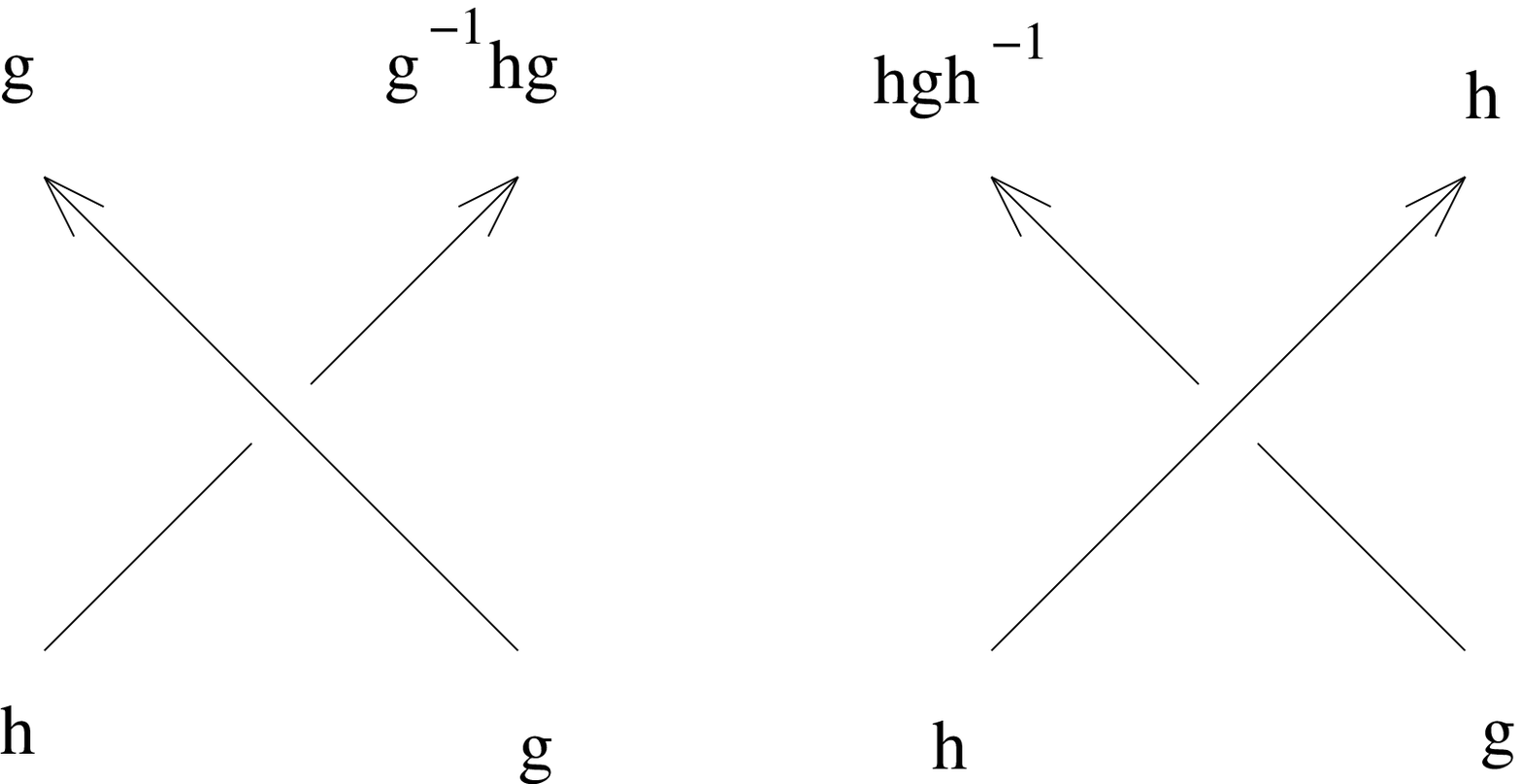}$$

\bigskip

\begin{center}
Fig. \ref{ex}.1
\end{center} .

Let $\sigma$ be an n-strand braid. To describe the 
image of $(g_1,...,g_n)$ under $\sigma$ we use braid diagrams with edges 
labelled by elements of $G$ as in Figure \ref{ex}.1. This figure 
describes how to label the two outgoing half-edges in a crossing given 
labels on the incoming ones. Label the bottom edges of $\sigma$ 
by $(g_1,...,g_n)$. Tracing through the braid, the rules 
in Figure \ref{ex}.1 then determine another $n$-tuple of 
$G$-elements $\sigma(g_1,...,g_n)$ on the top edges of the braid. This 
defines the braid action.

Let $L$ be the oriented link in $S^3$ obtained as the closure of the braid $\sigma$. From the above description it
is easy to see that the fixed point set 
$$J_G(L): = Fix(\sigma) = \{\rho\in G^{\times n} |\sigma(\rho) = \rho  \},$$ 
can be identified with the space $Hom(\pi_1(S^3 \setminus L), 
G)$ of homomorphisms 
of the fundamental group of the complement of $L$ into $G$. 
(This observation implies, in particular, that the space $J_G(L)$ is a link invariant.)
Indeed, a necessary and sufficient condition for the labelling of the diagram of $\sigma$
by group elements satisying the rules of Figure \ref{ex}.1 to extend over the 
``closing arcs'' of the closure of $\sigma$ is that $\sigma(\rho) = \rho$. Then the labelling 
describes a homomorphism of the Wirtinger presentation of 
the link complement into $G$ (see \cite{R1}, Lemma 4.6). Let us call such a labelling admissible.

All of the above arguments continue to hold if we replace 
the group $G$ by some specified conjugacy class $C$ in $G$. 
In this case $J_C(L):= Fix(\sigma) = \{\rho\in C^{\times n} |\sigma(\rho) = \rho  \} $ can be identified with $Hom(\pi_1(S^3 \setminus L), 
G)_C$ - the space of representations which send the positively oriented link meridians of $L$  
into the conjugacy class $C$. For the rest of this chapter, 
unless explicitly stated otherwise, we will be concerned 
with the case $G = SU(2)$ and $C = \{A \in SU(2) | tr(A) = 0\}$. 

Before going on to the examples, let us recall the behaviour
of the conjugation action of $C$ on itself. Viewing 
$SU(2)$ as the unit $3$-sphere in the space of quaternions, 
$C$ is the unit $2$-sphere in the space of pure quaternions 
$\mathbb{R}^3$. If $g$ and $h$ are in $C$, the element 
$g^{-1}hg$ is obtained by rotating $h$ by the angle $\pi$ 
around $g$. 

Equivalently, $\,g$ acts on $h$ by the involutive isometry 
$i_g$ (reflection in $g$) associated to $g$ if we regard $C 
\cong S²$ as a symmetric Riemannian manifold.

\subsection{$(2,n)$-torus links}

A $(2,n)$-torus link $T = T_{2,n}$ is the closure of a 2-strand braid
$\sigma = \sigma_1^{n}$ (where $\sigma_1$ is the elementary braid on $2$ strands generating  the braid group ${\mathcal B}_2$ ). If $n$ is odd, $T_{2,n}$ is a knot, if $n$ is 
even it is a link of two components. The space $J_{C}(T_{2,n})$ is 
the space of fixed points of the action of the braid $\sigma$ on $C \times C$.

Let us compute $\sigma(a,b)$ for two arbitrary points $a$,$b$ on the 
$2$-sphere $C$. Choose a great circle $S$ containing $a$ and $b$. If $a$ and 
$b$ are neither equal nor antipodal, there is a unique choice. Recall that 
the isometry $a \mapsto i_b(a)$ associated to $b$ is given by 
reflecting $a$ in $b$ and that for the braid generator $\sigma _1$ 
we have $$\sigma_1(a,b) = (b, i_b(a)).$$ We parametrise $S$ starting at $b$, 
with $a$ at the angle $\alpha$. Thus $(a,b)=(\alpha , 0)$. Then it is easy to see that 
$$\sigma(\alpha,0) = (-(n-1) \alpha , - n \alpha).$$ Hence, $(a,b)$ 
is a fix-point of $\sigma$ iff $n \alpha \cong 0 $ (mod $2\pi$), 
that is iff $\alpha = \pm \frac{2\pi j}{n}$, $0 \leq j \leq \frac{n-1}{2}$.   

Consider the projection on the second factor $p: J_C(T_{2,n}) 
\rightarrow C$. For $n$ odd, the fibre over $b$ is the disjoint 
union of $b$ and $\frac{n-1}{2}$ circles, each in a plane perpendicular 
to $b$, and consisting of the $a$ which have distance $\frac{2\pi j}{n}$ 
to $b$. It is clear from this description that $J_C(T_{2,n})$ 
is the disjoint union of an $S^2$ and $\frac{n-1}{2}$ copies 
of the unit tangent bundle to $S^2$ (which is the same as $\mathbb{R}P^3$). 

If $n$ is even, the fibre is the disjoint union of two points 
(corresponding to the angles $0$ and $\pi$, respectively) and 
$\frac{n-2}{2}$ circles. Thus in this case, $J_C(T_{2,n})$
is diffeomorphic to the union of two copies of $S^2$ and 
$\frac{n-2}{2}$ copies of $\mathbb{R}P^3$ unit tangent bundle 
of $S^2$. 

In other words, for $(2,n)$-torus links we have

\begin{equation}\label{ex1}
J_C(T_{2,n}) = S^2 \cup S^2 \cup \bigcup_{j=1}^{\frac{n-1}{2}} \mathbb{R}P^3,
\end{equation} 
where one of the spheres appears only if $n$ is even.

The Khovanov homology $Kh^{i,j}(T_{2,n})$ was computed already in the 
seminal paper \cite{K}. 
The non-zero groups are as follows.
First, there is a pair of groups each isomorphic to the integers. 
\begin{align*}
&Kh^{0,-n}(T_{2,n}) = \mathbb{Z} \\
&Kh^{0,2-n}(T_{2,n}) = \mathbb{Z} 
\end{align*} 
Then there is a number of ``triples'' of groups, one 
triple for each $1 \leq k \leq \frac{n-1}{2}$:
\begin{align*}
&Kh^{-2k-1, -4k-2-n}(T_{2,n}) = \mathbb{Z} \\
&Kh^{-2k, -4k-n}(T_{2,n}) = \mathbb{Z}_2 \\
&Kh^{-2k,-4k+2-n}(T_{2,n}) = \mathbb{Z}
\end{align*}
If $n$ is odd and these are placed on the $(i,2j+1)$-lattice 
(if $n$ is even, the $(i, 2j)$-lattice), they form a chess-board 
``knight move'' pattern (for more on this see below, Remark \ref{knight}).  

Finally, for even $n$, there is another pair 

\begin{align*}
&Kh^{-n ,-3n}(T_{2,n}) = \mathbb{Z} \\
&Kh^{-n,-2-3n}(T_{2,n}) = \mathbb{Z}
\end{align*}

If we collapse the bigrading in these formulas to a single grading 
$k = i - j$, the first and last groups 
become copies of the cohomology of the 
2-sphere $S^2$, with grading shifted by $n-2$ and $2n-2$ respectively. 
Similarily, each ``triple'' as above becomes a copy of the cohomology 
of real projective space $\mathbb{R}P^3$, shifted by $2k - 2 + n$. 
Calling the resulting groups $Kh^*$, we obtain 
$$
Kh^*(T_{2,n}) \cong H^*(S^2)\{n-2\} \oplus H^*(S^2)\{2n-2\} \oplus 
\bigoplus_{k=1}^{\frac{n-1}{2}} H^*(\mathbb{R}P^3) \{2k - 2 + n\},
$$
where the second term appears only when $n$ is even. We note that 
exactly the same spaces  as in (\ref{ex1}) appear here, only shifted by certain 
integers. All the cohomology here and below is with coefficients 
in the integers $\mathbb{Z}$. 

\subsection{The figure eight knot $4_1$} The invariant space for 
the figure eight knot was computed in \cite{R1}. The result is 
$$J_{C}(4_1) = S^2 \cup \mathbb{R}P^3 \cup 
\mathbb{R}P^3.$$ Consulting JavaKh, \cite{KT}, 
collapsed Khovanov homology for the figure eight knot is 
$$Kh^*(4_1) = H^*(S^2)\{-1\} \oplus H^*(\mathbb{R}P^3)\{-3\} 
\oplus H^*(\mathbb{R}P^3)\{0\}.$$

\subsection{The knot $5_2$} For the five crossing knot which is 
not a torus knot, $$J_C(5_2) = S^2 \cup \mathbb{R}P^3 \cup \mathbb{R}P^3 
\cup \mathbb{R}P^3.$$ 
This can be easily obtained by using the techniques of \cite{R1}, Section 5. The details are left to the reader. JavaKh tells us that Khovanov homology of the knot  $5_2$ is

$$ H^*(S^2)\{1\} \oplus H^*(\mathbb{R}P^3)\{5\} 
\oplus H^*(\mathbb{R}P^3)\{7\} \oplus H^*(\mathbb{R}P^3)\{11\} 
$$

\subsection{The knot $6_1$} This knot has 
$$J_C(6_1) = S^2 \cup \mathbb{R}P^3 \cup \mathbb{R}P^3 
\cup \mathbb{R}P^3 \cup \mathbb{R}P^3$$ and Khovanov homology, again 
by JavaKh, 
$$ H^*(S^2)\{- 1\} \oplus H^*(\mathbb{R}P^3)\{-3\} 
\oplus H^*(\mathbb{R}P^3)\{-1\} \oplus H^*(\mathbb{R}P^3)\{0\}
\oplus H^*(\mathbb{R}P^3)\{2\}.$$

\begin{rem} 
\label{rem942} The same pattern as above recurs for all seven crossing knots
and for at least some eight-crossing knots. The invariant spaces of those knots can be computed by techniques of \cite{R1}, Section 5. 
The limits are set by our 
capacity for solving the equations for the invariant spaces, which become 
very hard to solve when the number of labels used is more than three. We 
have been able to calculate one example which does not produce the 
Khovanov homology groups, namely the knot $9_{42}$, in Rolfsen's table.
The calculation is presented in Section \ref{counter}.
\end{rem}

\begin{rem} 
\label{knight}
In the bigraded Khovanov homology of knots, there is the 
well-known pattern of ``knight moves''; two $\mathbb{Z}$:s and one 
$\mathbb{Z}_2$ arranged in the chess board constellation of a knight 
move. In our examples, these are the source of the $H^*(\mathbb{R}P^3)$:s after 
collapse of the bigrading. There are also two additional copies of $\mathbb{Z}$ sitting 
in $i$-degree $0$. They correspond to  $H^*(S^2)$ after collapse. 
It has been conjectured by Alexander Shumakovitch \cite{Sh1}, Conjecture 3, that for alternating knots all the homology groups are arranged in
knight moves and a "sphere pair". (The same is conjectured for
non-split alternating links, and even more generally for H-thin links,
except that there are more sphere pairs. Compare for example
$T_{2,n}$-links with even $n$ above.) With rational coefficients this
conjecture (originally due to Bar-Natan, Garoufalidis and Khovanov \cite{BN}, \cite{G2}) is a
theorem due to Eun Soo Lee \cite{L3}. 
\end{rem}

\subsection{Reduced homology}

Consider $J_C(L) = Fix(\sigma)$ as the set of admissible labellings 
of the link diagram of $L$ by elements from $C$ as above. Then the 
choice of an arc in the link diagram determines a map from $J_C(L)$ 
to $C$. This map is a fibration. One may consider the 
cohomology of its fibre for a fixed point in $C$. In all the 
examples above (and all the ones we have computed), 
the fibre projection is the identity on the sphere components and the 
natural projection on the $\mathbb{R}P^3$:s. We get the following set of 
examples in which we see that this gives the ``reduced'' \cite{K2} (collapsed) 
Khovanov homology. 

\subsection{$T_{2,n}$-torus links - reduced case} The fibre projection to $C$
is the identity on the $S^2$-components and the natural projection on 
the $\mathbb{R}P^3$:s. Thus the fibre is $$\{*\} \cup \{*\} \cup 
\bigcup_{i=0}^{\frac{n-1}{2}} S^1\, \, , $$ where 
the second term appears only if $n$ is even. Reduced collapsed 
Khovanov homology is 

$$H^*(*)\{n-1\} \oplus H^*(*)\{2n-1\} \oplus \bigoplus_{i=0}^{\frac{n-1}{2}} 
H^*(S^1)\{2i - 1 + n\}.$$

\subsection{Two-bridge links} All of the above examples are special cases of 
two-bridge links. For two-bridge knots and with rational coefficients the observation (\ref{intro1}) of the Introduction holds true. Indeed, we have

\begin{thm}\label{exTh3} Let $\, K\,$ be a two-bridge knot and let $\, \{ J^r_C(K) \}_{r\in R}\,$ be the set of connected components of the representation variety   $J_C(K)$. Then there are integers $\, N_r=N_r(K), \,\, r\in R,\,$ such that  
\begin{equation*}
Kh^*(K; \mathbb{Q}) = \bigoplus_{r\in R} H^*(J_C^r(K);\mathbb{Q})\{N_r\}.
\end{equation*}
\end{thm} 

\begin{proof} Let $G = SU(2)$. A two-bridge knot group $\pi = \pi_1(S^3 \setminus K)$ can be 
presented by two meridional generators $a$,$b$ and one relation. For any representation 
in $J_C(K) = Hom(\pi, G)_C$, the images of the generators will lie on some 
great circle, and the single equation in the angle $\alpha$ between them 
can be solved just as in the case of $(2,n)$-torus links. 
The result is that $J_C(K)$ is the disjoint union of a copy of $S^2$ 
and a number of $\mathbb{R}P^3$:s. These connected components are exactly 
the $G$-conjugacy classes of these representations. 

Let $$S^1_A = \{a + bi: a^2 + b^2 = 1\} \subset G$$ and 
$$S^1_B = \{cj + dk: c^2 + d^2 =1\} \subset C \subset G.$$ 
The disjoint union of these circles in $G = SU(2)$ is called 
the binary dihedral group $N$, see \cite{Kl1}. 

The inclusion of $(N, S^1_B)$ into $(G,C)$ induces an obvious map
$$
Hom(\pi,N)_{S^1_B}/N \xrightarrow{\cong} Hom(\pi,G)_C/G.
$$

This map is surjective since any representation in $Hom(\pi,G)_C$ can be $G$-conjuga- ted into $N$. 
It is also injective. Indeed, let $\psi$ and $\phi$ be two representations in $Hom(\pi,N)_{S^1_B}$ 
and suppose they are conjugate modulo $G$. This amounts to saying that the pair 
$\phi(a), \phi(b)$ of points in $S^1_B$ can be rotated and reflected into $\psi(a), \psi(b)$. 
But in fact any rotation or reflection of $S^1_B$ can already be effected by conjugation by 
elements of $N$. 

According to \cite{Kl1}, Theorem 10, the number $n$ of non-abelian conjugacy classes 
of binary dihedral representations is given by the formula $$n = \frac{|\Delta_{-1}(K)| - 1}{2}$$ 
where $\Delta_{-1}(K)$ is the Alexander polynomial evaluated at $-1$. Furthermore, all the 
non-abelian representations are contained in $Hom(\pi, N)_{S^1_B}$.  

A representation in $Hom(\pi,N)_{S^1_B}$ is non-abelian if and only if 
it maps the two generators to distinct and non-antipodal elements, 
that is, if it lies in one of the $\mathbb{R}P^3$ components of 
$Hom(\pi,G)_C$. Thus $n$ is also 
the number of $\mathbb{R}P^3$:s. Taking into account also the one component 
of abelian representations (the sphere), we get $$dim(H^*(J_C(K);\mathbb{Q})) 
= 2(n + 1),$$ since every connected component contributes two dimensions 
to cohomology. 

Now, $\Delta_{-1}(K)$ is known as the determinant of $K$, and it 
is a classical theorem (see e.g. \cite{BZ}) that for alternating knots $L$ 
the determinant is the number of spanning trees in the Tait 
(checkerboard) graph of $L$. Furthermore, a theorem of Champanerkar 
and Kofman \cite{CK} says that the rank of reduced Khovanov homology 
for an alternating knot is exactly this number of spanning trees. 
Finally, the ranks of unreduced and reduced Khovanov homologies 
are related by 
$$dim(Kh(K);\mathbb{Q}) = dim(Kh_{red}(K);\mathbb{Q}) + 1.$$ 
It follows that 
$$ dim(H^*(J_C(K); \mathbb{Q})) = 2(n+1) = |\Delta_{-1}(K)| + 1 
= dim(Kh(K;\mathbb{Q})).$$ 
\end{proof}

\begin{rem} Theorem \ref{exTh3} is probably true also with integral coefficients. 
For example, it would follow from the conjecture mentioned in Remark \ref{knight} 
since all two-bridge knots are alternating.
\end{rem}

\subsection{Connected sum of knots - The square knot} It is a result of the second author (unpublished) that if $K_1$ and $K_2$ are knots  and we form their connected sum $K_1 \# K_2$ then the invariant space $J_C(K_1 \# K_2)= J_C(K_1) \times _C J_C(K_2)$ is the pullback of the fibrations $J_C(K_1)\rightarrow C$ and  $J_C(K_2)\rightarrow C$ over the base arcs. (This result holds for any conjugation quandle $Q$.)

 In the case of the 
connected sum of the trefoil knot with its mirror image, 
a.k.a. the square knot $Sq$, we get   

$$
J_C(Sq) = S^2 \cup \mathbb{R}P^3 \cup \mathbb{R}P^3 \cup \mathbb{R}P^3 
\times S^1.
$$

Again, a consultation of JavaKh gives 

$$
Kh^*(Sq) = H^*(S^2)\{-1\} \oplus H^*(\mathbb{R}P^3)\{-4\} 
\oplus H^*(\mathbb{R}P^3)\{1\} \oplus H^*(\mathbb{R}P^3 \times S^1)\{-2\}. 
$$

\vskip1truecm

\section{An example with non-isomorphic cohomology groups}\label{counter} 

\bigskip
In this section we calculate, as promised in Remark \ref{rem942}, the 
representation variety $J_C(L)$ of the knot $L=9_{42}$ and show that its cohomology 
groups do not coincide with the Khovanov homology groups. Cognoscenti 
will note that the knot in question is non-alternating and that it 
has all the possible ``card symbols'' in Bar-Natan's table \cite{BN}.

\begin{prop}\label{counterP1}
$\, J_C(9_{42}) = S^2 \cup \bigcup_{i=1}^{7} \mathbb{R}P^3\,$.

\end{prop}

\begin{proof}
We consider the knot diagram of $9_{42}$ pictured in Figure 
\ref{counter}.1 and regard the representation variety $J_C(9_{42})$ 
as the space of admissible colorings of the diagram in that figure, \cite{R1}. 
As we can see we need 
to use three colors $a, b, c$. For simplicity of notation let us use the functional notation for conjugation, so that $$ a(b) = a^{-1}ba\,\,\,  ( = aba^{-1}).$$

\bigskip

$$\includegraphics[width=10cm, height=7cm]{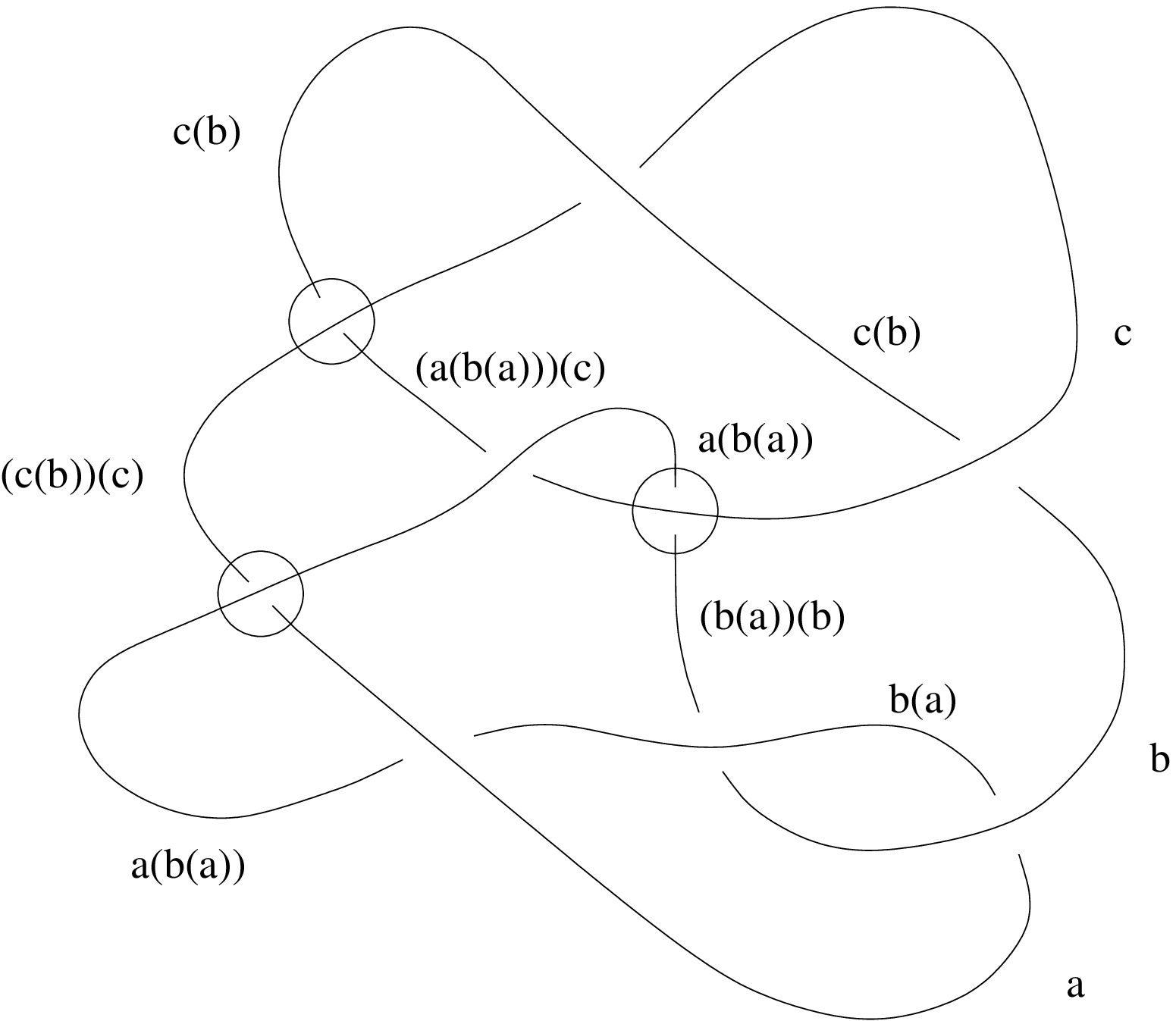}$$

\bigskip

\begin{center}
Fig. \ref{counter}.1
\end{center} .

The three circles in the figure encircle the crossings which give the three equations to be satisfied by the colors $\, a, b, c\,$.

It is easy to verify that the equations become 
\begin{displaymath}\tag{1}
c(b(a(b))) = a(b(a)),
\end{displaymath}
\begin{displaymath}\tag{2} 
a(b(a(b(a)))) = c(b(c))
\end{displaymath}
 and 
\begin{displaymath}\tag{3}
c(b(c(b)))= a(b(a(b(a(c))))).
\end{displaymath}

\subsection{Case 1: $a = \pm b$} If $a$ and $b$ coincide, then the first 
equation says that that $c(a) = a$. Thus, $c = \pm a$. The second equation, 
however, then becomes $a = c$, so only the solution $a = b = c$ survives. 
It trivially satisfies the last equation, too. The solutions of this form 
constitute the diagonal $S^2$ in $S^2 \times S^2 \times S^2$. If, instead, 
$a$ and $b$ are antipodal, then equation 1 says that $c$ should take the 
antipode of $a$ to $a$. Thus $c$ should be on the equator midway between 
$a$ and $b=-a$. Then $b(c) = -c$, so that equation 2 becomes $a = -c$, 
which is a contradiction. Thus, there are no solutions of the form $(a,-a,c)$.
This case contributes one copy of $S^2$ to $J_C(9_{42})$. 

\subsection{Case 2: $a \neq \pm b$} If $a$ and $b$ are not parallel then 
they lie on a unique great circle $\Gamma \subset S^2$. We will parametrize 
$\Gamma$ by arc-length starting at $a$ with $b$ corresponding to 
angle $\theta$. Consider the 
expressions $\alpha = a(b(a))$ and $\beta = b(a(b))$. They are both 
certain points on $\Gamma$. Equation 1 says that reflection in $c$ 
takes $\alpha$ to $\beta$. We will consider different subcases, based 
on the relative position of $\alpha$ and $\beta$. 

\subsubsection{Case 2a: $\alpha = \beta$} In this case, Equation 1 says 
that $c$ fixes $\alpha$, that is, $c = \pm \alpha$. In terms of the angle 
$\theta$, it is easy to see that $\alpha = \beta$ means $5\theta = 0$ $mod 
$ $2\pi$, so that $a$ and $b$ are vertices of a regular pentagon 
$P$ on $\Gamma$. Clearly $\alpha$ is also on the pentagon. 
Since $c(\alpha) = \alpha$, $c$ must be a point on the regular 
decagon $D$ containing $P$. (So far there
are two antipodal possibilities for c.) Next, note that the left hand side of 
Equation 2 collapses under the given assumption: 
$$ a(b(a(b(a)))) = a(b(b(a(b)))) = a(a(b)) = b. $$ Hence Equation 2 becomes
$$ b = c(b(c)),$$ which has solutions $3\vartheta = 0$ $ mod $ $2\pi$ in terms 
of the angle $\vartheta$ between $b$ and $c$. If $\vartheta \neq 0$ this is 
clearly incompatible with $c$ lying on $D$. If $\vartheta = 0$, 
so that $b = c$, then it is easy to see that Equation 1 reduces 
to $a = b$, so that we are back in Case 1. Hence there are no 
new solutions from this case.

\subsubsection{Case 2b: $\alpha = - \beta$} In this case, Equation 1 says that 
$c$ takes $\alpha$ to its antipode. Hence $c$ will lie on the great equatorial 
circle $\Delta$ all of whose points are equidistant from $\alpha$ and $\beta$. 
In terms of the angle $\theta$ between $a$ and $b$, Equation 1 says 
$5\theta + \pi = 0$ mod $2\pi$. So, $b$ and the {\em antipode of} $a$ are 
vertices of a regular pentagon $P'$ on $\Gamma$. The regular pentagon $P$ 
with $a$ as a vertex is the complement of $P'$ in a regular decagon $D$ on 
$\Gamma$. Since reflection in a vertex of $D$ preserves $P$ and $P'$ 
respectively, it follows that $\alpha$ is in $P$ and $\beta$ is in $P'$.   
Next we observe that the right hand side of Equation 3 equals  
$\alpha(c)$, which is $-c$, since $\alpha$ is perpendicular to $c$. Denote 
by $\tau$ the angle between $b$ and $c$ on the great circle $E$ that 
they span. ($\Delta$ contains no vertices of $D$, hence $b \neq \pm c$.) 
The equation 3 then has the solutions $\tau = \frac{\pi}{3} + 
\frac{2 \pi k}{3}$, that is $c$ is a vertex of an equilateral triangle 
on $E$, which also contains as a vertex the antipode of $b$. Assume 
first that $c$ is antipodal to $b$. In particular, $c$ is a vertex of 
$P$. Hence, $c(\alpha)$ is again in $P$. But by Equation 1, $c(\alpha) 
= \beta$ which is in P'. Contradiction. Assume instead that $c$ is not 
antipodal to $b$. Then there are exactly two possibilities for $c$: the two 
points on $\Delta$ at a distance $\frac{\pi}{3}$ from $b$. That these two 
points exist is clear from the relative position of $\Gamma$ and $\Delta$: 
$b$ cannot be $\beta$, since it would mean that $baba^{-1}b^{-1}=b$ which is equuivalent to $a(b)=b$ and, hence, to 
$a = \pm b$, which is excluded by assumption. 
Therefore the distance of $b$ to $\Delta$ along $\Gamma$ is at 
most $\frac{3\pi}{10} < \frac{\pi}{3}$. 

One readily checks from this geometric description that both 
choices of $c$ give solutions $(a,b,c)$ also to Equation 2. 
Indeed, one checks that $c(b(c))=-b$ while $a(b(a(b(a))))$ is the point on $\Gamma $ corresponding to the angle $-4\theta \equiv \theta + \pi\,\,  \text{mod} \,\, 2\pi\,$ which also gives $-b$. 

Hence 
given $a$, there are two circles worth of choices of $b$ at 
angle $\theta$ from $a$, and then exactly two choices of $c$ 
for each such pair $a$,$b$. 
These solutions give a compact subspace of the space of all solutions which  is a union of some connected components of that space. 
The two choices of $c$ give points $(a,b,c)$ which lie in different connected components. To see that we observe that the continuous function $\det (a, b, c)$, where  $a,\, b,\, c$ are interpreted as vectors in $\R ^3$, does not vanish on the subspace of solutions considered here and that it differs in sign for the two different choices of $c$, given the pair $a$,$b$.
Thus, points $(a, b, c_1)$ and  $(a, b, c_2)$ with $c_1\ne c_2\,$ lie in different connected components.
It follows that the solutions of 
this form constitute $4$ copies of $\mathbb{R}P^3$.  

\subsubsection{Case 2c: $\alpha \neq \pm \beta$} In the last case, 
where $\alpha$ and $\beta$ are non-parallel, Equation 1 forces $c$ to be 
on $\Gamma$. Denote, as before, the angle between $a$ and $b$ by $\theta$.
Then $\alpha$ corresponds to $-2\theta$ and $\beta$ to $3\theta$. Equation $1$ 
says that $c$ must be on the line bisecting the angle 
between $\alpha$ and $\beta$, that is at an angle $\vartheta = \frac{\theta}{2} 
+ k\pi$ from $a$. With this notation we have $\theta = 2 \vartheta$ and it is straightforward 
to compute that the equations 2 and 3 are both equivalent to 
$7\vartheta = 0$ $mod$ $2\pi$, where we should disregard the trivial 
solution $\vartheta =0$ since it belongs to a case already studied. Denote the
regular heptagon containing $a$ and $c$ by $H$. The equality  $\theta = 2 \,\vartheta$ tells us that $b$ is uniquely determined given $a$ and $c$. Once $a$ is given, there are three circles worth of choices of $c$ at angle $\vartheta$ from $a$ and then exactly one $b$ for every such pair $a$, $c$. The solutions of this form constitute $3$ copies of  $\mathbb{R}P^3.$ 

Thus Case $1$ contributes the component $S^2$, Case $2a$ nothing, Case $2b$ four copies of $\R P^3$ and Case $2c$ three copies of  $\R P^3$  to the solution space $J_C(9_{42})$. Proposition \ref{counterP1} follows.
\end{proof}

\bigskip

\noindent
{\bf  Comparison to Khovanov homology}

\bigskip

According to Proposition \ref{counterP1}, 
\begin{equation}\label{counter1}
J_C(9_{42}) = S^2 \cup 
\bigcup_{i=1}^{7} \mathbb{R}P^3.
\end{equation}

\medskip

\noindent 
On the other hand, consultation of JavaKh gives for collapsed Khovanov 
homology,
\smallskip 
\begin{displaymath}
\begin{split}
Kh^*(9_{42})&
=H^*(S^2)\{-1\} \oplus H^*(\mathbb{R}P^3)\{-5\} \oplus 
H^*(\mathbb{R}P^3)\{-3\}\oplus\\
&\qquad \qquad \qquad \qquad \qquad \qquad  \oplus H^*(\mathbb{R}P^3)\{-2\} \oplus 
H^*(\mathbb{R}P^3)\{0\}.
\end{split}
\end{displaymath}

\bigskip

\noindent
Hence the two homology groups differ 
by three copies of $\mathbb{R}P^3$.

\vskip1.2truecm

\appendix
\section{On the structure of the singularities of $K_{2n}$}

\bigskip

In this appendix we shall give a proof of Proposition 2.7 describing structure of the singularities  of the spaces $\, K_{2n}\,$ and  $\, \widetilde{K}_{2n}\,$.

As pointed out in Subsection 2.2 the singularities of  $\, K_{2n}\,$  and of  $\, \widetilde{K}_{2n}\,$ are contained in the subsets $\, \Sigma\,$ and $\, \widetilde{\Sigma}\,$ respectively. The sets  $\, \Sigma\,$ and $\, \widetilde{\Sigma}\,$ are disjoint unions of subsets $\, \Delta_{\epsilon _1, ... ,\epsilon _{2n}}\,$ and there is a smooth action of the group $\, (\Z /2\Z)^{2n}\,$ on $\, P_{2n}\,$ mapping $\, K_{2n} \cup \widetilde{K}_{2n}\,$ to itself and acting transitively on the family of all subsets   $\, \Delta_{\epsilon _1, ... ,\epsilon _{2n}}\,$. Hence the structure of  $\, K_{2n}\,$  and of  $\, \widetilde{K}_{2n}\,$ in a neighbourhood of each  $\, \Delta_{\epsilon _1, ... ,\epsilon _{2n}}\,$ is the same and it is enough to study it for just one $\, \Delta =  \Delta_{\epsilon _1, ... ,\epsilon _{2n}}\,$. We choose  $\, \Delta =  \Delta_{\epsilon _1, ... ,\epsilon _{2n}}\,$ with $\, \epsilon _j=(-1)^{j-1}\,$.

Thus $\, \Delta = \{ \, (A, -A, A, -A, ..., A, -A) \in  P_{2n} \,\, | \,\, A\in  {\mathcal S}\, \}\,$. Recall that $\,  {\mathcal S} \subset SU(2)\,$ consists of the matrices $\, \left( \smallmatrix it&z\\-{\bar z}&-it \endsmallmatrix\right)\,$ with $\, t\in\R, z\in \C\,$ and $\, t^2+|z|^2=1\,$.

We shall now equip $\, P_{2n}\,$ with a Riemannian metric.

Let us first choose a left- and right-invariant Riemannian metric on the Lie group $\, SU(2)\,$. The geodesics $\, \gamma \,$ in $\, SU(2)\,$  with $\, \gamma (0)=I\,$  are precisely the one-parameter subgroups of   $\, SU(2)\,$, see \cite{M2} , Lemma 21.2. 

Let $\, U\,$ be the linear subspace of the Lie algebra $\, {\mathfrak su}(2)\,$ consisting of the matrices  $\,\widetilde{u} =  \left( \smallmatrix 0&u\\-{\bar u}&0 \endsmallmatrix\right), \,\, u\in\C\,$.
The tangent space $\, T_J  {\mathcal S}\,$ to the submanifold  $\,  {\mathcal S}\,$ of $\, SU(2)\,$ at the point $\, J= \left( \smallmatrix i&0\\ 0&-i\endsmallmatrix\right)\,$ consists of the matrices  $\, \left( \smallmatrix 0&z\\-{\bar z}&0 \endsmallmatrix\right)=  \left( \smallmatrix 0&iz\\i{\bar z}&0 \endsmallmatrix\right) \cdot  J, \,\, z\in \C ,\,$ and, hence, it is equal to $\, U J = \{ \, \widetilde{u}J \, | \,\widetilde{u}\in U \, \} \,$.  

Whenever convenient we shall identify $\, U\,$ with $\, \C\,$ via $\, \widetilde{u}\leftrightarrow u \,$. Under this identification the restriction of the Riemannian metric to $\, U\,$ is a positive multiple of the standard real inner product on $\, \C = \R ^2\,$. 

For any matrix $\, \widetilde{u} =  \left( \smallmatrix 0&u\\-{\bar u}&0 \endsmallmatrix\right), \,\, u\in\C\,$, one has $\, \widetilde{u} ^2=-|u|^2 I\,$. Therefore 
\begin{displaymath}
\begin{split}
e^{ \widetilde{u}}&=\left(\sum\limits_{j=0}^{\infty}\frac {(-|u|^2)^j}{(2j)!}\right)
I
+ \left(\sum\limits_{j=0}^{\infty}\frac {(-|u|^2)^j}{(2j+1)!}\right)
\widetilde{u}
=\\
&=\cos (|u|) \, I + \frac {\sin (|u|)}{|u|} \,\, \, \widetilde{u}\, ,
\end{split}
\end{displaymath} 
where $\,  \frac {\sin (z)}{z}\,$ is to be interpreted as an entire analytic function. It follows that, if $\, z\in \C , \,\, |z|=1,\,\, B= \left( \smallmatrix 0&iz\\i{\bar z}&0 \endsmallmatrix\right)\,$, then
\begin{equation}
\gamma (t) = e^{tB}\cdot J = \left(\begin{array}{rr}
i\, \cos(t)& \sin (t) \, z\\
- \sin(t) \, {\bar z}& -i\, \cos (t) 
\end{array} \right)\,\, , \qquad t\in \R,
\end{equation} 
is the geodesic in $\, SU(2)\,$ with the starting point at $\, J\,$ and corresponding to the tangent vector $\,\left( \smallmatrix 0&z\\-{\bar z}&0 \endsmallmatrix\right)=  \left( \smallmatrix 0&iz\\i{\bar z}&0 \endsmallmatrix\right) \cdot  J \in  T_J  {\mathcal S}\,$. 

Observe also that, since $\, \text{Tr}( e^{tB}\cdot J)=0\,$,  one has  $\,\gamma (t) = e^{tB}\cdot J \in  {\mathcal S}\,$ for all $\, t\in \R\,$, showing that $\,  {\mathcal S}\,$ is a totally-geodesic submanifold of $\,  SU(2)\,$. 

We equip  $\,  {\mathcal S}\,$ with the Riemannian metric being  the restriction of that of  $\,  SU(2)\,$ and give $\, P_{2n}=\prod\limits_1^{2n}  {\mathcal S}\,$ the product Riemannian metric. Geodesics in  $\, P_{2n}\,$ are then products of geodesics in the factors. The group $\,  SU(2)\,$ acts on $\, P_{2n}\,$ through isometries. (The Riemannian metric we have introduced on  $\,  {\mathcal S}\,$ is, of course, nothing but the ordinary metric on the $2$-dimensional sphere. It will however be covenient to phrase it that way.)  

Recall that we have the $\,  SU(2)$-equivariant  mapping $\, p_{2n}: P_{2n}\rightarrow SU(2)\,$, $\, p_{2n}(A_1, ... , A_{2n})=A_1\cdot ... \cdot A_{2n}  \,$ and that $\, K_{2n}= p_{2n}^{-1}(I)\,$ and  $\, \widetilde{K}_{2n}= p_{2n}^{-1}(-I)\,$. From now on we denote $\, p_{2n}\,$ by $\,  p\,$.   

Let us consider the point $\,{\mathcal J}=(J, -J, J, -J, ... , J, -J)\in \Delta \subset P_{2n} \,$. The orbit of  $\,{\mathcal J}\,$ under the action of  $\,  SU(2)\,$ is equal to $\, \Delta\,$. The isotropy subgroup of  $\,{\mathcal J}\,$ is the maximal torus $\, {\mathbb T}\,$ consisting of the diagonal matrices in  $\,  SU(2)\,$.
We identify 
 the tangent space $\, T_{\mathcal J} P_{2n} \,$ to  $\, P_{2n}\,$ at the point  $\,{\mathcal J}\,$ with  $\, U^{2n} = \bigoplus\limits _{j=1}^{2n} \, U \cong \C ^{2n}\,$
 through  $\, T_{\mathcal J} P_{2n} \ni     ( \widetilde{u}_1 \cdot J, \,\, \widetilde{u}_2 \cdot (-J), ... ,  \widetilde{u}_{2n} \cdot (-J)) \leftrightarrow 
( \widetilde{u}_1,  \widetilde{u}_2, ... ,  \widetilde{u}_{2n})\in  U^{2n}   \,$. 

We have $\, p({\mathcal J})=I\,$. The derivative of $\, p\,$ at $\,{\mathcal J}\,$ is a $\, {\mathbb T}$-equivariant linear map    $\, dp_{\mathcal J}  :  T_{\mathcal J} P_{2n}\rightarrow {\mathfrak su}(2)\,$ given by 
\begin{displaymath}
dp_ {\mathcal J}: U^{2n}\rightarrow {\mathfrak su}(2),
\end{displaymath}
\begin{equation}
dp_{\mathcal J} (\widetilde{u}_1,  \widetilde{u}_2, ... ,  \widetilde{u}_{2n}) = \sum\limits_{j=1}^{2n}\, (-1)^{j-1} \widetilde{u}_j .
\end{equation}
The last claim follows from the fact that the conjugation by $\, J\,$ acts on $\, U\,$ as multiplication by $\, -1\,$. The image of $\,dp_{\mathcal J}\,$ is equal to $\, U\,$. 

The tangent space to $\, \Delta \,$ at the point  $\,{\mathcal J}\,$ is a subspace of  $\,  T_{\mathcal J} P_{2n}\,$ which corresponds to the diagonal $\, \widetilde{\Delta}\,$ of $\,  U^{2n}\,$, $\,  \widetilde{\Delta}=\{ (\widetilde{u}, \widetilde{u}, ... ,\widetilde{u})\in   U^{2n} \, \mid \, \widetilde{u}\in U \,\}\,$. It is a subspace of $\,\text{Ker} ( dp_{\mathcal J} )\,$ (as follows also directly from the fact that $\, p(\Delta )= \{ I\}\,$). As the circle  $\, {\mathbb T}\,$ acts on $\, \Delta \,$ with $\,{\mathcal J}\,$ being a fixed point, the tangent space  $\,  \widetilde{\Delta}= T_{\mathcal J}\Delta \,$ is an invariant subspace of  $\,  T_{\mathcal J} P_{2n}\,$. Let $\, V\,$ be its orthogonal complement in $\,  T_{\mathcal J} P_{2n}\,$.   
Since $\, {\mathbb T}$ acts on  $\,  T_{\mathcal J}  P_{2n}\,$ through isometries, $\, V\,$ is a $\,  {\mathbb T}$-invariant subspace as well.

We shall denote $\, SU(2)\,$ by $\, G\,$ to simplify notation.

The $\, G$-submanifold $\, \Delta \,$  of $\, P_{2n}\,$  has an equivariant regular  neighbourhood $\, G$-diffeomorphic to a neighbourhood of the zero-section in the normal bundle $\, G\times _{\mathbb T} V\,$ of $\, \Delta\,$ in $\, P_{2n}\,$ (see \cite {K1} or \cite{J1}, Thm.1.3, p.3). The $\,G$-diffeomorphism between them can be obtained as follows.

Given a tangent vector $\, v\in V\,$ let $\, \gamma _v\,$ be the geodesic in $\, P_{2n}\,$ starting at the point  $\,{\mathcal J}\,$ and such that $\, \frac {d  \gamma _v}{dt} (0) = v\,$. We define $\, \phi : G\times  _{\mathbb T} V \rightarrow P_{2n}\,$ by $\, \phi ([g, v])= g( \gamma _v (1))\,$. The mapping $\, \phi \,$ is $G$-invariant and, when restricted to a small enough neighbourhood $\, Y\,$ of the zero-section of the normal bundle $\,  G\times  _{\mathbb T} V\,$, it is a diffeomorphism of  $\, Y\,$ onto a neighbourhood of $\, \Delta\,$ in $\, P_{2n}\,$, \cite{J1}, p.3-4. 

Hence, in order to describe the singularities of $\, K_{2n}=p^{-1}(I)\,$ along $\, \Delta\,$, it is enough to look at the restriction of the mapping $\, p\,$ to the neighbourhood $\, \phi ( G\times _{\mathbb T} V)\,$ of $\, \Delta\,$ or, in other words, to study the singularities of the composed mapping $\, \widetilde{p}=p\circ \phi :  G\times _{\mathbb T} V \rightarrow SU(2)\,$ along the zero-section of $\, G\times _{\mathbb T} V\,$. Since the mapping $\,\widetilde{p}\,$ is $\, G$-equivariant, it will be enough to study the singularity of the restricted mapping $\, \widehat{p}: V\rightarrow SU(2), \,\,  \widehat{p}(v)=\widetilde{p} ([ I,  v]) = p(\gamma _v (1))\,$ at the origin.

Since $\, \phi\,$ is nothing but a restriction of the exponential map of the Riemannian manifold  $\, P_{2n}\,$, its derivative at the zero-section of $\,  G\times _{\mathbb T} V \,$ is equal to the identity. Therefore the derivative $\, d\widehat{p}_0:  V\rightarrow {\mathfrak su}(2)\,$ of $\,\widehat{p}\,$ at the point $\, 0\,$ is given by the restriction of (A.2) to $\, V\,$.  We want to determine the structure of the subspace $\, \widehat{p}^{-1}(I)\,$ of $\, V\,$ locally in a neighbourhood of the point  $\, 0\,$.

Since $\, V\,$ is an orthogonal complement in $\,  T_{\mathcal J} P_{2n}\,$ 
of a subspace contained in  $\,\text{Ker} ( dp_{\mathcal J} )\,$ it follows from (A.2) that the image of the derivative $\,  d\widehat{p}_0 \,$  is equal to $\, U\,$. Thus  $\,  d\widehat{p}_0 \,$ is of rank $2$ and $\, 0\in V\,$ is a singular point of $\, \widehat{p}: V \rightarrow SU(2)\,$. 

The circle $\,  {\mathbb T}\,$ consisting of the diagonal matrices is a submanifold of $\, SU(2)\,$. Its tangent space $\, T_I  {\mathbb T}\,$ at $\, I\in  {\mathbb T}\,$ is transversal (and complementary) to $\, U\,$ in $\, {\mathfrak su}(2)\,$. Therefore there is an open neighbourhood $\, \widetilde{V}\,$ of $\, 0\,$ in $\, V\,$  such that the restriction of $\,\widehat{p}\,$ to  $\, \widetilde{V}\,$, $\,\widehat{p} \, | _{\widetilde{V}} : \widetilde{V} \rightarrow SU(2)\,$, is transversal to the submanifold   $\,  {\mathbb T}\,$. 

Let $\, M\,$ be the pre-image of  $\,  {\mathbb T}\,$ under  $\,\widehat{p} \, | _{\widetilde{V}}\,$, $\, M= \widehat{p} ^{-1}({\mathbb T} ) \cap \widetilde{V}\,$. It is a smooth submanifold of $\, \widetilde{V}\,$ of codimension $2$ containing the point $\, 0\,$. The tangent space $\, T_0 M \,$ is equal to $\,( d\widehat{p}_0)^{-1}(T_I {\mathbb T})\,$. As $\,  d\widehat{p}_0(V) = U\,$ and $\, U\cap T_I {\mathbb T} = \{0\}\,$, it follows that  $\, T_0 M \,$ is also equal to the kernel of $\, d\widehat{p}_0:  V\rightarrow {\mathfrak su}(2)\,$.

Let us consider the restriction of $\, \widehat{p}\,$ to the submanifold  $\, M\,$, $\,\widehat{f}=\widehat{p} \,|_M : M\rightarrow  {\mathbb T}\,$. The space   $\,\widehat{p} ^{-1}(I ) \cap \widetilde{V}\,$ which we want to study is a subspace of $\, M\,$ equal to $\,  \widehat{f} ^{-1}(I )\,$. 

Let $\, W\,$ be a neighbourhood of $\, I\,$ in $\, SU(2)\,$ which is a domain of the chart $\, \text{Log} : W\rightarrow  {\mathfrak su}(2)\,$, (  $\, \text{Log}\,$ being the inverse  of the exponential mapping $\, \text{exp} :   {\mathfrak su}(2)\rightarrow SU(2)\,$). We can assume that $\, \widehat{p} (\widetilde{V}) \subset W\,$.  

Let $\, T\,$ be the linear subspace  of $\, {\mathfrak su}(2)\,$ which is the tangent space of the circle $\, {\mathbb T}\,$ at $\, I\,$. Then   $\, T\,$ consists of diagonal matrices in  $\, {\mathfrak su}(2)\,$ and it is the orthogonal complement of $\, U\,$. Let $\, \pi :  {\mathfrak su}(2)\rightarrow  T\,$ be the orthogonal projection of  $\, {\mathfrak su}(2)\,$ onto  $\, T\,$ and let $\, \iota :T \rightarrow \R \,$ be the isomorphism $\, \iota \left(\smallmatrix is & 0\\ 0& -is\endsmallmatrix\right) =s\,$. Finally, define mappings 
\begin{displaymath}
\overline{p}: \widetilde{V}\rightarrow T \qquad \text{and}  \qquad F: \widetilde{V}\rightarrow \R 
\end{displaymath}
to be $\,\overline{p}=\pi \circ \text{Log} \circ (\widehat{p} \, | _{ \widetilde{V}})\,$ and $\, F=\iota \circ \overline{p}\,$.

Observe that since $\, \widehat{p}(M)\subset {\mathbb T}\,$ and $\, \text{Log}(  {\mathbb T}) \subset T\,$, the restriction of  $\,\overline{p}\,$ to $\, M\,$ is equal to $\, \text{Log} \circ \widehat{f}\,$, $\, \overline{p}\,|_M= \text{Log} \circ \widehat{f}\,$. Therefore $\,\widehat{p} ^{-1}(I ) \cap \widetilde{V}= \widehat{f} ^{-1}(I ) = (F \, |_M)^{-1}(0)\,$.

As the derivative of $\, \text{Log}\,$ at $\, I\,$ is the identity  map, we have $\, d\overline{p}_0=\pi \circ d\widehat{p}_0\,$ and $\, dF_0= \iota \circ \pi \circ d\widehat{p}_0\,$. Since the image of $\, d\widehat{p}_0\,$ is equal to $\, U\,$ which is the kernel of $\, \pi \,$, we get that  $\, d\overline{p}_0 = 0\,$ and $\, dF_0= 0\,$. Thus the origin $\, 0\in \widetilde{V}\,$ is a critical point of the function $\, F\,$. 

Let us denote by $\, f\,$ the restriction of $\, F\,$ to $\, M\,$,
\begin{equation}
f=F \, | _M = \iota \circ \overline{p} \, |_M = \iota \circ \text{Log}\circ \widehat{f}\, :\, M\rightarrow \R \,\, .
\end{equation}
The origin $\, 0\in M\,$ is a critical point of $\, f\,$.

We want to determine the Hessian (quadratic form) of $\, f\,$ at the point $\, 0\,$. Let $\, \widehat{N}\,$ be the kernel of $\, d\widehat{p}_0 : V\rightarrow  {\mathfrak su}(2)\,$. It is a (linear) submanifold of codimension $2$ in $\, V\,$ containing $\, 0\,$. The tangent spaces of $\, M\,$ and $\,\widehat{N}\,$ at $\, 0\,$ coincide. Let $\, N=\widehat{N}\cap \widetilde{V}\,$ and $\, g=F \,|_N : N\rightarrow \R\,$. The origin $\, 0\,$ is a critical point of $\, g\,$ and we have

\begin{lma}
The Hessian quadratic forms of $\, g\,$ and of $\, f\,$ at $\, 0\,$ coincide,
\begin{displaymath}
Hess_0(f) = Hess_0(g) \,\, .
\end{displaymath} 
\end{lma}

\begin{proof} 
The equality of both Hessians depends on the fact that $\, 0\in \widetilde{V}\,$ is a critical point of the funktion $\, F\,$ (and it would not need to hold otherwise).

Let $\, X_0, Y_0 \in T_0M=T_0N\,$ be two tangent vectors to $\, M\,$, and hence to $\, N\,$, at $\, 0\,$.

Let us  extend $\, Y_0\,$ to a vector field $\, Y\,$ on $\, M\,$ and then extend $\,Y\,$ to a vector field $\, \widetilde{Y} \,$ on $\, \widetilde{V}\,$. Let $\, Hess_0(F)\,$ be the Hessian bilinear form of $\, F : \widetilde{V}\rightarrow \R \,$ at $\, 0\,$. Then $\,  \widetilde{Y}(F) \, |_M = Y(f)\,$ and
\begin{displaymath}
\begin{split}
Hess_0(f)(X_0, Y_0 )&= X_0(Y(f))=X_0( \widetilde{Y}(F) \, |_M ) = X_0( \widetilde{Y}(F))=\\
&= Hess_0(F)(X_0, Y_0) \,\, .
\end{split}
\end{displaymath}
(See \cite{M2}, p.4, for the definition of the Hessian bilinear form.)

Since $\, g=F \,|_N\,$, we get in a similar way that  $\, Hess_0(g)(X_0, Y_0 )= Hess_0(F)(X_0, Y_0)\,$. Therefore $\, Hess_0(f)(X_0, Y_0 )=Hess_0(g)(X_0, Y_0 )\,$, proving the Lemma.
\end{proof}

Hence, in order to describe the Hessian of $\, f\,$ at $\, 0\,$, it is enough to describe the Hessian of $\, g\,$ at $\, 0\,$, $\, g\,$ being the restriction of $\, F\,$ to the {\it linear} submanifold $\, N\,$  of $\, V\,$.

The identification of $\, U\,$ with $\, \C\,$ chosen before gives $\, U\,$ a structure of a complex vector space. The action of the circle $\, {\mathbb T}\,$ on $\, U\,$ preserves that structure. The derivative $\, dp_{\mathcal J }: U^{2n}\rightarrow U\,$ given by (A.2) is a complex linear map.  As $\, \widetilde{\Delta}=T_{\mathcal J } \Delta\,$ is a $\C $-linear subspace of $\, U^{2n}\,$ (being the diagonal), so is its orthogonal complement $\, V\,$. The map $\, d\widehat{p}_0: V\rightarrow U\,$ is $\C$-linear as a restriction of $\, dp_{\mathcal J }\,$.  

The tangent space  $\, T_0N\,$ to $\, N\,$ at $\, O\,$ is equal to $\, \widehat{N}=\text{ker}(d\widehat{p}_0:V\rightarrow U)\,$. Since $\, d\widehat{p}_0\,$ is $\C$-linear and surjective, $\, \widehat{N}\,$ is a $\C$-linear subspace of $\, V\,$ of complex codimension $1$. As $\, V\subset U^{2n}\,$ is the orthogonal complement of the diagonal $\, \widetilde{\Delta} \subset U^{2n}\,$, a $\C$-basis of $\, V\,$ is given by vectors $\, c_j=e_j-e_{j+1}, \,\, j=1, ... , 2n-1,\,$ where $\, e_1, ... ,e_{2n}\,$ is the standard $\C$-basis of $\, U^{2n}=\C ^{2n}\,$.  Since, according to (A.2), $\,  d\widehat{p}_0(c_j)=  d\widehat{p}_0(e_j-e_{j+1})=(-1)^{j-1}-(-1)^j=(-1)^{j-1}\cdot 2 \in U=\C\,$, it follows that the vectors \begin{equation}
w_j=c_j+c_{j+1}=e_j-e_{j+2} \in  U^{2n}=\C ^{2n} \, , \quad j=1, ... ,2n-2,
\end{equation}
form a $\C$-basis of the tangent space  $\, T_0N = \text{ker}(d\widehat{p}_0)\,$ and an $\R$-basis of  $\, T_0N\,$ is given by the vectors $\, v_k, \,\, k=1, ... ,4n-4\,$,
\begin{equation}
v_k=\left\{\begin{array}{lll}
w_j& \text{if}& k=2j-1,\\
iw_j&\text{if}& k=2j,
\end{array}\right. \qquad \quad  j=1, ... , 2n-2.
\end{equation}

Let $\, (x_1, ... , x_{4n-4})\,$ be the system of coordinates on $\, N= T_0N \cap \widetilde{V}\,$ given by the $\R$-basis $\,\{ v_k\}\,$. We shall compute the Hessian matrix of the function $\, g\,$  at $\, 0\,$ w.r.t. these coordinates.

Let $\, 1\le k \le l\le 4n-4\,$. We consider the mapping $\, g_{k,l}: \R ^2 \rightarrow \R ,\,\,  g_{k,l}(s, t )= g(sv_k+tv_l)\,$, in order to compute $\, \dfrac{\partial ^2 g}{\partial x_k \partial x_l} (\vec 0\, ) =  \dfrac{\partial ^2 g_{k,l}}{\partial t \, \partial s} (0, 0)\,$.

Let $\, \text{exp}:V\rightarrow P_{2n}\,$ be the exponential mapping of the Riemannian manifold $\, P_{2n}\,$ restricted to the normal space $\, V\,$ of the orbit $\, \Delta \,$. We have $\, g=F \,|_N=\iota \circ \pi \circ \text{Log} \circ p \circ  \text{exp}\, | _N\,$. The mappings $\, g\,$ and $\, g_{k,l}\,$ can be factorized through  the mappings $\, \overline{g}= \text{Log} \circ p \circ  \text{exp}\, | _N : N\rightarrow  {\mathfrak su}(2)\,$ and   $\, \overline{g}_{k,l}: \R ^2 \rightarrow  {\mathfrak su}(2),\,\,  \overline{g} _{k,l}(s, t )=  \overline{g}(sv_k+tv_l)\,$ respectively. 

We shall first compute $\,  \dfrac{\partial ^2  \overline{g}_{k,l}}{\partial t \, \partial s} (0, 0)\in {\mathfrak su}(2) \,$. 
Let us recall a result of the Lie theory, \cite{V1},

\begin{thm} Let $\, H\,$ be a Lie group and $\,  {\mathfrak h}\,$ its Lie algebra. There is a neighbourhood $\, W$ of $\, 0\,$ in  $\,  {\mathfrak h}\,$ such that for $\, X, Y \in W \subset  {\mathfrak h}\,$ one has 

(i) $\, \text{exp}\, X \,\cdot \text{exp}\, Y =  \text{exp} ( X + Y + \tfrac 12 [ X, Y ] + O(\rho ^3))\,$;

(ii)  $\, \text{exp}\, X \,\cdot \text{exp}\, Y \, \cdot  \text{exp}\, (-X) \,\cdot \text{exp}\,(-Y)   =  \text{exp}\, (\, [ X, Y ] + O(\rho ^3))\,$. 
\end{thm}
\noindent
(Here $\,exp\,$ is the exponential mapping of the Lie group $\, H\,$ and $\, \rho \,$ is a length of the vector $\, (X,Y)\,$ in $\, {\mathfrak h} \times {\mathfrak h}\,$.)   

Recall also that for $\, (u_1, ... , u_{2n})\in V \subset U^{2n}\,$ we have 
\begin{displaymath}
\begin{split}
\overline{g}(u_1, ... , u_{2n})&=\text{Log} \left( p(\text{exp}(u_1, ... , u_{2n}))\right) = \\ 
&=\text{Log} ( p( e^{\widetilde{u_1}} J, \,  e^{\widetilde{u_2}} (-J), ... ,  e^{\widetilde{u_{2n}}} (-J))) =\\ 
&= \text{Log} ( e^{\widetilde{u_1}} J\cdot  e^{\widetilde{u_2}} (-J)\cdot \cdot \cdot   e^{\widetilde{u_{2n}}} (-J)) = \\
&= \text{Log} ( e^{\widetilde{u_1}} \cdot  e^{-\widetilde{u_2}} \cdot \cdot \cdot  e^{\widetilde{u_{2n-1}}}   \cdot  e^{-\widetilde{u_{2n}}})\, .
\end{split}
\end{displaymath}

In the following $\, r=\sqrt{s^2+t^2}\,$. We consider several cases.

\medskip

(1) If $\, k=l=2j-1\,$ then
\begin{displaymath}
\begin{split}
 \overline{g}_{k,l}(s,t)&= \overline{g}(s v_k + t v_k)= \overline{g}((s+t) w_j)= \overline{g}((s+t)(e_j-e_{j+2}))=\\
&=\text{Log}\left(e^{(-1)^{j-1}\widetilde{(s+t)}} \cdot  e^{-(-1)^{j+1}\widetilde{(s+t)}}  \right) = \text{Log}\,  (I) =\\
&=0 \, .
\end{split}
\end{displaymath} 

We get the same result in cases when  

(2) $k=2j-1$ and $l=k+1,k+2, k+4$ or $l\ge k+6$, 

(3) $\, k=l=2j\,$,

(4) $k=2j$ and $l\ge k+4$.

In all those cases one gets, of course, $\,   \dfrac{\partial ^2  \overline{g}_{k,l}}{\partial t \, \partial s} (0, 0)=0 \,$.

\medskip

(5) If $\, k=2j-1, \, l=k+3  \,$ then
\begin{displaymath}
\begin{split}
 \overline{g}_{k,l}(s, t)&= \overline{g}( sv_k + t v_{k+3})= \overline{g}(s(e_j-e_{j+2}) + t(i e_{j+1} - i   e_{j+3}))=\\
&= \text{Log} \left(e^{(-1)^{j-1}\widetilde{s}} \cdot e^{(-1)^j\widetilde{(ti)}} \cdot e^{- (-1)^{j+1}\widetilde{s}}  \cdot e^{- (-1)^{j+2}\widetilde{(ti)}}  \right)=\\
&=[\,\, (-1)^{j-1}\widetilde{s}\, , \,\, (-1)^j\widetilde{(ti)}\,\, ] + O(r^3) =\qquad ( \text{by Theorem A.2 (ii)})\\
&= -st\, [\widetilde{1}\, , \,\, \widetilde{i} \,\, ] +  O(r^3) =\\
&= -2st \left( \smallmatrix i&0\\0&-i \endsmallmatrix \right) +  O(r^3)  \quad .
\end{split}
\end{displaymath} 
The last equality follows from $\, [\widetilde{1}\, , \,\, \widetilde{i} \,\, ]= \left[ \left( \smallmatrix 0&1\\-1&0 \endsmallmatrix \right)\, ,\,\, \left( \smallmatrix 0&i\\i&0 \endsmallmatrix \right)\right] = 2  \left( \smallmatrix i&0\\0&-i \endsmallmatrix \right) \,$.

Consequently,   $\,   \dfrac{\partial ^2  \overline{g}_{k,l}}{\partial t \, \partial s} (0, 0)= -2 \left( \smallmatrix i&0\\0&-i \endsmallmatrix \right)  \,$ in that case.

\medskip

(6) If  $\, k=2j-1, \, l=k+5  \,$ then 
\begin{displaymath}
\begin{split}
 \overline{g}_{k,l}(s, t)&= \overline{g}( sv_k + t v_{k+5})= \overline{g}(s(e_j-e_{j+2}) + t(ie_{j+2}-ie_{j+4}))=\\
&=\text{Log} \left(e^{(-1)^{j-1}\widetilde{s}} \cdot e^{(-1)^{j+1}\widetilde{(-s+ti)}} \cdot e^{- (-1)^{j+3}\widetilde{(ti)}}   \right)=\\
&=[\, \text{by applying Theorem A.2 (i) twice  }]=\\
&=\left( (-1)^{j-1} \widetilde{s} + (-1)^{j+1} \widetilde{(-s+ti)}+\tfrac12 [\, \widetilde{s} \, , \, \, \widetilde{(-s+ti)}\, ]\right) - (-1)^{j+3}\widetilde{(ti)} +\\
& \qquad\quad  + \tfrac 12 \, \left[  (-1)^{j-1} \widetilde{s} + (-1)^{j+1} \widetilde{(-s+ti)} \, , \,\,  - (-1)^{j+3}\widetilde{(ti)}  \right] +  O(r^3) =\\
&= \tfrac 12 \, [\, \widetilde{s} \, , \,\, \widetilde{(ti)} \,] +  O(r^3) = \tfrac 12 \,  st\, [\, \widetilde{1} \, , \,\, \widetilde{i} \,] +  O(r^3)= \\
&= st \left( \smallmatrix i&0\\0&-i \endsmallmatrix \right) +  O(r^3) \,\, . 
\end{split}
\end{displaymath}

Therefore   $\,   \dfrac{\partial ^2  \overline{g}_{k,l}}{\partial t \, \partial s} (0, 0)= \left( \smallmatrix i&0\\0&-i \endsmallmatrix \right)  \,$ in that case. 

\medskip

(7) If  $\, k=2j, \,\,  l=k+1  \,$ we get 
\begin{displaymath}
 \overline{g}_{k,l}(s, t)= 2  st \left( \smallmatrix i&0\\0&-i \endsmallmatrix \right) +  O(r^3)\, ,
\end{displaymath} 
and, consequently,  $\, \, \,   \dfrac{\partial ^2  \overline{g}_{k,l}}{\partial t \, \partial s} (0, 0)= 2 \left( \smallmatrix i&0\\0&-i \endsmallmatrix \right)  \,$. The proof is the same as in the case (5).

\medskip
Finally,

(8) if  $\, k=2j, \,\,  l=k+3  \,$ we get  
\begin{displaymath}
 \overline{g}_{k,l}(s, t)= -  st \left( \smallmatrix i&0\\0&-i \endsmallmatrix \right) +  O(r^3)\, ,
\end{displaymath} 
and, hence,  $\,   \dfrac{\partial ^2  \overline{g}_{k,l}}{\partial t \, \partial s} (0, 0)= - \left( \smallmatrix i&0\\0&-i \endsmallmatrix \right)  \,$. The proof is the same as in the case (6).

\bigskip

The results of the cases (1) - (8) give a proof of 

\begin{lma}
If $\, 1\le k \le l \le 4n-4\,$, then 
\begin{displaymath}
 \dfrac{\partial ^2  \overline{g}_{k,l}}{\partial t \, \partial s} (0, 0)= \left\{  \begin{array}{lll}
-2 \left( \smallmatrix i&0\\0&-i \endsmallmatrix \right)& \text{for} & k=2j-1, \,\,l= k+3,\\
\left( \smallmatrix i&0\\0&-i \endsmallmatrix \right)& \text{for} & k=2j-1, \,\,l= k+5,\\
2 \left( \smallmatrix i&0\\0&-i \endsmallmatrix \right)& \text{for} & k=2j, \,\,l= k+1,\\
- \left( \smallmatrix i&0\\0&-i \endsmallmatrix \right)& \text{for} & k=2j, \,\,l= k+3,\\
0&&\text{otherwise}.
\end{array} 
\right.
\end{displaymath} 
\end{lma}

As an immediate conesquence of Lemma A.3 we get

\bigskip

\begin{cor}
The Hessian matrix  $\, H\,$ of the function $\, g : N \rightarrow \R\,$ at the point $\, \vec 0 \,$ is given by 
\begin{displaymath}
H_{k,l}=  \dfrac{\partial ^2  g}{\partial x_k \, \partial x_l} (\vec 0)= \left\{  \begin{array}{rll}
-2 & \text{for} & k=\text{odd}, \,\,l= k+3,\\
 1& \text{for} & k= \text{odd} , \,\,l= k+5,\\
2 & \text{for} & k=\text{even}  , \,\,l= k+1,\\
-1& \text{for} & k=\text{even}   , \,\,l= k+3,\\
0&&\text{otherwise},
\end{array} 
\right.
\end{displaymath} 
for $\, 1\le k\le l \le 4n-4\,$.
\end{cor}

\begin{proof}
The function  $\, g : N \rightarrow \R\,$ is a composition of the mapping $\, \overline{g}: N\rightarrow  {\mathfrak su}(2)\,$ and of the linear function $\, \iota \circ \pi: {\mathfrak su}(2) \rightarrow \R\,$. Since  $\, \iota \circ \pi\,$ is linear, 
\begin{displaymath}
\dfrac{\partial ^2  g}{\partial x_k \, \partial x_l} (\vec 0) = (\iota \circ \pi) \left( \dfrac{\partial ^2  \overline{g}}{\partial x_k \, \partial x_l} (\vec 0)\right)= (\iota \circ \pi) \left(  \dfrac{\partial ^2  \overline{g}_{k,l}}{\partial t \, \partial s} (0, 0)    \right)\,\, . 
\end{displaymath}
Corollary follows now from Lemma A.3 since $\, ( \iota \circ \pi) \left( \smallmatrix i&0\\0&-i \endsmallmatrix \right) = 1\,$.
\end{proof}

Thus the (symmetric) Hessian matrix $\, H\,$ of $\, g\,$ (and hence of $\, f\,$) at $\, \vec 0\,$ is a $\, (4n-4)\times (4n-4)$-matrix in a block form
\begin{equation}
 H = \left(   
\begin{array}{cccccccc}
C&A&&&&&&\\
B&C&A&&&&&\\
&B&C&A&&&&\\
&&B&.&&&&\\
&&&&.&&&\\
&&&&&.&&\\
&&&&&&C&A\\
&&&&&&B&C
\end{array}
\right) 
\end{equation} 
with all blocks of the size $4\times 4\,$ and of the form 
\begin{equation}
C=\left(\begin{array}{rrrr}0&0&0&-2\\
0&0&2&0\\
0&2&0&0\\
-2&0&0&0
  \end{array}  \right) \, , \,\,\,\, A=\left(\begin{array}{rrrr}
0&1&0&0\\
-1&0&0&0\\
0&2&0&-1\\
-2&0&1&0
  \end{array}  \right) \,\, \text{and} \,\,\,\,  B=A^T \,\, .
\end{equation}

Recall that the signature of a real symmetric matrix is defined as the number of its  positive eigenvalues minus the number of the negative eigenvalues. 

\begin{prop}
The Hessian matrix $\, H\,$ is invertible and its signature is equal to $\, 0\,$. 
\end{prop}

\begin{proof} Let $\, P\,$ be the $(4n-4)\times (4n-4)$ permutation matrix of the permutation exchanging the first coordinate with  the second one, the third coordinate with the fourth one, and so on, i.e. of the permutation $\, ((1,2)(3,4) ... (4n-5,4n-4))\,$. One checks directly that $\, P^{-1} H P = PHP= - H \,$. This  implies that the number of positive egevalues of $\, H\,$ is equal to the number of its negative eigenvalues. Thus the signature of $\, H\,$ is equal to $0$.

It remains to show that the matrix $\, H\,$ is invertible.

It follows  from Corollary A.4 that $\, H_{k,l}\ne 0 \,$ implies that $\, k-l\,$ is odd. 

Let us consider the matrix $\, \widetilde{H}=PH\,$ which is obtained from $\, H\,$ by exchanging the first row with the second one, the third row with the fourth one, and so on. Then $\,\widetilde{H}_{k,l}\ne 0\,$ implies that $\, k-l\,$ is even. The matrix $\,\widetilde{H}\,$ is skew-symmetric. Indeed, $\,\widetilde{H} \, ^T = (PH)^T = H^T P^T = HP= P(P^{-1}HP)=P(-H)=-\widetilde{H}\,$.  It has the block form
\begin{equation}
\widetilde{H} = \left(   
\begin{array}{ccccccc}
\widetilde{C}&\widetilde{A}&&&&&\\
\widetilde{B}&\widetilde{C} & \widetilde{A}&&&&\\
&\widetilde{B}  &&&&&\\
&&.&&&&\\
&&&.&&&\\
&&&&.&&\\
&&&&&\widetilde{C} &\widetilde{A} \\
&&&&&\widetilde{B} &\widetilde{C}
\end{array}
\right) 
\end{equation} 
with 
\begin{equation}
\widetilde{C}=\left(\begin{array}{rrrr}
0&0&2&0\\0&0&0&-2\\
-2&0&0&0\\0&2&0&0
  \end{array}  \right) \, , \,\,\,\, \widetilde{A} =\left(\begin{array}{rrrr}
-1&0&0&0\\0&1&0&0\\
-2&0&1&0\\0&2&0&-1
  \end{array}  \right) \,\, \text{and} \,\,\,\,  \widetilde{B}  =-\widetilde{A}  \,^T \,\, .
\end{equation}

 Denote by $\, H'\,$ the $\, (2n-2)\times (2n-2)$-matrix obtained from  $\, \widetilde{H}\,$ by removing all rows and all columns with an even index. Similarily, denote by  $\, H''\,$ the $\, (2n-2)\times (2n-2)$-matrix obtained from  $\, \widetilde{H}\,$ by removing all rows and all columns with an odd index. Then 
\begin{equation}
\det(H) =\det( \widetilde{H} )= \det (H' ) \cdot  \det (H''  ) \,\, 
\end{equation}
and both $\, H'\,$ and $\, H''\,$ are skew-symmetric matrices. Moreover, on inspection of (A.9) we see that $\, H'' = -H'\,$.\,\,  Thus $\, \det (H')=  \det (H'')\,$ and $\, \det(H)=  \det (H') ^2\,$.

The skew-symmetric  $\, (2n-2)\times (2n-2)\,$-matrix $\, H'=H'(n)\,$ has the form 
\begin{equation}
H'(n) = \left(   
\begin{array}{rrrrrrrrrr}
0&2&-1&&&&&&&\\
-2&0&-2&1&&&&&&\\
1&2&0&2&-1&&&&&\\
&-1&-2&0&-2&1&&&&\\
&&&&.&&&&&\\
&&&&&.&&&&\\
&&&&&&.&&&\\
&&&&&&&0&-2&1\\
&&&&&&&2&0&2\\
&&&&&&&-1&-2&0
\end{array}
\right) 
\end{equation} 

Furthermore, $\,  \det (H'(n)) = (Pf(H'(n)))^2\,$, where $\, Pf(H'(n))\,$ is the Pfaffian of $\, H'(n) \,$, see \cite{A1}, Thm.3.27. 

If $\, A\,$ is a skew-symmetric matrix let us denote by $\, A_{\widehat{1i}}, \,\, i\ge 2,\,$ the matrix obtained from A by removing the first and the $i$-th row and the first and the $i$-th column. Observe that $\, H'(n+1)_{\widehat{12}} = H'(n)\,$.

Let us denote $\, Pf(n)= Pf(H'(n)), \,\, n\ge 2\,$. We have

\begin{lma}
$\, Pf(n+2)=2 Pf(n+1) + Pf(n)\,$ for $\, n\ge 2\,$.
\end{lma}

\begin{proof} If $\, A=(a_{ij})\,$ is a skew-symmetric $\, 2m\times 2m$-matrix then 
\begin{equation}
Pf(A)= \sum\limits_{i=2}^{2m} (-1)^i a_{1i} Pf(A_{\widehat{1i}}) \,\, ,
\end{equation}
see \cite{A1}, p.142. Applying this to the matrix $\, H'(n+2)\,$ and using $\, H'(n+2)_{\widehat{12}} = H'(n+1)\,$ we get 
\begin{equation}
Pf( H'(n+2))= 2\, Pf( H'(n+1)) +  Pf( H'(n+2)_{\widehat{13}} ) \,\, ,
\end{equation}
where 
\begin{displaymath}
 H'(n+2)_{\widehat{13}} = \left(  
\begin{array}{rrrrrrrr}
0&1&0&0&&&&\\
-1&0&-2&1&&&&\\
0&2&0&2&-1&&&\\
0&-1&-2&0&-2&&&\\
&&1&2&&&&\\
&&&&.&&&1\\
&&&&&.&0&2\\
&&&&&-1&-2&0
\end{array}
\right)
\end{displaymath}
Observe that  $\, ( H'(n+2)_{\widehat{13}} \, )\,_{\widehat{12}}= H'(n)\,$. Using (A.12) again we get 
\begin{displaymath}
  Pf( H'(n+2)_{\widehat{13}} )=  Pf( ( H'(n+2)_{\widehat{13}} \, )\,_{\widehat{12}} )=  Pf( H'(n)) \, \, ,
\end{displaymath}
and, finally, 
\begin{displaymath}
Pf( H'(n+2))= 2\, Pf( H'(n+1)) +  Pf( H'(n) )\,\,  ,
\end{displaymath}
as claimed in the Lemma.
\end{proof}

Now, if $\, A=(a_{ij})\,$ is  a skew-symmetric $2m\times 2m$-matrix then

$\,
\qquad  \text{for} \,\,\, m=1, \,\, \,  Pf(A)=a_{12} \,\, ;
\,$

$\,
\qquad  \text{for} \,\,\, m=2, \,\, \,  Pf(A)=a_{12} \, a_{34}+a_{13} \, a_{42}+ a_{14} \, a_{23}\,\, 
\,$,

\noindent
see \cite{A1}, p.142. 

Since
\begin{displaymath}  
H'(2)= \left(\begin{array}{rr}
0&2\\
-2&0
\end{array}  \right) \qquad \text{and}  \qquad  H'(3)= \left(\begin{array}{rrrr}
0&2&-1&0\\
-2&0&-2&1\\
1&2&0&2\\
0&-1&-2&0
\end{array}  \right)\,\, , 
\end{displaymath}
we get $\, Pf( H'(2) )= 2\,$ and  $\, Pf( H'(3) )= 5\,$.

It follows now from Lemma A.6, by induction on $n$, that $\, Pf(H'(n))> 0\,$ for all $\, n\ge 2\,$.
Consequently, $\, \det ( H'(n)) =( Pf(H'(n)))^2 > 0\,$ and
\begin{displaymath}
\det (H)= (\det (H'(n)))^2 > 0 \qquad \text{for \,\,\, all}\,\,\, n\ge 2\, .
\end{displaymath}
The Hessian matrix $\, H\,$ is therefore invertible. That completes the proof of Proposition A.5.
\end{proof}

Let us recall that  $\, p:P_{2n}\rightarrow SU(2)\,$ is the mapping $\, p(A_1, ... ,A_{2n})= A_1\cdot \cdot \cdot A_{2n}\,$, that $\, K_{2n}=p^{-1}(I)\subset P_{2n}\,$ and that  $\, \Delta = \{ (A,-A, A, -A,...,A,-A)\in P_{2n}\,| \, A\in SU(2), \,\, Tr(A)=0\,\} \subset K_{2n}\,$.

\begin{thm} 
 There is a neighbourhood $\, U\,$ of  $\, \Delta \,$ in $\,K_{2n}\,$ which is homeomorphic to the bundle $\,  G\times _{\mathbb T} C(  S^{2n-3}\times  S^{2n-3})   \,$, where $\, C(  S^{2n-3}\times  S^{2n-3}) \,$ is the cone over the product of two $\,(2n-3)$-dimensional spheres.  The homeomorphism identifies $\, \Delta\,$ with the section of the bundle given by the vertex of the cone.  

\end{thm}

\begin{proof}
Proposition A.5 and Lemma A.1 imply that the function $\, f: M\rightarrow \R\,$ has at $\, 0\in M\,$ a non-degenerate critical point of index equal to $\,\frac 12 \dim M= 2n-2\,$. It follows  that there is a $\, {\mathbb T}$-equivariant neighbourhood $\, W\,$ of the point $\, 0\,$ in $\, M\,$ such that the intersection of the level set $\, f^{-1}(0)\,$ with  $\, W\,$ is $\,  {\mathbb T}$-homeomorphic to the cone  $\, C(  S^{2n-3}\times  S^{2n-3}) \,$
over the product  $\, S^{2n-3}\times  S^{2n-3}  $\, S  of two spheres of dimension $\, 2n-3\,$.

Recall that the level set $\, f^{-1}(0)\,$ is equal to $\, \widehat{p} ^{-1}(I) \cap \widetilde{V}\,$, where $\, \widehat{p} : V \rightarrow SU(2)\,$ is the restriction of the mapping $\, p: P_{2n} \rightarrow SU(2)\,$ to the normal subspace of the orbit $\, \Delta \,$ in $\, P_{2n}\,$ and $\, \widetilde{V}\,$ is a small neighbourhood of $\, 0\,$ in $\, V\,$.  (Choosing $\, \widetilde{V}\,$  small enough we can assume that $\, W=M\,$.) 

In a neighbourhood $\, Y\,$ of $\, \Delta \,$ in $\, P_{2n}\,$, the mapping $\, p\,$ is equivalent to the $\, G$-equivariant extension $\,\widetilde{p}: G\times _{\mathbb T}V \rightarrow SU(2)\,$ of the $\, {\mathbb T}$-equivariant mapping  $\, \widehat{p} : V \rightarrow SU(2)\,$. Therefore,  $\,K_{2n} \cap Y =  p^{-1}(I)\cap Y\,$ is homeomorphic to $\,  G\times _{\mathbb T} f^{-1}(0)\,$ and, hence, to  $\,  G\times _{\mathbb T} C(  S^{2n-3}\times  S^{2n-3})   \,$.  Under this homeomorphism, the inclusion of $\, \Delta \,$ into  $\, p^{-1}(I)\cap Y\,$ corresponds to the inclusion of the ``zero-section'' $\, \Delta = G\times _{\mathbb T}\{{\vec 0}\} \cong G/{\mathbb T}  \,$  into the bundle   $\,  G\times _{\mathbb T} C(  S^{2n-3}\times  S^{2n-3}) \,$, the point $\,{\vec 0}\,$ being the vertex of the cone. Finally, put $\,U = K_{2n} \cap Y\,$. 
\end{proof}    
 
\medskip

\section{Evaluation of the first Chern class I}\label{almI}

\bigskip

In this appendix we shall give a proof of Theorem \ref{almChTh1}.

Recall that in Section \ref{almCh} for every integer $k$ such that $1\le k \le 2n-1\,$ and for $\epsilon = \pm 1$ we have defined the mappings $\gamma _{k, \epsilon}:S^2\rightarrow K$ by 
\begin{displaymath}
\gamma _{k, \epsilon }(A)= (J, J, ... ,J, A, \epsilon A, J, ...,(-1)^n \epsilon J ), \qquad \text{for} \,\,\,  A\in S^2,
\end{displaymath}
where, as in Section 2, $\, S^2\,$ has been identified with  
$\,{\mathcal S}\,$, \,\, $\,  J=  \left( 
\begin{array}{rr} i&0\\ 0&-i 
\end{array}\right)\,$  and the first factor of  $\, A\,$ on the RHS is in the 
$k$-th place. In the case when $\, k=2n-1\,$ the sign $\, (-1)^n \epsilon \,$ is placed at the first factor of $\, J\,$. The mappings $\,\gamma _{k, \epsilon } \,$ are embeddings of $\, S^2\,$  into $\, K\,$ and, hence, \, into $\,  \mathscr{M}\,$ and $\, P_{2n}\,$. 

Let $\ [ \gamma _{k, \epsilon } ]  \in \text{H}_2( \mathscr{M} ; \Z )\,$ be the homology classes represented by  the corresponding mappings.

\medskip

\begin{thm}\label{almITh1}
For all integers $\, k\,$ such that $ \, 1\le k \le 2n-1\,$ and $\, \epsilon =\pm 1\,$, the evaluation of the first Chern class $\, c_1(\mathscr{M})\,$ on the homology classes
 $\ [ \gamma _{k, \epsilon } ]\,$ is equal to $\, 0\,$,
\begin{displaymath}
\langle\, c_1(\mathscr{M}) \, | \,  [ \gamma _{k, \epsilon } ]\, \rangle = 0. 
\end{displaymath}
\end{thm}

\medskip

\begin{proof}
Let $\, G_{2n}\,$ be a product of $2n$ copies of $\, G\,$. $\, P_{2n}\,$ is a submanifold of $\, G_{2n}\,$. The formulas (3.2) and (3.3) make sense and define a smooth $2$-form on $\, G_{2n}\,$. We denote this form by $\, \widehat{\omega_c}\,$, it is an extension of the $2$-form $\, \omega_c\,$ on  $\, P_{2n}\,$.

For an integer $\, k\,$ such that $\,  1\le k \le 2n-1\,$ let  $\, G_{2n}^k\,$ be the submanifold of  $\, G_{2n}\,$ defined by 
\begin{displaymath}
 G_{2n}^k = \{ (g_1, ... ,g_{2n}) \in  G_{2n} \,\, |\,\, g_j\in  {\mathcal S} \,\, \text{for} \,\, j\ne k, k+1 \,\}.
\end{displaymath}
Then  $\, P_{2n}\,$ is a submanifold  of  $\, G_{2n}^k\,$ and, allowing for an abuse of notation, we consider  $\,\gamma _{k, \epsilon }\,$ also as an embedding of $\, S^2\,$ into  $\, G_{2n}^k\,$, $\, \gamma _{k, \epsilon}: S^2\rightarrow  G_{2n}^k\,$. 

Let $\, \xi_{k, \epsilon}\,$ be the restriction of the tangent bundle $\, T G_{2n}^k\,$ to the embedded sphere $\, \gamma _{k, \epsilon }( {\mathcal S}), 
\,\,   \xi_{k, \epsilon} = T G_{2n}^k |_{ \gamma _{k, \epsilon }( {\mathcal S}) }\,$.

Let $\, \widetilde{\omega_c}\,$ be the restriction of the $2$-form $\, \widehat{\omega_c}\,$ to the bundle  $\, \xi_{k, \epsilon}\,$. The vector bundle $\, \tau_{k,\epsilon} = TP_{2n}\, |\, _ { \gamma _{k, \epsilon }( {\mathcal S}) }\,$ is a subbundle of  $\, \xi_{k, \epsilon}\,$ and the restriction of $\, \widetilde{\omega_c}\,$ to this subbundle is equal to the restriction of $\, {\omega_c}\,$, \,\,  $\,\widetilde{\omega_c}\, |\, _{ \tau_{k,\epsilon}} =  \omega_c \,  |\, _{\tau_{k,\epsilon}}   \,$.

The factorization of  $\, G_{2n}^k\,$ as the product $\,  G_{2n}^k = \prod\limits_{j=1}^{2n} \, X_j\,$ with $\, X_k= X_{k+1}=G\,$ and $\, X_j=  {\mathcal S}\,$ for $\, j\ne k, k+1\,$ gives us a splitting of the tangent bundle $\, T G_{2n}^k\,$  as a direct sum 
\begin{displaymath}
 T G_{2n}^k = \bigoplus\limits_{j=1}^{2n} \, \pi _j^* (TX_j)\,\, ,
\end{displaymath}
where $\, \pi _j: G_{2n}^k \rightarrow X_j\,$ is the projection onto the $j$th factor.
When restricted to the embedded  sphere $\,  \gamma _{k, \epsilon }( {\mathcal S})\subset  G_{2n}^k  \,$ it gives us a corresponding splitting 
\begin{equation}\label{almI1}
\xi_{k, \epsilon} = \bigoplus\limits_{j=1}^{2n} \, \zeta _j \,\, ,
\end{equation} 
with 
\begin{equation}\label{almI2}
\zeta _j =  \pi _j^* (TX_j)\, | \,  _{\gamma _{k, \epsilon }( {\mathcal S})}=\left\{\begin{array}{ll} 
\gamma _{k, \epsilon }( {\mathcal S}) \times T_J({\mathcal S}) & \quad \text{if}\,\,\, j\ne k, k+1,\\
\gamma _{k, \epsilon }( {\mathcal S}) \times {\mathfrak g} & \quad \text{if}\,\,\, j=k, k+1\, .
\end{array}\right.
\end{equation} 
Here the identification of $\,\zeta _k\,$ and $\,\zeta _{k+1}\,$ with $\, \gamma _{k, \epsilon }( {\mathcal S}) \times {\mathfrak g}  \,$ is done via the identification of $\, TG\,$  with $\, G \times {\mathfrak g}\,$ through the left translations in $\, G\,$.  That means that every element of $\, TG\,$ is written uniquely in the form $\, g\cdot v\,$ with $\, g\in G\,$ and $\, v\in {\mathfrak g}\,$. (In cases when we shall need to use the identification of $\, TG\,$  with $\, G \times {\mathfrak g}\,$ through the right translations we shall indicate that.)

Let us extend the invariant inner product $\, \bullet\,$ from the Lie algebra $\, {\mathfrak g}\,$ to the tangent bundle $\, TG\,$ by translation. Let $\, s\,$ be a normal vector field to $\, {\mathcal S}\,$ in $\, G\,$ of norm 1 w.r.t. that inner product. One can, for example, take $\, s(J)=  \bigl( 
\begin{smallmatrix} -1&0\\ 0&-1 
\end{smallmatrix}\bigr)=   J\cdot  \bigl( 
\begin{smallmatrix} i&0\\ 0&-i 
\end{smallmatrix}\bigr)\,$, extend it to the whole sphere $\,  {\mathcal S}\,$ by conjugation (this will give a well-defined non-vanishing normal field) and then normalise it. The normal field $\, s\,$ on  $\,  {\mathcal S}\,$ gives us two linearly independent sections $\, s_1\,$ and $ \, s_2\,$ of the bundle $\, \xi_{k, \epsilon}\,$,  the section $\, s_1\,$ in the direct summand $\, \zeta _k\,$ defined by $\, s_1(A_1, ... , A_{2n}) = s(A_k)\,$ and the section $\, s_2\,$ in $\, \zeta _{k+1}\,$ defined by  $\, s_2(A_1, ... , A_{2n}) = s(A_{k+1})\,$.

\medskip

Let $\, A\in {\mathcal S}\,$ and let $\, v_i\,$ and $\, w_m\,$  belong to the fibers of the bundles $\, \zeta _i\,$ and  $\, \zeta _m\,$ over the point $\, \gamma _{k, \epsilon }(A)\,$ respectively, $\, v_i \in \zeta _i^{-1} (\gamma _{k, \epsilon }(A)), \,\, \, w_m \in \zeta _m^{-1} (\gamma _{k, \epsilon }(A)) \,$. 
We shall now compute the values $\, \widetilde{\omega _c}(v_i, w_m )\,$ of the $2$-form $\, \widetilde{\omega _c}\,$ on these elements. As  $\, \widetilde{\omega _c}\,$ is anti-symmetric, we can assume that $\, 1\le i \le m \le 2n\,$.

Recall that, according to (3.3),
\begin{equation}\label{almI3}
\widetilde{\omega _c} = - \sum\limits _{j=1}^{2n-1} \,\, \omega_{\, [z_1...z_j\, |\, z_{j+1} ]} \,\, .
\end{equation} 

Let $\, f_{z_1...z_j}:G_{2n}^k\rightarrow G\,$ be a mapping defined by $\,  f_{z_1...z_j}(A_1, ... , A_{2n}) = A_1\cdot A_2 \cdot ...\cdot A_j\,$ and let  $\, f_{z_{j+1}}:G_{2n}^k\rightarrow G\,$ be the projection onto the $\, (j+1)$-st coordinate.
For any two tangent vectors $\, x\,$ and $\, y\,$ in the same fiber of the tangent bundle $\, TG_{2n}^k\,$ we have 
\begin{equation}\label{almI4}
\begin{split}
\omega_{\, [z_1...z_j\, |\, z_{j+1} ]}\, (x,y)&= \Omega ((df_{z_1...z_j}(x), 
\, df_{z_{j+1}}(x)), \, (df_{z_1...z_j}(y), 
\, df_{z_{j+1}}(y)))\\
&= \frac 12 (\omega (df_{z_1...z_j}(x)) \bullet \bar{\omega} (df_{z_{j+1}}(y)) -\\ &\qquad \qquad \qquad  - \omega (df_{z_1...z_j}(y)) \bullet \bar{\omega} (df_{z_{j+1}}(x))).
\end{split}
\end{equation}

\medskip

\begin{lma}\label{almIL2} 
Let $\, 1\le i \le m\le 2n\,$.

(i) If $\, j+1\ne m\,$ then $\, \omega_{\,[z_1...z_j\, |\, z_{j+1} ]}\, (v_i, \, w_m ) = 0\,$.

(ii) If $\, i=m\,$ then $\, \widetilde{\omega _c} \, (v_m, w_m ) = 0\,$.
\end{lma}

\begin{proof}
Since $\, df_{z_{j+1}} (v_i) = 0\,$ unless $\, j+1=i\,$ and  $\, df_{z_{j+1}} (w_m) = 0\,$ unless $\, j+1=m\,$, it follows from (\ref{almI4}) that $\, \omega_{[z_1...z_j\, |\, z_{j+1} ]}\, (v_i,\, w_m)=0\,$ for $\, j+1\ne i, m\,$.
If $\, j+1=i\,$ then $\, df_{z_1...z_j}(v_i)=  df_{z_1...z_j}(w_m)=0\,$ and again $\, \omega_{[z_1...z_j\, |\, z_{j+1} ]}\, (v_i,\, w_m)=0\,$. The claims of the Lemma follow. 
\end{proof}

Thus, if $\, 1\le i\le m\le 2n\,$, then 
\begin{equation}\label{almI5}
\begin{split}
\widetilde{\omega _c} \, (v_i, w_m )
&= -\, \omega_{\,[z_1...z_{m-1}\, |\, z_m ]}\, (v_i,\, w_m)= \\
&= -\frac 12 \,\omega (df_{z_1...z_{m-1}}(v_i)) 
\bullet \bar{\omega} (df_{z_m}(w_m)).
\end{split}
\end{equation}

Let $\, U=\{ \bigl( \begin{smallmatrix}
0&u\\ -\bar{u} &0
\end{smallmatrix}\bigr) \in {\mathfrak g}\,\, | \,\, u\in \C \, \}\,$. Then $\, T_J( {\mathcal S}) = J\cdot U\,$. 

Let $\, \chi = \epsilon (-1)^n\,$ if $\, k=2n-1\,$ and $\, \chi =1\,$ if $\, k\ne 2n-1\,$ and let $\, \alpha =  \epsilon (-1)^n\,$ if $\, k<2n-1\,$ and $\, \alpha =1\,$ if $\, k=2n-1\,$. 

Since $\, \gamma _{k, \epsilon }(A)= ( \chi J, ... , J, A, \epsilon A, J, ... , \alpha J )\,$, we have, according to the identifications (\ref{almI2}), 
\begin{equation}\label{almI6}
v_i = \left\{ \begin{array}{lll} ( \gamma _{k, \epsilon }(A), \,\, \mu J\cdot v \,) & \quad \text{with} \,\, v\in U &\quad \text{if} \,\, i\ne k, k+1 \\
( \gamma _{k, \epsilon }(A), \,\, A\cdot v \,) & \quad \text{with} \,\, v\in {\mathfrak g} &\quad \text{if} \,\, i=k \\
( \gamma _{k, \epsilon }(A), \,\,\epsilon  A\cdot v \,) & \quad \text{with} \,\, v\in {\mathfrak g} &\quad \text{if} \,\, i=k+1, 
\end{array}\right.
\end{equation} 
and
\begin{equation}\label{almI7}
w_m = \left\{ \begin{array}{lll} ( \gamma _{k, \epsilon }(A), \,\, \lambda J\cdot w \,) & \quad \text{with} \,\, w\in U &\quad \text{if} \,\, m\ne k, k+1 \\
( \gamma _{k, \epsilon }(A), \,\, A\cdot w \,) & \quad \text{with} \,\, w\in {\mathfrak g} &\quad \text{if} \,\, m=k \\
( \gamma _{k, \epsilon }(A), \,\,\epsilon  A\cdot w \,) & \quad \text{with} \,\, w\in {\mathfrak g} &\quad \text{if} \,\, m=k+1. 
\end{array}\right.
\end{equation} 
Here $\, \mu= \epsilon (-1)^n\,$ if $\, i=1, \, k=2n-1\,$ or $\, i=2n, \, k<2n-1\,$ and $\, \mu =1\,$ otherwise, $\, \lambda = \epsilon (-1)^n\,$  if $\, m=2n, \, k<2n-1\,$ or $\, m=1, \, k=2n-1\,$ and $\, \lambda =1\,$ otherwise.

\medskip

\begin{lma}\label{almIL3}
Let $\, 1\le i \le m\le 2n\,$. Then
\begin{displaymath}
\widetilde{\omega _c} \, (v_i, w_m ) = \bigl(-\frac 12\bigr)\cdot \left\{
\begin{array}{ll}
0& \quad \text{if}\,\, \, i=m,\\
(-1)^{m-i-1}\, v\bullet (Ad(A)(w))&\quad \text{if} \,\,\, i<m=k,\\
(Ad(A)(v))\bullet w &\quad \text{if} \,\,\, i=k, \, m=k+1,\\
-(Ad(J^{m-i-2})(Ad(A)(v)))\bullet w &\quad \text{if} \,\,\, i=k<m-1,\\
-(Ad(J^{m-i-1})(v))\bullet w &\quad \text{otherwise.}
\end{array}\right.
\end{displaymath}
Here $\, v\,$ and $\, w\,$ are as in (\ref{almI6}) resp. (\ref{almI7}).
\end{lma}

\begin{proof}
The claim of the Lemma follows from (\ref{almI5}) by a direct check of nine cases. Observe that $\, Ad(J)\,$ acts on $\, U\subset {\mathfrak g}\,$ by multiplication by $\, -1\,$.

\smallskip

\noindent
Case 1: if $\, i=m\,$ then $\, \widetilde{\omega _c} \, (v_i, w_m ) = 0\,$ by Lemma \ref{almIL2} (ii).

\smallskip

\noindent
Case 2: if $\, i<m<k\,$ then $\, v_i=(\gamma _{k, \epsilon }(A), \,\, \mu J\cdot v \,)\,$ with $\, v\in U\subset {\mathfrak g}\,$ and $\, df_{z_1...z_{m-1}}(v_i)=\chi J^{m-1}\cdot (Ad(J^{m-1-i})^{-1}(v))=\chi J^{m-1}\cdot (Ad(J^{m-1-i})(v))\,$. Similarly, $\, w_m=(\gamma _{k, \epsilon }(A), \,\, \lambda J\cdot w \,)\,$ with $\, w\in U\,$  and $\, df_{z_m}(w_m) = \lambda J\cdot w = (Ad(J)(w))\cdot \lambda J = (-w)\cdot J\,$. (Here we use the right translation by $ J$, and the fact that $\, \lambda =1\,$ in this case.) It follows that $\, \widetilde{\omega _c} \, (v_i, w_m ) = -\frac 12\, \omega (df_{z_1...z_{m-1}}(v_i))\bullet \bar{\omega }( df_{z_m}(w_m)) = -\frac 12 (-1)^{m-1-i}\, v\bullet (-w)=-\frac 12 (-1)(Ad(J^{m-1-i})(v))\bullet w\,$, as claimed.

\smallskip

\noindent
Case 3: if $\, i<m=k\,$ then  $\, df_{z_1...z_{m-1}}(v_i)=\chi J^{m-1}\cdot (Ad(J^{m-1-i})^{-1}(v))=\chi J^{m-1}\cdot ((-1)^{m-1-i} \, v)\,$ with $\, v\in U\,$ just as in Case 2, while  $\, w_m=(\gamma _{k, \epsilon }(A), \,\, A\cdot w \,)\,$ with $\, w\in {\mathfrak g}\,$,  $\, df_{z_m}(w_m) = A\cdot w = (Ad(A)(w))\cdot A\,$\,\, (again, here we use the right translation by $\, A\,$) and $\, \widetilde{\omega _c} \, (v_i, w_m ) = -\frac 12\, \omega (\chi J^{m-1}\cdot ((-1)^{m-1-i} \, v))\bullet \bar{\omega }((Ad(A)(w))\cdot A)=(-\frac 12)(-1)^{m-1-i}\, v\bullet (Ad(A)(w))\,$ as claimed.

\smallskip

\noindent
Case 4: if $\, i<k, \, m=k+1\,$ then  $\, df_{z_1...z_{m-1}}(v_i)=\chi J^{m-2}A \cdot (Ad(J^{m-2-i}A)^{-1}(v))\,$ with $\, v\in U\,$ and  $\, df_{z_m}(w_m) = \epsilon A\cdot w = (Ad(A)(w))\cdot (\epsilon A)\,$ with $\, w\in {\mathfrak g}\,$. Consequently, 
\begin{equation*}
\begin{split}
\widetilde{\omega _c} \, (v_i, w_m )&=(-\frac 12 )(Ad(J^{m-2-i}A)^{-1}(v))\bullet (Ad(A)(w))= \\ &=(-\frac 12 )(Ad(-A)Ad(J^{m-2-i})^{-1}(v))\bullet (Ad(A)(w))=\\
&= (-\frac 12 )(Ad(A)Ad(J^{m-2-i})(v))\bullet (Ad(A)(w))=\\
&= (-\frac 12 )(-1) (Ad(J^{m-1-i})(v))\bullet w),
\end{split}
\end{equation*}
as claimed. The last equality follows by invariance of the inner product $\, \bullet \,$, the previous one from the fact that $\, Ad(A)=Ad(-A)=Ad(A^{-1})\,$ for all $\, A\in {\mathcal S}\,$.

\smallskip

\noindent
Case 5: if $\, i=k, \, \, m=k+1\,$ then  $\, df_{z_1...z_{m-1}}(v_i)=\chi J^{m-2}A \cdot v\,$ with $\, v\in {\mathfrak g}\,$ and  $\, df_{z_m}(w_m) = \epsilon A\cdot w = (Ad(A)(w))\cdot (\epsilon A)\,$ with $\, w\in {\mathfrak g}\,$. Thus 
\begin{equation*}
\begin{split}
\widetilde{\omega _c} \, (v_i, w_m )&=(-\frac 12 )\omega (\chi J^{m-2} A\cdot v) \bullet \bar{\omega }((Ad(A)(w))\cdot (\epsilon A))= \\ 
&=(-\frac 12 )\, v\bullet Ad(A)(w)=\\
&= (-\frac 12 )(Ad(A^{-1})(v)\bullet w)  \qquad (\, \text{by invariance of}\,\,\, \bullet \,\, )\\
&= (-\frac 12 )(Ad(A)(v)\bullet w)  \qquad (\, \text{since}\,\,\,  Ad(A^{-1})=Ad(-A)=Ad(A)\, ),
\end{split}
\end{equation*}
as claimed.

\smallskip

\noindent
Case 6: if $\, i<k<m-1\,$ then  $\, df_{z_1...z_{m-1}}(v_i)=\chi \epsilon J^{m-1} \cdot (Ad(J^{m-1-i})(v)), \,\, v\in U\,$ and  $\, df_{z_m}(w_m) = \lambda J\cdot w = Ad(J)(w) \cdot (\lambda J) = (-w)\cdot (\lambda J), \,\, w\in U\,$. Thus 
\begin{displaymath}
\widetilde{\omega _c} \, (v_i, w_m )=(-\frac 12 )\,( Ad(J^{m-1-i})(v)) \bullet (-w)
\end{displaymath}
as claimed.

\smallskip

\noindent
Case 7: if $\, i=k<m-1\,$  then  $\, df_{z_1...z_{m-1}}(v_i)=\chi \epsilon J^{m-1} \cdot (Ad(AJ^{m-2-i})^{-1}(v)), \,\, v\in {\mathfrak g}\,$ and  $\, df_{z_m}(w_m) = (-w)\cdot (\lambda J), \,\, w\in U\,$. Thus
\begin{equation*}
\begin{split}
\widetilde{\omega _c} \, (v_i, w_m )&=(-\frac 12 )(Ad(AJ^{m-2-i})^{-1}(v))\bullet (-w)=\\ &= (-\frac 12 )(-1)(Ad(J^{m-2-i})^{-1}(Ad(A^{-1})(v)))\bullet w = \\ 
&=(-\frac 12 )(-1)(Ad(J^{m-2-i})(Ad(A)(v)))\bullet w,
\end{split}
\end{equation*}
as claimed.

\smallskip

\noindent
Case 8: if $\, i=k+1<m\,$ then $\, df_{z_1...z_{m-1}}(v_i)=\epsilon J^{m-1} \cdot (Ad(J^{m-1-i})(v)), \,\, v\in {\mathfrak g}\,$ and  $\, df_{z_m}(w_m) = \lambda J\cdot w= (-w)\cdot (\lambda J), \,\, w\in U\,$. Therefore 
\begin{displaymath}
\widetilde{\omega _c} \, (v_i, w_m )= (-\frac 12 )(Ad(J^{m-1-i})(v))\bullet (-w),
\end{displaymath}
as claimed.

\smallskip

\noindent
Case 9: if $\, k+1<i<m\,$ then $\, df_{z_1...z_{m-1}}(v_i)=\epsilon J^{m-1} \cdot (Ad(J^{m-1-i})(v)), \,\, v\in U\,$  and  $\, df_{z_m}(w_m) = \lambda J\cdot w= (-w)\cdot (\lambda J), \,\, w\in U\,$. Therefore 
\begin{displaymath}
\widetilde{\omega _c} \, (v_i, w_m )= (-\frac 12 )(Ad(J^{m-1-i})(v))\bullet (-w),
\end{displaymath}
again as claimed.

This concludes the proof of Lemma \ref{almIL3}.

\end{proof}

\medskip

Let $\, W\,$ be the fibre of the vector bundle $\, \xi _{k,\epsilon }\,$ over the point $\, \gamma _{k,\epsilon }(J)\,$. Then $\, W=\bigoplus\limits _{j=1}^{2n} V_j\,$ is the direct sum of real vector spaces
\begin{equation}\label{almI8}
V_j=\left\{ \begin{array}{ll}
T_J({\mathcal S}) & \qquad j\ne k, k+1,\\
{\mathfrak g}& \qquad j=k, k+1.
\end{array}\right.
\end{equation}
The identities (\ref{almI1}) and (\ref{almI2}) yield a trivialization of the bundle  $\, \xi _{k,\epsilon }\,$
\begin{displaymath}
\psi _{k,\epsilon }:  \xi _{k,\epsilon } \xrightarrow{\cong }\gamma _{k,\epsilon }({\mathcal S})\times W
\end{displaymath}
which is the identity on the fiber over  $\, \gamma _{k,\epsilon }(J)\,$.

Let $\, \widetilde{\Phi}_k\,$ be the automorphism of the bundle $\, \gamma _{k,\epsilon }({\mathcal S})\times W\,$ given by 
\begin{displaymath}
 \widetilde{\Phi}_k (\gamma _{k,\epsilon }(A), \,\, \bigoplus\limits _{j=1}^{2n} v_j \,) = (\gamma _{k,\epsilon }(A), \,\, \bigoplus\limits _{j=1}^{2n} \widetilde{\phi}_k (A, \, v_j \,) )
\end{displaymath}
with 
\begin{displaymath}
 \widetilde{\phi}_k (A, \, v_j \,)  = \left\{\begin{array}{cl}
v_j & \qquad \text{if}\quad j\ne k,\\
Ad(-JA)(v_k) & \qquad  \text{if}\quad j= k,
\end{array}
\right.
\end{displaymath}
for $\, A\in {\mathcal S}, \, \, v_j\in V_j\,$.

Denote by  $\, {\Phi}_k\,$ the corresponding automorphism of the bundle  $\, \xi _{k,\epsilon }\,$.

$\, \widetilde{\omega _c}\,$ is a $2$-form on the bundle   $\, \xi _{k,\epsilon }\,$. Denote by $\, \omega _0\,$ its restriction to the fibre $\, W\,$ over the point $\,  \gamma _{k,\epsilon }(J)\,$.  The form  $\, \omega _0\,$ gives us the product form  $\, \widehat{\omega} _0\,$ on the trivial bundle $\, \gamma _{k,\epsilon }({\mathcal S})\times W\,$ and, via the isomorphism $\, \psi _{k,\epsilon }\,$, a $2$-form $\, \widetilde{\omega _0}\,$ on the bundle  $\, \xi _{k,\epsilon }\,$.

Let $\, ( \xi _{k,\epsilon }, \,\, \widetilde{\omega _c})\,$ and  $\, ( \xi _{k,\epsilon }, \,\, \widetilde{\omega _0})\,$ denote the real vector bundle  $\, \xi _{k,\epsilon }\,$ equipped with the $2$-forms  $\, \widetilde{\omega _c}\,$ and  $\, \widetilde{\omega _0}\,$ respectively.

\smallskip

\begin{lma}\label{almIL4}

(i) The automorphism $\,  {\Phi}_k: \xi _{k,\epsilon }\rightarrow  \xi _{k,\epsilon } \,$ of the real vector bundle  $\, \xi _{k,\epsilon }\,$  

\qquad \qquad \qquad \,  satisfies $\,  {\Phi}_k^*(  \widetilde{\omega _0}) = \widetilde{\omega _c}\,$.

\qquad \quad \quad\,  (ii) $\,  \widetilde{\omega _0}\,$ is a non-degenerate $2$-form on  $\, \xi _{k,\epsilon }\,$.

\qquad \quad \quad  (iii)  $\,  \widetilde{\omega _c}\,$ is a non-degenerate $2$-form on  $\, \xi _{k,\epsilon }\,$.
\end{lma}

\begin{proof}
(i) Observe that when $\, A=J\,$ the expressions of Lemma \ref{almIL3} give us a formula for the $2$-form  $\, {\omega _0}\,$. Moreover, Lemma \ref{almIL3} shows that for any $\, A\in {\mathcal S}\,$ and $\, v,w\in  \xi _{k,\epsilon }^{-1} ( \gamma _{k,\epsilon }(A))\,$ one has
\begin{displaymath}
\widetilde{\omega _c} (v,w) =  \widetilde{\omega _0} ({\Phi}_k(v),\, {\Phi}_k(w)). 
\end{displaymath}
This shows that $\,  {\Phi}_k^*(  \widetilde{\omega _0}) = \widetilde{\omega _c}\,$.

\smallskip

(ii) The vector bundle  $\, \xi _{k,\epsilon }\,$  contains as a subbundle of real codimension $2$ the restriction $\, \tau_{k,\epsilon } = TP_{2n} \, | \, _{\gamma _{k,\epsilon }({\mathcal S}) }\,$ of the tangent bundle $\, TP_{2n}\,$ of $\, P_{2n}\,$ to the sphere $\,\gamma _{k,\epsilon }({\mathcal S})\,$. The restriction of the $2$-form $\, \widetilde{\omega _c} \,$ to $\,  \tau_{k,\epsilon }\,$ is equal to the restriction of the symplectic form   $\,{\omega _c} \,$ to 
 $\,  \tau_{k,\epsilon }\,$ and, hence, is non-degenerate by Theorem 3.1 and Corollary 3.4.

Let $\, s_1(\gamma _{k,\epsilon }(J)), \, s_2(\gamma _{k,\epsilon }(J))\in W\,$ be the values of the sections $\, s_1, \, s_2\,$ at the point $\, \gamma _{k,\epsilon }(J)\,$. In the decomposition (\ref{almI8}) the value $\, x_1= s_1(\gamma _{k,\epsilon }(J))\,$ is contained in $\, V_k={\mathfrak g}\,$ while  $\, x_2= s_2(\gamma _{k,\epsilon }(J))\,$ is contained in $\, V_{k+1}={\mathfrak g}\,$ and both are orthogonal to $\, U\subset {\mathfrak g}\,$ w.r.t. the inner product $\, \bullet \,$.

All the elements of the fibre of $\,  \tau_{k,\epsilon }\,$ over the point $\,  \gamma _{k,\epsilon }(J)\,$ belong to 
$\, \bigoplus\limits _{j=1}^{2n} U \,$ in the decomposition (\ref{almI8}). Moreover, the subspace $\, U\subset {\mathfrak g}\,$ is invariant under the action of $\, Ad(J)\,$ and is orthogonal to $\, s(J)\,$ w.r.t. the inner product $\, \bullet \,$. It follows then from the formulas of Lemma \ref{almIL3} that $\, x_1\,$ and $\, x_2\,$ are orthogonal w.r.t. the $2$-form $\, \widetilde{\omega _c}\,$ (and hence  w.r.t. the $2$-form $\, \widetilde{\omega _0}\,$) to the fibre of $\,   \tau_{k,\epsilon }\,$ over 
 $\, \gamma _{k,\epsilon }(J)\,$. Furthermore, again from Lemma \ref{almIL3}, case $\, i=k, \, m=k+1\,$, we have 
\begin{displaymath}
 \widetilde{\omega _0}(x_1,\, x_2)=  \widetilde{\omega _c}(x_1,\, x_2)= (Ad(J)(x_1))\bullet x_2 \,\, .
\end{displaymath}
Since $\, s(J)\,$ is normal to $\, T_J({\mathcal S})\,$, we have $\, Ad(J)(x_1) = x_1\,$. Moreover, by the choice of the sections $\, s_1, s_2\,$, we have $\, x_1 \bullet x_2 =1\,$. Therefore
\begin{displaymath}
 \widetilde{\omega _0}(x_1,\, x_2)= x_1\bullet x_2 = 1\,\, .
\end{displaymath}
It follows that the $2$-form $\, \omega _0\,$ on the space $\, W\,$ is an orthogonal direct sum of the non-degenerate form  $\, \omega _c\,$ on the fibre of 
 $\,   \tau_{k,\epsilon }\,$ and a hyperbolic plane spanned by $\, x_1\,$ and $\, x_2\,$. Therefore $\, \omega _0\,$ is non-degenerate on  $\, W\,$. It follows that the product $2$-form $\,  \widetilde{\omega _0}\,$ is non-degenerate on $\, \xi _{k,\epsilon }\,$. That proves (ii).

\smallskip

(iii) Since, according to (i), $\,  \widetilde{\omega _c}\,$ is isomorphic to  $\,  \widetilde{\omega _0}\,$ through $\, \Phi _k\,$, it follows from (ii) that 
 $\,  \widetilde{\omega _c}\,$ is non-degenerate as well.
\end{proof}

\medskip

As a consequence of Lemma \ref{almIL4} we see that   $\, ( \xi _{k,\epsilon }, \,\, \widetilde{\omega _c})\,$ and  $\, ( \xi _{k,\epsilon }, \,\, \widetilde{\omega _0})\,$ are isomorphic symplectic vector bundles. Moreover, there is a trivialization of  $\,  \xi _{k,\epsilon }\,$ in which the $2$-form $\,  \widetilde{\omega _0}\,$ is constant i.e. $\, ( \xi _{k,\epsilon }, \,\, \widetilde{\omega _0})\,$ is a trivial symplectic vector bundle. It follows that the first Chern class of both these bundles is the same and that it vanishes,
\begin{equation}\label{almI9}
c_1 ( \xi _{k,\epsilon }, \,\, \widetilde{\omega _c}) = c_1 ( \xi _{k,\epsilon }, \,\, \widetilde{\omega _0}) = 0\, .
\end{equation}

The bundle $\, \tau _{k,\epsilon }\,$ is a vector subbundle of  $\, ( \xi _{k,\epsilon }, \,\, \widetilde{\omega _c})\,$ of real codimension $2$ and $\, \widetilde{\omega _c}\, |\, _{ \tau _{k,\epsilon } } = {\omega _c}\, |\, _{ \tau _{k,\epsilon } }\,$ is a non-degenerate $2$-form on $\, \tau _{k,\epsilon }\,$ by Theorem 3.1 and Corollary 3.4. Hence  $\, \tau _{k,\epsilon }\,$ is a symplectic subbundle of  $\, ( \xi _{k,\epsilon }, \,\, \widetilde{\omega _c})\,$. Let $\, \eta\,$ be the symplecto-orthogonal complement of  $\, \tau _{k,\epsilon }\,$ in  $\, ( \xi _{k,\epsilon }, \,\, \widetilde{\omega _c})\,$. Then $\, \eta\,$ is a symplectic vector bundle of real dimension $2$ and there is an isomorphism of symplectic vector bundles between   $\, ( \xi _{k,\epsilon }, \,\, \widetilde{\omega _c})\,$ and the orthogonal direct sum of symplectic bundles $\, \tau _{k,\epsilon } \oplus \eta \,$, 
\begin{equation}\label{almI10}
\xi _{k,\epsilon }\cong \tau _{k,\epsilon } \oplus \eta 
\end{equation}

Let us choose complex structures on  $\, \tau _{k,\epsilon }\,$ and on $\, \eta \,$ compatible with their symplectic forms. The isomorphism (\ref{almI10}) equipps then  
$\, \xi _{k,\epsilon }\,$ with a complex structure compatible with the $2$-form $\,\widetilde{\omega _c}\,$ and  
\begin{equation}\label{almI11}
c_1(\xi _{k,\epsilon }, \,\, \widetilde{\omega _c})= c_1(  \tau _{k,\epsilon } \oplus \eta ) =    c_1(\tau _{k,\epsilon }) + c_1( \eta ) \, .
\end{equation}
On the other hand the values of the section $\, s_1\,$ of $\, \xi _{k,\epsilon }\,$ never lie in  $\, \tau _{k,\epsilon }\,$. Therefore the decomposition of  $\, s_1\,$ in the direct sum (\ref{almI10}) gives us a nowhere vanishing section $\, \widetilde{s}\,$ of  $\, \eta \,$. Since  $\, \eta \,$ is a $1$-dimensional complex bundle it follows  that as a complex bundle  $\, \eta \,$  is trivial and $\, c_1( \eta ) = 0\,$. Therefore, from (\ref{almI9}) and (\ref{almI11}) we get
\begin{displaymath}
 c_1(\tau _{k,\epsilon }) = c_1(\xi _{k,\epsilon }, \,\, \widetilde{\omega _c}) -  c_1( \eta ) = 0 \, .
\end{displaymath}
As  $\, \langle c_1 \, | \,  [ \gamma _{k, \epsilon } ] \rangle =  c_1(\tau _{k,\epsilon })  \,$,   that proves Theorem \ref{almITh1}.

\end{proof}

The subspace $\, K_{2n}\,$ of $\, P_{2n}\,$ is $\, G$-invariant, thus $\, G=SU(2)\,$ acts on $\, K_{2n}\,$. Let us consider the projection onto the orbit space $\, q: K_{2n}\rightarrow  K_{2n}/G\,$.

\begin{lma}\label{almIL5}
Let $\, n\ge 2\,$. The composed mappings $\, q\circ \gamma _{k,\epsilon} : S^2 \rightarrow K_{2n}/G\,$ are contractible  for all $\, 1\le k\le 2n-1\,$ and $\, \epsilon =\pm 1\,$. 
\end{lma}
  
\begin{proof}
The mapping $\, \gamma _{k,\epsilon}\,$ is an embedding of the 2-dimensional sphere $\ S^2\,$  into  $\, K_{2n}\,$ with the image 
\begin{displaymath}
\gamma _{k,\epsilon} (S^2)=\{ \, (J, ... , J, A, \epsilon A, J, ... ,(-1)^n\epsilon J) \in P_{2n} \, \, \vert \, \, A\in {\mathcal S} \, \} \, .
\end{displaymath}
Let $\, {\mathbb T}\,$ be the $1$-dimensional torus in $\, SU(2)\,$ which is the isotropy subgroup of the point $\, J\,$ in $\, G\,$ under the conjugacy action of $\, G\,$ on itself. ($\, {\mathbb T}\,$ consists of the diagonal matrices in $\,  SU(2)\,$.) The torus $\, {\mathbb T}\,$ is the largest subgroup of $\, G\,$ keeping the subset $\,\gamma _{k,\epsilon} (S^2)\,$ invariant and no element of $\, G- {\mathbb T}\,$  maps any point of  $\,\gamma _{k,\epsilon} (S^2)\,$ into $\,\gamma _{k,\epsilon} (S^2)\,$. Hence, the image of   $\,\gamma _{k,\epsilon} (S^2)\,$ in  $\,  K_{2n}/G\,$ under the quotient projection $\, q\,$  is equal to  $\,\gamma _{k,\epsilon} (S^2)/{\mathbb T} \,$. But, 
$\,\gamma _{k,\epsilon} (S^2)/{\mathbb T} \,$ is homeomorphic to 
$\, {\mathcal S}/{\mathbb T} \,$. The space $\, {\mathcal S} /{\mathbb T} \,$ is homeomorphic to an interval and hence is contractible. Thus the image of  $\,\gamma _{k,\epsilon} (S^2)\,$ under $\, q\,$ is a contractible subspace of $\, K_{2n}/G\,$. Consequently, the composition  $\, q\circ \gamma _{k,\epsilon} : S^2 \rightarrow K_{2n}/G\,$ is a contractible mapping.
\end{proof}
\medskip

\section{Evaluation of the first Chern class II}\label{almII}

\bigskip

Let $\, \widetilde{{\mathcal S}}\,$ be the smooth manifold  diffeomorphic to the $2$-dimensional sphere defined in Section \ref{almCh} (see Figure \ref{almCh}.1). Let $\, f=f_n: \widetilde{{\mathcal S}} \rightarrow K_{2n}\,$ be the mapping defined by (\ref{almCh2}) and   (\ref{almCh3}).

Consider the homomorphism  $\, f_*: H_2( \widetilde{{\mathcal S}}, \, \Z )\rightarrow  H_2( K_{2n}, \, \Z )\,$ induced by $\, f\,$ on the second homology groups. The orientation of $\,  \widetilde{{\mathcal S}}\,$ gives us $\, [ \widetilde{{\mathcal S}}]\in H_2( \widetilde{{\mathcal S}}, \, \Z )\,$. Let us consider the element $\, f_*  [ \widetilde{{\mathcal S}}] \in  H_2( K_{2n}, \, \Z )\,$.

Let $\, {\mathscr M} \subset P_{2n}\,$ be the open neighbourhood of $\, K_{2n}\,$ given by Theorem 3.1.  $\, {\mathscr M}\,$ is a symplectic manifold with the symplectic form $\, \omega_{\mathcal C}\,$. We assume that  $\, {\mathscr M}\,$ has been choosen so that it is homotopy equivalent to $\, K_{2n}\,$. Let $\, c_1 (  {\mathscr M}) \in H^2( {\mathscr M}, \Z )=H^2 (K_{2n}, \Z )\,$ be the first Chern class of the symplectic structure on  $\, {\mathscr M}\,$.

\medskip

\begin{thm}\label{almIITh1}
For every $\, n\ge 2\,$
\begin{displaymath}
\langle \,  c_1 (  {\mathscr M})   \, \vert \,  f_* [ \widetilde{{\mathcal S}}]     \, \rangle = -2.
\end{displaymath}
\end{thm}

\begin{proof}
Let us denote by $\, K'_4\,$ the subspace of $\, K_{2n}\,$ consisting of points $\, (A_1, ... , A_{2n})\in K_{2n}\,$ such that $\, A_j=(-1)^j J\,$ for $\, j\ge 5\,$.  Observe that such points satisfy the condition 
\begin{equation}\label{7.3}
A_1\cdot A_2\cdot A_3\cdot A_4 = I.
\end{equation}
By the definition of the mapping $\, f\,$, its image is contained in $\,  K'_4\,$.

Let $\, \xi\,$ be the restriction of the tangent bundle of $\, P_{2n}\,$ to  $\, K'_4\,$, $\, \xi = T P_{2n} \vert _ {K'_4} =  T {\mathscr M} \vert _ {K'_4}\,$. It is a symplectic vector bundle over  $\, K'_4\,$ with the symplectic form $\, \omega _{\mathcal C} = \omega _c\,$ (see Corollary 3.4). The product decomposition of $\, P_{2n} = \prod {\mathcal S}\,$ gives us a canonical identification of vector bundles 
\begin{equation}\label{7.4}
\xi \cong \bigoplus\limits _{j=1}^{2n} \, p^*_j(T{\mathcal S}) \vert _ {K'_4}\,\,\, . 
\end{equation}
We denote by  $\,dp_j: \xi \rightarrow  p^*_j(T{\mathcal S}) \vert _ {K'_4}\,$ the projection onto the $j$-th summand.

Let us consider the vector subbundle $\, \xi _1 =  \bigoplus\limits _{j=5}^{2n} \, p^*_j(T{\mathcal S}) \vert _ {K'_4}\,$ of $\, \xi \,$. Since $\, p_j(x) = (-1)^j J\,$ for every $\, j\ge 5\,$ and $\, x\in K'_4\,$, we have a trivialization
\begin{equation}\label{7.5}
\Phi : \xi _1 \rightarrow K'_4 \times  \bigoplus\limits _{j=5}^{2n} \, V_j  \,\, , 
\end{equation}
where $\, V_j = T_{(-1)^j J}{\mathcal S}\,$ is the tangent space to $\, {\mathcal S}\,$ at the point $\, (-1)^j J\,$. 

Let us denote $\, V = \bigoplus\limits _{j=5}^{2n} \, V_j \,$. We identify $\, \xi _1\,$ with $\,  K'_4 \times V\,$ through $\, \Phi\,$.

\begin{lma}
The restriction of $\, \omega _c\,$ to $\, \xi _1\,$ is a non-degenerate $2$-form on  $\, \xi _1\,$ and it does not depend on the first coordinate $\, x\in  K'_4\,$ in $\, \xi _1 =  K'_4 \times V\,$.
\end{lma}

\begin{proof}
There is an obvious identification of the vector space $\, V\,$ with the tangent space $\, T_YP_{2n-4}\,$ to $\, P_{2n-4}\,$  at the point $\, Y=(-J, J, -J, ... ,-J, J)\in P_{2n-4} \,$. It follows from the definition (3.3) of the $2$-form $\, \omega _c = - \sum \omega _{[z_1...z_j\vert z_{j+1}]}\,$ and from the identity (\ref{7.3}) that, for every point  $\,  x\in K'_4\,$, this identification maps the restriction of $\, \omega _c\,$ to the fiber of $\, \xi _1\,$ over $\, x\,$ to the $2$-form $\, \omega _c\,$ on $\, T_YP_{2n-4}\,$. Therefore the restriction of $\, \omega _c\,$ to $\, \xi _1\,$ does not depend on $\, x\in  K'_4\,$. Moreover, since the $2$-form $\, \omega _c\,$ on  $\, T_YP_{2n-4}\,$ is non-degenerate (Theorem 3.1), so is the restriction of   $\, \omega _c\,$ to $\, \xi _1\,$.   
\end{proof}

It follows that $\, \xi _1\,$ is a symplectic subbundle of  $\, \xi \,$ and, as a symplectic bundle, $\, \xi _1\,$ is trivial.  Therefore 
\begin{equation}\label{7.6}
c_1(\xi _1) = 0 \, .
\end{equation}

Let $\, \xi _2\,$ be the symplectic orthogonal complement of $\, \xi _1\,$ in 
 $\, \xi \,$. Thus  $\, \xi _2\,$ is a symplectic subbundle of  $\, \xi \,$, of real dimension $8$, and
\begin{equation}\label{7.7}
\xi = \xi _2 \oplus \xi _1
\end{equation} 
as symplectic bundles. Therefore 
\begin{equation}\label{7.8}
c_1(\xi ) = c_1(\xi _2) + c_1(\xi _1) =  c_1(\xi _2)\, .
\end{equation}

Let $\, \zeta \,$ be the vector subbundle of $\, \xi\,$ given by 
\begin{equation}\label{7.9}
\zeta = \bigoplus\limits _{j=1}^4 \, p^{*}_j (T{\mathcal S})\vert _{K'_4}
\end{equation} 
in the decomposition (\ref{7.4}).  Then $\, \xi = \zeta \oplus \xi _1\,$ as vector bundles. Let $\, \varphi : \zeta \rightarrow \xi _2\,$ be the projection of $\, \zeta \,$ onto $\, \xi _2\,$ in the decomposition (\ref{7.7}). It is an isomorphism of vector bundles. See Figure \ref{almII}.1. 

\vskip1truecm
$$\includegraphics[width=7cm, height=3cm]{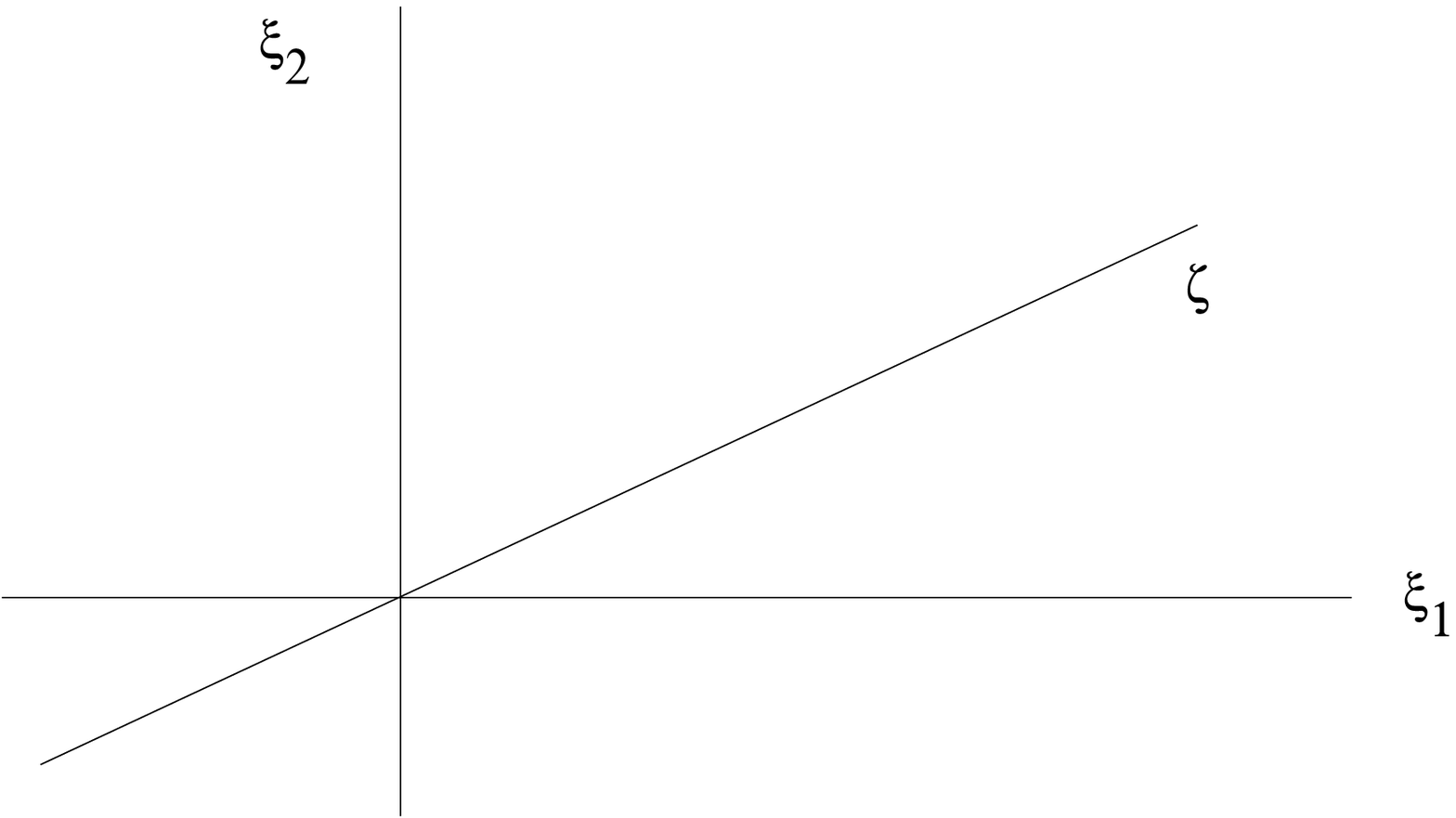}$$

\bigskip
\begin{center}
Fig. \ref{almII}.1
\end{center}

\smallskip

We shall now describe the isomorphism $\, \varphi\,$ explicitly. 
Let  $\, U\,$ be the subspace of  $\,{\mathfrak g}= su(2)\,$  consisting of matrices of the form $\, \bigl(\begin{smallmatrix} 0&u\\-\overline{u}&0\end{smallmatrix}\bigr), \,\, u\in \C\,$, and let $\, \text{proj}:{\mathfrak g}\rightarrow U\,$  be the orthogonal projection of  $\,{\mathfrak g}\,$ onto  $\, U\,$ w.r.t. the inner product $\, \bullet\,$. As  $\, A\bullet B = -\frac 12 \, \alpha \,\text{tr}(AB)\,$ for some $\, \alpha \in \R, \,\, \alpha >0\,$, we see that $\,\text{proj} \bigl(\begin{smallmatrix} it&z\\-\overline{z}&-it\end{smallmatrix}\bigr)= \bigl(\begin{smallmatrix} 0&z\\-\overline{z}&0\end{smallmatrix}\bigr)\,$ for $\,t\in \R, \,\, z\in \C\,$.

Let $\, x=(A_1, A_2, A_3, A_4, -J, J, ... , -J, J)\,$ be an arbitrary point in $\, K'_4\,$, with $\, A_i\in {\mathcal S}, \,\, A_1A_2A_3A_4 = I\,$, and let  
\begin{displaymath}
{\mathbf v}= {\mathbf v}_1\oplus {\mathbf v}_2\oplus {\mathbf v_3\oplus }{\mathbf v}_4 \in \zeta _x
\end{displaymath} 
be an arbitrary vector in the fiber $\, \zeta _x\,$ of $\, \zeta \,$  over $\ x \,$. Here $\, {\mathbf v}_j\in T_{A_j}{\mathcal S}\,$ and therefore $\,  {\mathbf v}_j= v_j\cdot A_j\,$ for some $\, v_j\in {\mathfrak g}, \,\, j=1, ... ,4\,$. 

Since $\, \varphi : \zeta \rightarrow \xi _2 \subset \xi = \bigoplus\limits _{j=1}^{2n} \, p^*_j(T{\mathcal S}) \vert _{K'_4}\,$,\,\, we have 
\begin{displaymath}
\varphi ({\mathbf v}) = \varphi ({\mathbf v}_1\oplus {\mathbf v}_2\oplus {\mathbf v_3\oplus }{\mathbf v}_4) = \bigoplus\limits_{j=1}^{2n} {\mathbf w}_j \, \, ,
\end{displaymath}
with $\,  {\mathbf w}_j= w_j\cdot A_j \in T_{A_j}{\mathcal S}\,$, some $\, w_j\in {\mathfrak g}, \,\, j=1, ... , 2n\,$  (and where $\, A_j=(-1)^jJ\,$ for $\, j\ge 5\,$).

The explicit description of the isomorphism $\, \varphi\,$ is given by
\begin{lma}{\label{L7.3}}
\begin{displaymath}
w_j = \begin{cases}
\quad v_j & \qquad {\rm if} \,\,\, j=1, ... ,4,\\
-\sum\limits _{i=1}^4 \, {\rm proj} \bigl( Ad(A_4A_3...A_i)(v_i)\bigr) & \qquad  {\rm if} \,\,\,j=5, ... ,2n.
\end{cases}
\end{displaymath}
\end{lma}

\begin{rmk}\label{R7.4}
Observe that $\, w_5=w_6= ... = w_{2n} \in U\,$.
\end{rmk}

\begin{proof}
Let $\, w_j\in  {\mathfrak g}\,$ be given by the formulas of Lemma \ref{L7.3}. Since 
$\, w_5,w_6, ... , w_{2n} \in U\,$, we have $\, {\mathbf w}_j = w_j \cdot (-1)^j J \in T_{ (-1)^j J}({\mathcal S})\,$ for $\, j=5, ... ,2n\,$. Therefore $\, \bigoplus\limits _{j=5}^{2n} {\mathbf w}_j \in \xi _1\,$.

For $\, i=1, ..., 4\,$ we have $\, w_i=v_i\,$ and, hence, $\, {\mathbf w}_i = w_i \cdot A_i =  v_i \cdot A_i =  {\mathbf v}_i \in T_{A_i}({\mathcal S})\,$.

It follows that 
\begin{displaymath}
 \bigoplus\limits _{j=1}^{2n} {\mathbf w}_j- {\mathbf v} =  \bigoplus\limits _{j=1}^{2n} {\mathbf w}_j-  \bigoplus\limits _{i=1}^4 {\mathbf v}_i =  \bigoplus\limits _{j=5}^{2n} {\mathbf w}_j \in \xi _1 \,\, .
\end{displaymath}

To prove that $\, \varphi ({\mathbf v}) = \bigoplus\limits _{j=1}^{2n} {\mathbf w}_j\,$ it is therefore enough  to show that $\, \bigoplus\limits _{j=1}^{2n} {\mathbf w}_j\in \xi _2\,$ i.e. to show that  $\, \bigoplus\limits _{j=1}^{2n} {\mathbf w}_j\,$ is symplectic orthogonal to $\, \xi _1\,$ w.r.t. the $2$-form $\,\omega _c\,$.

Let $\, u\in U\,$ and let $\, {\mathbf u}_k = u\cdot (-1)^k J \in T_{A_k}({\mathcal S}), \,\, k=5, ..., 2n\,$. Then, for $\, i=1, ... ,4\,$, 
\begin{equation}{\label{7.10}}
\begin{split}
\omega _c({\mathbf v}_i, \, {\mathbf u}_k ) &= -\omega _{[z_1...z_{k-1} \vert z_k]}( {\mathbf v}_i, \, {\mathbf u}_k ) = -\tfrac 12 \,\,\omega(df_{z_1...z_{k-1}}( {\mathbf v}_i)) \bullet \overline{\omega} (df_k( {\mathbf u}_k)) = \\
&= -\tfrac 12\,\, \omega (A_1...A_{i-1}\cdot v_i\cdot A_i ... A_4(-1)^5J...(-1)^{k-1}J) \bullet  \overline{\omega} ( u\cdot (-1)^k J ) = \\
& = -\tfrac 12\,\, \omega (A_1...A_4(-1)^5J...(-1)^{k-1}J\cdot Ad((A_i...A_4J^{k-5})^{-1})(v_i))\bullet \\ 
&\qquad \qquad \qquad \qquad \qquad \qquad \qquad \qquad \qquad \qquad \qquad  \bullet  \overline{\omega} ( u\cdot (-1)^k J ) = \\
&=-\tfrac 12\,\, ( Ad(J^{k-5}A_4 ... A_i)(v_i))\bullet u = \\
&=-\tfrac 12\,\, ( Ad(A_4 ... A_i)(v_i))\bullet (Ad(J^{k-5})(u))=\\
&=(-1)^k\, \tfrac 12\,\,  Ad(A_4 ... A_i)(v_i))\bullet u=\\
&=(-1)^k\, \tfrac 12\,\, \bigl( \text{proj} ( Ad(A_4 ... A_i)(v_i))\bullet u \,\bigr)
\end{split}
\end{equation}

Furthermore, for $\, 5\le j \le k\,$ and $\, {\mathbf w}_j = w_j \cdot A_j =  w_j \cdot (-1)^jJ\,$,
\begin{equation*}
\begin{split}
\omega _c({\mathbf w}_j, \, {\mathbf u}_k ) &= -\omega _{[z_1...z_{k-1} \vert z_k]}( {\mathbf w}_j, \, {\mathbf u}_k ) =\\
&= \begin{cases}
-\tfrac 12 \,\,\omega(df_{z_1...z_{k-1}}( {\mathbf w}_j)) \bullet \overline{\omega} (df_k( {\mathbf u}_k))& \qquad \text{if} \,\, 5\le j <k,\\
\quad 0 &\qquad \text{if} \,\, j=k,
\end{cases}\\
&= \begin{cases}
-\tfrac 12 \,\,\omega(A_1...A_4(-1)^{k-1}J^{k-5}\cdot Ad(J^{k-j})(w_j))\, \bullet &\\
\qquad \qquad \qquad \qquad \qquad \qquad \bullet \,\,\overline{\omega} (u\cdot (-1)^k J))&\qquad \text{if} \,\, 5\le j <k,\\
\quad 0 &\qquad \text{if} \,\, j=k,
\end{cases}\\
&= \begin{cases}
-\tfrac 12 \,\, Ad(J^{k-j})(w_j))\, \bullet 
\, u&\qquad \text{if} \,\, 5\le j <k,\\
\quad 0 &\qquad \text{if} \,\, j=k,
\end{cases}\\
&= \begin{cases}
\tfrac 12 \,(-1)^{k-j-1}\,w_j\, \bullet 
\, u&\quad \qquad \,\,\, \text{if} \,\, 5\le j <k,\\
\quad 0 &\quad \qquad \,\,\, \text{if} \,\, j=k.
\end{cases}
\end{split}
\end{equation*}

Exchanging the roles of $\, {\mathbf w}_j\,$ and $\, {\mathbf u}_k\,$ we get, for $\, 5\le j,k \le 2n\,$ and $\, {\mathbf w}_j = w_j \cdot  (-1)^j J\,$,
\begin{equation*}
\omega _c({\mathbf w}_j, \, {\mathbf u}_k )= \begin{cases}
\tfrac 12 \,(-1)^{k-j-1}\,w_j\, \bullet 
\, u&\quad \qquad \,\,\, \text{if} \,\, j <k,\\
\quad 0 & \,\,\, \quad \qquad \text{if} \,\, j=k,\\
\tfrac 12 \,(-1)^{k-j}\,w_j\, \bullet 
\, u&\quad \qquad \,\,\, \text{if} \,\, k <j .
\end{cases}
\end{equation*}
Hence, with $\, w_5=w_6= ... =w_{2n}\,$, we have 
\begin{equation}{\label{7.11}}
\begin{split}
\omega _c ( \sum\limits _{j=5}^{2n} {\mathbf w}_j , \,  {\mathbf u}_k )&= \sum\limits _{j=5}^{2n} \omega _c ( {\mathbf w}_j , \,  {\mathbf u}_k )= 
 \tfrac 12 \bigl(\sum\limits _{j=5}^{k-1}(-1)^{k-j-1} +  \sum\limits _{j=k+1}^{2n}(-1)^{k-j} \bigr) (w_5 \bullet u )=\\
&= \tfrac 12 (-1)^k \,w_5 \bullet u\, . 
\end{split}
\end{equation}
Since $\,  w_5=w_6= ... =w_{2n}= -\sum\limits _{i=1}^4 \text{proj}(Ad(A_4A_3..A_i)(v_i))\,$, the identities (\ref{7.10}) and (\ref{7.11}) imply
\begin{equation*}
\begin{split}
\omega _c ( \sum\limits _{j=1}^{2n} {\mathbf w}_j , \,  {\mathbf u}_k )&=\sum\limits _{i=1}^4 \omega _c ( {\mathbf v}_i , \,  {\mathbf u}_k ) +  \sum\limits _{j=5}^{2n} \omega _c ( {\mathbf w}_j , \,  {\mathbf u}_k )=\\
&= \tfrac 12 \, (-1)^k \bigl(\,\, \sum\limits _{i=1}^4 \text{proj}(Ad(A_4A_3..A_i)(v_i)) + w_5\,\,\bigr) \bullet \, u =\\
&= 0 
\end{split}
\end{equation*}
for all $\, {\mathbf u}_k \in T_{A_k}({\mathcal S})\,$ with $\, A_k=(-1)^k J\,$ and $\, k=5, ... ,2n\,$.  Thus $\, \sum\limits _{j=1}^{2n} {\mathbf w}_j   \,$ is symplectic orthogonal to $\, \xi _1\,$ w.r.t. the $2$-form $\, \omega _c\,$. That proves Lemma {\ref{L7.3}}.
\end{proof}

\medskip

\begin{cor}\label{C7.4}
There exists an isotropic subbundle $\, \eta \,$ of $\, \xi _1\,$ such that for every $\, {\mathbf v}\in \zeta\,$
\begin{displaymath}
\varphi ( {\mathbf v}) -  {\mathbf v} \in \eta \,\, .
\end{displaymath}
\end{cor}

\begin{proof}
Let $\, \eta = \{ \sum\limits _{j=5}^{2n} {\mathbf w}_j \in \xi _1 \, \vert \, {\mathbf w}_j = w_j \cdot (-1)^j J \,\,\, \text{with} \,\,\, w_5=w_6= ... =w_{2n}\in U \,\, \}\,$. Then, according to Lemma \ref{L7.3} and Remark \ref{R7.4},
\begin{displaymath}
\varphi ( {\mathbf v}) -  {\mathbf v} \in \eta 
\end{displaymath}
for all  $\, {\mathbf v}\in \zeta\,$. Furthermore, if $\,  \sum\limits _{j=5}^{2n} {\mathbf w}_j, \,\,  \sum\limits _{k=5}^{2n} {\mathbf u}_k \in \eta\,$ with $\, {\mathbf w}_j = w_j \cdot (-1)^j J, {\mathbf u}_k = u_k \cdot (-1)^k J\,$ 
and 
$\, w_j, u_k \in U\,$ then $\, w_5= ... =w_{2n}, \,\, u_5 = ... = u_{2n}\,$ and, as follows from (\ref{7.11}),
\begin{equation*}
\begin{split}
\omega _c \bigl( \sum\limits _{j=5}^{2n} {\mathbf w}_j , \,  \sum\limits _{k=5}^{2n} {\mathbf u}_k \, \bigr)&=  \sum\limits _{k=5}^{2n} \, \omega _c \bigl( \sum\limits _{j=5}^{2n} {\mathbf w}_j , \, {\mathbf u}_k \, \bigr)= \tfrac 12 \, 
 \sum\limits _{k=5}^{2n}(-1)^k ( w_5 \bullet u_k ) =\\
&=\tfrac 12 \, \bigl(\,\,
 \sum\limits _{k=5}^{2n}(-1)^k \,\, \bigr) ( w_5 \bullet u_5 )=\\
&= 0\,\,.
\end{split}
\end{equation*}  
Thus $\,\eta \,$ is an isotropic  subbundle of $\, \xi _1\,$.
\end{proof}

\medskip

\begin{cor}\label{C7.5}
The vector bundle isomorphism $\, \varphi :\zeta \rightarrow \xi _2\,$ satisfies 
\begin{displaymath}
\omega _c (\, \varphi ({\mathbf v}), \,  \varphi ({\mathbf w})\, ) = \omega _c (\, {\mathbf v}, \,  {\mathbf w}\, )
\end{displaymath}
for all $\,  {\mathbf v},   {\mathbf w} \in \zeta _x\,$ and all $\, x\in K'_4\,$. Thus $\, \zeta\,$ and $\, \xi _2\,$ are isomorphic symplectic bundles.
\end{cor}

\begin{proof}
Let $\, x\in K'_4 \,$ and let $\, {\mathbf v},   {\mathbf w} \in \zeta _x\,$. Let $\, \eta \,$ be the isotropic subbundle of Corollary \ref{C7.4}. Then  
$\, \varphi ( {\mathbf v}) -  {\mathbf v} \in \eta _x, \,\,\,  \varphi ( {\mathbf w}) -  {\mathbf w} \in \eta _x\,$ and, therefore, 
\begin{displaymath} 
 \omega _c ({\mathbf v}   -\, \varphi ({\mathbf v}), \, {\mathbf w}  -  \varphi ({\mathbf w})\, ) = 0\,.
\end{displaymath}
Moreover, 
\begin{displaymath}  \omega _c ({\mathbf v}   -\, \varphi ({\mathbf v}), \,  \varphi ({\mathbf w})\, ) = 0 =  \omega _c (\, \varphi ({\mathbf v}), \, {\mathbf w}  -  \varphi ({\mathbf w})\, )\,
\end{displaymath}
since $\,  \varphi ({\mathbf v}), \, \varphi ({\mathbf w}) \in \xi _{2_x}\,$ and $\, {\mathbf v}   -\, \varphi ({\mathbf v}), \, {\mathbf w}  -  \varphi ({\mathbf w})\in \xi _{1_x}\,$ and $\,  \xi _{2_x}\,$ and  $\,  \xi _{1_x}\,$ are symplectic orthogonal w.r.t. $\, \omega _c\,$ to each other  by the definition of $\, \xi _2\,$.

Consequently,
\begin{equation*}
\begin{split}
\omega _c (\, {\mathbf v}, \,  {\mathbf w}\, )& = \omega _c ( ({\mathbf v}   -\, \varphi ({\mathbf v})) +\varphi ({\mathbf v}) , \, ({\mathbf w}  -  \varphi ({\mathbf w})) +  \varphi ({\mathbf w})\, )=\\
&=\omega _c ({\mathbf v}   -\, \varphi ({\mathbf v}), \, {\mathbf w}  -  \varphi ({\mathbf w})\, )+   \omega _c ({\mathbf v}   -\, \varphi ({\mathbf v}), \,  \varphi ({\mathbf w})\, )+ \\
&\qquad \qquad \qquad \qquad \qquad + \omega _c (\, \varphi ({\mathbf v}), \, {\mathbf w}  -  \varphi ({\mathbf w})\, )+ \omega _c (\, \varphi ({\mathbf v}), \,  \varphi ({\mathbf w})\, )=\\
&= \omega _c (\, \varphi ({\mathbf v}), \,  \varphi ({\mathbf w})\, )\,\, .
\end{split}
\end{equation*}
\end{proof}

Corollary \ref{C7.5} implies that $\, c_1(\xi _2 ) = c_1( \zeta )\,$ and, following (\ref{7.8}), that
\begin{equation}\label{7.12}
c_1(\xi ) = c_1(\xi _2) = c_1(\zeta ) \,\, .
\end{equation}

\medskip

We identify the space $\, K_4\,$ with $\, K'_4\,$ by identifying points $\, (A_1, A_2, A_3, A_4) \in K_4\,$ with points  $\, (A_1, A_2, A_3, A_4, -J, J, ... , -J, J) \in K'_4\,$. The bundle $\,\zeta \,$  over    $\, K'_4\,$ is then identified with the restriction $\, TP_4 \vert _{  K_4}\,$ of the tangent bundle of $\, P_4\,$ to  $\, K_4\,$. We denote even this bundle by $\, \zeta \,$.

Recall that the mapping $\, f:\widetilde{\mathcal S} \rightarrow K_{2n} \,$ factorizes as 
\begin{displaymath} 
\begin{CD}
\widetilde{\mathcal S}@>f>> K_{2n}\\
@VfVV   @AA\cup A\\
K_4 @>=>> K'_4 
\end{CD}
\end{displaymath}

\smallskip

We shall now study the symplectic bundle $\, f^*(\zeta )\,$ over $\, \widetilde{\mathcal S}\,$.

Let $\, B\,$ be a 1-dimensional manifold diffeomorphic  to a circle and described as a union of three subspaces $\, B_1, B_2\,$ and $\, B_3\,$, as in Fig.\ref{almII}.2

$$\includegraphics[width=7cm, height=2cm]{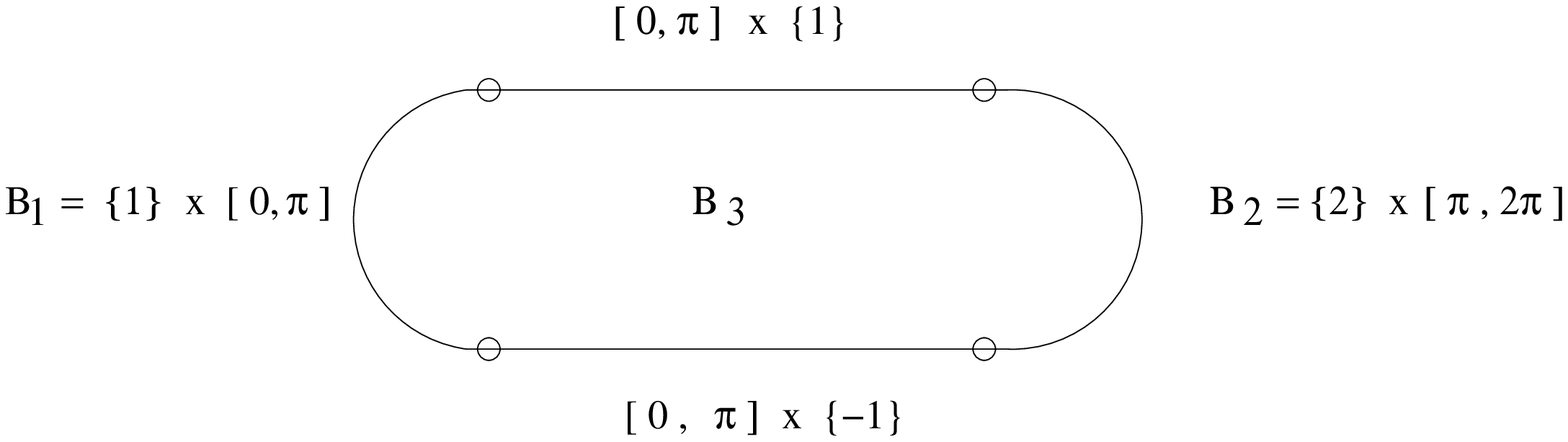}$$

\begin{center}
Fig. \ref{almII}.2
\end{center} 

\smallskip

\noindent
with $\, B_1\,$ and $\, B_2\,$ being arcs and $\, B_3\,$ being a disjoint union of two intervals. We think of $\, B_1\,$ as the interval $\, \{1\}\times [0, \pi ]\,$, of $\, B_2\,$ as the interval $\, \{2\}\times [\pi , 2\pi ]\,$ and of $\, B_3\,$ as the product $\, [0, \pi ]\times \{\pm 1\}\,$ and we identify $\, (1,0)\in B_1=\{1\}\times [0, \pi ]\,$ with $\, (0,1)\in B_3\,$, $\, (1, \pi )\in B_1 =\{1\}\times [0, \pi ]\,$ with $\, (0, -1)\in B_3\,$, $\, (2, \pi )\in B_2=\{2\}\times [\pi , 2\pi ]\,$ with $\, (\pi , -1)\in B_3\,$ and $\, (2, 2\pi )\in B_2\,$ with $\, (\pi , 1)\in B_3\,$, see Fig. \ref{almII}.2.

Let $\, D=B\times I, \,\, I=[0, \pi]\,$, and $\, D_j=B_j\times I, \,\, j=1,2,3\,$, see Fig. \ref{almII}.3.
\vskip0.2truecm
$$\includegraphics[width=7cm, height=2cm]{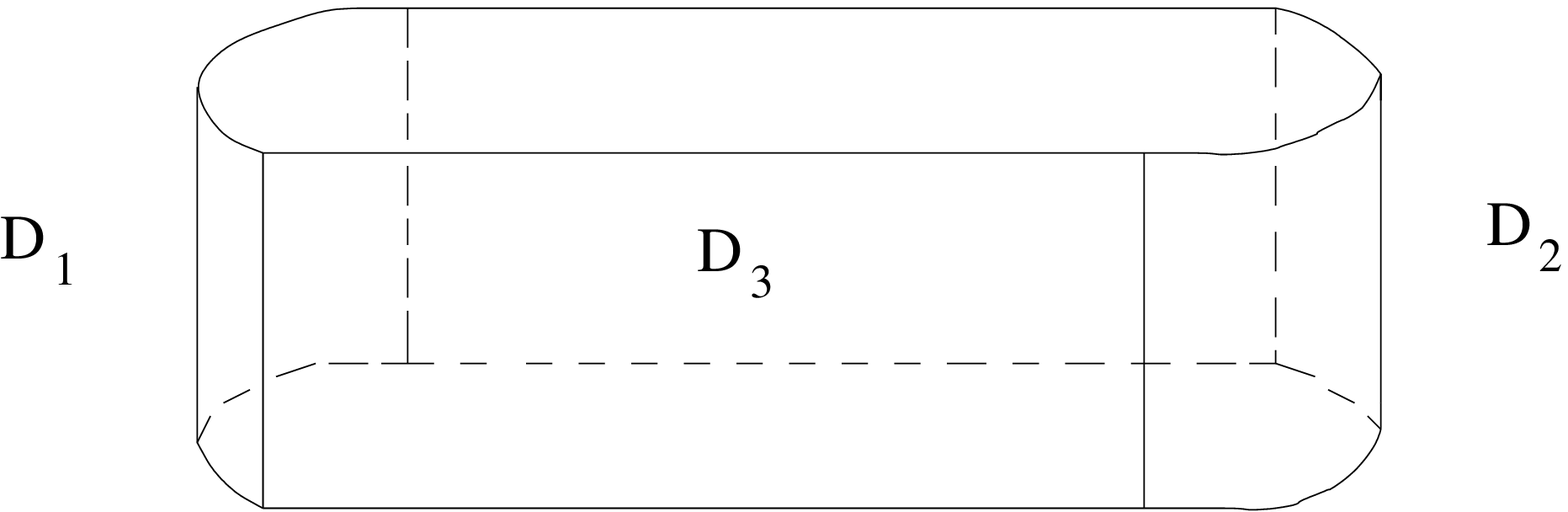}$$
\begin{center}
Fig. \ref{almII}.3
\end{center}

\smallskip
$\, D\,$ is  a $2$-manifold with the boundary consisting of two circles $\, C_0\,$ and $\, C_1\,$, $\, \partial D = C_0\cup C_1\,$, where $\, C_0=B\times\{ 0\}\,$ and $\, C_1=B\times \{ \pi \}\,$.

We define a smooth mapping $\, h:D\rightarrow  \widetilde{\mathcal S}\,$ such that $\, h(D_j)\subset U_j\,$ by
\begin{displaymath}
h(j, \theta , t) = \left(\begin{array}{cc}
i\cos (t)& \sin (t) e^{i\theta }\\
-\sin (t) e^{-i\theta }& -i\cos (t)
\end{array}
\right)\in U_j
\end{displaymath}
if $\, j=1\,$ and $\, (j, \theta , t)\in D_1 = \{ 1\}\times [0, \pi ]\times I\,$ or  $\, j=2\,$ and $\, (j, \theta , t)\in D_2 = \{ 2\}\times [\pi, 2\pi ]\times I\,$, and define
\begin{displaymath}
h(\theta , \epsilon , t ) = (\theta , [\epsilon t ])\in U_3
\end{displaymath}
if $\, j=3\,$ and $\,(\theta , \epsilon , t )\in D_3 = [0, \pi ]\times \{\pm 1\} \times I\,$. 

Observe that $\, h\,$ maps the boundary circles $\, C_k, \,\, k=0,1,\,$ to the subsets $\, \widetilde{C}_k\,$ of $\, U_3\,$ consisting of the points $\, (\theta , [ k\pi ]), \,\, \theta \in [0, \pi ]\,$. The subsets $\, \widetilde{C}_0\,$ and  $\, \widetilde{C}_1\,$ are disjoint and each one is diffeomorphic to an interval, see Fig. \ref{almII}.4.

\vskip0.4truecm
$$\includegraphics[width=5cm, height=3cm]{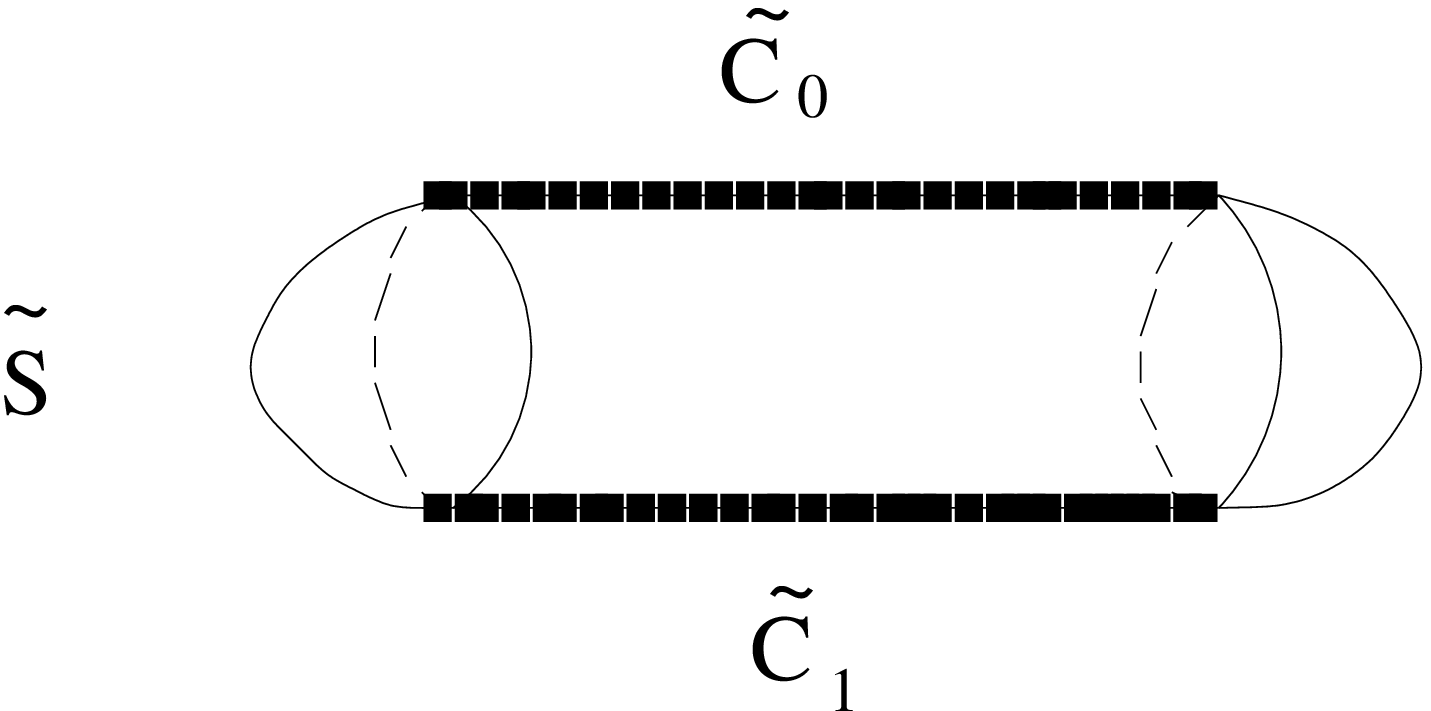}$$

\begin{center}
Fig. \ref{almII}.4
\end{center}

\smallskip

The composition $\, f\circ h\,$ maps $\, D\,$ into $\, K_4\,$. By the definition of $\, f\,$, for every point of $\, U_3\,$ of the form $\, (\theta, [k\pi ]), \,\, k=0,1,\,$ we have
\begin{equation}\label{7.13}
f(\theta, [k\pi ]) = \left( \left( \begin{array}{cc}
i\cos \theta &\sin \theta \\
-\sin \theta & -i\cos \theta
\end{array}\right),\,\,  J,\,\,  (-1)^k J,\,\, (-1)^k\left( \begin{array}{cc}
i\cos \theta &\sin \theta \\
-\sin \theta & -i\cos \theta
\end{array}\right) \right).
\end{equation}

Let us denote by $\, \widehat{C}_k, \,\, k=0,1,\,$ the subsets of $\, K_4\,$ consisting of all the points of the form (\ref{7.13}) with $\, 0\le \theta \le \pi\,$. Again, both $\,  \widehat{C}_0\,$ and $\,  \widehat{C}_1\,$ are diffeomorphic to a compact interval and they are disjoint. The composition $\, f\circ h : D\rightarrow K_4\,$ maps the boundary circles $\, C_k\,$ into the subsets $\, \widehat{C}_k,  \,\, k=0,1\,$.

We shall construct a trivialization of the symplectic bundles $\, \zeta \vert _{ \widehat{C}_k }, \,\, k=0,1,\,$ and exibit a compatible complex structure on them.

Let us first consider the restriction of the bundle $\, \zeta \,$ to $\, \widehat{C}_0\,$. Since $\, \zeta = \bigoplus\limits _{j=1}^4 p^*_j(T{\mathcal S})\,$, at a point
\begin{displaymath}
x=(A_1, A_2, A_3, A_4)=\left( \left( \begin{array}{cc}
i\cos \theta &\sin \theta \\
-\sin \theta & -i\cos \theta
\end{array}\right),\,\,  J,\,\,   J,\,\, \left( \begin{array}{cc}
i\cos \theta &\sin \theta \\
-\sin \theta & -i\cos \theta
\end{array}\right) \right)
\end{displaymath}
of $\,\widehat{C}_0\,$ one has a basis $\, {\mathbf v}_{j,k}=  {\mathbf v}_{j,k}(x), \,\, j=1, ...,4, \,\, k=1,2,\,$ of the real $8$-dimensional fiber $\,\zeta _x\,$ given by 
\begin{equation}\label{7.14}
 {\mathbf v}_{j,k}(x)=v_{j,k} \cdot A_j \in p^*_j(T_{A_j}{\mathcal S})\qquad \text{with} \qquad v_{j,k}=v_{j,k}(x)\in {\mathfrak g} \, ,
\end{equation} 
where
\begin{equation}\label{7.15}
\begin{aligned}
v_{j,1}& = \left( \begin{array}{cc}
0&i\\
i&0
\end{array}\right) \qquad j=1, ... , 4,\\
v_{j,2}&=\begin{cases}
\left( \begin{array}{cc}
i\sin \theta &-\cos \theta \\
\cos \theta & -i\sin \theta
\end{array}\right) & \qquad j=1, 4,\\
\left( \begin{array}{cc}
0&-1\\
1&0
\end{array}\right)& \qquad j=2,3.
\end{cases}
\end{aligned}
\end{equation} 
The basis vectors are ordered lexicographically in $\, (j,k)\,$. The basis $\,\{ {\mathbf v}_{j,k} \}\,$ gives a trivialization of the restriction of the bundle  $\, \zeta \,$ to  $\,\widehat{C}_0\,$.

Since  $\,\widehat{C}_0\subset K_4\,$, we have $\, \omega _{\mathcal C} = \omega _c = -\sum\limits _{j=1}^3 \omega _{[z_1...z_j\vert z_{j+1}]}\,$ on $\, \zeta\,$ (see Corollary 3.4). It follows that $\, \omega _c( {\mathbf v}_{j,k_1}, \,\,  {\mathbf v}_{j,k_2} ) = 0 \,$ for  $\, k_1,k_2 = 1,2,\,$ and, if $\, j_1 < j_2\,$, then 
\begin{equation}\label{7.16}
\begin{split}
\omega _c( {\mathbf v}_{j_1,k_1}, \,\,  {\mathbf v}_{j_2,k_2} )&= 
-\omega _{[z_1...z_{{j_2}-1}\vert z_{j_2}]}( {\mathbf v}_{j_1,k_1}, \,\,  {\mathbf v}_{j_2,k_2} )=\\
&=-\tfrac 12 \omega \bigl( df_{z_1...z_{j_2-1}}( {\mathbf v}_{j_1,k_1})\bigr) \bullet \overline{\omega}\bigl( df_{z_{j_2}}( {\mathbf v}_{j_2,k_2})\bigr)=\\
&=-\tfrac 12 \omega \bigl( A_1...A_{j_2-1}\cdot Ad\bigl((A_{j_1}...A_{j_2-1})^{-1}\bigr)( v_{j_1,k_1})\bigr) \bullet  \overline{\omega}\bigl( v_{j_2,k_2}\cdot A_{j_2})\bigr)= \\
&=-\tfrac 12 \bigl(  Ad\bigl((A_{j_1}...A_{j_2-1})^{-1}\bigr)( v_{j_1,k_1})\bigr) \bullet \,  v_{j_2,k_2} \, .
\end{split}
\end{equation}

\medskip
Let us, from now on, choose the inner product  $\,\bullet \,$ to be $\, a \bullet b = -\frac 12 \, \text{tr} (ab)\,$ for $\, a,b\in {\mathfrak g}= su(2)\,$.
A direct calculation, left to the reader, using (\ref{7.16}), shows that the matrix of the symplectic form   $\, \omega _{\mathcal C}\,$ on $\, \zeta _x\,$ in the basis   $\,\{ {\mathbf v}_{j,k} \}\,$ is equal to 

\begin{equation}\label{7.17}
{\mathcal A}_x = {\mathcal A}(\theta )=\frac 12 \left(\begin{array}{cccccccc}
0&0&1&0&-1&0&1&0\\
0&0&0&c&0&-c&0&1\\
-1&0&0&0&1&0&-1&0\\
0&-c&0&0&0&1&0&-c\\
1&0&-1&0&0&0&1&0\\
0&c&0&-1&0&0&0&c\\
-1&0&1&0&-1&0&0&0\\
0&-1&0&c&0&-c&0&0
\end{array}  \right) \quad \text{with} \quad c=\cos \theta \,\, .
\end{equation}
For example, the entry in the second row and the fourth column is equal to 
\begin{equation*}
\begin{split}
\omega _c( {\mathbf v}_{1,2}, \,\,  {\mathbf v}_{2,2} )&= -\omega _{[z_1 \vert z_2 ]}( {\mathbf v}_{1,2}, \,\,  {\mathbf v}_{2,2} )=\\
&=  -\tfrac 12 \omega \,\bigl( \left(\begin{array}{cc}
i\sin \theta &-\cos \theta\\
\cos \theta & -i\sin \theta
\end{array}  \right) \cdot A_1 \bigr) \bullet \overline{\omega}\,\bigl(\left(\begin{array}{cc}
0 &-1\\
1& 0
\end{array}  \right) \cdot A_2 \bigr)=\\
& =  -\tfrac 12 \omega \,\bigl( A_1 \cdot \left(\begin{array}{cc}
-i\sin \theta &\cos \theta\\
-\cos \theta & i\sin \theta
\end{array}  \right) \bigr) \bullet \overline{\omega}\,\bigl(\left(\begin{array}{cc}
0 &-1\\
1& 0
\end{array}  \right) \cdot A_2 \bigr)=\\
&= -\tfrac 12 ( -\tfrac 12) \text{tr} \bigl( \left(\begin{array}{cc}
-i\sin \theta &\cos \theta\\
-\cos \theta & i\sin \theta
\end{array}  \right)\left(\begin{array}{cc}
0 &-1\\
1& 0
\end{array}  \right)   \bigr) = \tfrac 12 \cos \theta =\\
&= \,\, \,  \tfrac 12 \, c\,\, .
\end{split}
\end{equation*}

\medskip

We shall now construct a complex structure in the bundle $\, \zeta \vert _{\widehat{C}_0}\,$ compatible with the symplectic form  $\, \omega _{\mathcal C}\,$.

Let us consider the subbundle $\, W'\,$ of $\, \zeta \vert _{\widehat{C}_0}\,$ spanned by $\, {\mathbf v}_{1,1}, \,{\mathbf v}_{1,2}, \, {\mathbf v}_{4,1} \,$ and $\, {\mathbf v}_{4,2} \,$. According to (\ref{7.17})
\begin{equation}
\begin{array}{l}
\omega _{\mathcal C}({\mathbf v}_{1,1}, \,  {\mathbf v}_{4,1}) = \omega _{\mathcal C}({\mathbf v}_{1,2}, \,  {\mathbf v}_{4,2})= 1\\
\omega _{\mathcal C}({\mathbf v}_{1,1}, \,  {\mathbf v}_{1,2}) = \omega _{\mathcal C}({\mathbf v}_{1,1}, \,  {\mathbf v}_{4,2})=  \omega _{\mathcal C}({\mathbf v}_{1,2}, \,  {\mathbf v}_{4,1}) = 0.
\end{array}
\end{equation} 
Thus $\, {\mathbf v}_{1,1}, \,{\mathbf v}_{1,2}, \, {\mathbf v}_{4,1}, \, {\mathbf v}_{4,2} \,$ is a symplectic basis of $\, W'\,$.

Let $\, W''\,$ be the symplectic orthogonal complement of $\, W'\,$ in  $\, \zeta \vert _{\widehat{C}_0}\,$. We define sections $\,{\mathbf u}_1, \, {\mathbf u}_2, \,   {\mathbf u}_3 \,$ and $\, {\mathbf u}_4 \,$ of  $\, \zeta \vert _{\widehat{C}_0}\,$ to be 
\begin{equation}\label{7.19}
\begin{array}{ll}
{\mathbf u}_1= {\mathbf v}_{1,1}+ {\mathbf v}_{2,1}- {\mathbf v}_{4,1},& \qquad  {\mathbf u}_2= c\, {\mathbf v}_{1,2}+ {\mathbf v}_{2,2}- c\, {\mathbf v}_{4,2},\\
 {\mathbf u}_3= -{\mathbf v}_{1,1}+ {\mathbf v}_{3,1}+ {\mathbf v}_{4,1}, & \qquad  {\mathbf u}_4=- c\, {\mathbf v}_{1,2}+ {\mathbf v}_{3,2}+ c\, {\mathbf v}_{4,2},
\end{array}
\end{equation}
where $\, c=\cos \theta\,$. A direct check, using (\ref{7.17}) shows that the sections  $\,{\mathbf u}_1, \, {\mathbf u}_2, \,   {\mathbf u}_3, \, {\mathbf u}_4 \,$ belong to and form a symplectic basis of  $\, W''\,$.

It follows that an isomorphism $\, {\mathcal J}: \zeta \vert _{\widehat{C}_0}\rightarrow  \zeta \vert _{\widehat{C}_0} \,$ defined by 
\begin{equation}\label{7.20}
\begin{array}{llll}
{\mathcal J} (  {\mathbf v}_{1,1}) = {\mathbf v}_{4,1},& {\mathcal J} (  {\mathbf v}_{4,1}) = -{\mathbf v}_{1,1},&{\mathcal J} (  {\mathbf v}_{1,2}) = {\mathbf v}_{4,2},&{\mathcal J} (  {\mathbf v}_{4,2}) = -{\mathbf v}_{1,2},\\
{\mathcal J} (  {\mathbf u}_1) = {\mathbf u}_3,& {\mathcal J} (  {\mathbf u}_3) = -{\mathbf u}_1,&{\mathcal J} (  {\mathbf u}_2) = {\mathbf u}_4,&{\mathcal J} (  {\mathbf u}_4) = -{\mathbf u}_2,
\end{array}
\end{equation}
is a complex structure on $\,  \zeta \vert _{\widehat{C}_0}  \,$ compatible with the $2$-form $\,  \omega _{\mathcal C}  \,$. It follows also that the sections
\begin{equation}\label{7.21}
{\mathbf w}_1= {\mathbf v}_{1,1}, \,\,\, {\mathbf w}_2= {\mathbf v}_{1,2}, \,\,\, {\mathbf w}_3 = {\mathbf u}_1, \,\,\,  {\mathbf w}_4 = {\mathbf u}_2,
\end{equation}
of  $\,  \zeta \vert _{\widehat{C}_0}  \,$ form a basis over $\, \C\,$ of the fibres  of  $\,  \zeta \vert _{\widehat{C}_0}  \,$ w.r.t. this complex structure.

In terms of that basis one has
\begin{equation}\label{7.22}
\begin{array}{ll}
{\mathbf v}_{1,1} = {\mathbf w}_1, & {\mathbf v}_{1,2} = {\mathbf w}_2,\\
{\mathbf v}_{2,1} = (-1+{\mathcal J})  {\mathbf w}_1+  {\mathbf w}_3, & 
{\mathbf v}_{2,2} = c \, (-1+{\mathcal J})  {\mathbf w}_2+  {\mathbf w}_4,\\
{\mathbf v}_{3,1} = (1-{\mathcal J})  {\mathbf w}_1+{\mathcal J} {\mathbf w}_3, & 
{\mathbf v}_{3,2} = c \, (1-{\mathcal J})  {\mathbf w}_2+ {\mathcal J} {\mathbf w}_4,\\
{\mathbf v}_{4,1} ={\mathcal J} {\mathbf w}_1, & {\mathbf v}_{4,2} ={\mathcal J} {\mathbf w}_2.
\end{array}
\end{equation}

A similar discussion concerning  the restriction of the bundle $\, \zeta\,$ to the subspace $\, \widehat{C}_1 \,$ differs from the case of  $\, \widehat{C}_0 \,$ only in the fact that $\, A_3=-J\,$ and $\,  A_4=-\left(\begin{array}{cc}
i\cos \theta & \sin \theta \\
-\sin \theta & -i \cos \theta
\end{array} \right)\,$. Everything else follows through in exactly the same way and with the same results.

We now choose and fix for the rest of the paper a complex  structure  $\, {\mathcal J}\,$ on  the bundle  $\, \zeta\,$ compatible with the symplectic $2$-form 
 $\,  \omega _{\mathcal C}  \,$ and such that it coincides  over $\, \widehat{C}_0 \cup \widehat{C}_1\,$ with the complex structure $\, {\mathcal J}\,$ described in (\ref{7.20}).

We have $\, f^*(\zeta )= f^*(TP_4) = f^* (T {\mathscr M}_4)\,$, where $\,  {\mathscr M}_4 ={\mathscr M} \subset  P_4\,$ is an open neigbourhood of $\, K_4\,$ given by Theorem 3.1 in case $\, n=2\,$. 

Let us consider the mapping $\, \widetilde{r}:{\mathscr M}_4 \rightarrow {\mathfrak g}\,$ given by $\,\widetilde{r}(A_1, A_2, A_3, A_4) = \exp ^{-1}(A_1\newline A_2 A_3 A_4)\,$ (see Section 3). According to \cite{GHJW}, Theorem 8.12, $\, \widetilde{r}\,$ is a momentum mapping for the symplectic manifold $\, {\mathscr M}\,$ equipped with the $2$-form  $\,  \omega _{\mathcal C}  \,$ and w.r.t. the diagonal action of $\, G=SU(2)\,$ on  $\, {\mathscr M}\,$ by conjugation.

The subspace $\, K_4 =  \widetilde{r}^{-1} (0)\,$ is a $\, G$-invariant subspace  of $\, P_4\,$. As observed in Section 2, $\, K_4\,$ contains four exceptional  orbits $\, \Delta _{1,1,1,1}, \,  \Delta _{1,-1,1,-1}, \, \Delta _{1,1,-1,-1}, \,   \,$ and  $\, \Delta _{1,-1,-1,1}\,$ of the action of $\, G\,$. All these orbits are diffeomorphic to the $2$-sphere $\, {\mathcal S}\,$ and the isotropy subgroups at their points are isomorphic to circles.

Let us denote by $\, L_4\,$ the complement  of these four orbits in $\, K_4\,$. It is a smooth $\, G$-invariant manifold. The isotropy subgroups  at the points of $\, L_4\,$ are all equal to $\, \{\pm I\}\,$ and the group $\, SO(3)=SU(2)/\{ \pm I\}\,$ acts freely on $\, L_4\,$. 

The matrices $\, X_1= \left(\begin{array}{cc}
i&0\\
0&-i
\end{array}  \right), \,\, X_2=\left(\begin{array}{cc}
0&i\\
i&0
\end{array}  \right), \,\, X_3=\left(\begin{array}{cc}
0&-1\\
1&0
\end{array}  \right)\,$ form a basis of $\, {\mathfrak g}\,$. Through the action of $\, G\,$  they induce tangent vector fields $\, \widetilde{X}_1, \,  \widetilde{X}_2, \, \widetilde{X}_3\,$ on $\, P_4\,$ which, at every point of $\, L_4\,$, are tangent to $\, L_4\,$ and linearly independent. As $\, f^*(\zeta )\,$ is a pullback of the tangent bundle $\, TP_4\,$ the vector fields  $\, \widetilde{X}_1, \,  \widetilde{X}_2, \, \widetilde{X}_3\,$ induce sections $\, 
 \widetilde{Y}_1, \, \widetilde{Y}_2, \, \widetilde{Y}_3\,$ of the bundle $\, f^*(\zeta )\,$.

Let us now consider the pull-back $\, \vartheta = (f \circ h)^*(\zeta ) \,$ of $\, \zeta \,$ to $\, D\,$. It is a symplectic vector bundle of the real fiber dimension $8$ equipped with a compatible complex structure $\, {\mathcal J}\,$. The boundary $\, \partial D\,$ is the  union of the two circles $\, C_0\,$ and $\, C_1\,$ and $\, (f \circ h)(C_i) = \widehat{C_i}\,$ in $\, K_4\,$. Therefore the restriction of $\,\vartheta \,$ to the boundary  $\, \partial D\,$ is equipped with  the basic sections $\, {\mathbf v}_{j,k}\, , \,\, j=1, ... ,4,\,\, k=1,2,\,$ over $\, \R\,$ defined by (\ref{7.14} - \ref{7.15}), with the complex structure   $\, {\mathcal J}\,$ defined by (\ref{7.20}) and also with the basic sections $\, {\mathbf w}_1, ... ,  {\mathbf w}_4\,$ over $\, \C\,$ choosen as in (\ref{7.21}). Since $\, \vartheta =h^* f^* (TP_4)\,$ is a pull-back of $\, TP_4\,$, the vector fields  $\, \widetilde{X}_1, \,  \widetilde{X}_2, \, \widetilde{X}_3\,$ induce sections  $\, 
 Y_1, \,  Y_2, \, Y_3\,$ of  $\, \vartheta \,$.  As the composition $\, f \circ h\,$ maps the interior $\, D-\partial D\,$ of $\, D\,$ into $\, L_4\,$, the sections  $\, 
 Y_1, \,  Y_2, \, Y_3\,$ are linearly independent over every point of  $\, D-\partial D\,$.

Let again $\, x\,$ be a point of $\,  \widehat{C_k}, \, k=0,1\,$,
\begin{displaymath}
x = \left( \left( \begin{array}{cc}
i\cos \theta &\sin \theta \\
-\sin \theta & -i\cos \theta
\end{array}\right),\,\,  J,\,\,  (-1)^k J,\,\, (-1)^k\left( \begin{array}{cc}
i\cos \theta &\sin \theta \\
-\sin \theta & -i\cos \theta
\end{array}\right) \right)
\end{displaymath}
with $ \, 0\le \theta \le\pi\,$. We shall determine the values of the vector fields
 $\, \widetilde{X}_1, \,  \widetilde{X}_2, \, \widetilde{X}_3\,$ at $\, x\,$.

Let $\, X\in {\mathfrak g}\,$ and let $\, g(t)\,$ be a $1$-parameter  subgroup of $\, G=SU(2)\,$ such that $\, g'(0)=X\,$. Let $\, G\,$ act on itself by conjugation and let $\, \widetilde{X}\,$  be the vector field on $\, G\,$ induced by $\, X\,$ through this action. If $\, A\in G\,$ then the value of   $\, \widetilde{X}\,$ at  $\, A\,$ is 
\begin{equation}\label{7.23}
 \widetilde{X}(A) = \frac d{dt}( g(t) A\, g(t)^{-1}) \vert _{t=0} = XA-AX = (X-Ad(A)X )\cdot A \,\, .
\end{equation}
If $\, A=\left( \begin{array}{cc}
i\cos \theta &\sin \theta \\
-\sin \theta & -i\cos \theta
\end{array}\right)\,$ then 
\begin{equation}\label{7.24}
\begin{array}{l}
X_1-Ad(A)X_1 = 2\sin \theta \left( \begin{array}{cc}
i\sin \theta &-\cos \theta \\
\cos \theta & -i\sin \theta
\end{array}\right),\\
X_2-Ad(A)X_2 = 2 \left( \begin{array}{cc}
0 & i \\
i & 0
\end{array}\right),\\
X_3-Ad(A)X_3 = 2\cos \theta \left( \begin{array}{cc}
i\sin \theta &-\cos \theta \\
\cos \theta & -i\sin \theta
\end{array}\right).
\end{array}
\end{equation}
Applying now (\ref{7.24}) to each one of the four coordinates of $\, x\,$ and using (\ref{7.14}),(\ref{7.15}) and (\ref{7.22}) we obtain 
\begin{equation}{\label{7.25}}
\begin{array}{l}
 \widetilde{X}_1(x) = 2 \sin \theta \cdot \, ({\mathbf v}_{1,2} + {\mathbf v}_{4,2}) =  2 \sin \theta \cdot \, (1+{\mathcal J}) {\mathbf w}_2, \\
 \widetilde{X}_2(x) = 2 ({\mathbf v}_{1,1} + {\mathbf v}_{2,1} + {\mathbf v}_{3,1} + {\mathbf v}_{4,1}) = 2 (1+{\mathcal J}) ({\mathbf w}_1 + {\mathbf w}_3 ), \\
 \widetilde{X}_3(x) = 2(\cos \theta \cdot {\mathbf v}_{1,2} + {\mathbf v}_{2,2}  + {\mathbf v}_{3,2} + \cos \theta \cdot {\mathbf v}_{4,2}) = 2 (1+{\mathcal J})(c \, {\mathbf w}_2 + {\mathbf w}_4).
\end{array}
\end{equation}
That gives also values of the sections $\, \widetilde{Y}_1, \widetilde{Y}_2, \widetilde{Y}_3\,$ at the points of $\,\widetilde{C}_0\,$ and  $\,\widetilde{C}_1\,$.

We shall now define a fourth section of the the bundle $\, \vartheta = (f\circ h)^*(\zeta )\,$ over $\, D\,$. Let $\, X_4\,$ be the tangent field on $\, D=B\times I\,$ parallel to the factor $\, I\,$, $\, X_4(b,t) =-\frac d{dt}\,$, see Fig.\ref{almII}.5.

\vskip0.5truecm

$$\includegraphics[width=7cm, height=2cm]{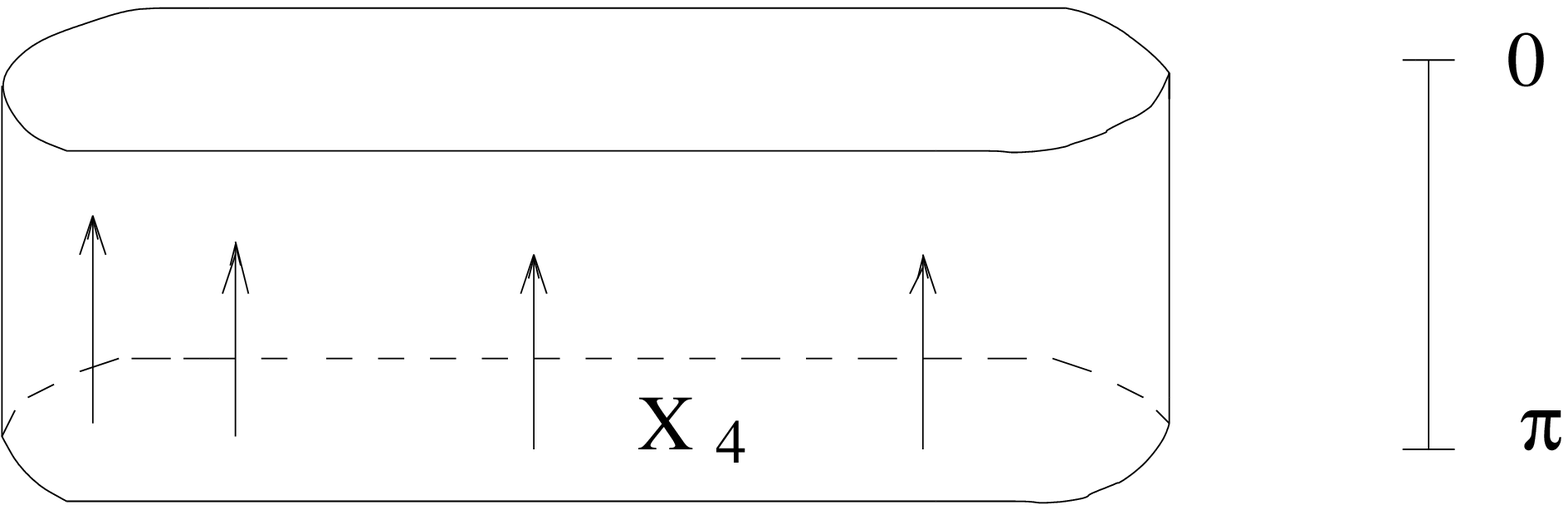}$$

\begin{center}
Fig. \ref{almII}.5
\end{center}

The derivative $\, d(f\circ h)\,$ maps the tangent bundle $\, TD\,$ of $\, D\,$ into $\, TP_4\vert _{K_4} = \zeta \,$. Hence the tangent field  $\, X_4\,$ induces a section $\, Y_4\,$ of the bundle $\, \vartheta = (f\circ h)^*(\zeta )\,$.

\begin{lma}\label{L7.6}
Values of the  sections $\, Y_1, Y_2, Y_3, Y_4\,$ of the bundle  $\, \vartheta\,$  are $\, \R$-linearly independent  over every point of the interior $\, D-\partial D\,$. 
\end{lma}

\begin{proof}
Suppose that at a point $\, (b,t)\in D=B\times I\,$ with $\, 0<t<\pi\,$ (i.e. belonging to the interior of $\, D\,$) one has 
\begin{equation}\label{7.26}
\sum\limits _{i=1}^4 \lambda _i \, Y_i(b,t) = 0
\end{equation} 
for some $\, \lambda _i\in\R, \,\, i=1, ...,4.\,$

Let $\, p_j:P_4\rightarrow {\mathcal S}, \,\, j=1, ...,4,\,$ be the projection on the $j$-th coordinate. As $\, p_2(f(y)) = J\in {\mathcal S}\,$ for all $\, y\in \widetilde{\mathcal S}\,$, we get that $\, p_2(f\circ h)(z)=J\in {\mathcal S}\,$ for all $\, z\in D\,$. Therefore the derivative $\, d(p_2\circ f\circ h)\,$ maps $\, X_4(b,t)\,$ to $\, 0\,$. It follows that $\, \varphi = dp_2\circ f\circ h\,$ maps $\, Y_4(b,t)\,$  to $\, 0\,$ and that it maps $\,\sum\limits _{i=1}^4 \lambda _i \, Y_i(b,t)\,$ to $\,\sum\limits _{i=1}^3 \lambda _i \, \varphi ( Y_i(b,t))= \sum\limits _{i=1}^3 \lambda _i \, dp_2 ( \widetilde{X}_i(x))\,$ with $\, x=f(h(b,t))\,$. Hence (\ref{7.26}) implies 
\begin{displaymath}
 \sum\limits _{i=1}^3 \lambda _i \, dp_2 ( \widetilde{X}_i(x))=0\,\, .
\end{displaymath} 
However, at the point $\, p_2(x) = p_2(f\circ h)(b,t) = J\,$ we have $\, dp_2(\widetilde{X}_1(x))= (X_1-Ad(J) X_1)\cdot J = 0\,$ while $\, dp_2(\widetilde{X}_2(x))= (X_2-Ad(J) X_2)\cdot J= 2\bigl(\begin{smallmatrix}
0&i\\i&0
\end{smallmatrix}\bigr) \cdot J\,$ and $\, dp_2(\widetilde{X}_3(x))= (X_3-Ad(J) X_3)\cdot J= 2\bigl(\begin{smallmatrix}
0&-1\\1&0
\end{smallmatrix}\bigr) \cdot J \,$ are linearly independent in $\, T_J{\mathcal S}\,$. (See (\ref{7.24})). Therefore $\, \lambda _2 = \lambda _3=0\,$.

Furthermore, for every point $\, (b,t)\in D-\partial D\,$ one has $\, (p_3\circ f\circ h)(b,t) \ne \pm J\,$. It follows then directly from the definition of the tangent field $\, X_4\,$ on $\, D\,$ and of the map $\, f\circ h\,$ that $\,( dp_3\circ f \circ h)(Y_4(b,t))\,$ is a non-zero tangent vector to $\,{\mathcal S}\,$ at the point $\, A=(p_3\circ f\circ h)(b,t)\,$ and that this vector is tangent to the great circle joining $\, A\,$ and $\, J\,$ in $\,{\mathcal S}\,$.  On the other hand, the tangent vector $\, ( dp_3\circ f \circ h)(Y_1(b,t)) = dp_3(\widetilde{X}_1((f\circ h)(b,t)))= (X_1-Ad(A)X_1)\cdot A\,$ is transversal to that great circle at the point $\, A\,$. Therefore  $\, ( dp_3\circ f \circ h)(Y_1(b,t))\,$ and  $\,( dp_3\circ f \circ h)(Y_4(b,t))\,$ are $\, \R$-linearly independent vectors in $\, T_A{\mathcal S}\,$. Consequently, $\, \lambda _1=\lambda _4 =0\,$. It follows that the values  of $\, Y_1, Y_2, Y_3, Y_4\,$ are $\, \R$-linearly independent over every point $\, (b,y)\,$ in the interior of $\, D\,$.
\end{proof}

\begin{lma}\label{L7.7}
For every point $\, y\in D-\partial D\,$ the values $\, Y_j(y), \,\, j=1,...,4,\,$ spann a Lagrangian subspace of the fiber $\, \vartheta _y\,$. 
\end{lma}

\begin{proof}
The statement of the Lemma is equivalent to the claim that, for $\, x=(f\circ h)(y)\,$, the values $\, \widetilde{X}_1(x), \widetilde{X}_2(x),  \widetilde{X}_3(x)\,$ and $\,  \widetilde{X}_4(x)=d(f\circ h)(X_4(y))\,$ spann a Lagrangian subspace of $\, \zeta _x\,$.  Now, $\, \zeta _x = T_xP_4= T_x{\mathscr M}_4\,$ and  $\, \widetilde{r}:{\mathscr M}_4\rightarrow {\mathfrak g}, \,\,  \widetilde{r}(A_1, A_2, A_3,\newline A_4 )= \exp ^{-1}( A_1  A_2  A_3 A_4 ),\,$ is a momentum mapping for the symplectic manifold $\, {\mathscr M}_4 \,$ with the diagonal action of $\, G\,$ by the conjugation.

As the vector fields $\, \widetilde{X}_1,  \widetilde{X}_2,  \widetilde{X}_3\,$ are induced by $\, X_1, X_2, X_3 \in {\mathfrak g}\,$ under that action of $\, G\,$ and the point $\, x=(f\circ h)(y)\,$ belongs to the $\, G$-invariant subspace $\, K_4=\widetilde{r} ^{-1}(0)\,$, we have
\begin{displaymath}
 \widetilde{X}_i (x) \in \text{Ker} (\, d\, \widetilde{r}:T_x{\mathscr M}_4 \rightarrow  {\mathfrak g}\, ) \,\,\, , \quad i=1, 2, 3.
\end{displaymath}

Similarly, since $\, f\circ h\,$ maps $\, D\,$ into $\,  K_4=\widetilde{r} ^{-1}(0)\,$, we have 
\begin{displaymath}
 \widetilde{X}_4 (x) \in \text{Ker} (\, d\, \widetilde{r}:T_x{\mathscr M}_4 \rightarrow  {\mathfrak g}\, ) \,\, .
\end{displaymath}
Therefore, by the momentum map property,
\begin{displaymath}
 \omega _{\mathcal C}(\, \widetilde{X}_j (x), \,\widetilde{X}_k (x)\, )=  (\, d\, \widetilde{r} ( \,\widetilde{X}_k (x)\, )\, ) \bullet \, X_j = 0\, \bullet  X_j = 0 \, 
\end{displaymath}
for $\, j=1,2,3\,$ and $\, k=1, ..., 4\,$.

Hence $\, \widetilde{X}_1 (x), ... , \widetilde{X}_4 (x)\,$ spann an isotropic subspace of $\, \zeta _x\,$. By Lemma \ref{L7.6} that subspace has real dimension $4$ and, therefore, is Lagrangian. 
\end{proof}

Let us recall that the bundle $\, \vartheta\,$ is equipped with the complex structure $\, {\mathcal J}\,$ compatible with $\,  \omega _{\mathcal C}\,$.

\begin{cor}\label{C7.9}
For every point $\, y\in D-\partial D\,$ the values $\, Y_j(y), \,\, j=1, ... , 4,\,$ form a basis of the complex vector space $\, \vartheta _y\,$.
\end{cor}

\begin{proof}
The values  $\, Y_j(y), \,\, j=1, ... , 4,\,$ are linearly independent over $\, \R\,$ by Lemma \ref{L7.6} and belong to a Lagrangian subspace of $\, \vartheta _x\,$ by Lemma \ref{L7.7}. Hence, they are linearly independent over $\, \C\,$. As $\, \dim _{\R}\vartheta _x = 8\,$, Corollary \ref{C7.9} follows.
\end{proof}

In order to compute $\,c_1(\zeta ) \,$  we shall now study the behaviour of the sections $\, Y_1, ... , Y_4\,$ in a neighbourhood of the boundary $\, \partial D = C_0 \sqcup C_1\,$. We shall deform continuosly some of the sections.

The mapping $\, f\circ h\,$ maps $\, C_k\,$ onto $\, \widehat{C}_k \subset K_4, \, \, k=0,1\,$. By the definition the values  of $\, Y_1, Y_2, Y_3\,$ at a point $\, y\in \partial D\,$ are given by $\, \widetilde{X}_1(x),\,  \widetilde{X}_2(x),\,  \widetilde{X}_3(x)\in \zeta _x\,$ with $\, x = (f\circ h)(y)\,$. According to (\ref{7.25}), both $\, Y_2(y)= \widetilde{X}_2(x)\,$ and  $\, Y_3(y)= \widetilde{X}_3(x)\,$ are non-zero (and actually $\, \C$-linarily independent) at all points $\, y\in \partial D\,$. By Corollary \ref{C7.9} the same is true for $\, y\in D-\partial D\,$, hence for all $\, y\in D\,$. We shall leave the sections $\, Y_2\,$ and  $\, Y_3\,$ unchanged.

In order to describe the deformations of $\, Y_1\,$ and $\, Y_4\,$ it will be useful to see the sphere $\,{\mathcal S}\,$ as embedded in $\, G=SU(2)\,$, to see the tangent bundle $\, T{\mathcal S}\,$ of the sphere as embedded in the tangent bundle $\, TG\,$ of $\, G\,$, to identify  $\, TG\,$ with $\, {\mathfrak g}\times G\,$ through the {\it right} translations, in that way to look upon every fiber of the bundle  $\, T{\mathcal S}\,$ as embedded in $\, {\mathfrak g}\,$ and, finally, to look upon every fiber of the bundle $\, TP_4 = ( T{\mathcal S})^4\,$ as embedded in the vector space $\,  {\mathfrak g}^4= {\mathfrak g}\oplus {\mathfrak g}\oplus{\mathfrak g}\oplus {\mathfrak g}\,$. Thus a tangent vector $\, ( {\mathbf v}_1,  {\mathbf v}_2,  {\mathbf v}_3,  {\mathbf v}_4 )\in T_{(A_1, A_2, A_3, A_4)}P_4\,$ with $\,  {\mathbf v}_j= v_j\cdot A_j, \,\, A_j\in {\mathcal S}, \, v_j\in {\mathfrak g}, \, \, j=1,...,4,\,$ is embedded as the vector $\, v_1\oplus v_2\oplus v_3\oplus v_4\,$ in  $\,  {\mathfrak g}^4\,$. 

The value $\, Y_1((b,t))\,$ of the section $\, Y_1\,$ at the point $\, (b,t)\in D\,$ can be identified with the value of the section $\, \widetilde{Y}_1= \widetilde{X}_1 \circ f\,$ of the bundle $\, f^*(\zeta )\,$ at the point $\, h(b,t)\in \widetilde{\mathcal S}\,$. According to (\ref{7.25}) and Lemma \ref{L7.6} the section $\, \widetilde{Y}_1  \,$ is non-zero at all except four points of $\,  \widetilde{\mathcal S}\,$. The four exceptional points are: $\, z_1=(0,0)\in U_3=[0, \pi ] \times S^1\,$ identified with $\, J\in U_1, \,\, z_2=(\pi , 0)\in U_3\,$ identified with 
 $\, J\in U_2, \,\, z_3=(0, \pi )\in U_3\,$ identified with  
$\, -J\in U_1\,$ and $\, z_4=(\pi , \pi )\in U_3\,$ identified with $\, -J\in U_2\,$. Of these four points $\, z_1, z_2 \in \widetilde{C}_0\,$ while  $\, z_3, z_4 \in \widetilde{C}_1\,$.

We shall now calculate directional derivatives of $\, \widetilde{Y}_1= \widetilde{X}_1 \circ f\,$ at the points $\, z_1, ... , z_4\,$. The values of  $\, \widetilde{Y}_1\,$ are seen as elements of the vector space  $\,  {\mathfrak g}^4\,$ and so will be their derivatives. Observe that since the section  $\, \widetilde{Y}_1= \widetilde{X}_1 \circ f\,$ vanishes at all $\, z_i\,$, to calculate its directional derivatives amounts to changing its length and taking a limit. Hence, when pulled back over $\, D\,$, the section can be continuously deformed, by changing only its length over $\, D-\partial D\,$, to a section which has the directional derivatives as its values along parts of $\ \partial D\,$.

\bigskip

{\it The directional derivatives of  $\, \widetilde{Y}_1\,$ at $\, z_1\,$ and $\, z_2\,$}

\bigskip

$$\includegraphics[width=9cm, height=3.5cm]{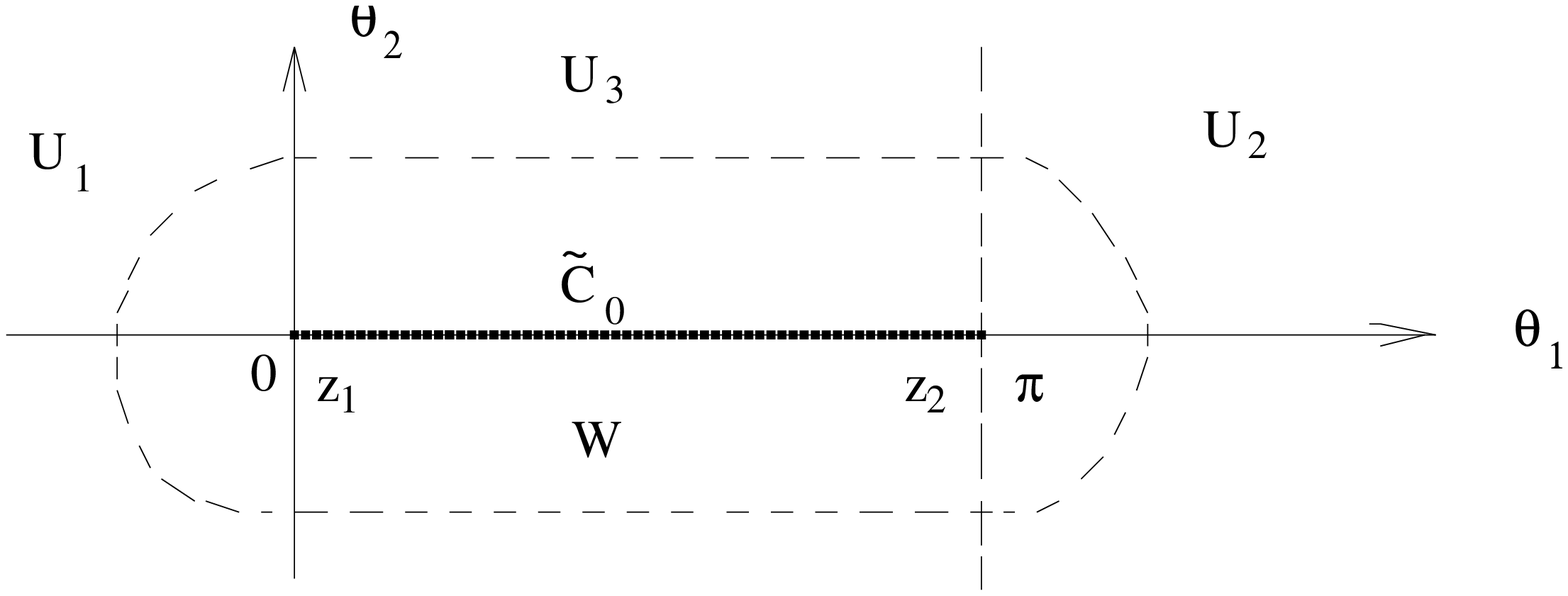}$$
\begin{center}
Fig. \ref{almII}.6
\end{center}

Figure \ref{almII}.6 describes a neighbourhood of $\, \widetilde{C}_0\,$ in $\, \widetilde{\mathcal S}\,$. It is a union of a neighbourhood of  $\, \widetilde{C}_0\,$ in $\, U_3\,$ (the middle strip $\, 0\le \theta _1\le \pi\,$ ), of  a neighbourhood of $\, z_1\,$ in $\, U_1\,$ (in the left half-plane $\, \theta _1\le 0\,$) and of a neighbourhood of $\, z_2\,$ in $\, U_2\,$ (in the right half-plane $\, \theta _1\ge \pi\,$). In the middle strip
\begin{multline*}
f(\theta _1, \, \theta _2 ) =\left( \left(\begin{array}{rr}
i\cos \theta _1& \sin \theta _1\\
-\sin \theta _1&-i\cos \theta _1
\end{array}
\right), \,\, J, \left(\begin{array}{rr}
i\cos \theta _2& \sin \theta _2\\
-\sin \theta _2&-i\cos \theta _2
\end{array}
\right), \,\,\right. \\
\left.\left(\begin{array}{rr}
i\cos (\theta _1 + \theta _2)& \sin (\theta _1+\theta _2)\\
-\sin (\theta _1+\theta _2)&-i\cos (\theta _1 + \theta _2)
\end{array}
\right)\right)\in P_4
\end{multline*}
(see the definition of $\, f\,$). Hence, according to (\ref{7.24}),
\begin{multline*}
\widetilde{Y}_1(\theta _1, \, \theta _2 )=( \widetilde{X}_1 \circ f)(\theta _1, \, \theta _2 )  =\\
=\left( 2\sin \theta _1\left(\begin{array}{rr}
i\sin \theta _1&- \cos \theta _1\\
\cos \theta _1&-i\sin \theta _1
\end{array}
\right), \,\, 0,\,\, 2 \sin \theta _2\left(\begin{array}{rr}
i\sin \theta _2& -\cos \theta _2\\
\cos \theta _2&-i\sin \theta _2
\end{array}
\right), \,\,\right. \\
\left.2 \sin (\theta _1+\theta _2)  \left(\begin{array}{rr}
i\sin (\theta _1 + \theta _2)& -\cos (\theta _1+\theta _2)\\
\cos (\theta _1+\theta _2)&-i\sin (\theta _1 + \theta _2)
\end{array}
\right)\right)\in {\mathfrak g}^4.
\end{multline*}
Thus the directional derivative of $\, \widetilde{Y}_1\,$ at the point $\, (\theta _1, \, \theta _2 )= (0,0)\,$ in the direction of the vector $\, {\mathbf u}=(-\sin (t), \, \cos (t) ), \,\, \pi \le t \le 2\pi,\,$ is
\begin{multline*}
D_{ {\mathbf u}}( \widetilde{Y}_1)(0,0) = 2\left( -\sin (t)\left(\begin{array}{rr}
0&-1\\
1&0
\end{array}
\right), \,\, 0, \,\, \cos(t)\left(\begin{array}{rr}
0&-1\\
1&0
\end{array}
\right), \right. \\ \left. (\cos (t) - \sin (t))\left(\begin{array}{rr}
0&-1\\
1&0
\end{array}
\right)\right) \in {\mathfrak g}^4,
\end{multline*}
which, according to (\ref{7.15}) and (\ref{7.22}), is equivalent to 
\begin{equation}\label{7.27}
\begin{split} 
D_{ {\mathbf u}}( \widetilde{Y}_1)(0,0)& =-2\sin (t)\, ({\mathbf v}_{1,2} + {\mathbf v}_{4,2}) + 2\cos (t)\, ({\mathbf v}_{3,2} + {\mathbf v}_{4,2})=\\
&=-2\sin (t)\,(1+{\mathcal J}){\mathbf w}_2 + 2\cos (t)\, (c(1-{\mathcal J}) + {\mathcal J}){\mathbf w}_2 + {\mathcal J}{\mathbf w}_4)=\\
&=-2\sin (t)\, (1+{\mathcal J}){\mathbf w}_2 + 2\cos (t)\,  ({\mathbf w}_2 + {\mathcal J}{\mathbf w}_4).
\end{split}
\end{equation} 
The last equality holds because $\, c=\cos \theta _1 = 1\,$ at $\, (\theta _1, \, \theta _2 )=(0,0)\,$.

Let now $\, k(t)\,$ be the real analytic function $\, k(t) = \dfrac {\sin (t)}t= \sum\limits _{i=0}^{\infty} (-1)^i \dfrac {t^{2i}}{(2i+1)!}\,$ and let $\, r=\sqrt{\theta _1^2 + \theta _2^2}\,$. Define $\, k( \theta _1, \, \theta _2 )\,$ to be   $\, k( \theta _1, \, \theta _2 )= k(r) = \sum\limits _{i=0}^{\infty} (-1)^i \dfrac {(\theta _1^2 + \theta _2^2)^i}{(2i+1)!}\,$.  We have $\, k(0,0)=1\,$ and $\, \dfrac{\partial k}{\partial \theta _1} (0,0) = 0 =  \dfrac{\partial k}{\partial \theta _2} (0,0)\,$.

To describe $\, f\,$ in the left half-plane we can choose coordinates $\,( \theta _1, \, \theta _2 )\,$ there in such a way that 
\begin{equation}\label{7.28}
f(\theta _1, \, \theta _2 ) = (J, J, A, A)
\end{equation} 
with $\, A=\left(\begin{array}{rr}
is&v\\
-\overline{v}&-is
\end{array}
\right), \,\, s=\cos (r) \in \R \,$ and $\, v=v(\theta _1, \, \theta _2 )=  k( \theta _1, \, \theta _2 )\cdot (\theta _2 - i \theta _1 ) \in \C\,$.

According to (\ref{7.23}) we get 
\begin{displaymath}
\widetilde {Y}_1 (\theta _1, \, \theta _2 ) = (B, B, C, C)
\end{displaymath}
with $\, B=X_1-Ad(J)X_1= 0\,$ and 
\begin{multline*}
C=X_1-Ad(A)X_1=\left(\begin{array}{rr}
i&0\\
0&-i
\end{array}
\right) - A\left(\begin{array}{rr}
i&0\\
0&-i
\end{array}
\right) A^{-1}=\\= \left(\begin{array}{cc}
i(1+|v|^2-s^2)&-2sv\\
2s\overline{v}&-i(1+|v|^2-s^2)
\end{array}
\right).
\end{multline*}
We have $\, v(0,0)=0, \, \, \, \, \dfrac {\partial v}{\partial \theta _1}(0,0) = -i, \,\, \,\,  \dfrac {\partial v}{\partial \theta _2}(0,0) = 1, \,\, \,\, s(0,0)=1, \,\,   \,\,  \dfrac {\partial s}{\partial \theta _1}(0,0) = 0 =\dfrac {\partial s}{\partial \theta _2}(0,0)\,$. Thus 
\begin{displaymath}
\biggl. \frac {\partial C}{\partial \theta _1}\biggr\vert _{(\theta _1, \, \theta _2 )=(0,0)}= 2\left(\begin{array}{rr}
0&i\\
i&0
\end{array}
\right)\qquad \text{and} \qquad \biggl. \frac {\partial C}{\partial \theta _2}\biggr\vert _{(\theta _1, \, \theta _2 )=(0,0)}= 2\left(\begin{array}{rr}
0&-1\\
1&0
\end{array}
\right)\,\, .
\end{displaymath}

It follows that the directional derivative of $\, \widetilde{Y}_1\,$ at the point  $\, (\theta _1, \, \theta _2 )= (0,0)\,$ in the direction of the vector $\, {\mathbf u}=(-\sin (t), \, \cos (t) ), \,\, 0\le t \le \pi,\,$ is
\begin{equation}\label{7.29}
\begin{split}
D_{\mathbf u} (\widetilde{Y}_1)(0,0)& = -\sin (t)\,\left( 0, \,\, 0, \,\, 2 \left(\begin{array}{rr}
0&i\\
i&0
\end{array}
\right), \,\, 2\left(\begin{array}{rr}
0&i\\
i&0
\end{array}
\right) \right) + \\
& \quad  + \cos (t) \,\,  \left( 0, \,\, 0, \,\, 2 \left(\begin{array}{rr}
0&-1\\
1&0
\end{array}
\right), \,\, 2\left(\begin{array}{rr}
0&-1\\
1&0
\end{array}
\right) \right)=\\
&= -2\sin (t)\,\, ({\mathbf v}_{3,1} +{\mathbf v}_{4,1}) +  2\cos (t) \,\,  ({\mathbf v}_{3,2} +{\mathbf v}_{4,2})=\\
&= -2\sin (t)\,\, ({\mathbf w}_1 +{\mathcal J}{\mathbf w}_3) +  2\cos (t) \,\,   ({\mathbf w}_2 +{\mathcal J}{\mathbf w}_4).
\end{split}
\end{equation}
Thus, for  $\, {\mathbf u}=(-\sin (t), \, \cos (t) )\,$, we have at the point 
$\, (\theta _1, \, \theta _2 )=(0,0)\,$
\begin{equation}\label{7.30}
D_{\mathbf u} (\widetilde{Y}_1)(0,0)= \begin{cases}
-2\sin (t)\,\,({\mathbf w}_1 +{\mathcal J}{\mathbf w}_3)  +  2\cos (t) \,\,   ({\mathbf w}_2 +{\mathcal J}{\mathbf w}_4) & \qquad \text{if} \quad 0\le t \le \pi ,\\
-2\sin (t)\,\, ( 1 +{\mathcal J}){\mathbf w}_2 +  2\cos (t) \,\,   ({\mathbf w}_2 +{\mathcal J}{\mathbf w}_4) & \qquad \text{if} \quad \pi \le t \le 2\pi . 
\end{cases}
\end{equation}

A similar calculation performed at the point $\, (\theta _1, \, \theta _2 )=(\pi ,0)\,$ in the strip $\, 0\le \theta _1 \le \pi\,$ shows that the directional derivative of $\, \widetilde{Y}_1\,$ at  $\, (\theta _1, \, \theta _2 )=(\pi ,0)\,$ in a direction of  $\, {\mathbf u}=(-\sin (t), \, \cos (t) ), \,\, 0\le t\le \pi ,\,$ is equal to 
\begin{equation*}
\begin{split}
D_{\mathbf u} (\widetilde{Y}_1)(\pi ,0)& = -\sin (t)\,\left( 2 \left(\begin{array}{rr}
0&-1\\
1&0
\end{array}
\right) , \,\, 0, \,\, 0, \,\,  2 \left(\begin{array}{rr}
0&-1\\
1&0
\end{array}
\right) \right) + \\
& \quad  + \cos (t) \,\,  \left( 0, \,\, 0, \,\, 2 \left(\begin{array}{rr}
0&-1\\
1&0
\end{array}
\right), \,\, 2\left(\begin{array}{rr}
0&-1\\
1&0
\end{array}
\right) \right),
\end{split}
\end{equation*}
which, according to (\ref{7.15}), is equivalent at $\, \theta =\theta _1 = \pi \,$ to  
\begin{equation*}
\begin{split}
D_{\mathbf u} (\widetilde{Y}_1)(\pi ,0)& =  2\sin (t)\,\, ({\mathbf v}_{1,2} +{\mathbf v}_{4,2}) +  2\cos (t) \,\,  ({\mathbf v}_{3,2} -{\mathbf v}_{4,2})=\\
&= 2\sin (t)\,\, (1+{\mathcal J}){\mathbf w}_2 +  2\cos (t) \,\, ((c(1-{\mathcal J})- {\mathcal J}) {\mathbf w}_2 +{\mathcal J}{\mathbf w}_4)=\\
& =  2\sin (t)\,\, (1+{\mathcal J}){\mathbf w}_2 +  2\cos (t) \,\, (-{\mathbf w}_2 + {\mathcal J}{\mathbf w}_4).
\end{split}
\end{equation*}
(Note that at $\, \theta =\theta _1 = \pi \,$ one has $\, c=\cos \theta _1 = -1\,$.)

A calculation which follows exactly the case considered in the left half-plane $\, \theta _1\le 0\,$ shows that in the right half-plane $\, \theta _1 \ge \pi\,$ and for  $\, {\mathbf u}=(-\sin (t), \, \cos (t) )\,$ with $\, \pi \le t\le 2\pi\,$ one has 
\begin{equation*}
\begin{split}
D_{\mathbf u} (\widetilde{Y}_1)(\pi ,0)& = -\sin (t)\,\left( 0, \,\, 0, \,\, 2 \left(\begin{array}{rr}
0&i\\
i&0
\end{array}
\right) , \,\,  2 \left(\begin{array}{rr}
0&i\\
i&0
\end{array}
\right) \right) + \\
& \quad  + \cos (t) \,\,  \left( 0, \,\, 0, \,\, 2 \left(\begin{array}{rr}
0&-1\\
1&0
\end{array}
\right), \,\, 2\left(\begin{array}{rr}
0&-1\\
1&0
\end{array}
\right) \right),
\end{split}
\end{equation*}
which, according to (\ref{7.15}), is equivalent at $\, \theta =\theta _1 = \pi \,$ to 
\begin{equation*}
\begin{split}
D_{\mathbf u} (\widetilde{Y}_1)(\pi ,0)& = - 2\sin (t)\,\, ({\mathbf v}_{3,1} +{\mathbf v}_{4,1}) +  2\cos (t) \,\,  ({\mathbf v}_{3,2} -{\mathbf v}_{4,2})=\\
&=- 2\sin (t)\,\, ({\mathbf w}_1+ {\mathcal J}{\mathbf w}_3) +  2\cos (t) \,\, ((-{\mathbf w}_2 +{\mathcal J}{\mathbf w}_4).
\end{split}
\end{equation*} 
Thus, for  $\, {\mathbf u}=(-\sin (t), \, \cos (t) )\,$, we have at the point 
 $\, (\theta _1, \, \theta _2 )=(\pi ,0)\,$
\begin{equation}\label{7.31}
D_{\mathbf u} (\widetilde{Y}_1)(\pi ,0)= \begin{cases}
\quad 2\sin (t)\,\,(1 +{\mathcal J}){\mathbf w}_2  +  2\cos (t) \,\,   (-{\mathbf w}_2 +{\mathcal J}{\mathbf w}_4) & \qquad \text{if} \quad 0\le t \le \pi ,\\
-2\sin (t)\,\, ({\mathbf w}_1  +{\mathcal J}{\mathbf w}_3) +  2\cos (t) \,\,   (-{\mathbf w}_2 +{\mathcal J}{\mathbf w}_4) & \qquad \text{if} \quad \pi \le t \le 2\pi . 
\end{cases}
\end{equation}

\medskip

Let $\, W\,$ be an open contractible neighbourhood of the subset $\, \widetilde{C}_0\,$ in $\, \widetilde{\mathcal S}\,$ as indicated in Figure \ref{almII}.7.

$$\includegraphics[width=11cm, height=4cm]{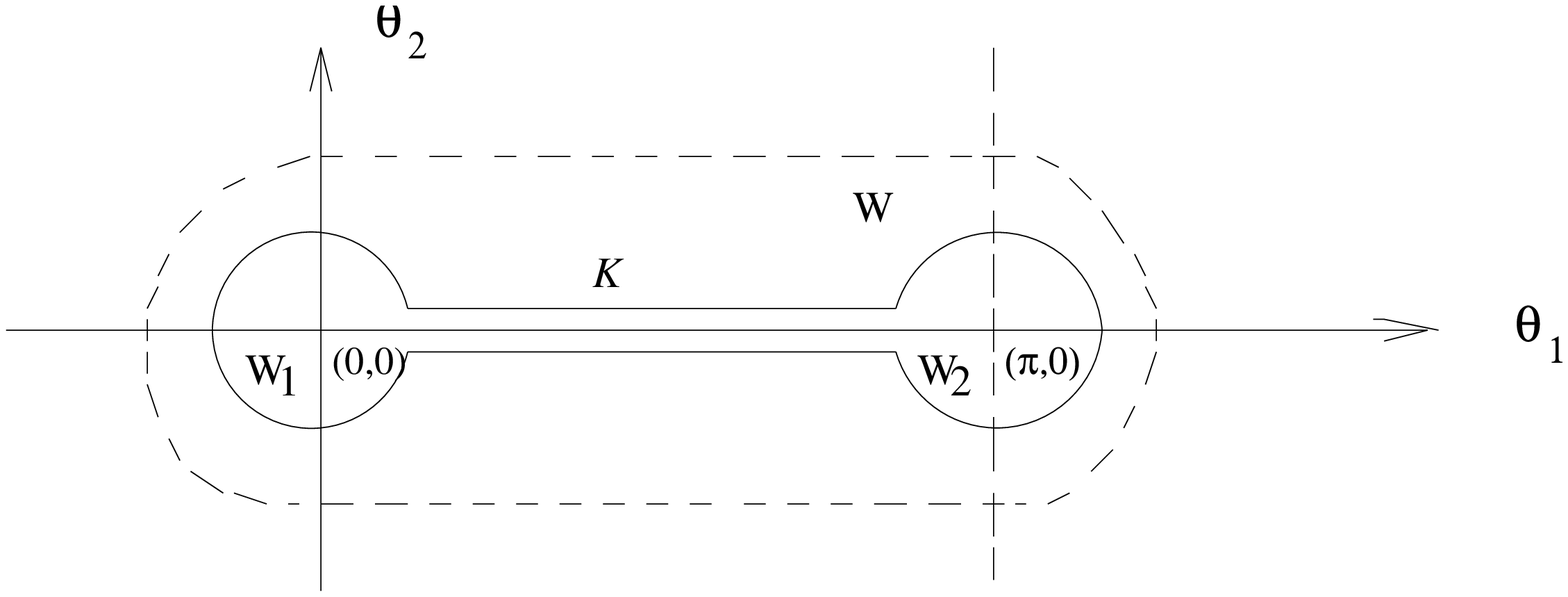}$$

\begin{center}
Fig. \ref{almII}.7
\end{center} 

Let $\, W_1\,$ and $\, W_2\,$ be two small open discs around the point 
 $\, (\theta _1, \, \theta _2 )=(0 ,0)\,$ and around  $\, (\theta _1, \, \theta _2 )=(\pi ,0)\,$ respectively. We introduce a cut in the  $\, (\theta _1, \, \theta _2 )$-plane from  $\, W_1\,$ to  $\, W_2\,$ along the  $\, \theta _1$-axis. Let $\, {\mathcal K}\,$ be the contour in the $\, (\theta _1, \, \theta _2 )$-plane consisting of four paths: $\, \tau _1 = \partial W_1, \,\, \tau _2 =$ lower edge of the cut, $\, \tau _3=\partial W_2\,$ and $\, \tau _4 =$ upper edge of the cut. The contour  $\, {\mathcal K}\,$ is oriented counter-clock-wise. See Figure \ref{almII}.7. 

What we have done above shows that the section  $\, \widetilde{Y}_1\,$ of the bundle $\, f^*(\zeta )\,$ can be deformed continuously on a compact subset of $\, W\,$, by changing  only its length and taking limits, to a section $\, \widetilde{Y}_1'\,$ defined in the complement  of $\,  W_1 \cup W_2\,$  and such that the values  of  $\, \widetilde{Y}_1'\,$ on $\, \tau _1\,$ are given by (\ref{7.30}),\,  its values on  $\, \tau _3\,$ are given by (\ref{7.31}) and that on the both edges  $\, \tau _2, \, \tau _4 \,$  of the cut   $\, \widetilde{Y}_1'\,$ is equal to 
\begin{equation}\label{7.32}
 \widetilde{Y}_1'= 2(1+{\mathcal J})\, {\mathbf w}_2
\end{equation}  
(see (\ref{7.25})). We change the length of the section $\,  \widetilde{Y}_1 =  \widetilde{X}_1 \circ f\,$ along the cut by dividing it by $\, \sin (\theta _1)\,$. The reason for introducing the cut will be explained below.

We shall now consider the values of the section $\, Y_4\,$ over the points of the contour $\,{\mathcal K}\,$. Let us denote by $\, \widetilde{Y}_4\,$ the section of the bundle $\, f^*(\zeta )\,$  which corresponds to  $\, Y_4\,$ over $\, W\,$ with $\, W_1\cup W_2\,$ removed and with the cut between $\, W_1\,$ and   $\, W_2\,$. The reason for introducing the cut is the fact that the values of  $\, \widetilde{Y}_4\,$ over both its edges $\, \tau _3\,$ and  $\, \tau _4\,$ are different.

By the definition of the vector field  $\, X_4\,$  on $\, D\,$ and of the mapping $\,  f\circ h\,$  we get that the values  of  $\, \widetilde{Y}_4\,$ are tangent vectors  to the  curves $\, f(j(s))\,$ in $\, P_4\,$, where the curves $\, j(s)\,$ are given as in the Figure \ref{almII}.8.

$$\includegraphics[width=9cm, height=4cm]{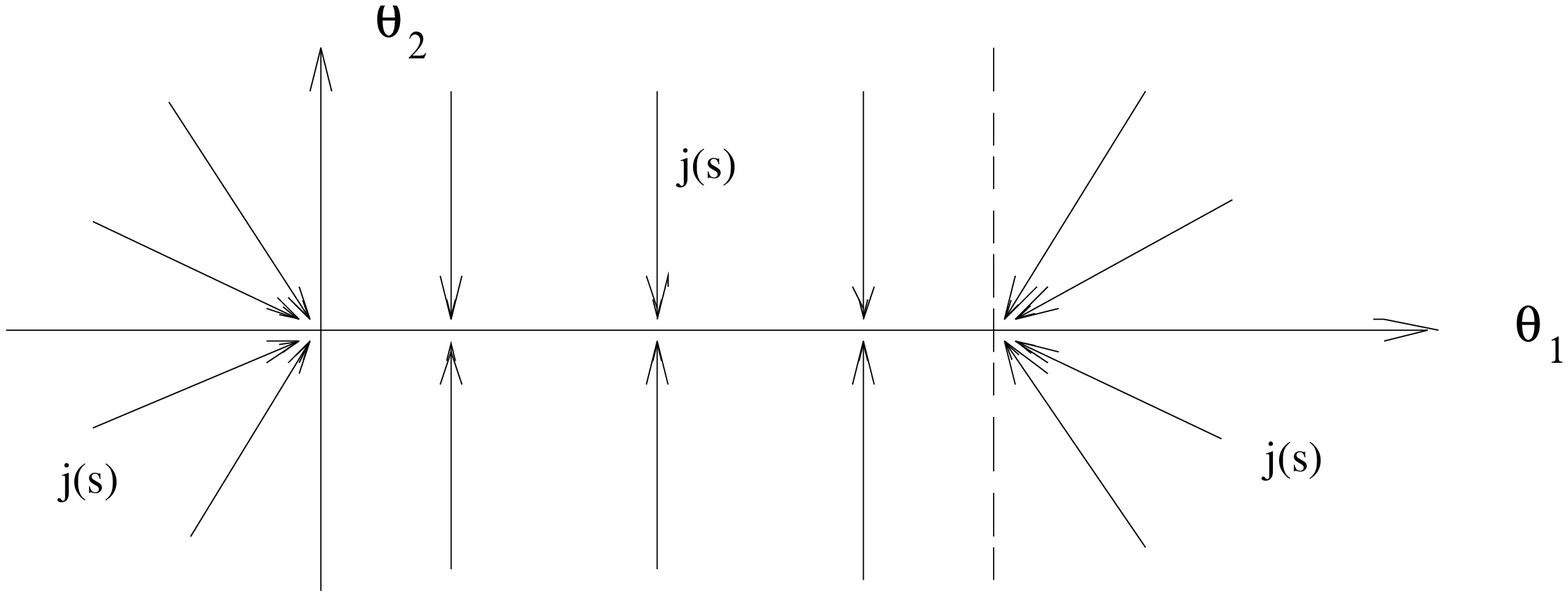}$$

\begin{center}
Fig. \ref{almII}.8
\end{center}
 \medskip
The curves  $\, j(s)\,$  are half-lines parametrised linearly, orthogonal to $\,\theta _1$-axis in the strip $\, 0\le \theta _1 \le \pi \,$ and converging radially to $\, (\theta _1, \, \theta _2 )= (0,0)\,$ in the half-plane $\, \theta _1\le 0\,$ and to  $\, (\theta _1, \, \theta _2 )= (\pi ,0)\,$ in the half-plane $\, \theta _1\ge \pi\,$. All the curves  $\, j(s)\,$ are oriented towards the segment $\, \{\,  0\le \theta _1 \le \pi, \,\, \theta _2=0\, \}\,$.

To find the values of  $\, \widetilde{Y}_4\,$ at the point  $\, (\theta _1, \, 0_+ )\,$ of the upper edge $\, \tau _4\,$ of the cut we take $\, j(s)=(\theta _1, \, -s ), \,\, -\infty <s\le 0\,$. Then
\begin{multline*}
f(j(s))= \left( \left(\begin{array}{cc}
i\cos \theta _1& \sin \theta _1\\
-\sin \theta _1&-i\cos \theta _1
\end{array}
\right), \,\, J, \left(\begin{array}{cc}
i\cos (s)& -\sin (s)\\
\sin (s)&-i\cos (s)
\end{array}
\right), \,\,\right. \\
\left.\left(\begin{array}{cc}
i\cos (\theta _1 - s)& \sin (\theta _1 - s)\\
-\sin (\theta _1 - s)&-i\cos (\theta _1 - s)
\end{array}
\right)\right)\in P_4
\end{multline*}
and 
\begin{multline*}
\widetilde{Y}_4(\theta _1, \, 0_+) =\biggl. \frac d{ds}\bigl(f(j(s))\bigr)\biggr\vert _{s=0}= \left( 0, \,\, 0, \,\, \left(\begin{array}{cc}
0& -1\\
1&0
\end{array}
\right), \,\,\right. \\
\left.\left(\begin{array}{cc}
i\sin \theta _1 & -\cos \theta _1\\
\cos \theta _1&-i\sin \theta _1
\end{array}
\right)\right)\in TP_4.
\end{multline*}
At the point 
\begin{displaymath}
f(\theta _1, \, 0_+) =  \left( \left(\begin{array}{cc}
i\cos \theta _1& \sin \theta _1\\
-\sin \theta _1&-i\cos \theta _1
\end{array}
\right), \,\, J, \,\, J, \,\,\left(\begin{array}{cc}
i\cos \theta _1 & \sin \theta _1 \\
-\sin \theta _1 &-i\cos \theta _1 
\end{array}
\right)\right)
\end{displaymath}
of $\, \widehat{C}_0\,$ the last expression, according to (\ref{7.14}-\ref{7.15}), is equivalent to 
\begin{displaymath}
\widetilde{Y}_4(\theta _1, \, 0_+) = -{\mathbf v}_{3,1} -{\mathbf v}_{4,1}.
\end{displaymath}  
By ( \ref{7.20}-\ref{7.21}) we get 
\begin{equation}\label{7.33}
\widetilde{Y}_4(\theta _1, \, 0_+) = -({\mathbf w}_1 +{\mathcal J}{\mathbf w}_3)\, .
\end{equation}

A similar calculation at the lower edge $\, \tau _2\,$ of the cut gives 
\begin{equation}\label{7.34}
\widetilde{Y}_4(\theta _1, \, 0_-) = {\mathbf v}_{3,1} +{\mathbf v}_{4,1}  = {\mathbf w}_1 +{\mathcal J}{\mathbf w}_3\, .
\end{equation}
Observe that (\ref{7.33}) gives also values of $\, \widetilde{Y}_4\,$ at those parts of $\, \partial W_1\,$ and $\, \partial W_2\,$ which lie in the domain $\, 0\le \theta _1\le \pi , \,\, \theta _2\ge 0,\,$ while (\ref{7.34}) gives values of  $\, \widetilde{Y}_4\,$ at the parts of  $\, \partial W_1\cup \partial W_2\,$ lying in the domain $\, 0\le \theta _1\le \pi , \,\, \theta _2\le 0\,$. Observe also that, since the RHS of  (\ref{7.33}) and (\ref{7.34}) are non-vanishing vectors at all points of $\, \widetilde{C}_0\,$ (including end points), we do not take directional derivatives here.

To find the value of  $\, \widetilde{Y}_4\,$ at a point $\, ({\mathbf u}(t), 0, 0 )\,$ of $\, \partial W_1\,$ lying in the domain $\, \theta _1\le 0\,$ and corresponding to the vector $\, {\mathbf u}(t)= (-\sin(t), \cos (t)), \,\, 0\le t \le \pi,\,$ we take $\, j(s)=(s  \sin(t),-s  \cos (t))\,$ and get 
\begin{displaymath}
f(j(s))= (J, J, A, A) \in P_4
\end{displaymath}
with
\begin{equation}\label{7.35}
A=\left(\begin{array}{cc}
i\cos (s) & v \\
-\overline{v} &-i\cos (s) 
\end{array}
\right), \,\,\, v=-s\, k(-s \, \sin (t), s \, \cos (t))\cdot (\cos (t) + i \sin (t)) \in \C\, .
\end{equation}
(see ( \ref{7.28} ). Thus 
\begin{multline*}
\widetilde{Y}_4({\mathbf u}(t), 0, 0 ) = \biggl. \frac d{ds}\bigl(f(j(s))\bigr)\biggr\vert _{s=0}=\\= \left( 0, \,\, 0, \,\, \left(\begin{array}{cc}
0& z\\
-\overline{z}&0
\end{array}
\right), \,\, \left(\begin{array}{cc}
0& z\\
-\overline{z}&0
\end{array}
\right) \right)\in T_{(J,J,J,J)}P_4
\end{multline*}
with $\, z=-\cos (t)-i \sin (t) \in \C\,$. According to (\ref{7.14}-\ref{7.15}) this equality, at the point $\, (\theta _1, \theta _2)=(0,0)\,$ i.e. with $\, \theta =0\,$, is equivalent to 
\begin{equation*}
\begin{split}
\widetilde{Y}_4({\mathbf u}(t), 0, 0 )&= -\cos (t) \left( 0, \,\, 0, \,\, \left(\begin{array}{cc}
0& 1\\
-1&0
\end{array}
\right), \,\, \left(\begin{array}{cc}
0& 1\\
-1&0
\end{array}
\right) \right)-\\
&\qquad \qquad \qquad -\sin (t) \left( 0, \,\, 0, \,\, \left(\begin{array}{cc}
0& i\\
i&0
\end{array}
\right), \,\, \left(\begin{array}{cc}
0& i\\
i&0
\end{array}
\right) \right)=\\
&= -\cos (t)\, ({\mathbf v}_{3,1} +{\mathbf v}_{4,1}) - \sin (t)\, ({\mathbf v}_{3,2} +{\mathbf v}_{4,2}), \quad 0\le t \le \pi.
\end{split}
\end{equation*} 
That in turn, by (\ref{7.21})-(\ref{7.22}), with $\, c= \cos (0)=1\,$ is equaivalent to 
\begin{equation}\label{7.36}
\widetilde{Y}_4({\mathbf u}(t), 0, 0 )= -\cos (t) \, ( {\mathbf w}_1 + {\mathcal J} {\mathbf w}_3) - \sin (t) \, ( {\mathbf w}_2 + {\mathcal J} {\mathbf w}_4) \,.
\end{equation}

For a similar  calculation for points  $\, ({\mathbf u}(t), \pi , 0 )\,$ of $\, \partial W_2\,$ lying in the domain $\, \theta _1\ge \pi\,$ and corresponding to the vector $\, {\mathbf u}(t)= (-\sin(t), \cos (t))\,$ with $ \, \pi \le t \le 2\pi,\,$ we take $\, j(s)=(\pi + s  \sin(t),-s  \cos (t))\,$ and get 
\begin{displaymath}
f(j(s))= (-J, J, A, -A) \in P_4
\end{displaymath}
with $\, A\,$ as in (\ref{7.35}). Consequently, 
\begin{multline*}
\widetilde{Y}_4({\mathbf u}(t), \pi , 0 ) = \biggl. \frac d{ds}\bigl(f(j(s))\bigr)\biggr\vert _{s=0}=\\= \left( 0, \,\, 0, \,\, \left(\begin{array}{cc}
0& z\\
-\overline{z}&0
\end{array}
\right), \,\, \left(\begin{array}{cc}
0& -z\\
\overline{z}&0
\end{array}
\right) \right)\in T_{(-J,J,J,-J)}P_4
\end{multline*}
with $\, z=-\cos (t)-i \sin (t) \in \C\,$. According to (\ref{7.14}-\ref{7.15}) this equality, at the point $\, (\theta _1, \theta _2)=(\pi ,0)\,$ i.e. with $\, \theta =\pi \,$, is equivalent to 
\begin{equation*}
\begin{split}
\widetilde{Y}_4({\mathbf u}(t), \pi , 0 )&= -\cos (t) \left( 0, \,\, 0, \,\, \left(\begin{array}{cc}
0& 1\\
-1&0
\end{array}
\right), \,\, \left(\begin{array}{cc}
0& -1\\
1&0
\end{array}
\right) \right)-\\
&\qquad \qquad \qquad -\sin (t) \left( 0, \,\, 0, \,\, \left(\begin{array}{cc}
0& i\\
i&0
\end{array}
\right), \,\, \left(\begin{array}{cc}
0& -i\\
-i&0
\end{array}
\right) \right)=\\
&= -\cos (t)\, ({\mathbf v}_{3,1} +{\mathbf v}_{4,1}) - \sin (t)\, ({\mathbf v}_{3,2} - {\mathbf v}_{4,2}), \quad \pi \le t \le 2\pi.
\end{split}
\end{equation*} 
 By (\ref{7.21})-(\ref{7.22}) with $\, c= \cos (\pi )=-1\,$, we get 
\begin{equation}\label{7.37}
\widetilde{Y}_4({\mathbf u}(t), \pi , 0 )= -\cos (t) \, ( {\mathbf w}_1 + {\mathcal J} {\mathbf w}_3) + \sin (t) \, ( {\mathbf w}_2 - {\mathcal J} {\mathbf w}_4) \,.
\end{equation}  
for $\, \pi \le t \le 2\pi\,$.

Let us now calculate the determinant 
\begin{displaymath}
{\mathcal D} =\left\vert \begin{array}{c}
-\,\,\, \widetilde{Y}_1'\,\,\, -\\
-\,\,\, \widetilde{Y}_2\,\,\,  -\\
-\,\,\, \widetilde{Y}_3\,\,\,  -\\
-\,\,\, \widetilde{Y}_4\,\,\,  -
\end{array} \right\vert
\end{displaymath}
of the sections $\, \widetilde{Y}_1', \widetilde{Y}_2, \widetilde{Y}_3, \widetilde{Y}_4\,$ in the $\, \C$-basis $\, {\mathbf w}_1, ... , {\mathbf w}_4\,$ over points of the contour $\, {\mathcal K}\,$. We first subdivide  $\, {\mathcal K}\,$ into paths as in Figure \ref{almII}.9.

\medskip

$$\includegraphics[width=9cm, height=2.5cm]{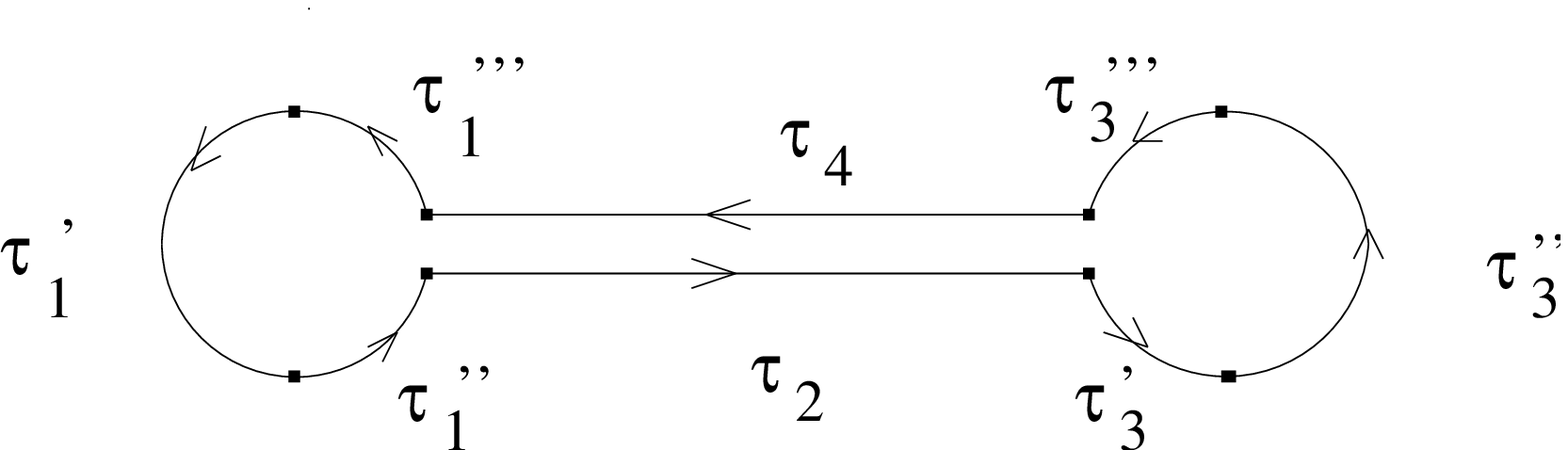}$$

\begin{center}
Fig. \ref{almII}.9 
\end{center}

Over $\, \tau _1'\,$ we have, according to (\ref{7.25}), (\ref{7.30}) and (\ref{7.36}),
\begin{equation}\label{7.38}
\begin{array}{l}
\widetilde{Y}_1'({\mathbf u}(t), 0, 0)= -2\sin (t) \, ({\mathbf  w}_1 + {\mathcal J}{\mathbf  w}_3) + 2 \cos (t) \,  ({\mathbf  w}_2 + {\mathcal J}{\mathbf  w}_4),\\
\widetilde{Y}_2({\mathbf u}(t), 0, 0)= 2 (1+{\mathcal J})\,  ({\mathbf  w}_1 + {\mathbf  w}_3),\\
\widetilde{Y}_3({\mathbf u}(t), 0, 0)= 2 (1+{\mathcal J})\,  ({\mathbf  w}_2 + {\mathbf  w}_4),\\
\widetilde{Y}_4({\mathbf u}(t), 0, 0)= -\cos (t) \, ({\mathbf  w}_1 + {\mathcal J}{\mathbf  w}_3) - \sin (t) \,  ({\mathbf  w}_2 + {\mathcal J}{\mathbf  w}_4),
\end{array}
\end{equation}
with $\, 0\le t \le \pi\,$. Hence
\begin{equation}\label{7.39}
{\mathcal D}= {\mathcal D}(t) = \left\vert \begin{array}{cccc}
-2\sin (t)& 2\cos (t)& -2i\sin (t)&2i \cos (t)\\
2(1+i)&0&2(1+i)&0\\
0&2(1+i)&0&2(1+i)\\
-\cos (t)& -\sin (t)& -i \cos (t)& -i \sin (t)
\end{array}\right\vert = 2^5
\end{equation}
for $\, 0\le t\le \pi\,$ over the path $\, \tau _1'\,$.

Over $\, \tau _1''\,$, according to (\ref{7.25}), (\ref{7.30}) and (\ref{7.34}) together with the  remark following it, one has 
\begin{equation}\label{7.40}
\begin{array}{l}
\widetilde{Y}_1'({\mathbf u}(t), 0, 0)= -2\sin (t) \, (1 + {\mathcal J}){\mathbf  w}_2 + 2 \cos (t) \,  ({\mathbf  w}_2 + {\mathcal J}{\mathbf  w}_4),\\
\widetilde{Y}_2({\mathbf u}(t), 0, 0)= 2 (1+{\mathcal J})\,  ({\mathbf  w}_1 + {\mathbf  w}_3),\\
\widetilde{Y}_3({\mathbf u}(t), 0, 0)= 2 (1+{\mathcal J})\,  ({\mathbf  w}_2 + {\mathbf  w}_4),\\
\widetilde{Y}_4({\mathbf u}(t), 0, 0)= {\mathbf  w}_1 + {\mathcal J}{\mathbf  w}_3,
\end{array}
\end{equation}
for $\, \pi \le t \le \frac 32 \pi\,$. Hence
\begin{equation}\label{7.41}
\begin{split}
{\mathcal D}&= {\mathcal D}(t) = \left\vert \begin{array}{cccc}
0 & 2(\cos (t)-\sin (t) (1+i))&0&2i \cos (t)\\
2(1+i)&0&2(1+i)&0\\
0&2(1+i)&0&2(1+i)\\
1&0& i&0
\end{array}\right\vert =\\
&= -2^5 (\cos (t) -i \sin (t)),
\end{split}
\end{equation}
for $\, \pi \le t\le \frac 32 \pi\,$ over $\, \tau _1''\,$.

Over $\, \tau _1'''\,$ the only difference compared to (\ref{7.40}) is that, according to (\ref{7.33}) and the remark following it, $\,\widetilde{Y}_4({\mathbf u}(t), 0, 0)=-( {\mathbf  w}_1 + {\mathcal J}{\mathbf  w}_3)\,$. Hence 
\begin{equation}\label{7.42}
{\mathcal D}= {\mathcal D}(t) = 2^5 (\cos (t) -i \sin (t)) \qquad \qquad \text{for}\quad \tfrac 32 \pi \le t \le 2\pi
\end{equation} 
over  $\, \tau _1'''\,$.

Over $\, \tau _2\,$, according to (\ref{7.25}), (\ref{7.32}) and  (\ref{7.34}), one has
\begin{equation}\label{7.43}
\begin{array}{l}
\widetilde{Y}_1'(\theta _1, 0)= 2(1 + {\mathcal J}){\mathbf  w}_2,\\
\widetilde{Y}_2(\theta _1  , 0)= 2 (1+{\mathcal J})\,  ({\mathbf  w}_1 + {\mathbf  w}_3),\\
\widetilde{Y}_3(\theta _1  , 0)= 2 (1+{\mathcal J})\,  (c \,{\mathbf  w}_2 + {\mathbf  w}_4),\\
\widetilde{Y}_4(\theta _1  , 0)= {\mathbf  w}_1 + {\mathcal J}{\mathbf  w}_3,
\end{array}
\end{equation}
with $\, 0\le \theta _1 \le \pi\,$ and $\, c=\cos (\theta _1 )\,$. Therefore
\begin{equation}\label{7.44}
{\mathcal D}= {\mathcal D}(\theta _1) = \left\vert \begin{array}{cccc}
0& 2(1+i)&0&0\\
2(1+i)&0&2(1+i)&0\\
0&2c\, (1+i)&0&2(1+i)\\
1&0& i & 0
\end{array}\right\vert = -2^5 i
\end{equation}
for $\, 0\le \theta _1\le \pi\,$ over the path $\, \tau _2\,$.

The case of $\, \tau _4\,$ differs from that of $\, \tau _2\,$ only in $\, 
\widetilde{Y}_4(\theta _1  , 0)= -( {\mathbf  w}_1 + {\mathcal J}{\mathbf  w}_3)\,$. Thus, over  $\, \tau _4\,$ we have 
\begin{equation}\label{7.45}
{\mathcal D}= {\mathcal D}(\theta _1) = 2^5 i
\end{equation}
for $\, 0\le \theta _1\le \pi\,$.

Over $\, \tau _3'\,$, according to (\ref{7.25}), ( \ref{7.31}) and (\ref{7.34}), one has 
\begin{equation}\label{7.46}
\begin{array}{l}
\widetilde{Y}_1'({\mathbf u}(t), \pi , 0)= 2\sin (t) \, (1 + {\mathcal J}){\mathbf  w}_2 + 2 \cos (t) \,  (- {\mathbf  w}_2 + {\mathcal J}{\mathbf  w}_4),\\
\widetilde{Y}_2({\mathbf u}(t), \pi , 0)= 2 (1+{\mathcal J})\,  ({\mathbf  w}_1 + {\mathbf  w}_3),\\
\widetilde{Y}_3({\mathbf u}(t), \pi , 0)= 2 (1+{\mathcal J})\,  (- {\mathbf  w}_2 + {\mathbf  w}_4),\\
\widetilde{Y}_4({\mathbf u}(t), \pi , 0)= {\mathbf  w}_1 + {\mathcal J}{\mathbf  w}_3,
\end{array}
\end{equation}
with $\, \frac {\pi}2 \le t \le \pi\,$. Therefore 
\begin{equation}\label{7.47}
\begin{split}
{\mathcal D}&= {\mathcal D}(t) = \left\vert \begin{array}{cccc}
0 & 2(-\cos (t)+\sin (t) (1+i))&0&2i \cos (t)\\
2(1+i)&0&2(1+i)&0\\
0&-2(1+i)&0&2(1+i)\\
1&0& i&0
\end{array}\right\vert =\\
&= 2^5 (\cos (t) -i \sin (t)),
\end{split}
\end{equation}
for $\,\frac 12 \pi \le t\le \pi\,$ over $\, \tau _3'\,$.

Over $\, \tau _3'''\,$ the only difference compared to (\ref{7.46}) is that, according to (\ref{7.33}), $\, \widetilde{Y}_4({\mathbf u}(t), \pi , 0)= -( {\mathbf  w}_1 + {\mathcal J}{\mathbf  w}_3)\,$. Hence 
\begin{equation}\label{7.48}
{\mathcal D}= {\mathcal D}(t) = -2^5 (\cos (t) -i \sin (t)) \qquad \qquad \text{for}\quad 0 \le t \le \tfrac {\pi}2
\end{equation} 
over  $\, \tau _3'''\,$.

Finally, over  $\, \tau _3''\,$, according to  (\ref{7.25}), ( \ref{7.31}) and (\ref{7.37}), one has 
\begin{equation*}
\begin{array}{l}
\widetilde{Y}_1'({\mathbf u}(t), \pi , 0)=- 2\sin (t) \, ({\mathbf  w}_1 + {\mathcal J}{\mathbf  w}_3) + 2 \cos (t) \,  (- {\mathbf  w}_2 + {\mathcal J}{\mathbf  w}_4),\\
\widetilde{Y}_2({\mathbf u}(t), \pi , 0)= 2 (1+{\mathcal J})\,  ({\mathbf  w}_1 + {\mathbf  w}_3),\\
\widetilde{Y}_3({\mathbf u}(t), \pi , 0)= 2 (1+{\mathcal J})\,  (- {\mathbf  w}_2 + {\mathbf  w}_4),\\
\widetilde{Y}_4({\mathbf u}(t), \pi , 0)= -\cos (t)\, ({\mathbf  w}_1 + {\mathcal J}{\mathbf  w}_3) + \sin (t)\, ({\mathbf  w}_2 - {\mathcal J}{\mathbf  w}_4) 
\end{array}
\end{equation*}
for $\, \pi \le t\le 2\pi\,$. Hence 
\begin{equation}\label{7.49}
\begin{split}
{\mathcal D}&= {\mathcal D}(t) = \left\vert \begin{array}{cccc}
-2\sin (t) & -2\cos (t)&-2 i \sin (t)&2i \cos (t)\\
2(1+i)&0&2(1+i)&0\\
0&-2(1+i)&0&2(1+i)\\
-\cos (t)&\sin (t)& -i \cos (t)&-i \sin (t)
\end{array}\right\vert  = -2^5,
\end{split}
\end{equation}
for $\,\pi \le t\le 2\pi\,$ over $\, \tau _3''\,$.

It follows from (\ref{7.39}) - (\ref{7.49}) that the values of $\, {\mathcal D}\,$ along the contour $\, {\mathcal K}\,$ lie in $\, \C -\{0\}\,$ and transverse once the circle of radius $\, 2^5\,$ with center at $\, 0\,$ in the clock-wise direction. Hence, the winding number of  $\, {\mathcal D}\,$ with respect to  $\, 0\,$ along  $\, {\mathcal K}\,$ is equal to $\, -1\,$,
\begin{equation}\label{7.50}
\text{Ind} _ {\mathcal K} ( \, {\mathcal D}, \, 0\, ) = -1
\end{equation} 
(see Figure \ref{almII}.10). 
$$\includegraphics[width=7cm, height=5.7cm]{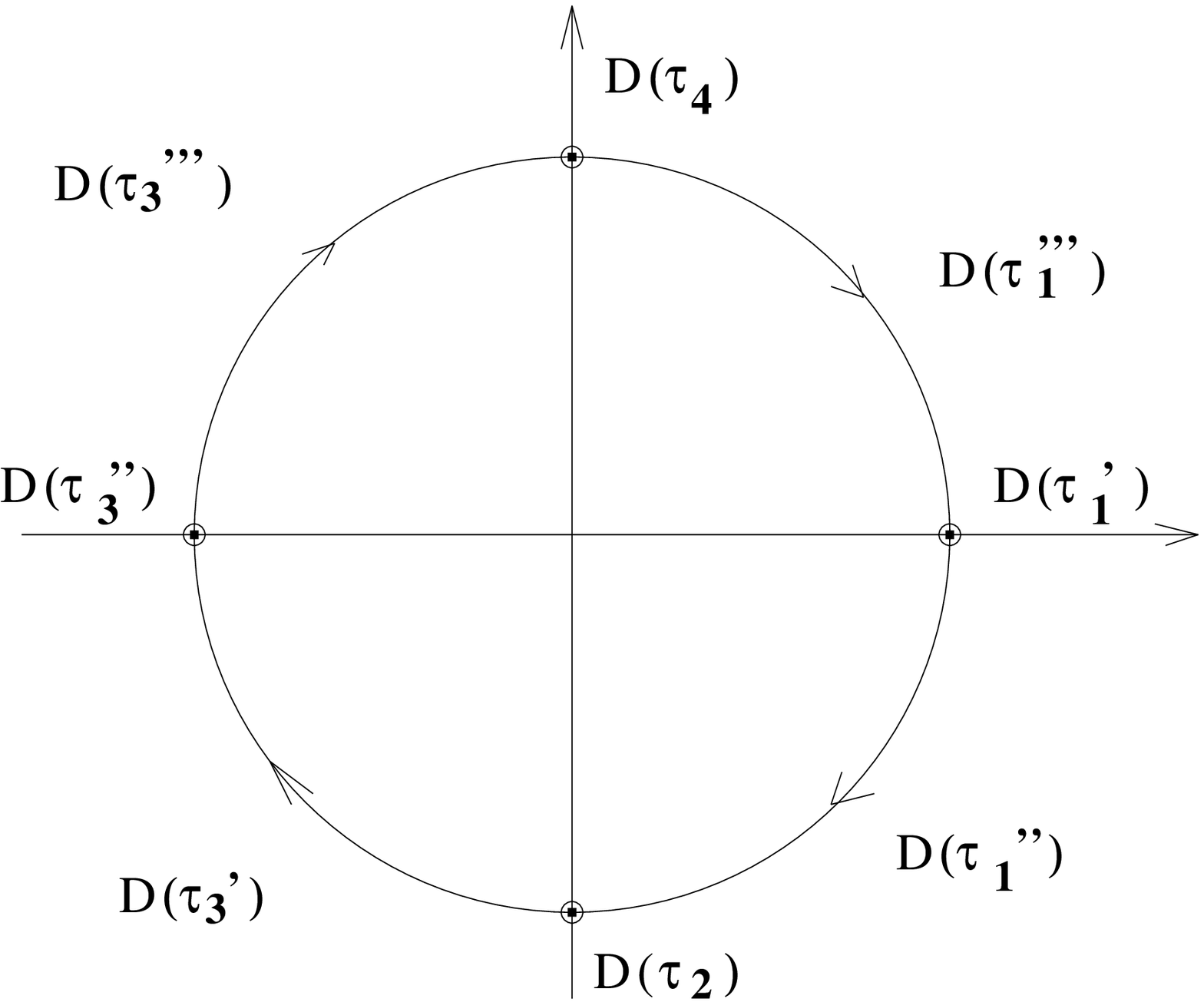}$$
\begin{center}
Fig. \ref{almII}.10
\end{center}

Let us recall that the sphere $\, \widetilde{\mathcal S}\,$ is oriented in such a way that the orientation of $\, U_3= [0, \pi ]\times  \R/2\pi \Z\,$ is equal to the product of the standard orientations of  $\,  [0, \pi ]\,$ and of $\,  \R/2\pi \Z\,$. The contour $\, {\mathcal K}\,$ is then positively oriented.

We need also to study the sections $\, \widetilde{Y}_1, ... ,\widetilde{Y}_4\,$ in a neighbourhood of the subset $\, \widetilde{C}_1\,$. Let $\, W'\,$ be a neighbourhood of  $\, \widetilde{C}_1\,$ in $\, \widetilde{\mathcal S}\,$ analogous to the neighbourhood $\, W\,$ of  $\, \widetilde{C}_0\,$. Let  $\, {\mathcal K}'\,$ be a contour in  $\, W'\,$ analogous to the contour  $\, {\mathcal K}\,$ in  $\, W\,$. Let  $\, {\mathcal K}'\,$ be oriented in the same way with respect to the orientation of  $\, \widetilde{\mathcal S}\,$  as  $\, {\mathcal K}\,$ is. Let $\,  {\mathcal D}'\,$ be the corresponding determinant function. By calculations analogous to those in the case of  $\,  {\mathcal D}\,$ one proves that also in that case 
\begin{equation}\label{7.51}
\text{Ind} _ {{\mathcal K}'} ( \, {\mathcal D}', \, 0\, ) = -1\, .
\end{equation} 
(Actually, if $\, {\mathcal K}'\,$ is obtained from  $\, {\mathcal K}\,$ through a translation along the $\, \theta _2$-axis then we get $\,{\mathcal D}'=-{\mathcal D}\,$. That proves (\ref{7.51}).)

As  $\, \widetilde{Y}_1, ... ,\widetilde{Y}_4\,$ are sections of  the bundle $\, f^*(\zeta )\,$ which are  $\, \C$-linearly independent over   $\, \widetilde{\mathcal S} - ( \widetilde{C}_0 \cup \widetilde{C}_1) \,$ we get 
\begin{equation*}
\langle \,\, [ \widetilde{\mathcal S}]\, , \, c_1(f^*(\zeta )) \, \rangle = \text{Ind} _ {\mathcal K} ( \, {\mathcal D}, \, 0\, ) + \text{Ind} _ {{\mathcal K}'} ( \, {\mathcal D}', \, 0\, ) = -2\, .
\end{equation*} 
Since, according to (\ref{7.12}), $\, c_1(\zeta ) = c_1 ( \xi )\,$, one has finally 
\begin{equation*}
\begin{split}
\langle \,\,f_* [ \widetilde{\mathcal S}]\, , \, c_1 \, \rangle &= \langle \,\,f_* [ \widetilde{\mathcal S}]\, , \, c_1 ( {\mathscr M} ) \, \rangle = 
\langle \,\,f_* [ \widetilde{\mathcal S}]\, , \, c_1 ( \xi ) \, \rangle =\\
& = \langle \,\,f_* [ \widetilde{\mathcal S}]\, , \, c_1 ( \zeta ) \, \rangle =
\langle \,\, [ \widetilde{\mathcal S}]\, , \, c_1(f^*(\zeta )) \, \rangle =\\
&= -2 \, .
\end{split}
\end{equation*}
That proves Theorem \ref{almII}.1.
\end{proof}

\vskip1truecm 
\end{document}